\documentclass[psamsfonts, 11pt]{amsbook}

\newtheorem{theorem}{Theorem}[subsection]
\newtheorem{lemma}[theorem]{Lemma}
\newtheorem{proposition}[theorem]{Proposition}
\newtheorem{corollary}[theorem]{Corollary}

\theoremstyle{definition}

\newtheorem{definition}[theorem]{Definition}
\newtheorem{example}[theorem]{Example}

\newtheorem{problem}[theorem]{Problem}

\newtheorem{remark}[theorem]{Remark}

\newtheorem{situation}[theorem]{Situation}

\newtheorem{problemintro}[]{Problem}

\theoremstyle{remark}

\numberwithin{equation}{section}

\newcommand{\NN}{\mathbb{N}}
\newcommand{\ZZ}{\mathbb{Z}}
\newcommand{\QQ}{\mathbb{Q}}
\newcommand{\RR}{\mathbb{R}}
\newcommand{\CC}{\mathbb{C}}

\newcommand{\PP}{\mathbb{P}}
\renewcommand{\AA}{\mathbb{A}}

\newcommand  {\shE}     {\mathcal{E}}
\newcommand  {\shF}     {\mathcal{F}}
\newcommand  {\shG}     {\mathcal{G}}
\newcommand  {\shH}     {\mathcal{H}}

\newcommand  {\shM}     {\mathcal{M}}

\newcommand  {\shN}     {\mathcal{N}}
\newcommand  {\shL}     {\mathcal{L}}
\newcommand  {\shR}     {\mathcal{R}}
\newcommand  {\shS}     {\mathcal{S}}
\newcommand  {\shT}     {\mathcal{T}}

\newcommand  {\shQ}     {\mathcal{Q}}

\newcommand  {\foa}     {\mathfrak{a}}

\newcommand  {\fom}     {\mathfrak{m}}

\newcommand  {\fop}     {\mathfrak{p}}

\newcommand  {\foq}     {\mathfrak{q}}

\newcommand  {\Char}    {\operatorname{char}}
\newcommand  {\cd}      {\operatorname{cd}}
\newcommand  {\cht}     {\operatorname{cht}}

\newcommand  {\Det}    {\operatorname{Det}}

\newcommand  {\dual}    {\vee}

\newcommand  {\Ext}     {\operatorname{Ext}}

\newcommand  {\Hom}     {\operatorname{Hom}}

\newcommand  {\id}      {\operatorname{id}}
\newcommand  {\im}      {\operatorname{im}}

\newcommand  {\lra}     {\longrightarrow}

\newcommand  {\mult}    {\operatorname{mult}}

\renewcommand{\O}       {\mathcal{O}}

\newcommand  {\Proj}    {\operatorname{Proj}}

\newcommand  {\ra}      {\rightarrow}

\newcommand  {\rk}    {\operatorname{rk}}

\newcommand  {\Rel}     {\operatorname{Rel}}

\newcommand  {\Spec}    {\operatorname{Spec}}

\newcommand  {\Supp}    {\operatorname{Supp}}

\newcommand {\soclo}{\star}
\newcommand {\pasoclo}{\star}

\newcommand {\mindeg}{\rho}
\newcommand {\maxdeg}{\lambda}
\newcommand {\gr} {+ \rm gr}

\usepackage{amscd}
\usepackage{amssymb}

\begin{document}

\setlength{\unitlength}{1cm}

\setlength{\parindent}{0cm}

\begin{titlepage}

\centerline{}
\bigskip

\bigskip

\LARGE
\centerline{The Theory of Tight Closure}
\centerline{from the Viewpoint of Vector Bundles}

\vspace{5cm}
\Large
\centerline{Habilitationsschrift}
\smallskip
\large

\vspace{5cm}
\centerline{}
\smallskip
\Large
\centerline{Holger Brenner}
\smallskip

\vspace{3cm}
\large
\centerline{Bochum 2003}

\end{titlepage}

\newpage

\pagenumbering{arabic}
\setcounter{page}{0}
\thispagestyle{empty}

\setcounter{tocdepth}{3}

\tableofcontents

\newpage

\thispagestyle{empty}

\markboth{Introduction}{Introduction}

\section*{Introduction}

In this habilitation thesis we link
together two mathematical subjects which are unrelated so far:
the theory of tight closure on one hand,
and the theory of vector bundles on the other.
The aim is to translate algebraic problems in
tight closure theory into problems about
vector bundles and projective bundles in order to
attack them with the help of the geometric tools then
available.

A vector bundle is a geometric object over a base space
which is locally isomorphic to a standard bundle
(like $\RR^n$, $\CC^n$, $\AA^n_k$, depending on the category).
Vector bundles, and also their relatives,
locally free sheaves, projective bundles and affine-linear bundles
play a leading role in different parts of
mathematics: in algebraic geometry,
in complex analysis, in differential geometry,
in algebraic and topological K-theory, in mathematical physics,
in gauge theory and in the theory of moduli spaces, to mention a few.
The importance and the ubiquity of bundles
in all these branches is obvious and deserves no further
explanation.

The theory of tight closure on the other hand
is a rather new and fascinating development in commutative algebra.
It has been developed since the end of the eighties
by M. Hochster and C. Huneke starting with
\cite{hochsterhuneketightclosure}, \cite{hochsterhunekebriancon}.
One motivation for its creation
was to provide a systematic foundation
for the technique called ``reduction to positive characteristic'',
which has been applied with success to
questions of homological algebra,
to invariant theory, to birational geometry
and to the characterization of singularities.
Reduction to positive characteristic means
to prove a statement first in positive characteristic
and then to deduce from this its validity in characteristic zero.
By reduction to positive characteristic $p$,
the Frobenius homomorphism $f \mapsto f^p$ becomes available
which is ultimately what gives the method its power.

We shall give a survey about the theory of tight closure in greater detail
in section \ref{tightsection}
(consult also \cite{brunstight}, \cite{hochstertightsolid},
\cite{hunekeapplications}, \cite{huneketightparameter} and
\cite{smithtightintroduction}).
In this introduction we shall only recall briefly
its definition and emphasize some important problems of this theory
which will be of interest for us.

\medskip

\medskip

\medskip
{\em Some problems in tight closure theory}

\medskip
Tight closure is, as the name suggests,
a closure operation. To every ideal $I$
in a commutative Noetherian ring containing
a field it assigns its tight closure $I^* \supseteq I$.
If the field has positive characteristic $p>0$, then
this ideal is defined with the help of the Frobenius homomorphism.
Suppose that $I=(f_1, \ldots ,f_n)$,
then the containment $f \in I^*$ is defined by the
property that
there exists an element $c \in R$, not contained in any minimal prime,
such that $cf^q \in (f_1^q, \ldots,f_n^q)$ holds for almost all
powers $q=p^{e}$.

If the field has characteristic zero,
then the notion of tight closure
is defined within a relative setting by reduction to positive characteristic.
We shall restrict our description of tight closure
to positive characteristic in this introduction.

It is the interplay between inclusion and exclusion properties of tight closure
which makes it a strong tool in commutative algebra.
The main exclusion property is that in a regular ring every ideal is
tightly closed, that is $I=I^*$.
The main inclusion properties are the contraction property,
the colon-capturing property and the Brian\c{c}on-Skoda property.
These properties together with the persistence of tight closure
are somehow the axioms of this theory.
In spite of these general and useful inclusion and exclusion results
the following is a difficult problem.

\begin{problemintro}
Let $I$ denote an ideal in a commutative Noetherian ring $R$
over a field and let
$f \in R$. Give criteria to decide whether $f \in I^*$ holds or not.
\end{problemintro}

The definition of tight closure mentioned above is not well suited
for a decision procedure, since
we have to check infinitely many conditions.
About the possibility to compute the tight closure of an ideal
or to answer the question whether $f \in I^*$ holds, C. Huneke writes
``Tight closure is very difficult to compute; indeed that is necessarily
the case. It contains a great deal of information concerning subtle
properties of the ring and the ideal''
(\cite[Basic Notions]{hunekeapplications}).

\smallskip
We have a closer look at two problems which are special cases
of problem 1.
First we consider the contraction property mentioned above.
Suppose that $R$ is a Noetherian domain over a field of positive
characteristic. Suppose that $I \subseteq R$ is an ideal
and that $R \subseteq S$ is a finite extension of domains such that
$f \in IS$ holds. Then $f \in I^*$ holds, that is
the contraction from a finite extension belongs to
the tight closure.
The set of elements $f \in R$ for which there exists such a finite extension
with $f \in IS$ form an ideal which is called the plus closure of $I$,
written $I^+$.
One of the most difficult questions of tight closure theory
is whether the inclusion $I^+ \subseteq I^*$ is, in fact,
an equality.
Hochster calls this a ``tantalizing question'' (\cite{hochstersolid}).

\begin{problemintro}
Let $R$ denote a Noetherian domain over a field of positive
characteristic.
Is it true that $I^* =I^+$?
\end{problemintro}

A positive answer to this question would solve
the localization problem of tight closure. This is
the question whether $(I^*)R_F =(IR_F)^*$ holds
for any multiplicatively closed subset
$F \subseteq R$.

The main result with respect to this problem is the
theorem of K. Smith which asserts that the equality
holds for a parameter ideal in a
locally excellent Noetherian domain
over a field of positive characteristic (\cite{smithparameter}).
A parameter ideal in a local ring
is an ideal primary to the maximal ideal
and generated by $n$ elements, where $n$ is the dimension of the ring.

Beside this deep and striking result there is not much known about
this problem. ``There is no non-trivial class of rings in which this open problem
has been solved'', K. Smith writes in \cite{smithtightintroduction}.
In particular this problem is completely open for two-dimensional normal
graded $K$-domains (for a regular ring there is no problem, since
in a regular ring every ideal is tightly closed).

\smallskip
To formulate another interesting problem we
suppose now that $R$ is a standard-graded domain
over a field $K$,
that is $R= \oplus_{d \in \NN}R_d$, $R_0 =K$
and $R$ is finitely generated
by elements of $R_1$. Suppose that
$I = (f_1, \ldots, f_n)$ is a homogeneous, $R_+$-primary ideal 
generated by homogeneous forms $f_i$.
It is quite easy to see that $I^*$ is again a homogeneous
ideal. What can we say about the homogeneous
components $(I^*)_d$, for $d \in \NN$?
Again there exists a striking answer for parameter ideals:
if $R$ is normal, Cohen-Macaulay and non-singular outside
the maximal ideal $R_+$, and if the $f_i$ are homogeneous parameters,
then the so-called Strong Vanishing Theorem, proved by Hara
in \cite{harafrobeniusrational},
gives the following numerical characterization
of the tight closure
(it holds for $p=0$ or $p \gg 0$ sufficiently large;
this last condition is typical for tight closure and
we will explain this in section \ref{tightsection} and \ref{relsituation}).
$$ (f_1, \ldots, f_n)^* =(  f_1, \ldots, f_n)
+ R_{\geq \deg (f_1) + \ldots + \deg(f_n)} \, .$$
This theorem implies the Kodaira Vanishing Theorem for ample invertible
sheaves on projective varieties
(\cite{hunekesmithkodaira}), hence the name.
Hara's result gives a strong bound on the degree of an element
to be in the tight closure.
The degree bound
$N= \sum_{i=1}^n \deg (f_i)$ has the property
that above this bound every homogeneous element belongs to
the tight closure, whereas a homogeneous element
of degree strictly smaller than this bound belongs to the tight closure only if it
belongs to the ideal itself.
Note that the Strong Vanishing Theorem needs essentially
that the ideal is generated by parameters. This naturally leads to
the following problem by dropping this condition.

\begin{problemintro}
Let $R$ denote a normal standard-graded domain
and let $(f_1, \ldots ,f_n)$ be a homogeneous $R_+$-primary ideal.
Give degree bounds for the inclusion and the exclusion for the tight closure
$(f_1, \ldots, f_n)^*$.
Under which conditions does there exist a strong bound theorem.
Can one describe the bounds in terms of the degrees of the $f_i$?
\end{problemintro}

It is clear that we cannot expect a result without any further
condition.
For example, if we take one of the ideal generators twice,
then this does not change the ideal nor its tight closure,
but it does change a numerical expression in terms of
the degrees.
The situation here is similar to the problem
concerning the plus closure.
For homogeneous parameter ideals the Strong Vanishing Theorem
gives a complete answer,
but for arbitrary primary ideals very little is known.
There exist some general (inclusion and exclusion) bounds
due to K. Smith (see \ref{smithinclusionbound} and \ref{exclusionbound} below),
but they are quite coarse for non-parameter ideals,
even in dimension two. In fact the only non-trivial computation
for a non-parameter ideal
is that the containment $xyz \in (x^2,y^2,z^2)^*$ holds for the
ring $K[x,y,z]/(x^3+y^3+z^3)$ (\cite{singhcomputation}).
Our attack on these problems is via forcing bundles,
which we shall describe now.

\medskip

\medskip

\medskip
{\em Via forcing algebras to projective bundles}

\medskip
We shall now describe briefly the main idea of the 
construction, which
assigns to ``forcing data'' for tight closure various geometric objects.
Here forcing data consist of a set
of ideal generators $I=(f_1, \ldots ,f_n) \subseteq R$ in
a Noetherian ring and another element $f_0 \in R$.
Suppose that $I$ is primary to a maximal ideal $\fom$ of
height $d $, i.e. $V(I)=V(\fom)$.

The starting point for our geometric interpretation
is the result of Hochster (see \cite{hochstersolid}
and Proposition \ref{solidcohomology} below), that for a local complete
$K$-domain $(R,\fom)$ of dimension $d$ over a field of positive
characteristic the containment
$f \in (f_1, \ldots, f_n)^*$
is equivalent to the cohomological property that $H^d_\fom(A) \neq 0$,
where
$$A = R[T_1, \ldots, T_n]/(f_1T_1+ \ldots +f_nT_n+f)$$
is the forcing algebra for the forcing data $f_1, \ldots ,f_n;f$.
The spectrum of this forcing algebra yields
over the punctured spectrum $D(\fom)= \Spec\, R- \{ \fom \}$
an affine-linear bundle
$$B= \Spec \, R[T_1, \ldots, T_n]/(f_1T_1+ \ldots +f_nT_n+f)|_{D(\fom)} \, .$$
This means that it looks locally
like an affine space over the base and that the transition mappings are
affine-linear.
Thus we assign to forcing data $f_1, \ldots ,f_n;f$
an affine-linear bundle.

The scheme
$V = \Spec \, R[T_1, \ldots, T_n]/(f_1T_1+ \ldots +f_nT_n)|_{D(\fom)}$
is a vector bundle over $D(\fom)$ of rank $n-1$.
Its sheaf of sections is the sheaf of relations
$\shR= \Rel(f_1, \ldots ,f_n)$ for the elements
$f_1, \ldots, f_n$ (Proposition \ref{relationbundle}).
These sections are also the
translations for the affine-linear bundle $B$, i.e. $V$ acts on $B$
and $B$ is a principal fiber bundle (or a torsor) with
structure group $V$ (Proposition \ref{relationaction}).

The containment $ f \in (f_1, \ldots ,f_n)^*$
is now (in positive characteristic, $d \geq 2$) equivalent with
$H^{d-1}(B, \O_B) \neq 0$, and this means
that the cohomological dimension of $B$ equals $d-1$,
which is also the cohomological dimension of the
base $D(\fom)$ (Proposition \ref{solidcd}).
The cohomological dimension of a scheme $X$ is the maximal
number $i$ such that there exists a coherent sheaf
on $X$ with $H^{i}(X, \shF) \neq 0$
(see section \ref{cohodimsection}).
This reinterpretation of tight closure
leads to the following more general problem.

\begin{problemintro}
Let $Y$ denote a scheme and let $V \ra Y$ be a vector bundle acting
on an affine-linear bundle $B \ra Y$.
What can we say about the cohomological dimension of $B$,
in particular compared with the cohomological dimension of
the base scheme $Y$?
\end{problemintro}

This is in general a quite difficult and far reaching problem.
The vector bundle $V$ itself (the trivial affine-linear bundle for $V$)
has the zero-section $Y \ra V$, hence the cohomological
dimensions of $Y$ and $V$ coincide.

In order to study the cohomological properties of such an
affine-linear bundle and to attack the last problem
it is helpful to embed the situation into
a projective setting.
The affine-linear bundle corresponds to
a $\rm\check{C}$ech-cohomology class
$c \in H^1(Y, \shS)$, where $\shS$ is the sheaf of sections in $V$
(Proposition \ref{cohotorsor}).
This class $c \in H^1(Y, \shS)= \Ext^1( \O_Y, \shS)$
defines also an extension $0 \ra \shS \ra \shS' \ra \O_Y \ra 0$
of locally free sheaves.
The corresponding inclusion of vector bundles $V \hookrightarrow V'$
yields an inclusion of projective bundles
$\PP(V) \hookrightarrow \PP(V')$. Then we find that
$B \cong \PP(V') -\PP(V)$, i.e. $B$ is the complement
of a projective subbundle of codimension one inside a
projective bundle (Proposition \ref{torsorprojective}).
In particular $B$ is the complement of an irreducible divisor,
which we call the forcing divisor.

\smallskip
The cohomology class $c$
corresponding to the affine-linear bundle given by forcing data
$f_1, \ldots ,f_n;f_0 \in R$ in the forcing situation has another description:
if we consider the short exact sequence
$$0 \lra \shR=Rel (f_1, \ldots ,f_n) \lra
\O_U^n \stackrel{f_1, \ldots ,f_n}{\lra} \O_U \ra 0 $$
of locally
free sheaves on $U=D(\fom)$, then
$c= \delta (f_0) \in H^1(D(\fom), \shR)$,
where $\delta$ is the connecting homomorphism
(Proposition \ref{forcingbundletorsortorsor}).
The extension defined by this class has also an easy description,
namely
$V'= \Spec \, R[T_0,T_1, \ldots, T_n]/(f_1T_1+ \ldots +f_nT_n+f_0T_0)|_{D(\fom)}$,
and the inclusion of vector bundles $V \subset V'$ is given by
$T_0 =0$ (Proposition \ref{forcingsequence1}).

\medskip

\medskip

\medskip
{\em The graded situation}

\medskip
The realization of an affine-linear bundle in terms
of projective bundles and subbundles
is in particular useful in the graded situation.
Let $R$ be a normal standard-graded $K$-domain, let $f_1, \ldots, f_n$
denote homogeneous, $R_+$-primary elements and let $f_0$ denote another
homogeneous element. Then we can do the same constructions
to get bundles on the projective variety $\Proj \,R$
(sections \ref{gradedsection} and \ref{gradedforcingsection}).
In particular we assign to homogeneous forcing data
$f_1, \ldots ,f_n;f_0$ a projective bundle over $\Proj \, R$
together with a projective subbundle
(caller the forcing subbundle or the forcing divisor) in such a way that
the cohomological properties of its complement 
characterize $f_0 \in (f_1, \ldots ,f_n)^*$ (Proposition
\ref{solidgraded}).
Then we are in an entirely projective setting
and we can use all the methods and tools of
projective algebraic geometry to attack the underlying tight closure problem.
Moreover, this projective situation
is smooth under the condition that the graded ring
$R$ has an isolated singularity.

\smallskip
We may also interpret the containment in the plus closure in terms
of the geometric situation.
If $R$ is a normal standard-graded $K$-domain, then the question whether
$f_0 \in (f_1, \ldots ,f_n)^{\gr}$ holds
(this is a graded version of the
plus closure)
is in our setting equivalent to the existence
of a projective subvariety inside
the corresponding affine-linear bundle $B$ on $Y=\Proj \, R$ of
maximal dimension $\dim \, Y$ (Proposition \ref{plus}).
Such a projective subvariety represents a finite solution for the tight
closure problem.
To paraphrase, if $B \cong \PP(V') - \PP(V)$, then
$f_0 \not\in (f_1, \ldots ,f_n)^{+{\rm gr}}$
is equivalent to the property that
every projective subvariety in $\PP(V')$ of dimension $\dim\, Y$
intersects the forcing divisor $\PP(V)$ positively.

Therefore the problem whether the tight
closure and the plus closure of an ideal
is the same has the following natural generalization
(both problems are relevant only in positive characteristic).

\begin{problemintro}
Let $Y$ denote a projective variety over a field of positive characteristic
and let $B \ra $Y denote an affine-linear bundle.
Suppose that the cohomological dimension of $B$ equals
the dimension of the base $Y$.
Does there exist a projective subvariety of dimension $\dim \, Y$
inside $B$?
\end{problemintro}

If $R$ is a two-dimensional normal domain,
then we are in a particularly manageable situation
(section \ref{dimensiontwosection}).
The containment $f \in I^*$ is then equivalent with the property
that the affine-linear bundle $B$ is not an affine scheme.
Recall  that a
scheme $X$ is called affine if it is isomorphic to the spectrum of a
commutative ring. If $X$ is of finite type over a field $K$, this means that
we can embed $X$ as a closed subscheme
into an affine space $\AA^N_K$.

If moreover $R$ is standard-graded,
then $\Proj\, R$ is a smooth projective curve and
the construction yields projective bundles and subbundles
over this curve, and again
the containment $f_0 \in (f_1, \ldots ,f_n)^*$ ($f_i$ homogeneous)
is equivalent with the non-affineness of the complement
of the forcing divisors $\PP(V) \subset \PP(V')$
(Theorem \ref{tightaffinecrit}).

Since a lot of problems of tight closure are still unsolved
in dimension two, it seems to be justified to concentrate
on the two-dimensional situation
(the corresponding projective situation is one-dimensional).
Therefore we shall deal to a great extent with the following special
case of problem 4.

\begin{problemintro}
Let $Y$ denote a smooth projective curve and let
$\PP(V')$ denote a projective bundle together
with a projective subbundle $\PP(V) \subset \PP(V')$
of codimension one.
When is the complement $\PP(V')- \PP(V)$ an affine scheme?
\end{problemintro}

Such a situation is equivalent to give a locally free sheaf $\shS$
(the sheaf of sections of $V$)
and a cohomology class $c \in H^1(Y, \shS)$.
In the forcing situation this class is given by
$\delta(f_0) \in H^1(Y, \shR(m))$, where $\shS=\shR(m)$
is the sheaf of relations for homogeneous elements $f_1, \ldots, f_n$
of total degree $m$ and where $f_0$ has degree $m$.

The containment to the graded plus closure is equivalent
to the question whe\-ther the bundle $B$ does contain projective curves
or not, or, equivalently, whether the forcing divisor intersects every
curve positively (Theorem \ref{pluscrit}).

\medskip

\medskip

\medskip
{\em Using algebraic geometry}

\medskip
What is the advantage of this geometric interpretation?
First of all it allows to use
the whole methods and tools of algebraic geometry 
and apply them to tight closure problems.
We will use among other things the following
items of algebraic geometry:
ample divisors and ample vector bundles,
big divisors, the notion of the slope of a vector bundle
and semistability, intersection theory on projective varieties,
in particular the top self intersection number of a hypersurface,
group schemes and torsors,
classification of vector bundles.
In fact rather basic properties of these topics suffice
to obtain new and interesting results for tight closure.
Algebraic geometry provides also the language to express
the complexity and the geometrical richness
which lie behind a tight closure problem.

Our geometric reformulation of tight closure theory
gives at once an important
numerical condition and a candidate for the bound degree
for which we asked in problem 3.
Heuristically, the affineness of the complement of a divisor
is a richness property (like being big or being ample)
and the behavior of the sections of the multiples of a divisor
is given (up to correcting higher cohomology terms)
by the top self intersection number of the divisor.
Hence a coarse numerical orientation for the tight closure
problem for a two-dimensional graded domain
is given by the top self intersection number
of the forcing divisor. We will see (Theorem \ref{topselfintersection}) that
$$\PP(V)^n =( d_1+ \ldots + d_n- d_0(n-1)) \deg (Y) \, ,$$
where $d_i= \deg (f_i)$ and $ \deg (Y)$ is the degree of
the ample invertible sheaf $\O_Y(1)$ on $Y$.
This number is positive if and only if
$\deg (f_0) < \frac{d_1+ \ldots +d_n}{n-1}$ holds,
a condition which will play a crucial role in the following.

The ampleness of the forcing divisor $\PP(V) \subset \PP(V')$
implies that its complement is affine, and this property implies
that the forcing divisor is big,
that is some multiple of it defines a rational mapping to a projective
space such that the image has maximal dimension.
The bigness property is a useful necessary condition for affineness
(Theorem \ref{slopemaxkrit} and our main inclusion criterion for tight closure
Theorem \ref{maxin} are based on it).
Bigness implies also that there exists a multiple of $\PP(V)$
which is linearly equivalent $a \PP(V) \sim H+F$, where
$H$ is another projective subbundle and where $F$ consists
of fiber components (Proposition \ref{bigaffine}).
Under this condition a projective curve which does not meet the forcing
divisor must lie on $H$. This gives a strong restriction for
the existence of such curves
and therefore for the existence of a finite solution
for a tight closure problem
(Propositions \ref{effectivevertreter} and \ref{nupositive}).

Regarding a tight closure problem $f_0 \in (f_1, \ldots, f_n)^*$
as a property of the corresponding cohomology class
$c=\delta(f_0) \in H^1(Y, \shR(m))$
provides the following feature,
which we will use several times:
to argue ``along a short exact sequence''. If
$c \in H^1(Y, \shS)$ and 
$0 \ra \shT \ra \shS \ra \shQ \ra 0$ is exact, then we may consider
the image $c'$ of the class in $H^1(Y,\shQ)$.
If the $\shQ$-bundle defined by $c'$
is affine, then so is the $\shS$-bundle defined by $c$ (\ref{monotonielemma}).
If $c'=0$, then $c= \varphi(e)$, and we may argue on $\shT$.
With this principle we may divide problems of tight closure
into problems of lower rank.

\medskip

\medskip

\medskip
{\em The parameter case in dimension two}

\medskip
The first case to test our geometric approach
to tight closure is for $n=2$ (in dimension two),
the parameter case.
The construction yields in this case a projective bundle of rank one,
that is a ruled surface over the
corresponding smooth projective curve together with a forcing section
(Corollary \ref{forcingrule}).
For this situation we cannot expect new results for tight
closure due to the simple reason that
``everything'' is known about the tight closure
of a parameter ideal.
However, our geometric interpretation yields new geometric proves
for known facts.

The numerical condition mentioned above
becomes $ \deg (f_0) < \deg (f_1) + \deg (f_2)$,
which we know from the Strong Vanishing Theorem.
The affineness of the complement of the forcing section
is equivalent to its ampleness and we recover
Strong Vanishing Theorem (Theorem \ref{ampleaffineruled} and Corollary
\ref{degreecrit1}).
For positive characteristic we recover also that
$(f_1,f_2)^* = (f_1,f_2)^{\gr }$ using the Artin-Schreier
sequence and \'{e}tale cohomology (Proposition \ref{frobeniusartinschreier}).

Already here we see
the geometrical richness behind a tight closure problem
(sections \ref{examplesparameter} and \ref{complex}), for example
we encounter Hirzebruch surfaces, incidence va\-rie\-ties, the
graph of a meromorphic function.
Moreover, tight closure provides
also analytically interesting examples of complements of sections
on a ruled surface.
One of the most prominent examples of tight closure,
that $z^2 \in  (x,y)^*$ holds in $K[x,y,z]/(x^3+y^3+z^3)$,
yields for $K= \CC$
a classical construction of Serre of a Stein but non-affine variety
(Corollary \ref{steinnonaffine}).

Moreover the tight closure of a parameter ideal provides
a new class of counter-examples to the hypersection problem
of complex analysis
(Proposition \ref{hypersectioncounterexample}):
suppose that $D \subseteq X$ is a hypersurface in a complex Stein space
$X$ of dimension $\geq 3$ and suppose that
for every analytic surface $S \subset X$ the intersection
$S \cap (X-D)$ is Stein. Is then $X-D$ Stein?
The first counter-example to this problem was given
by Coltoiu and Diederich (\cite{coltoiudiederich}).
The first instance of our class of counter-examples
(example \ref{hypersectionexample})
is given by
$$X=(\Spec \CC[x,y,z,T_1,T_2]/(x^4+y^4+z^4,xT_1+yT_2+z^3))^{\rm an}$$
together with the hypersurface $D=V(x,y) \subset X$.

\medskip

\medskip

\medskip
{\em Slope criteria for affineness}

\medskip
What can we say about the affineness of an affine-linear
bundle $B \cong \PP(V') - \PP(V)$
in general (problem 6 for arbitrary rank), where the situation
corresponds to a short exact sequence
$0 \ra \shS \ra \shS' \ra \O_Y \ra 0$ of locally free
sheaves on the smooth projective curve $Y$ over an
algebraically closed field
and to a cohomology class $c \in H^1(Y, \shS)$.

The ampleness of the divisor $\PP(V) \subset \PP(V')$
is by definition the ampleness of the locally
free sheaf $(\shS')^\dual$ (the dual sheaf of $\shS'$).
This property implies the affineness of the complement.
The ampleness of $\shS^\dual$ is equivalent to the
ampleness of the normal sheaf on $\PP(V)$. This property
does not depend on the cohomology class.
If $\shS^\dual$ is ample and $c \neq 0$, then also
$(\shS')^\dual$ is ample in characteristic zero
(and then also in characteristic $p \gg 0$)
and the complement is affine.
Therefore it is natural to look first for criteria for
ample vector bundles on a curve.

This is a well-studied notion and we are in the lucky position 
to use several classical results or to extend them for our needs. 
The numerical conditions for ampleness work best
with the notion of the slope of a bundle.

Recall that the slope of a locally free sheaf $\shE$ on a smooth
projective curve $Y$ is defined by $\mu (\shE) =\deg (\shE) /\rk (\shE)$.
We will use also the following variants
(see section \ref{sectionslope}): the minimal slope
$\mu_{\rm min} (\shE)$, which is defined as the minimal slope
of every quotient sheaf of $\shE$,
and the maximal slope, $\mu_{\rm max} (\shE)$, which is defined
as the maximal slope of every subsheaf of $\shE$.
A sheaf is called semistable if
the minimal and the maximal slope are the same.
In positive characteristic we will define also
$\bar{\mu}_{\rm min} (\shE)$ and $\bar{\mu}_{\rm max} (\shE)$,
which take into account also the behavior of subbundles and quotient bundles
under a Frobenius morphism.

The main numerical ampleness criterion, which is essentially due to Barton
(\cite[Theorem 2.1]{barton}), is that $\shE$ is ample
if and only if
$\bar{\mu}_{\min} (\shE) >0$ (Theorem \ref{amplekrit}).
From this we get for our question easily (in characteristic $0$ or $p \gg 0$)
that $\mu_{\rm max} (\shS) < 0$ and $c \neq 0$
imply that $\PP(V')- \PP(V)$ is affine (Theorem \ref{ampleaffine}).
With some more effort we get the following sufficient criterion:
suppose that $\varphi: \shS \ra \shT$
is a sheaf homomorphism such that $\shT$ is semistable
of negative slope $\mu (\shT) < 0$
and such that $\varphi(c) \neq 0$ holds in $H^1(Y, \shT)$,
then $\PP(V')- \PP(V)$ is affine (Theorem \ref{slopekritaffin}).
Important candidates for $\shT$ are invertible sheaves of negative
degree and the semistable quotient of minimal slope
in the Harder-Narasimhan filtration of $\shS$
(Corollaries \ref{semistablequotientaffine} and \ref{quotientnegdegree}).

Our main criterion
for $\PP(V') - \PP(V)$ to be not affine is
the condition that
$\bar{\mu}_{\rm min} (\shS) \geq 0$ or equivalently that
$\bar{\mu}_{\rm max} (\shG) \leq 0$ (Theorem \ref{slopemaxkrit}),
where $\shG= \shS^\dual$.
From this condition we deduce (Lemma \ref{muminpseudo})
that $\mu_{\rm max} (S^k \shG) \leq 0$
and therefore the forcing divisor is not big, hence its complement
is not affine.

If $\shS$ is semistable,
then these results together yield the numerical condition
that $\PP(V') -\PP(V)$ is affine if and only if
$\deg (\shS) <0 $ and $c \neq 0$
(Corollaries \ref{semistablepositiveaffine} and \ref{semistablenotaffine},
characteristic $0$ or $p \gg 0$).
If $\shS$ has a decomposition
$\shS= \shS_1 \oplus \ldots \oplus \shS_s$ with semistable summands
$\shS_j$, then we also get a complete numerical description:
$\PP(V') -\PP(V)$ is affine if and only if there exists $j$ such that
$\deg (\shS_j) <0 $ and $c_j \neq 0$, where $c_j$ is the component
of $c$ in $H^1(Y, \shS_j)$ (Theorem \ref{splittingsemistable}).
If the base curve is a projective line, then every locally free sheaf splits
into invertible sheaves and therefore we get the equivalence
that $\PP(V') - \PP(V)$ is affine if and only if $c \neq 0$
(Corollary \ref{projectivelineaffinetorsor}).
This corresponds to the property that every ideal in a regular ring
is tightly closed.

In general we have to work along the Harder-Narasimhan filtration of
$\shS$,
$$ 0 =\shS_0 \subset \shS_1 \subset \ldots \subset \shS_s =\shS \, ,$$
where the quotients $\shS_i /\shS_{i-1}$ are semistable.
Going through the short exact sequences corresponding to this filtration
we get the first steps of an algorithm to decide whether
$\PP(V') - \PP(V) $ is affine or not (section \ref{algorithmsection}).
If the Harder-Narasimhan filtration consists only of two terms
($s=2$),
that is if we have a short exact sequence
$0 \ra \shS_1 \ra \shS \ra \shS/ \shS_1 \ra 0$
with $\shS_1$ and $\shS/ \shS_1$ semistable, then we get a complete answer.
This applies in particular to the case where the rank of $\shS$ is two,
which corresponds to the tight closure of three elements.

\newpage

\medskip
{\em From slope criteria to inclusion and exclusion bounds}

\medskip
From these criteria for affineness and non-affineness we get now at once
criteria for the inclusion and exclusion for tight closure.
Suppose that $f_1, \ldots ,f_n \in R$ are
homogeneous $R_+$-primary elements in a two-dimensional normal
standard-graded domain over an algebraically closed field
$K$, and let $\shR(m)$ denote the sheaf of relations
of total degree $m$ on $Y= \Proj \, R$.
Another homogeneous element $f_0$ of degree $m$ yields
the cohomology class $c \in H^1(Y, \shR(m))$ and the
embedding $\PP(V) \subset \PP(V')$ as before.
We set $\mu_{\rm max} (f_1, \ldots ,f_n) :=
\mu_{\rm max} (\shR^\dual(0))$.
Note that its slope is
$\mu(\shR^\dual (0))= \deg \O_Y(1) (\sum_{i=1}^n \deg (f_i))/(n-1)$.

We obtain then the following degree bound for the inclusion
(Theorem \ref{maxin}):
if
$\deg \, (f_0) \geq $ $ \bar{\mu}_{\max} (f_1, \ldots ,f_n)/\deg \, (\O_Y(1))$,
then $f_0 \in (f_1, \ldots ,f_n)^*$.
It is worth noting that this kind of inclusion to tight closure
rests upon the non-bigness of the forcing divisor
and is not covered by the three standard containments
contraction, colon capturing and Brian\c{c}on-Skoda.
With additional conditions (see below) we get better estimates
for $\bar{\mu}_{\max} (f_1, \ldots ,f_n)$ and therefore better inclusion
bounds.
Without any further condition we may derive that
$R_m \subseteq (f_1, \ldots ,f_n)^*$, whenever
$m$ is greater or equal the sum of the two biggest degrees of the $f_i$
(Corollary \ref{inclusionbound}),
which improves slightly the bound $2 \max d_i$ of Smith given
in \cite[Proposition 3.1]{smithgraded}.

In the other direction we prove
(if the characteristic of the field is $0$ or $p  \gg  0$) that if
$\deg \, (f_0) < \mu_{\rm min} (f_1, \ldots ,f_n)/\deg \, (\O_Y(1))$ holds,
then $f_0 \in (f_1, \ldots, f_n)^*$
is only possible if already
$f_0 \in (f_1, \ldots ,f_n)$ holds (Theorem \ref{minex}).

If the sheaf of relations is semistable, then
the two slope bounds coincide
and we get the strong bound theorem (our main contribution to problem 3)
$$(f_1, \ldots ,f_n)^*=(f_1, \ldots ,f_n) + R_{\geq k} \, \, \, \, ,
k= \lceil \frac{\deg (f_1) + \ldots + \deg (f_n)}{n-1} \rceil $$
(Theorems \ref{semistablevanishing} and \ref{semistablevanishingp}).
Here we meet again the bound expected from the heuristic reasoning
about the top self intersection of the forcing divisor.
In particular, if the sheaf of sections is semistable,
there exists an easy numerical criterion for tight closure.
A complete numerical description of the tight closure
is also available if the sheaf of relations
is the direct sum of semistable sheaves.

The computation of the tight closure of three $R_+$-primary
homogeneous elements $f_1,f_2,f_3$
is already a very subtile problem
with a multitude of new phenomena, which become visible through
the geometric interpretation, see the examples in section \ref{correlations}.
However, due to the algorithm mentioned above
it is always possible to compute $(f_1,f_2,f_3)^*$,
at least if we are able to find the Harder-Narasimhan filtration
of the relation bundle (of rank two) for $f_1,f_2,f_3$.
If the sheaf of relations is semistable, then
the numerical criterion gives the answer. Otherwise
we have a short exact sequence $0 \ra \shL \ra \shR(m) \ra \shM \ra 0$
such that $\deg (\shL) > \deg (\shM)$ and we may decide
whether $f_0 \in (f_1,f_2,f_3)^*$ holds by looking
at the behavior of the class $c =\delta(f_0) \in H^1(Y, \shR(m))$
under the cohomology sequence corresponding to this short exact sequence
(Proposition \ref{exactsequencecrit} and Corollary \ref{exactsequencedecide}).
A similar reasoning yields also (but less complete) results on the plus closure
in positive characteristic.

Under the condition that the sheaf of relations for $f_1,f_2,f_3$
is indecomposable
we obtain degree bounds for inclusion and exclusion
which are quite near to the expected
number $k=(d_1+d_2+d_3)/2$ ($d_i = \deg (f_i)$):
the difference is at most
$(g-1)/\delta $, where $g$ is the genus of the curve
and $\delta= \deg (\O_Y(1))$ is its degree
(Corollary \ref{genusboundtight}).
In particular, if $R=K[x,y,z]/(F)$, where $\deg (F)= \delta$, and if
the sheaf of relations for $f_1,f_2,f_3 \in R$
is indecomposable, then
$R_m \subseteq (f_1,f_2,f_3)^*$ holds for $m \geq \frac{d_1+d_2+d_3}{2}
+ \frac{\delta -3}{2}$,
and $R_m \cap (f_1,f_2,f_3)^*= R_m \cap (f_1,f_2,f_3)$
for $m < \frac{d_1+d_2+d_3}{2} - \frac{\delta -3}{2}$
(Corollary \ref{degreebound} in characteristic zero, the results in positive
characteristic are slightly worse).
In positive characteristic we get under the same assumptions
the inclusion
$R_m \subset (f_1,f_2,f_3)^{\gr}$ for
$m \geq \frac{d_1+d_2+d_3}{2} + {\delta -3}$
(Proposition \ref{indecomposableplus}).

The existence of global relations of certain degree
has also consequences on the structure of $\shR(m)$
and hence on the corresponding tight closure problem
(section \ref{relations}).
This includes criteria for semistability, indecomposability
and ampleness in terms of the existence of global relations.
If we know that there does not exist
a global relation of degree
$k \leq \frac{d_1+d_2+d_3}{2} + \frac{g-1}{\delta}$,
then $R_m \subset (f_1,f_2,f_3)^*$ holds
for $m \geq   d_1+d_2+d_3-k + \frac{g-1}{\delta }$
(Corollary \ref{relationboundinclusion}).
On the other hand, if there does not exist a global relation of
degree $k$, then
$(I^*)_m =I_m$ holds for $m \leq k- \frac{2g}{\delta }$
(Proposition \ref{amplecritrelation}).
Other results require the condition that certain relations are
primary, which means that they define a subbundle
of $\shR(m)$.

\medskip

\medskip

\medskip
{\em Tight closure and plus closure over elliptic curves}

\medskip
Suppose now that the base curve
is an elliptic curve $Y$, that is a curve with genus $g(Y)=1$.
Then we are in a very favorable situation:
the vector bundles on an elliptic curve have been completely
classified by Atiyah (\cite{atiyahelliptic}).
This classification has led
to subsequent results on vector bundles over an elliptic curve
such as the numerical ampleness criterion of Hartshorne-Gieseker
(Theorem \ref{amplekritelliptic}),
the property that every indecomposable sheaf is semistable
(Proposition \ref{ellipticproperties})
and results about the behavior of cohomology classes under finite
morphisms like the theorem of Oda (Theorem \ref{odageneral}).
We extend this last theorem and give a converse to it in positive characteristic
which says the following (Lemma \ref{pvanish}):
if $c \in H^1(Y, \shS)$ is a cohomology class,
where $\shS$ is an indecomposable sheaf of degree
$\deg (\shS) \geq 0$,
then a sufficiently high power of the $p$-multiplication
(in the group structure of the abelian curve $Y$)
annihilates $c$.

From these results we deduce the following numerical criterion
for affineness for $\PP(V')- \PP(V)$ given by
$c \in H^1(Y, \shS)$: if
$\shS= \shS_1 \oplus \ldots \oplus \shS_s$
is the decomposition of $\shS$ into indecomposable
locally free sheaves, then $\PP(V')- \PP(V)$
is affine if and only if
there exists $j$ such that $\deg (\shS_j) < 0$ and
$c_j \neq 0$ (Theorem \ref{numkritaffinelliptic}).
The same numerical characterization holds in positive characteristic also for
the non-existence of projective curves
inside $\PP(V') - \PP(V)$ (Theorem \ref{numkritcurve}).
Therefore we get for an elliptic curve ($p >0$)
a positive answer to problem 5:
if $B \cong \PP(V')- \PP(V)$ is an affine-linear bundle,
then $B$ is an affine scheme if and only if it does not
contain any projective curve (Corollary \ref{geokritaffin}).

These results imply at once the following result concerning
problem 2: if $R=K[x,y,z]/(F)$ is the homogeneous coordinate ring
over an elliptic curve ($\deg (F)=3$, $R$ normal)
over a field of positive characteristic,
then for every $R_+$-primary homogeneous ideal $I$ we have the identity
$I^* = I^{\gr}$ (Theorem \ref{tightpluselliptic}).

\newpage

\thispagestyle{empty}

\setcounter{section}{0}

\markboth{1. Foundations}{1. Foundations}

\section{Foundations}

\bigskip
\subsection{A survey about the theory of tight closure}

\label{tightsection}
\markright{A survey about the theory of tight closure}
\

\bigskip
In this section we present an overview about the theory of tight closure.
This shall help the reader who is not acquainted with this
theory to understand the motivation for this habilitation thesis.
Other introducing literature for this subject are
\cite{brunstight}, \cite{hochstertightsolid},
\cite{hunekeapplications}, \cite{huneketightparameter} and
\cite{smithtightintroduction}, to which we also refer for proves.

The theory of tight closure has been developed
since the end of the eighties by Melvin Hochster and Craig Huneke.
It is one of the most interesting and fascinating
discoveries in commutative algebra over the last
15 years. Its fascinating and somehow mysterious flavor
is due to the fact that this theory is defined
at first stage only for commutative Noetherian rings
containing a field of positive characteristic with the help
of the Frobenius-Endomor\-phism. In a second stage it is then
also defined for rings containing a field of characteristic zero
by reduction to positive characteristic.

In fact it was an important motivation for developing
this theory to give a systematic foundation
for proves which work by ``reduction to positive characteristic'',
i. e. to establish a proposition
first in positive characteristic and to deduce from this fact that
it must hold in zero characteristic as well.

\medskip

\medskip

\medskip
{\em Definition of tight closure in positive characteristic
and basic properties}

\medskip
The theory of tight closure assigns to every ideal $I$ in a
Noetherian commutative ring $R$ over a field $K$ another
ideal $I^*$ which is called the tight closure of $I$.
We restrict first to the case where the
ground field $K$ has positive characteristic $p >0$.
Then the mapping $F:R \ra R$, $a \mapsto a^p$, is a ring-endomorphism,
the so-called {\em Frobenius}. Working in positive characteristic 
means studying and using the Frobenius morphism.

We denote the extended ideal under a Frobenius power
$F^{e}: R \ra R$, $a \mapsto a^q$, $q =p^{e}$, by $I^{[q]}$.
If the ideal $I$ is given by generators
$I=(f_1, \ldots ,f_n)$, then $I^{[q]}=(f_1^q, \ldots ,f_n^q)$.
With these Frobenius powers we define 
the tight closure of an ideal by the following condition.
$$ I^* = \{f \in R :\, \exists \, c \in R, \,
c \not\in \mbox{ minimal prime},\,
cf^q \in I^{[q]} \mbox{ for almost all }  q=p^{e}   \} \, .$$
If $R$ is reduced and $f \in I^*$, then we may skip
the word ``almost'' in the definition, i.e. we find then
also a $c$ (not in any minimal prime of $R$)
such that $cf^q \in I^{[q]}$ holds for all $q=p^{e}$.
It turns out that $I^*$ is again an ideal, which contains $I$ and
which has the property that $(I^*)^* = I^*$.
For $I \subseteq J$ we get the inclusion $I^* \subseteq J^*$.
Hence the name ``closure'' is justified.
Ideals which coincide with their tight closure are called
{\em tightly closed}.

We look first at such properties of this closure operation
which justify the name ``tight''.
An important
\--- and maybe at first glance disappointing \---
property of tight closure is the following theorem
($p$ is temporarily positive).

\begin{theorem}
\label{tightregular}
Let $R$ denote a regular ring over a field of characteristic
$p$. Then every ideal in $R$ is tightly closed.
\end{theorem}
\proof
See \cite[Theorem 1.3(e) and Theorem 1.5(3)]{hunekeapplications}.
\qed

\medskip
This proposition is due to the fact that for
regular rings the Frobenius morphism is flat.
Rings with the property that every ideal is tightly closed
are called {\em weakly F-regular} (F is for Frobenius)
and {\em F-regular}, if this property holds for every
localization (it is not known whether this is the same).
Hence theorem \ref{tightregular} tells us that regular rings are F-regular.
Theorem \ref{tightregular} gives rise to the question
whether there exist other F-regular rings and 
what kind of properties may be deduced from this property.
It also indicates that there should be a relation between
tight closure and singularity theory. We will discuss this
relationship below.

Theorem \ref{tightregular} has also consequences for non-regular rings
in connection with the following property which is called
the {\em persistence} of tight closure.

\begin{theorem}
\label{tightpersistence}
Let $R$ and $S$ denote Noetherian commutative rings over a
field of characteristic $p$
and let $\varphi: R \ra S$ denote a ring homomorphism.
Suppose that $R$ is essentially of finite type over
an excellent local ring or that $R_{\rm red}$ is F-finite.
Then $I^*S \subseteq (IS)^*$.
\end{theorem}
\proof
See \cite[Theorem 2.3]{hunekeapplications}.
\qed

\medskip
The proof of this statement is rather intricate (under the first condition)
and deserves the theory of test ideals, which we will see below.
The problem here is due to the fact that the element $c \in R$,
which we need in the definition to show
that an element belongs to the tight closure, may be mapped
under $\varphi: R \ra S$ to zero,
and then we cannot use the image of $c$ for this purpose anymore.

The property {\em F-finite} means that the Frobenius morphism is a
finite mapping. This is true for
$K$-algebras which are (essentially) of finite type
over a perfect field $K$.
The technical conditions in the theorem are also fulfilled
for complete local rings.

\smallskip
The relation of tight closure $I^*$ to the integral closure
$\bar{I}$ is given by the following corollary.
One way to characterize the containment $ f \in \bar{I}$
\--- among a lot of other equivalent possibilities \---
is by the property that for every discrete valuation domain
$\varphi : R \ra B$ we have $\varphi(f) \in IB$.
Let $\sqrt{I}$ denote the radical of $I$,
i.e. the intersection of all prime ideals which contain $I$.

\begin{corollary}
\label{tightintegral}
Let $R$ denote a Noetherian ring over a field
of characteristic $p$ and let $I \subseteq R$ denote an ideal.
Then
$$ I^* \subseteq \bar{I} \subseteq \sqrt{I} \, .$$
\end{corollary}
\proof
Since discrete valuation domains are regular, the first
inclusion follows from the persistence of tight closure
and from the F-regularity of regular rings.
\qed

\medskip
The statement \ref{tightintegral} tells us that the tight
closure of an ideal is much more closer to the
ideal than the integral closure,
in particular if we recall that in a regular ring
an ideal is not integrally closed in general.

\medskip
The strength of tight closure lies in the fact that it is
very tight as already shown but also wide enough to contain
interesting information.
So now we look for elements which belong
to the tight closure without belonging to the ideal itself.
Recall that elements $x_1, \ldots,x_d$ in a local Noetherian ring
$(R, \fom)$ of dimension $d$ are called {\em parameter}
if $V(x_1, \ldots, x_d) = V(\fom)$ holds.
An ideal which is generated by a system of parameters
is called a {\em parameter ideal}.
A local Noetherian ring is called Cohen-Macaulay
if one (and then every) system of parameters form a regular sequence,
i.e.  $x_k$ is not a zero divisor in
$R/(x_1, \ldots, x_{k-1})$ for $k=1, \ldots, d$.
This is equivalent with the property that
$$    ((x_1, \ldots,x_{k-1}) : x_k) \subseteq (x_1, \ldots,x_{k-1}) $$
holds in $R$ for $k=1, \ldots,d$.
The point is that
this statement holds without further conditions
if we replace the ideal on the right by its tight closure.

\begin{theorem}
\label{coloncapturing}
Let $x_1, \ldots, x_d$ denote parameters in a $d$-dimensional local
Noetherian ring $R$ over a field of characteristic $p$.
Suppose that $R$ is the homomorphic image of
a local Cohen-Macaulay ring.
Then we have the inclusion for $k=1, \ldots, d$
$$    ((x_1, \ldots,x_{k-1}) : x_k) \subseteq (x_1, \ldots,x_{k-1})^* $$
\end{theorem}
\proof
See \cite[Theorem 3.1]{hunekeapplications}.
\qed

\medskip
This property of tight closure is called ``colon-capturing''.
The condition that the ring is the image of a
Cohen-Macaulay ring is of technical nature and holds
for complete local rings and for rings which are essentially of finite type
over a field.
Hence tight closure provides a measure for the deviation
from the Cohen-Macaulay property.
Formulated in another way, we get the following.

\begin{corollary}
\label{regularcohenmacaulay}
A local weakly F-regular ring
which is the homomorphic image of a local Cohen-Macaulay ring
is Cohen-Macaulay.
\end{corollary}
\proof
This follows directly from \ref{coloncapturing}
\qed

\medskip
The properties of tight closure developed so far allow
us to prove interesting statements in which
tight closure is only a tool for proving,
but not a part of the statement itself.

\begin{theorem}
\label{directsummand}
A Noetherian ring $S$ which is a
direct summand of a local regular ring $R$ over a 
field of characteristic $p$ is {\rm (}F-regular and in particular{\rm)}
Cohen-Macaulay.
\end{theorem}
\proof
Let $S \hookrightarrow R$ be a direct summand and $I \subseteq S$
an ideal.
From $f \in I^* $ it follows due to the persistence \ref{tightpersistence}
(in fact we need here only a trivial case of persistence)
that $f \in (IR)^*$ and therefore $f \in IR$ due to
the (F-)regularity of $R$ (\ref{tightregular}).
Since $S$ is a direct summand it follows that $f \in I$.
(Since we may pass to the completion, the
technical conditions in the cited statements are granted.)
\qed

\medskip
An important special case of this statement is the following theorem of
Hochster and Roberts.
Its clear proof using the theory of tight closure is
already a justification for the effort!
The original proof in \cite{hochsterrobertsinvariant}
was a rather complicated reduction to the case of
positive characteristic.

\begin{theorem}
\label{tightinvariant}
Let $G$ denote a linearly reductive group, which acts on a polynomial
ring $K[X_1, \ldots, X_n]$ over a field $K$ of
characteristic $p$ as a group of linear automorphisms.
Then the ring of invariants $K[X_1, \ldots, X_n]^G$ is Cohen-Macaulay.
\end{theorem}
\proof
The Reynolds operator shows that the ring of invariants
is a direct summand in the polynomial ring.
Furthermore it is a finitely generated
graded $K$-algebra. Hence this follows from \ref{directsummand}.
\qed

\medskip
Note that the statement of \ref{tightinvariant} holds so far
only in positive characteristic, and therefore it applies
to rather few groups, namely to finite groups,
if their order is not divisible by the characteristic,
and to tori (the operations of tori correspond to multigraduations
on the polynomial ring).
When we will have developed the theory of tight closure also in
zero characteristic together with the basic properties
\ref{tightregular} - \ref{coloncapturing}, then
theorem \ref{tightinvariant} is also established in zero
characteristic and applies
also to the classical groups, the general linear group,
the special linear group, the orthogonal group and the
symplectic group.

\smallskip
The theorem of Brian\c{c}on-Skoda is another theorem
where tight closure provides both a simple proof and a more
general statement.

\begin{theorem}
\label{tightbrianconskoda}
Let $I \subseteq R$ denote an ideal in a Noetherian commutative
ring over a field of characteristic $p$.
Suppose that $I$ is generated by $k$ elements.
Then $\overline{I^{k+m}} \subseteq (I^{m+1})^*$.
In particular $\overline{I^{k}} \subseteq I^*$.
\end{theorem}
\proof
See \cite[Theorem 5.7]{hunekeapplications}.
\qed

\medskip
From this follows the theorem of Brian\c{c}on-Skoda
for regular rings.

\begin{theorem}
\label{brianconskoda}
Let $I \subseteq R$ denote an ideal in a regular
ring over a field of characteristic $p$, which is generated
by $k$ elements.
Then $\overline{I^{k+m}} \subseteq I^{m+1}$ for $m \geq 0$.
\end{theorem}

This theorem was proved by Brian\c{c}on and Skoda analytically
for the ring of convergent power series over $\CC$ (\cite{brianconskoda})
and later algebraically for arbitrary regular rings (even in mixed
characteristic)
by Lipman and Sathaye (\cite{lipmansathayebrianconskoda}) and by
Lipman and Teissier (\cite{lipmanteissierbrianconskoda}).
The theory of tight closure requires that the ring contains a field,
but on the other hand it yields also a statement for non-regular rings
without any further assumption.

\medskip

\medskip

\medskip
{\em The definition of tight closure in zero characteristic}

\medskip
We shall now explain how to get a theory of tight closure for Noetherian
rings which contain a field of zero characteristic
in such a way that the properties described in the last section
hold again.

The main idea is to consider a relative situation over a finitely generated
$\ZZ$-domain, where the fiber over the generic point has zero characteristic,
but the special fibers have positive characteristic.
We take then the behavior of the tight closure in almost all special fibers
as the leading principle for the definition of tight closure
in the generic fiber.
Here we encounter different possibilities to fix such a relative setting
and therefore there exist different possibilities
to develop a theory of tight closure in characteristic zero
by reduction to positive characteristic.
It is not known whether these different notions
lead really to different theories.
We consider here only the following definition, see
\cite[Definition 3.1]{hochsterappendix} and \cite{hochsterhuneketightzero}.

\begin{definition}
Let $R$ denote a Noetherian $\QQ$-algebra, $f \in R$ an element
and $J$ an ideal in $R$.
Then we declare that $f \in J^*$ if and only if
there exists a finitely generated $\ZZ$-subalgebra $S$ of $R$
which contains $f$ and such that for $I=J \cap R$ the following holds.

For almost all prime numbers $p \in \ZZ$
the containment $f_p \in I_p^*$ holds in $S_p$.
Here we set $S_p = S \otimes _\ZZ \ZZ/(p)= S/(pS)$, $f_p$ denotes
the image of $f$ in $S_p$ and $I_p$ is the extended ideal of $I$ in $S_p$.
\end{definition}

The theory in zero characteristic arising from this definition
is called {\em equational tight closure} or {\em $\QQ$-tight closure}.
We just call it tight closure in this survey.
To show that this conception yields a useful theory,
forces some amount of technical work, witnessed by the voluminous
manuscript \cite{hochsterhuneketightzero}.
In the definition above we may enrich the finitely generated $\ZZ$-algebra
$S$ and we may assume that it contains
a set of ideal generators for $J$.
Since the definition relies only on almost all prime numbers,
we may invert elements of $S$ if we want to.

With this definition the main results of tight closure hold also
in characteristic zero;
in particular the statements \ref{tightregular}, \ref{tightpersistence},
\ref{coloncapturing} and \ref{tightbrianconskoda}
together with their consequences are true for $p \geq 0$.

\smallskip
The notion of F-regularity for a $\QQ$-algebra $R$
deserves some care. It is defined as before by the property that
every ideal in $R$ is tightly closed.
There is however another similar notion called {\em of F-regular type}.
This means that we have a family of models for $R$
such that on an open non-empty subset all special fibers are F-regular.
A {\em family of models} for $R$ consists of a finitely generated
$\ZZ$-domain $A$ and a flat finitely generated $A$-algebra $S$
with $R= S \otimes_A K = (S \otimes_A Q(A) ) \otimes_{Q(A)} K$
for a suitable field $Q(A) \subseteq K$.

\smallskip
Is it possible to define a closure operation in zero characteristic
with similar properties as tight closure without reduction to positive
characteristic? What about mixed characteristic?
A serious attempt in this direction
is the notion of solid closure due to
Hochster (see \cite{hochstersolid}), which we will discuss in the next section
\ref{solidsection}. However, over a field of zero characteristic,
this notion does not have the right properties,
since in dimension $\geq 3$ not every ideal
in a regular ring is closed with respect to solid closure.
For a proposal how to rescue this notion see \cite{brennerproposal}.

\medskip

\medskip

\medskip
{\em Tight closure, plus closure and the localization problem}

\medskip
So far we have met three possibilities to force
an element to lie in the tight closure of an ideal:
through colon-capturing, through the theorem of Brian\c{c}on-Skoda and
through persistence. These possibilities exist
also in characteristic zero.
In this section however we explain a forth possibility
to find elements in the tight closure, which is useful
only in positive characteristic.

\begin{theorem}
\label{contraction}
Let $R \subseteq S$ be a finite extension of commutative
Noetherian domains over a field of characteristic
$p >0$. Let $I \subseteq R$ denote an ideal and $f \in R$ an element
such that $f \in IS$. Then $f \in I^*$ holds.
\end{theorem}

\proof
It is a rather easy to show
(see \cite[Key Lemma]{smithtightintroduction} or
\cite[Corollary 2.3]{hochstersolid})
that there exists an $R$-linear mapping
$\varphi: S \ra R$ such that $\varphi(1) =c \neq 0$
(this property means by definition that $S$ is a solid $R$-algebra, see section
\ref{solidsection}). Suppose that $f \in IS= (f_1, \ldots, f_n)S$.
Then also $f^q \in (f_1^q, \ldots, f_n^q)= I^{[q]}S$ holds,
say $ f^q = s_1f_1^q + \ldots + s_nf_n^q $.
Applying the $R$-linear mapping $\varphi$ to this equation we get
$ cf^q = \varphi(s_1)f_1^q + \ldots + \varphi(s_n)f_n^q $,
hence $cf^q \in (f_1^q, \ldots, f_n^q)$ holds in $R$.
\qed

\medskip
The property described in \ref{contraction}
is called the {\em contraction property}.
This property holds in every characteristic,
but it is not very useful in zero characteristic:
under the condition that $R$ is a normal domain,
the trace map $t : S \ra R$ shows that an
element $f \in IS$ must belong to $I$ itself.

The statement in \ref{contraction} leads to the definition of
the {\em plus closure}. For an ideal in a Noetherian domain
we set
$$ I^+ = \{ f \in R: \exists \, R \subseteq S
\mbox{ finite extension of domains such that } f \in IS \} \, .$$
Hence \ref{contraction}
tells us that $I^+ \subseteq I^*$.

Another way to define the plus closure is to consider
the absolute integral closure of $R$.
For a Noetherian domain $R$ let $Q=Q(R)$ denote its quotient
field and let $\overline{Q}$ denote an algebraic closure.
The integral closure of $R$ in $\overline{Q}$ is denoted by $R^+$
and is called the {\em absolute integral closure} of $R$.
With this we may also write $I^+ = IR^+ \cap R$.

One of the most difficult questions of tight closure theory
is whether \ref{contraction} has a converse,
i.e. does there exist for every $f \in I^*$ in a Noetherian domain $R$
over a field of positive characteristic
a finite extension $R \subseteq S$ such that $f \in IS$.
Hochster calls this a ``tantalizing question''
(\cite{hochstersolid}).

\begin{problem}
\label{tightplustantal}
Let $R$ denote a Noetherian domain over a field of positive characteristic.
Is it true for every ideal $I\subseteq R$
that the tight closure is the same as the plus closure, $I^*=I^+$?
\end{problem}

The most important result so far in this direction
is the following theorem of K. Smith.

\begin{theorem}
\label{tightparameter}
Let $R$ be an excellent local Noetherian domain of dimension $d$
over a field of positive characteristic.
Let $x_1, \ldots, x_d$ denote parameters in $R$.
Then 
$$(x_1, \ldots, x_d)^+ =(  x_1, \ldots, x_d)^* \, .$$
\end{theorem}
\proof
See \cite[Theorem 7.1]{huneketightparameter} or \cite{smithparameter}.
\qed

\medskip
Beside this beautiful result there is not much known about
this problem.
``There is no non-trivial class of rings in which this open problem
has been solved'', K. Smith writes in \cite{smithtightintroduction}.
In particular this problem is open for two-dimensional normal
graded $K$-domains.
One result of this habilitation thesis is that for a normal homogeneous
coordinate ring $R$ over an elliptic curve
and for every homogeneous $R_+$-primary ideal $I$ the equality
$I^* =I^+$ holds.

\medskip
What about the elements of tight closure which
arise by colon-capturing and by the theorem of Brian\c{c}on-Skoda,
do they belong also to the plus closure?
For the first kind of containment the following
theorem of Hochster-Huneke gives an answer.

\begin{theorem}
\label{plusbig}
Let $R$ denote an excellent local domain over a field
of positive characteristic.
Then the absolute integral closure $R^+$
is a big Cohen-Macaulay algebra.
\end{theorem}
\proof
See \cite[Theorem 7.1]{hunekeapplications}.
\qed

\medskip
A {\em big Cohen-Macaulay algebra} means that every system of parameters
of $R$ becomes a regular sequence in $R^+$
(but $R^+$ is not Noetherian anymore, hence the name ``big'').
From this property of the absolute integral closure $R^+$ it follows
that for every system of parameters $x_1, \ldots, x_d$ in $R$
the inclusion
$((x_1, \ldots, x_{k-1}): x_k) \subseteq (x_1, \ldots,x_{k-1})R^+$
holds in $R^+$.
Therefore the part of tight closure given by
colon-capturing belongs also to the plus closure.

This is also true for the part of tight closure coming from the
theorem of Brian\c{c}on-Skoda, see \cite[Theorem 12.8]{hunekeapplications}.

\medskip
The problem \ref{tightplustantal} is connected with another open
problem in tight closure, the localization problem.

\begin{problem}
\label{localization}
Let $R$ denote a Noetherian domain over a field
and let $F$ denote a multiplicatively closed system.
Does $(IR_F)^*= I^*R_F$ hold?
\end{problem}

Here the inclusion $\supseteq $ is a special case of
the persistence property and $\subseteq $ is the open problem.
A positive answer to \ref{tightplustantal}
would imply a positive answer to \ref{localization},
since it is easy to see that the plus closure
commutes with localization.

\medskip

\medskip

\medskip
{\em Computation of tight closure, degree bounds and vanishing theorems}

\medskip
About the possibility to compute the tight closure of an ideal
or to answer the question whether $f \in I^*$ holds, C. Huneke writes
``Tight closure is very difficult to compute; indeed that is necessarily
the case. It contains a great deal of information concerning subtle
properties of the ring and the ideal''
(\cite[Basic Notions]{hunekeapplications}).
The problem lies in the fact that due to the definition
we have to check infinitely many conditions.

Since we know often by the theory of test ideals
for which $c$ we may check the conditions $cf^q \in I^{[q]}$,
we get in case $f \not\in I^*$ after finitely many tests also a
negative answer.
In the case $f \in I^*$ however every single test has
a positive answer, but we cannot deduce the containment
after finitely many steps.

\medskip
Lets restrict in the following to the {\em graded situation},
i.e. let $R$ denote an $\NN$-graded $K$-algebra of finite type,
where $R_0 =K$ is a field of arbitrary characteristic.
The algebra is called {\em standard-graded}
if it is generated by finitely many
homogeneous elements of $R_1$.
For an {\em $R_+$-primary ideal} $I$ (i.e. $V(I)=V(R_+)$)
there exists a degree bound such that $R_{\geq N} \subseteq I$ holds
and there exists also a (in general smaller) degree bound for $I^*$.
Can we give estimates for such {\em inclusion bounds}
for a given homogeneous $R_+$-primary ideal, for example in terms of the degree
of some ideal generators?

The situation with this question is similar to the problem
about the plus closure explained in the previous section.
For parameter ideals where exists a satisfactory answer,
whereas for arbitrary primary ideals very few is known.
K. Smith has given the following two estimates.

\begin{theorem}
\label{smithinclusionbound}
Let $R$ denote an $\NN$-graded $K$-domain of dimension $d$ and let $I$
denote an $R_+$-primary ideal which is generated
by  homogeneous elements $f_1, \ldots,f_n$
of degree $d_i$. Then the following two estimates hold.

\renewcommand{\labelenumi}{(\roman{enumi})}
\begin{enumerate}

\item
$R_{\geq  \dim (R) (\max_i \{ d_i\}) } \subseteq I^*$

\item
$R_{\geq d_1+ \ldots + d_n} \subseteq I^*$

\end{enumerate}

\end{theorem}
\proof
See \cite[Proposition 3.1 and Proposition 3.3]{smithgraded}.
\qed

\medskip
These statements rely on the theorem \ref{tightbrianconskoda}
of Brian\c{c}on-Skoda and are therefore also true for the (graded)
plus closure.
If the ideal is a parameter ideal
(generated by $\dim R$ elements), then the second bound
is sharp, as we shall see in Theorem \ref{strongvanishing} below.

\medskip
Lets have a look at bounds from below, i. e. {\em exclusion bounds}.
Such bounds are not really exclusion bounds in the strong sense
(think of $0$),
but have the following property: if $\deg (f) \leq M $, then
$f \in I^*$ holds if and only if $f \in I$ holds.
A general result in this direction is the following.

\begin{theorem}
\label{exclusionbound}
Let $R$ denote a normal $\NN$-graded domain of finite type
over a perfect field $K$ of positive characteristic.
Let $I$ denote a homogeneous ideal and suppose
that $I$ is generated by elements of degree $\geq M$.
Then for a homogeneous element $f$ of degree $ \leq M$ we have
that $f \in I^*$ if and only if $f \in I$.
\end{theorem}

\proof
See \cite[Theorem 2.2]{smithgraded} for a proof
in terms of differential operators or proposition \ref{exclusionboundgeo}
below for a proof within our geometric interpretation.
\qed

\medskip
As we have already mentioned, the degree bound
in Theorem \ref{smithinclusionbound}(ii) is sharp for the parameter case.
This rests upon the following theorem,
which has been conjectured and proved for $\dim\, R=2$ in
\cite{hunekesmithkodaira} and was proved in general
by N. Hara in \cite{harafrobeniusrational}.
In the formulation of this theorem we encounter the condition
$p \gg 0$ for the characteristic. This means that we have a family of models
over a finitely generated $\ZZ$-domain $A$,
that the conditions are true in the generic point
and that the conclusion is true for the fibers
over an open non-empty subset.

\begin{theorem}
\label{strongvanishing}
Let $R$ denote a normal standard-graded $K$-domain
of dimension $d$
with an isolated singularity and
let $x_1, \ldots,x_d$ denote homogeneous parameters of degree $d_i$.
Set $N= \sum_i d_i$.
Suppose that the characteristic of $K$ is zero or $p \gg 0$.
Then the following hold.

\renewcommand{\labelenumi}{(\roman{enumi})}
\begin{enumerate}

\item
$$ (x_1, \ldots,x_d)^* =
\sum _i (x_1, \ldots, x_{i-1},x_{i+1}, \ldots,x_d)^* + R_{\geq N} \, .$$

\item
Moreover, if $R$ is Cohen-Macaulay, then
$$ (x_1, \ldots,x_d)^* =(  x_1, \ldots,x_d) + R_{\geq N} \, .$$
\end{enumerate}

\end{theorem}
\proof
See \cite[Theorems 5.15 and  6.1]{huneketightparameter},
\cite[Section 3]{smithtightintroduction} and \cite{harafrobeniusrational}.
\qed

\medskip
The Theorem \ref{strongvanishing} is called
``Strong vanishing theorem''. This deserves an explanation.
The statement (i) in \ref{strongvanishing}
is in positive characteristic equivalent
with the property that the Frobenius acts injectively
on the negative part of the local cohomology module
$H^d_{R_+}(R)$, where $d$ is the dimension  of $R$.

This last property is not only true for the dimension $d$,
but also more general for $H^{i}_{R_+}(R)$, $0 < i \leq d$ (for $p \gg 0$).
From this it follows at once that the negative
part of $H^{i}_{R_+}(R)$ is zero for $0 <i<d$ and for $p \gg 0$.
This property holds then also in zero characteristic.
Therefore we get an algebraic proof
of the famous Kodaira vanishing theorem from algebraic geometry,
if we apply this statement to the ring of global sections
of an ample invertible sheaf on a projective variety.

\begin{theorem}
\label{kodairavanishing}
Let $X$ denote a projective non-singular variety over a field of
characteristic zero and let $\shL$ be an ample invertible sheaf.
Then $H^{i}(X, \shL^{-1}) =0$ holds for $0 \leq i < \dim \, X$.
\end{theorem}

We emphasize that this proof works with positive characteristic
although the statement does not hold in every positive
characteristic.
For another algebraic proof of the vanishing theorem of Kodaira,
which work also by reduction to positive characteristic, but not
with tight closure, look at \cite{deligneillusie}.

\medskip
Theorem \ref{strongvanishing} means in particular that for
homogeneous parameters $f_1, \ldots, f_d$ in a graded
Cohen-Macaulay algebra there exists a degree bound,
namely $N=\sum_i \deg (f_i)$ such that 
every element of degree $\geq N$ belongs to $(f_1, \ldots, f_d)^*$
and such that every element of degree $< N$ which belongs
to $(f_1, \ldots, f_d)^*$ belongs already to the ideal itself.
This gives a completely satisfactory
numerical characterization for the tight closure of a parameter ideal
in a graded domain.

\medskip
Even in the first non-trivial case where
$(f_1, \ldots, f_n)$ is a homogeneous
$R_+$-primary ideal in a two-dimensional normal standard-graded
$K$-domain there is not much known apart
from the bounds mentioned in \ref{smithinclusionbound}
and in \ref{exclusionbound}.
In particular it is not clear
whether and under which conditions there may exist a common bound
which separates the inclusion to the tight closure
and the exclusion from the tight closure
like in the parameter case.

This habilitation thesis contains various results for
two-dimensional normal stan\-dard-graded domains
which give under some weak conditions
on the relation module for some ideal generators $f_1, \ldots, f_n$
much better estimates.
Moreover, if the sheaf of relations is
semistable, we also get a common bound for inclusion
and exclusion.
This bound is given by $(\deg \, (f_1) + \ldots +\deg\, (f_n))/(n-1)$
and has a geometric interpretation within our setting,
see theorems \ref{topselfintersection} and \ref{semistablevanishing}.

\medskip

\medskip

\medskip
{\em Singularities, test ideal and multiplier ideal}

\medskip
We cannot resist to mention the connection
between the theory of tight closure and singularities,
though we will not come back to this point in the following.

Regularity-properties defined by properties of the
Frobenius have been studied even before the rise of
tight closure, for example by
Fedder, Goto, Hochster, Kunz, Mehta, Ramanathan, J. Roberts, Watanabe
(see \cite{kunzpositive}, \cite{fedderwatanabe},
\cite{hochsterrobertsfrobenius},
\cite{mehtaramanathanfrobenius}).
These concepts have been systematized since the introduction of tight closure.
The work of both Hara and Smith yield to astonishing relations between
properties of singularities which are defined
by reduction to positive characteristic with the help of the Frobenius
morphism (so called $F$-singularities)
and with the test ideal on one hand,
and properties of singularities which
are defined with the help of a resolution
of singularities and the multiplier ideal on the other hand.

\medskip
We have already mentioned the F-regu\-larity defined by the property
that every ideal is tightly closed with the variant
being of F-regular type in characteristic zero.

\begin{definition}
A local Noetherian ring over a field $K$
is called {\em F-rational} if for every system $f_1, \ldots, f_d$
of parameters the identity
$(f_1, \ldots, f_d)= (f_1, \ldots, f_d)^*$ holds.
\end{definition}

The property F-rationality is in general a weaker notion than
F-regularity. For a Gorenstein ring however
both notions coincide,
see \cite[Proposition 5.1]{hochsterhuneketightclosure}.
Moreover, it is then even enough to know
that just one single parameter ideal is tightly closed.
An F-rational ring is Cohen-Macaulay
(this follows as in \ref{regularcohenmacaulay})
and also normal (\cite[Proposition 10.3.2]{brunsherzog}).

\begin{definition}
A Noetherian local ring over a field $K$ of positive
characteristic is called {\em F-pure}
if the Frobenius morphism is a pure mapping.
\end{definition}

The definition of $F$-rational carries over to zero characteristic
with the variant being {\em of $F$-rational type}.
In zero characteristic the notion of (dense) F-pure means
that we have a family of models such that the set of points is dense
for which the fibers are F-pure
(that we demand this property
only for a dense subset instead for an open non-empty
subset is due to the fact that for a cone over an elliptic curve
the property F-pure jumps from one prime number to another).

\medskip
For the definition of test elements and the test ideal
we have to go back to the definition
of tight closure in positive characteristic.

\begin{definition}
Let $R$ denote a Noetherian Ring containing a field of positive
characteristic $p >0$.
An element $c \in R$ is called {\em test element}
if  for every ideal $I \subseteq R$
and for every $f \in I^*$ the condition
$cf^q \in I^{[q]}$ holds for all $q=p^{e}$.
\end{definition}

The set of all test elements gives an ideal
which is called the {\em test ideal} of $R$.
We may take every element of the test ideal which is not
contained in any minimal prime of $R$
to test the containment $f \in I^*$.
The test ideal of a ring is the unit ideal if and only if the ring is
$F$-regular.
The main existence result about test elements is the following theorem.

\begin{theorem}
\label{testexistence}
Let $R$ be a Noetherian ring over a field of 
positive characteristic.
Suppose that $R$ is essentially of finite type over
an excellent local ring or that $R_{\rm red}$ is F-finite.
Let $0 \neq c \in R$ such that $R_c$ is regular.
Then there exists a power of $c$ which is a test element.
\end{theorem}
\proof
See \cite[Theorem 3.2]{huneketightparameter}.
\qed

\medskip
If $R,\fom$ has an isolated singularity it follows from this
existence result that the test ideal is primary to the maximal ideal $\fom$.
Test elements and test ideals play an important role
in the proof of persistence of tight closure.
The test ideal in characteristic zero is defined
by the elements in a family of models which
are test elements over an open non-empty subset.

\medskip
Let us now recall some properties and notions
of a variety $X$ (or of a finitely generated algebra $R$)
over a field of characteristic zero
which are defined with the help
of a {\em resolution of singularities}, see \cite{kollarpairs}.

Let $X=\Spec \, R$ denote an affine variety over a field $K$ of
characteristic zero and let $\pi: \tilde{X} \ra X$
denote a resolution of singularities.
We say that $X$ has {\em rational singularities}
if $X$ is normal and Cohen-Macaulay and if
$H^{i} (\tilde{X}, \O_{\tilde{X}})= 0$ for $i \geq 1$.
This property is independent of the resolution.

Now the following relationship between rationality and F-rationality holds.

\begin{theorem}
\label{rationalfrational}
Suppose that $R$ is a finitely generated domain
over a field $K$ of cha\-racteristic zero.
Then $R$ is of F-rational type if and only if it has
rational singularities.
\end{theorem}
\proof
See \cite[Theorem 3.2]{smithtightintroduction},
\cite{harafrobeniusrational} and \cite{smithrational}.
\qed

\medskip
Suppose further that the exceptional divisors
$E_1, \ldots, E_k$ of the resolution have simple
normal crossings and that $X$ is {\em $\QQ$-Gorenstein}, i.e.
the canonical module $\omega_X$ is a torsion element in the class group,
i.e. $(\omega_X)^{\otimes r} \cong \O_X$ holds for some $r \in \NN$.
Then it is possible to compare
the canonical divisor $K_{\tilde{X}}$ on the resolution
with the pullback of $K_X$ and to get an equation
of $\QQ$-Weil-divisors,
$$ K_{\tilde{X}} = \pi^* K_X + \sum a_i E_i \,,$$
where the $a_i$ are rational numbers.

\begin{definition}
The singularities of $X$ are called {\em log-terminal}
if $a_i > -1$ holds for all $i=1, \ldots, k$.

The singularities of $X$ are called {\em log-canonical},
if $a_i \geq -1$ holds for all $i=1, \ldots, k$.

The {\em multiplier ideal} of $X$ is defined by
$ \pi_* (\O_{\tilde{X}}(\lceil \sum a_i E_i \rceil ))$.
\end{definition}

These types of singularities arises in the study of birational geometry
and the minimal model program, see \cite{minimalmodel} or
\cite{kollarmori}.
The multiplier ideal was first defined analytically by Nadel and later
algebraically by Kollar, see \cite{kollarpairs},
\cite{einvanishing}, \cite{lazarsfeldpositive}.
The multiplier ideal is the unit ideal if and only if the singularities
are log-terminal.
Theorem \ref{rationalfrational} implies the following theorem.

\begin{theorem}
\label{fregularterminal}
Let $X$ denote a normal $\QQ$-Gorenstein variety over a field
of characteristic zero.
Then $X$ has $F$-regular type if and only if
$X$ has log-terminal singularities.
\end{theorem}
\proof
This follows from theorem \ref{rationalfrational} using the
canonical cover $X' \ra X$. This may be constructed
with the help of an identification
$(\omega_X)^{\otimes r} \cong \O_X$.
Now due to results of Watanabe and of Kawamata both
properties carry over from $X$ to $X'$ and back.
Since $X'$ Gorenstein, the statement follows.
\qed

\begin{problem}
Let $R$ denote a normal $\QQ$-Gorenstein domain of finite type
over a field of characteristic zero.
Is it true that $R$ is of (dense) F-pure type if and only
if the singularities of $R$ are log-canonical.
\end{problem}

This is a difficult problem with deep connections to number theory.
The answer is yes in dimension $\leq 2$
due to the classification of such singularities,
see \cite[Theorem 4.4 (4)]{smithvanishingsingularities}
and \cite{haraclassification}.

The best result so far about the relationship of
$F$-singularities and singularities defined by resolution
is the following theorem which was proved independently
by Hara and Smith. It implies at once
the theorem \ref{fregularterminal}.

\begin{theorem}
\label{testmultiplier}
Let $R$ denote a normal $\QQ$-Gorenstein local domain over a field
of characteristic zero.
Then the multipier ideal of $R$ and the test ideal of $R$ are the same.
\end{theorem}
\proof
See \cite[Theorem 3.1]{smithtestmultiplier} and \cite{haratestmultiplier}.
\qed

\bigskip

\subsection{Solid closure and forcing algebras}
\label{solidsection}

\markright{Solid closure and forcing algebras}
\

\bigskip
The theory of solid closure was an attempt by Hochster
(see \cite{hochstersolid} and \cite{hochstertightsolid})
to define directly without any reduction procedure
a closure operation for ideals in every Noetherian ring
(even in mixed characteristic) in the hope to get a general theory
with similar properties as tight closure.
For rings containing a field of positive characteristic it coincides with
tight closure under some mild conditions (see Theorem \ref{tightsolidpositive}
below).
It turned out however by an example due to P. Roberts
(see \cite{robertscomputation} and \ref{robertsexample} below)
that regular rings containing a field of characteristic zero
of dimension three do not have the expected property that
every ideal is solidly closed. This ``discouraging result''
(Hochster, \cite{hochstersolid}) shows that solid closure is too big.

This notion is anyway very important for us because of the following reasons.
Solid closure gives a cohomological characterization
of tight closure in positive characteristic.
It gives a good notion
with all expected properties for every Noetherian ring in dimension two,
the case we will mostly deal with.
Furthermore it provides the notion of forcing algebra which will play a crucial role
in the following.

The starting point of the definition of solid closure is the notion
of a solid module.

\begin{definition}
\label{solidmoduledefinition}
Let $R$ denote a commutative domain and let $M$ denote an $R$-module.
Then $M$ is called {\em solid} if there exists
a non-trivial $R$-module homomorphism
$\varphi: M \ra R $.
\end{definition}

Hence an $R$-module is solid if and only if its dual module
$ M^\dual = \Hom_R (M,R) \neq 0$.
This definition applies in particular to $R$-algebras.
The argument used in the proof of \ref{contraction}
says that a finite extension $R \subseteq S$ is a solid
$R$-algebra. This proof gives at once the following more general result:
if $I \subseteq R$ is an ideal in a domain over a field of positive
characteristic and if $R \ra S$ is a solid algebra such that
$f \in IS$, then $f \in I^soclo$.

The solidity of an $R$-module may also be characterized with the
help of local cohomology.
For local cohomology we refer to \cite{brodmannsharp},
\cite{brunsherzog}, \cite{grothendiecklocalcoho}
or \cite{SGA2}.

\begin{proposition}
\label{solidmodulecohomology}
Suppose that $R,\fom$ is a complete local Noetherian domain of
dimension $d$. Then an $R$-module is solid if and only if
$H^d_\fom (M) \neq 0$.
\end{proposition}
\proof
This is an application of Matlis duality,
see \cite[Corollary 2.4]{hochstersolid}.
\qed

\medskip
The notion of solid closure is now defined by the following condition
(see also \cite[Definition 1.2]{hochstersolid}).

\begin{definition}
Let $R$ denote a Noetherian ring, let $I \subseteq R$ be an ideal and
let $f \in R$ be an element.

\renewcommand{\labelenumi}{(\roman{enumi})}
\begin{enumerate}

\item
If $R$ is a complete local domain,
then $f$ belongs to the {\em solid closure} of $I$,
$f \in I^{\soclo}$, if and only if
there exists a solid $R$-algebra $S$ such that $f \in IS$.

\item
In general we declare that $f \in I^{\soclo}$
if $f \in (I \hat{R}_{\fom}/ \foq)^\soclo$ holds for every maximal ideal
$\fom$ of $R$ and every minimal prime $\foq$
of the completion $\hat{R}_{\fom}$ of $R_\fom$.
\end{enumerate}

\end{definition}

For properties of this closure operation we refer to
\cite{hochstersolid} and \cite{hochstertightsolid}.
Another important concept which arises naturally in the study of solid
closure is the notion of a forcing algebra
(Hochster calls this a generic forcing algebra)

\begin{definition}
Let $R$ denote a commutative ring and let $f_1, \ldots,f_n \in R$ and
$f_0 \in R$ be elements. The $R$-algebra
$$A=R[T_1,\ldots,T_n]/(f_1T_1+\ldots +f_nT_n+f_0) \,,$$
is called the {\em forcing algebra}
for the elements (the {\em forcing data}) $f_1, \ldots,f_n;f_0$.
\end{definition}

\begin{remark}
The forcing algebra forces that $f_0 \in (f_1, \ldots, f_n)A$.
Every other $R$-algebra $S$ with the property that
$f_0 \in (f_1, \ldots, f_n)S$ factors (not uniquely) through $A$:
if $- f_0 = s_1f_1+ \ldots + s_nf_n$ with $s_i \in S$, then we just
send $T_i \mapsto s_i$.
The forcing algebra has a section $\Spec \, R \ra \Spec A$
if and only if $f_0 \in (f_1, \ldots ,f_n)$.

If $A$ is the forcing algebra for elements $f_1, \ldots, f_n; f_0 \in R$
and if $ \varphi: R \ra R'$
is a ring homomorphism, then $A'=A \otimes_R R'$ is the forcing algebra
for the elements $\varphi(f_1), \ldots ,\varphi(f_n); \varphi(f_0)$.
\end{remark}

\medskip
The forcing algebra $A$ for forcing data $f_1, \ldots, f_n; f_0 \in R$
is solid if and only if there exists a solid $R$-algebra $S$ such
that $f_0 \in (f_1, \ldots, f_n)S$.
Hence in the definition of solid closure we only have to look
whether the forcing algebra is solid or not,
and therefore it is possible to
characterize the notion of solid closure in terms of
forcing algebras and local cohomology as follows.

\begin{proposition}
\label{solidcohomology}
Let $R$ be a Noetherian ring and let $f_1, \ldots,f_n, f_0 \in R$.
Then $f_0 \in (f_1, \ldots,f_n)^{\soclo}$ if and only if
for every local complete domain $R' =  \hat{R}_{\fom}/\foq $
{\rm(}where $\fom $ is a maximal ideal of $R$ and $\foq$ is a minimal prime
of $\hat{R}_{\fom}${\rm)} we have that
$H^d_{\fom'}(A') \neq 0$, where $d = \dim \, R'$ and $A'$ is the
forcing algebra over $R'$.
\end{proposition}
\proof
This follows from \ref{solidmodulecohomology}, see
also \cite[Proposition 5.3]{hochstersolid}.
\qed

\begin{remark}
\label{solidremarks}
This cohomological condition must only be checked for the maximal ideals
$\fom \supseteq (f_1, \ldots,f_n)$.
\end{remark}

The connection to tight closure in positive characteristic
is given by the following theorem.

\begin{theorem}
\label{tightsolidpositive}
Suppose that $R$ is a Noetherian ring containing
a field of characteristic $p >0$. Suppose furthermore that $R$
is essentially of finite type over an excellent local ring
or that the Frobenius endomorphism is finite.
Then $I^*= I^\soclo$.
\end{theorem}
\proof
See \cite[Theorem 8.6]{hochstersolid}.
\qed

\begin{remark}
In this habilitation thesis we shall deal with the notion of solid closure
rather than tight closure, and therefore we will also use the notation
$I^\soclo$.
As the previous theorem states this does not make any difference in
positive characteristic.
In characteristic zero and in dimension two it is not clear whether
solid closure coincides with any version of tight closure, but solid closure
has all the properties which one expect from a tight closure theory.
In dimension $\geq 3$ solid closure does not has the right properties,
as the following example of Roberts shows.
For a variant of solid closure with the right properties,
which is also defined without any reduction to positive characteristic,
see \cite{brennerproposal}.
\end{remark}

\begin{example}
\label{robertsexample}
A computation of P. Roberts in \cite{robertscomputation} shows
that regular rings over a field of characteristic zero
are not solidly closed in general.
Roberts considers the ideal $(x^3,y^3,z^3)$
in the polynomial ring $K[x,y,z]$, where $K$ is a field of characteristic
zero, and proves that $(xyz)^2 \in (x^3,y^3,z^3)^\soclo$
by showing that
$$H^3_{(x,y,z)}
(K[x,y,z][T_1,T_2,T_3]/(x^3T_1+y^3T_2+z^3T_3- x^2y^2z^2) )\neq 0 \, .$$
\end{example}

\subsection{Cohomological dimension}
\label{cohodimsection}

\markright{Cohomological dimension}

\ 

\bigskip
Recall that the {\em cohomological dimension} $ \cd \,X$
of a scheme $X$ is the maximal number $i$ such
that there exists a quasicoherent sheaf $\shF$ with $H^{i}(X, \shF) \neq 0$.
This notion was developed by R. Hartshorne,
see \cite{hartshornecohdim}. In the following proposition we
gather together some results on
cohomological dimension which we will need in the sequel.

\begin{proposition}
\label{cdproperties}
The cohomological dimension has the following properties.

\renewcommand{\labelenumi}{(\roman{enumi})}
\begin{enumerate}

\item 
Let $X$ denote a Noetherian scheme. Then $\cd\, X \leq \dim X$

\item
If $X' \ra X$ is an affine morphism of separated Noetherian schemes,
then $\cd \, X' \leq \cd \, X$.

\item
A Noetherian scheme $X$ is an affine scheme if and only if $\cd\, X \leq 0$
{\rm(}the cohomological dimension of the empty scheme is $-1${\rm)}.

\item
An irreducible separated variety $X$ is proper
if and only if $\cd \, X =\dim \, X$.

\item
Suppose that $R, \fom$ is a local Noetherian ring of dimension $d$.
Then the cohomological dimension of $D(\fom)$ is $d-1$.

\item
Suppose that $A$ is a graded $K$-algebra of finite type over a field $K$
and let $\foa$ denote a homogeneous ideal.
Then $\cd \, (D(\foa)) = \cd \, (D_+(\foa))$.
\end{enumerate}
\end{proposition}
\proof
(i) is another formulation of the vanishing theorem of Grothendieck,
see \cite[Theorem III.2.7]{haralg}.

(ii) follows from $\rm \check{C}$ech-cohomology, see \cite[III.4]{haralg},
in particular exc. 4.1.

(iii) This is the cohomological characterization of affine schemes
due to Serre,
see \cite[Theorem III.3.7]{haralg}.

(iv) is a theorem of Lichtenbaum, see \cite[III \S 3]{haramp}
or \cite{kleimanvanishing}.

(v) The estimate `$\leq$' follows from (i).
The estimate `$ \geq $' is clear for $d=0,1$ and follows for $d \geq 2$
from the Theorem of Grothendieck
that $H^{d-1}(D(\fom), \O_X) = H^d_{\fom} (R) \neq 0$,
see \cite[Theorem 3.5.7]{brunsherzog}.

(vi) The cone mapping $D(\foa) \ra D_+(\foa)$ is affine, hence
`$\leq$' follows from (i). On the other hand,
every quasicoherent sheaf on the open subset
$D_+(\foa) \subseteq \Proj\, A$
is quasicoherent extendible to $\Proj\, A$ and hence of type $\tilde{M}$,
where $M$ is a graded $A$-module
(\cite[Propositions II.5.8 and II.5.15]{haralg}).
\qed

\medskip
For an ideal $\foa \subseteq R$ we call the maximal number $j$
such that there exists an $R$-module $M$ with
$H^{j}_\foa (M) \neq 0$ the {\em cohomological height}, $\cht \, (\foa)$
(this is also called the {\em local cohomological dimension}).
For $\cd (D(\foa)) \geq 1$ we have $\cht \,(\foa) = \cd \,(D(\foa)) +1 $,
due to the long exact sequence relating local cohomology and global
cohomology.

We may express now the solid closure of an ideal $I=(f_1, \ldots, f_n )$
with the help of the forcing algebra and the cohomological dimension
in the following way.

\begin{proposition}
\label{solidcd}
Let $R $ be a normal excellent domain.
Let $f_1, \ldots,f_n \in R$ be elements which are
primary to a maximal ideal $\fom$
of height $d$, let $f_0 \in R$ and let
$A=R[T_1,\ldots,T_n]/(f_1T_1+\ldots +f_nT_n+f_0)$
denote the forcing algebra for this data. Then the following hold.

\renewcommand{\labelenumi}{(\roman{enumi})}
\begin{enumerate}

\item
$f_0 \in (f_1, \ldots, f_n)^\soclo$ if and only if the
cohomological height
of the extended ideal $\fom A$ is $d$.

\item
If $d \geq 2$, then $f_0 \in (f_1, \ldots, f_n)^\soclo$ if and only if
the cohomological dimension of $W=D(\fom A) \subset \Spec\, A$
is $d-1$.

\item
If $d=2$,
then $f_0 \in (f_1,\ldots, f_n)^\soclo $
if and only if $\, D(\fom A)$ is not an affine scheme.

\end{enumerate}
\end{proposition}
\proof
(i).
Since the completion of a normal and excellent domain is again a domain,
the condition $f_0 \in (f_1, \ldots,f_n)^\soclo$
is due to \ref{solidcohomology} equivalent to $H^d_{\fom '}(A') \neq 0$, where $R'$ is the completion of
$R_{\fom}$ and $A'=A \otimes_RR'$.
Since local cohomology commutes with completion
this is equivalent to $H^d_\fom (A) \neq 0$.
Since $H^d_\fom (A) = H^d_{\fom A} (A)$
this implies that $\cht \, (\fom A) \geq d$, and equality must hold
since the cohomological height of $\fom A$ can not be bigger
than $\cht \, (\fom) = d$ due to \ref{cdproperties}(ii),(v).
On the other hand, if $H^d_\fom (A) =0$, then this holds for every
$A$-module $M$, since we find a surjection
$A^{(J)} \ra M \ra 0$ and since $H^{d+1}_\fom(-)=0$.

(ii) follows from (i) by the long exact sequence of local and
global cohomology. (iii) follows from (ii) and
the cohomological characterization of affine schemes, \ref{cdproperties}(iii).
\qed

\begin{remark}
\label{affinefunctionsremark}
An open subset $U=D(\foa) \subseteq \Spec \, B$, $\foa =(a_1, \ldots ,a_k)$
is affine if and only if there exist
elements $q_1, \ldots, q_k \in \Gamma( U, \O_U)$ such that
$\sum a_i q_i=1$.
So if $f_0=1$, then the forcing equation gives
that $f_1T_1 + \ldots + f_nT_n= -1$, hence $D(f_1, \ldots ,f_n)$ is affine.
If $(f_1, \ldots ,f_n)$ is $\fom$-primary, then
$f_0 \not\in (f_1, \ldots, f_n)^\soclo)$.
\end{remark}

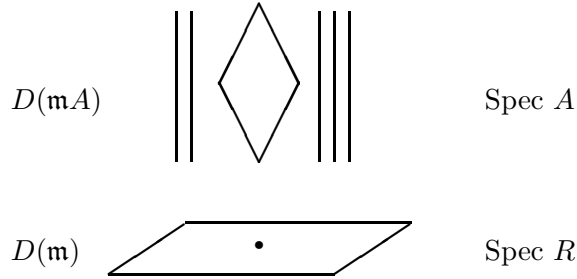
\unitlength1cm
\begin{picture}(12,5)

\put(6, 1.4){\circle*{0.1}}

\thicklines

\put(4,1){\line(1,0){3}}
\put(4,1){\line(3,2){1}}
\put(5.02,1.69){\line(1,0){3}}
\put(7.,1){\line(3,2){1}}

\put(2.7,1.2){$D(\fom)$}
\put(9,1.2){$\Spec\, R$}
\put(9,3.2){$ \Spec\, A$}
\put(2.7,3.2){$ D(\fom A)$}

\put(6,2.5){\line(1,2){0.53}}
\put(6,2.5){\line(-1,2){0.53}}
\put(6,4.6){\line(1,-2){0.53}}
\put(6,4.6){\line(-1,-2){0.53}}

\thinlines
\put(4.9,2.5){\line(0,1){2}}
\put(5.1,2.5){\line(0,1){2}}
\put(6.8,2.5){\line(0,1){2}}
\put(7.,2.5){\line(0,1){2}}
\put(7.2,2.5){\line(0,1){2}}
\end{picture}

For studying tight closure problems in terms of the forcing algebra
for forcing data $f_1, \ldots ,f_n;f_0$
we have to study the cohomo\-lo\-gical and the geometric properties of 
the open subset $D(\fom A) \subseteq \Spec \, A$ for a maximal
ideal $\fom \in \Spec\, R$ as shown in the picture.
In particular, for $\dim\, R=2$, we have to look for criteria
which imply that $D(\fom A)$ is an affine scheme or not.
Note that the cohomological height of
$\fom \subset R$ is the same as its height $d$.
If after building the forcing algebra the cohomological height of $mA$
is still $d$, then $f_0 \in I^\soclo$,
but if the cohomological height has dropped,
then $f_0 \not\in I^\soclo$.

\medskip
In studying the geometric properties of $D(\fom A)$
we will also use the notion of superheight.
Recall that the {\em superheight} of an ideal $\foa \subseteq R$
is the maximal height of $\foa R'$ in any Noetherian $R$-algebra $R'$.
The superheight of an ideal is less or equal to the cohomological
height due to \ref{cdproperties}(ii) and (v).

\begin{proposition}
\label{finiteextension}
Let $R$ be a Noetherian ring, 
let $f_1, \ldots,f_n \in R$ be elements primary to
a maximal ideal $\fom$ of height $d$.
Let $f_0 \in R$ denote another element
and let
$A=R[T_1,\ldots,T_n]/(f_1T_1+\ldots +f_nT_n+f_0)$ denote the
forcing algebra for these data.
Let $W=D(\fom A) \subset \Spec\, A$.

Suppose that there exists a local Noetherian ring $R'$ of dimension $d$
and a ring homomorphism $R \ra R'$ such that
$V(\fom R') =V(\fom_{R'})$ and
$f_0 \in (f_1,\ldots,f_n)R'$ hold.
Then the superheight of $\fom A$ is $d$ and hence also
the cohomological height is $d$
{\rm (}and $f_0 \in (f_1, \ldots,f_n)^\soclo$, if $R$ is normal and
excellent{\rm )}.
\end{proposition}

\proof
The morphism $\Spec \, R' \ra \Spec \, R$ lifts to a
morphism $\varphi: \Spec\, R' \ra \Spec\, A$ such that
$\varphi^{-1} ( D(\fom A)) = D(\fom_{R'})$.
Hence the superheight of $\fom A$ is $\geq d$ and hence equals $d$.
\qed

\medskip
In particular, if $I=(f_1, \ldots ,f_n) \subseteq R$
is primary to a maximal ideal of height $d$ and $R \subseteq S$
is a finite extension of Noetherian domains
such that $f_0 \in IS$, then the superheight
and the cohomological height of $mA$ is $d$.
The relation between the superheight and the cohomological height
of $\fom A$ in the forcing situation depends heavily on the characteristic of the base field.
In positive characteristic it is open whether
the two notions fall apart.
In zero characteristic however the two notions do fall apart
quite often.
In particular tight closure (solid closure) in zero characteristic
yields at once examples of ideals $\foa$ with superheight one,
but such that the open subset $D(\foa)$ is not affine.
For other examples see \cite{neeman} and \cite{brennersuperheight}.
We will apply this in section \ref{complex} to give new counterexamples
to the hypersection problem of complex analysis.

\begin{proposition}
\label{superhoehe}
Let $K$ be a field of characteristic zero and let
$R$ be a normal excellent $K$-domain.
Let $f_1,\ldots, f_n$ be primary to a maximal ideal $\fom$
of height $d$
and let $f_0 \in R$.
Suppose that
$f_0 \not\in (f_1, \ldots, f_n)$, but $f_0 \in (f_1, \ldots, f_n)^\soclo$.
Then the cohomological height of
$\fom A \subseteq A= R[T_1, \ldots,T_n]/ (\sum f_iT_i+f_0)$
is $d$, but its superheight is $< d$.
\end{proposition}

\proof
Let $R'$ denote a local normal Noetherian domain
of dimension $n \leq d$ and let
$ \varphi: A \ra R'$
be a homomorphism such that
$ V( (\fom A)R')= V(\fom_{R'})$.
This gives a homomorphism
$\psi: R \ra R'$ such that
$V(\fom R') = V(\fom_{R'})$.
If $n = d$, then $\hat{R_\fom} \ra \hat {R'}$ would be finite
(after enlarging the base field).
But a finite extension $R \subseteq S$
such that $f_0 \in (f_1, \ldots, f_n)S$
with $R$ normal forces that
$f_0 \in (f_1, \ldots, f_n)$ holds already in $R$
due to the existence of the trace map, see
\cite[Remarks 9.2.4]{brunsherzog}.
Hence the superheight is $< d$, but the cohomological height
is $d$ due to \ref{solidcd} (i).
\qed

\bigskip

\subsection{Vector bundles, locally free sheaves and projective bundles}

\markright{Vector bundles, locally free sheaves and projective bundles}

\

\bigskip
In this section we recall briefly the notion of vector bundles and its
relatives, locally free sheaves and projective bundles,
and we fix our notations.

Let $X$ denote a scheme.
A scheme $p: V \ra X$ is called a {\em geometric vector bundle} of rank $r$
if there exists an open affine covering
$X= \bigcup_i U_i$ and $U_i$-isomorphisms
$ \psi_i: V|_{U_i} = p^{-1}(U_i) \ra U_i \times \AA^r =\AA^r_{U_i}$
such that for every open affine subset
$U \subset U_i \cap U_j$
the transition mappings
$ \psi_j \circ \psi_i^{-1}: \AA_U^r \ra \AA_U^r$ 
are linear automorphisms, i.e. they are induced by an automorphism
of the polynomial ring $ \Gamma(U, \O_X)[T_1, \ldots, T_r]$
given by $T_i \mapsto \sum a_{ij} T_j$
(if $X$ is separated, then $U_i \cap U_j$ is affine and the condition has
to be checked only for such subsets).
$\AA^r_X$ together with its natural structure of a vector bundle
is called the trivial bundle over $X$ of rank $r$.

A vector bundle $V$ over $X$ is in particular a commutative group scheme over
$X$, where the addition $V \times_X V \ra V$ is locally given
by the identifications with $\AA_U^r$.
Due to the compatibility conditions this gives a well defined addition
(it is even
``un X-sch\'{e}ma en modules sur le X-sch\'{e}ma en anneaux X[T]'',
see \cite[9.4.13]{EGAI}).

\medskip
A vector bundle $V \ra X$ comes along with two locally free sheaves.
The {\em sheaf of sections} $\shS_V$ is given by
$$ \Gamma(U, \shS_V)= \{ s: U \ra V|_U :\,  p \circ s = \id_U \, \} \, .$$
The linear structure on the vector bundle makes
this into a quasicoherent sheaf. Since locally $V \cong \AA^r$,
the sheaf of sections is locally isomorphic to $\O_X^r$,
hence $\shS_V$ is a {\em locally free sheaf} of rank $r$ on $X$.
The group scheme $V \ra X$ represents the group functor
$X' \mapsto \Gamma(X', \varphi^*(\shS_V))$, where
$\varphi:X' \ra X$ is any scheme over $X$.

The {\em sheaf of linear forms} $\shF_V$ for a vector bundle
$V \ra X$ is given by
$$ \Gamma(U, \shF_V)
= \{ f: V|_U \ra \AA|_U \mbox{ is a linear morphism} \, \} \, .$$
Here a morphism $f:V \ra W$ of two vector bundles over $X$ is called linear
if for an open covering $X= \bigcup_i U_i$ with simultaneous identifications
$\psi_i: V|_{U_i} \ra \AA^r_{U_i}$
and $\varphi_i: W|_{U_i} \ra \AA^s_{U_i}$ (which are compatible with the
identifications which define the linear structures on the bundles)
the induced mappings
$$ \varphi_i \circ \psi_i^{-1} : \AA^r_{U_i} \lra \AA^s_{U_i} $$
are linear.
This is again a locally free sheaf of rank $r$, and we have the relationship
$\shS_V^\dual = \Hom_{\O_X}(\shS_V , \O_X) = \shF_V$.

Vector bundles and locally free sheaves are essentially the same objects.
For a locally free sheaf $\shF$ of rank $r$ on $X$
the symmetric algebra
$S(\shF)= \oplus_{k \geq 0} S^k(\shF)$ is a graded $\O_X$-algebra.
The associated spectrum of this algebra
(see \cite[9.4]{EGAI} or  \cite[1.7]{EGAII}),
$$V= \Spec \, S(\shF) = \Spec \, \oplus_{k \geq 0} S^k(\shF)$$
is a geometric vector bundle over $X$ with the sheaf of linear forms $\shF$.
This gives a correspondence of vector bundles and locally free sheaves.

We feel free to jump from one viewpoint to another,
but we will always make it clear whether we are dealing
with vector bundles or with locally free sheaves.
In particular we will always speak about the sheaf of sections
or the sheaf of linear forms and not about the sheaf associated to a
vector bundle.

There is one point however where the two concepts behave quite different.
A subbundle $W \subseteq V$ is a linear mapping
$W \ra V$ which gives for every point $x \in X$
a linear subspace $W_{\kappa(x)} \cong \AA^s_{\kappa(x)} \hookrightarrow
\AA^r_{\kappa(x)} \cong V_{\kappa(x)}$,
whereas a subsheaf $\shR \subseteq \shS$ means only that the
sheaf morphism $\shR \ra \shS$ is injective.
There may be points $x \in X$ such that
$\shR \otimes_{\O_X} \kappa(x) \ra \shS \otimes_{\O_X} \kappa(x)$
is not injective or is even zero.
A locally free subsheaf $\shR \subseteq \shS$ corresponds to a subbundle
if and only if the quotient sheaf $\shS/\shR$ is also locally free.

\medskip
Suppose that $\shF$ is a locally free sheaf on $X$.
Then the graded algebra $S(\shF) = \oplus_{k \geq 0} S^k(\shF)$ gives
also rise to a {\em projective bundle}
$$ \PP(\shF) = \Proj \oplus_{k \geq 0} S^k(\shF) \stackrel{\pi}{\lra} X\, ,$$
see \cite[4.1]{EGAII} or \cite[Proposition II.7.11]{haralg}.
$\PP(\shF)$ has a relatively ample invertible sheaf $\O_{\PP(\shF)} (1)$
with the property that
$\pi_* \O_{\PP(\shF)} (k) =S^k(\shF)$, see \cite[II.7.11]{haralg}.
Tensorizing the locally free sheaf $\shF$ with an invertible sheaf $\shL$
does not change the projective bundle,
i.e. there exists a canonical isomorphism
$i_\shL: \PP(\shF) \ra \PP(\shF \otimes \shL)$.
Under this isomorphism we get (\cite[Proposition 4.1.4]{EGAII}) that
$i_\shL^* ( \O_{\PP(\shF \otimes \shL)} (n))
=  \O_{\PP(\shF)}(n) \otimes \pi^* \shL^n $.

For a vector bundle $V \ra X$ with sheaf of linear forms $\shF$
we set $\PP(V) =\PP(\shF)$.
Hence a point in $\PP(V_{\kappa(x)})$ represents a line in $V_{\kappa(x)}$.
The natural cone mapping
$$V^\times  = V- \{\mbox{zero section} \} \lra \PP(V)$$
maps every closed point $v \neq 0$ to the line $[v]$.

We set
$\AA_{\PP(\shF)} (1)= \Spec \, \oplus_{k \geq 0} \O_{\PP(\shF)}(k)$.
This line bundle $\AA_{\PP(\shF)} (1)$ on $\PP(\shF)$ has therefore
the sheaf of linear forms $\O_{\PP(\shF)}(1)$ and
the sheaf of sections $\O_{\PP(\shF)}(-1)$.
We use this notation also for the corresponding line bundle
on a projective space.

If $\varphi: Y \ra X$ is a morphism, then a morphism $s:Y \ra \PP(\shF)$ over $X$
is equivalent with an invertible sheaf $\shL$ on $Y$
and a surjection $\varphi^* (\shF) \ra \shL \ra 0$,
where $\shL =s^* \O_{\PP(\shF)}(1)$, see \cite[Proposition 4.2.3]{EGAII} or
\cite[Proposition II.7.12]{haralg} for this correspondence.

A subbundle $W \subseteq V$ corresponds to a surjection
of linear forms $\shF \ra \shG \ra 0$. This induces a surjection
$\pi^* \shF \ra \pi^* \shG \ra \O_{\PP(\shG)} \ra 0$ and hence
a closed embedding $\PP(W) \subseteq \PP(V)$.

\newpage

\thispagestyle{empty}

\markboth{2. Geometric interpretation of tight closure via bundles}
{2. Geometric interpretation of tight closure via bundles}

\newsavebox{\viereck}

\savebox{\viereck}(5,4)[bl]{
\thinlines
\put(0,0){\line(1,0){3}}
\put(0,0){\line(3,2){1}}
\put(1.0,0.67){\line(1,0){3}}
\put(3.,0){\line(3,2){1}}}

\newsavebox{\viereckdick}

\savebox{\viereckdick}(5,4)[bl]{
\thicklines
\put(0,0){\line(1,0){3}}
\put(0,0){\line(3,2){1}}
\put(1.0,0.67){\line(1,0){3}}
\put(3.,0){\line(3,2){1}} }

\newsavebox{\strichedreistufev}

\savebox{\strichedreistufev}(5,4)[bl]{\thinlines
\put(0., 0){\line(0,1){1.2}}
\put(0.2, 0.15){\line(0,1){1.2}}
\put(0.4, 0.3){\line(0,1){1.2}}
}

\newsavebox{\bogenfett}

\savebox{\bogenfett}(5,2)[bl]{\thicklines
\bezier{300}(0,0)(1, -1)(3,0)}

\newsavebox{\strichedreistufe}

\savebox{\strichedreistufe}(5,4)[bl]{\thinlines
\put(0., 0){\line(0,1){2}}
\put(0.3, 0.2){\line(0,1){2}}
\put(0.6, 0.4){\line(0,1){2}}
}

\newsavebox{\strichefuenf}

\savebox{\strichefuenf}(5,4)[bl]{
\put(0., 0){\line(0,-1){2}}
\put(0.75, 0){\line(0,-1){2}}
\put(1.5, 0){\line(0,-1){2}}
\put(2.25, 0){\line(0,-1){2}}
\put(3,0){\line(0,-1){2}}
}

\newsavebox{\bogen}

\savebox{\bogen}(5,2)[bl]{
\bezier{300}(0,0)(1, -1)(3,0)}

\newsavebox{\krummbogen}

\savebox{\krummbogen}(5,4)[bl]{
\thicklines
\bezier{300}(0, 0.)( 0.5 ,0.4 )(1, 0)
\bezier{300}(1.,0.)(2, -0.7)(3,0.4)
}

\newsavebox{\ruledsurface}

\savebox{\ruledsurface}(6,5)[bl]{
\thicklines
\bezier{300}(0, 1)(1, 0)(3,1)
\thinlines
\bezier{300}(0,2)(1, 1)(3,2)
\bezier{300}(0,4)(1, 3)(3,4)
\put(0.1, 1.9){\line(0,1){2}}
\put(0.5, 1.65){\line(0,1){2}}
\put(0.9, 1.53){\line(0,1){2}}
\put(1.3, 1.5){\line(0,1){2}}
\put(1.7, 1.54){\line(0,1){2}}
\put(2.1, 1.64){\line(0,1){2}}
\put(2.5, 1.77){\line(0,1){2}}
\put(2.9, 1.94){\line(0,1){2}}
}

\newsavebox{\subbundle}

\savebox{\subbundle}(5,5)[bl]{

\put(0, 0.2){\usebox{\bogenfett}}
\put(0.7, 0.9){\usebox{\bogenfett}}

\thicklines

\put(0, 0.7){\line(1,1){0.7}}
\put(3., 0.7){\line(1,1){0.7}}

\thinlines
\put(0.29, .5){\line(1,1){0.7}}
\put(0.57, .34){\line(1,1){0.7}}
\put(0.9, .26){\line(1,1){0.7}}
\put(1.3, .24){\line(1,1){0.7}}
\put(1.71, .27){\line(1,1){0.7}}
\put(2.14, .38){\line(1,1){0.7}}
\put(2.55, .49){\line(1,1){0.7}}
}

\newsavebox{\viereckfaser}

\savebox{\viereckfaser}(5,4)[bl]{
\thinlines
\put(0,0){\line(0,1){2}}
\put(0,0){\line(1,1){0.7}}
\put(0.7, .7){\line(0,1){2}}
\put(0,2){\line(1,1){0.7}}
}

\newsavebox{\bundleranktwo}

\savebox{\bundleranktwo}(5,5)[bl]{

\put(0,0){\usebox{\bogenfett}}

\put(0,1.){\usebox{\bogen}}
\put(0,3.){\usebox{\bogen}}
\put(0.7, 1.7){\usebox{\bogen}}
\put(0.7, 3.7){\usebox{\bogen}}

\put(0, 1.5){\usebox{\viereckfaser}}

\put(1.3, 1.){\usebox{\viereckfaser}}

\put(3, 1.5){\usebox{\viereckfaser}}
}

\section{Geometric interpretation of tight closure via bundles}

\bigskip
In this chapter we shall introduce several geometric objects
associated to forcing data consisting of ideal
generators $f_1, \ldots ,f_n$ and another element $f_0$
in a commutative Noetherian ring $R$ in order to provide
a geometric interpretation of the property
$f_0 \in (f_1, \ldots ,f_n)^\soclo$.

We have mentioned already (section \ref{solidsection})
the forcing algebra corresponding to 
such data.
Ideal generators $f_1, \ldots, f_n$ themselves
give rise to the coherent sheaf of relations which is locally free over
$D(f_1, \ldots ,f_n)$.
This sheaf of relations is the sheaf of sections in the relation bundle,
which is a commutative group scheme over $\Spec \, R$
and which is a vector bundle
over $D(f_1, \ldots ,f_n)$ (section \ref{relationsection}).

This group scheme acts on the spectrum of the forcing algebra for
$f_1, \ldots, f_n;f_0$. This spectrum induces on $D(f_1, \ldots ,f_n)$
an affine-linear bundle and it is a geometric torsor
(or principal homogeneous space) for the relation bundle
 (section \ref{affinelinearsection}).
Hence the theory of tight closure leads us to study the
cohomological properties of geometric torsors for a vector bundle in general.
For this purpose it is helpful to provide a projective realization
of such a torsor as the complement of a projective
subbundle inside a projective bundle (sections \ref{torsorsection}
and \ref{projectivebundlessection}).

The objects which we construct and describe in this chapter have
also graded counterparts.
This provides an entirely projective setting consisting of a projective bundle
over a projective variety together with a projective
subbundle of codimension one
(sections \ref{gradedsection} and \ref{gradedforcingsection}).
This subbundle is called the forcing divisor. Properties
of this divisor such as ampleness and bigness are of interest
for the underlying tight closure problem (section \ref{amplesection}).
This relationship is especially useful
in dimension two (\ref{dimensiontwosection}).

\subsection{Relation bundles}
\label{relationsection}

\markright{Relation bundles}

\

\bigskip
Let $R$ denote a commutative ring. Elements $f_1, \ldots, f_n \in R$
define a vector bundle over the open subset
$D(f_1, \ldots,f_n) \subseteq \Spec\, R$
in the following way.

\begin{proposition}
\label{relationbundle}
Let $R$ denote a commutative ring and let $I \subseteq R$ denote an ideal
with ideal generators $f_1, \ldots, f_n$.
Set $U=D(I)=D(f_1, \ldots, f_n) \subseteq \Spec \, R$.
Then the following hold.

\renewcommand{\labelenumi}{(\roman{enumi})}
\begin{enumerate}
\item 
$\tilde{V} = \Spec \, R[T_1, \ldots , T_n]/(f_1T_1 + \ldots +f_nT_n)$
is a commutative group scheme over $\Spec\, R$.

\item
The sheaf of sections in $\tilde{V}$ is the quasicoherent
sheaf of relations for the elements $f_1, \ldots ,f_n$.

\item
The restriction $V= \tilde{V}|_U$ is a vector bundle on $U$
of rank $n-1$.

\end{enumerate}
\end{proposition}
\proof
(i). The coaddition is given by
$$R[T_1, \ldots , T_n]/(\sum f_iT_i )
\lra R[T_1, \ldots , T_n, S_1, \ldots, S_n]/
( \sum f_iT_i, \,\sum f_iS_i) \, ,$$
where $T_i \mapsto T_i + S_i$.
It is easy to see that this fulfills the axioms of a group scheme.

(ii). A section
$$s: \Spec \, S \ra
\Spec\, S[T_1, \ldots , T_n]/(f_1T_1 + \ldots +f_nT_n) \, $$
over an affine subset
$W=\Spec \, S \subseteq \Spec \,R$ is just
a tupel $(s_1, \ldots , s_n) \in S^n$ satisfying
$f_1s_1 + \ldots +f_ns_n=0$.

(iii).
The bundle $V$ is on $D(f_i)$, $i=1, \ldots,n$,
isomorphic to
$$\Spec R_{f_i}[T_1, \ldots,T_{i-1},T_{i+1}, \ldots,T_n] \,,$$
and the transition functions send
$$T_i \mapsto
-\frac{1}{f_i}(f_1T_1+ \ldots +f_{i-1}T_{i-1}+f_{i+1}T_{i+1} + \ldots +f_nT_n) \, ,$$
thus they are linear and $V$ is a vector bundle
on $U=\bigcup_{i=1}^n D(f_i)$.
\qed

\begin{definition}
We call
$\tilde{V} =\Spec \, R[T_1, \ldots , T_n]/(f_1T_1 + \ldots +f_nT_n)$
the {\em relation group scheme}
and its restriction $V= \tilde{V}|_U$ the {\em relation bundle}
corresponding to the ideal generators $f_1, \ldots, f_n$.
We will denote the sheaf of relations by
$\shR = Rel (f_1, \ldots,f_n) = {\rm Rel} (f_1, \ldots ,f_n)^{\tilde{}}$
and the sheaf of linear forms of $V$ by $\shF = \shR^\dual$.
\end{definition}

\begin{remark}
The bundle $\tilde{V}=\Spec \, R[T_1, \ldots , T_n]/(f_1T_1 + \ldots +f_nT_n)$
is not flat over $X$. The fiber over a point $x \in D(f_1, \ldots ,f_n)$
has dimension $n-1$ and the fiber over
a point $x \in V(f_1, \ldots ,f_n)$ has dimension $n$.
\end{remark}

\bigskip

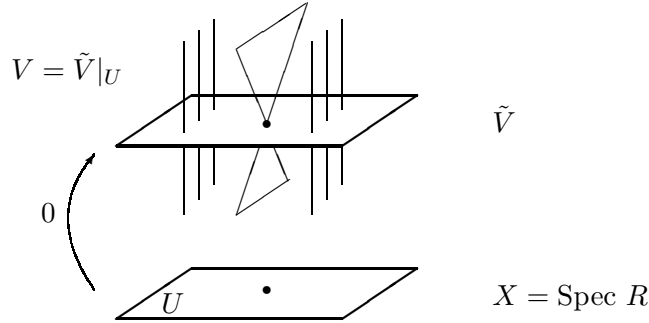
\setlength{\unitlength}{1cm}
\begin{picture}(12,5)

\put(6.5, 1.1){\circle*{0.1}}

\put(4.5,.7){\usebox{\viereckdick}}

\put(4.5,3){\usebox{\viereckdick}}

\put(5.1, .8){$U $}
\put(9.5, .9){$X =\Spec\, R$}
\put(9.5, 3.2){$ \tilde{V}$}
\put(3.1,3.9){$ V= \tilde{V}|_U$}

\put(6.5, 3.3){\circle*{0.1}}

\put(6.1, 2.1){\line(3,2){0.7}}
\put(6.1, 4.3){\line(3,2){0.92}}
\put(6.1, 2.1){\line(1,3){0.3}}
\put(6.5, 3.3){\line(1,3){0.54}}
\put(6.1, 4.3){\line(2,-5){0.4}}
\put(6.6, 3){\line(2,-5){0.18}}

\put(4.2 , 2.9){\vector(2,3){.01}}
\put(3.5,2){$0$}

\bezier{300}(4.2, 1.1)(3.5, 2.1)(4.2, 2.9)

\put(5.4,3){\line(0,-1){0.9}}
\put(5.6,3){\line(0,-1){0.7}}
\put(5.8,3){\line(0,-1){0.5}}

\put(5.4,3.2){\usebox{\strichedreistufev}}
\put(7.1,3){\line(0,-1){0.9}}
\put(7.3, 3){\line(0,-1){0.7}}
\put(7.5,3){\line(0,-1){0.5}}

\put(7.1,3.2){\usebox{\strichedreistufev}}
\end{picture}

\savebox{\strichedreistufev}{}

\begin{remark}
It is also possible to describe the relation algebra
as a symmetric algebra.
If the $R$-module $E$ is defined by
the presentation
$$R \stackrel{f_1, \ldots ,f_n}{\lra} R^n \lra E \lra 0, \,$$
then $S(E) \cong R[T_1, \ldots ,T_n]/(\sum f_iT_i)$,
see \cite[Introduction]{vasconcelos}.

\end{remark}

\begin{remark}
\label{presentingsequence}
The relation bundle $V=\tilde{V}|_U$ for elements $f_1, \ldots ,f_n \in R$
sits in the following short exact sequence
of vector bundles on $U=D(f_1, \ldots, f_n)$,
$$0 \lra V \lra  \AA^n_U \stackrel{f_1, \ldots, f_n}{\lra} \AA_U \lra 0 \, ,$$
where the mapping sends $(s_1, \ldots, s_n) \mapsto
\sum_i s_if_i$.
We call this sequence the {\em presenting sequence} or
{\em defining sequence} for the relation bundle.
We adopt the same name for the short exact sequences
of the sheaf of relations and the sheaf of linear forms,
$$0 \lra \shR
\lra  \O_U^n \stackrel{f_1, \ldots, f_n}{\lra} \O_U \lra 0 \, $$
and
$$0 \lra \O_U \stackrel{f_1, \ldots, f_n}{\lra}  \O_U^n \lra \shF \lra 0 \, .$$
\end{remark}

\begin{remark}
Suppose that we have an inclusion 
$(f_1, \ldots,f_n) \subseteq (g_1, \ldots, g_k)$ of ideals and let
$M=(a_{ij})$ be a matrix such that
$f_i = \sum_j a_{ij}g_j$.
This matrix defines a linear mapping
$$R[S_1, \ldots, S_k]/(\sum g_j S_j) \lra R[T_1, \ldots, T_n]/(\sum f_i T_i )$$
by sending $S_j \mapsto \sum_i a_{ij} T_i$.
This gives a morphism
$\tilde{V}(f_1, \ldots,f_n) \ra \tilde{V}(g_1, \ldots,g_k)$
of group schemes
and a morphism of vector bundles on $D(f_1, \ldots ,f_n)$.

If some ideal generator is superfluous,
say $f_n = \sum_{i=1}^{n-1} a_if_i$, then we have the
decomposition
$$Rel(f_1, \ldots, f_n) = Rel(f_1, \ldots,f_{n-1} ) \oplus \O_U$$
by sending
$(r_1, \ldots ,r_n) \mapsto ((r_1+a_1r_n, \ldots, r_{n-1}+a_{n-1} r_n), r_n)$.
\end{remark}

\smallskip

\subsection{Affine-linear bundles arising from forcing
algebras}
\label{affinelinearsection}

\markright{Affine-linear bundles arising from forcing algebras}

\

\bigskip
Suppose further that $R$ is a commutative ring and that
$f_1, \ldots ,f_n \in R$ are ideal generators of an ideal $I$.
Let $f_0 \in R$ denote another element so that we can build
the forcing algebra
$$A=R[T_1, \ldots, T_n]/(f_1 T_1 + \ldots +f_nT_n +f_0) \, .$$
We will see in a minute that its spectrum is not a vector bundle,
but an affine-linear bundle over $D(f_1, \ldots ,f_n)$.
We recall the definition.

\begin{definition}
A scheme $p: B \ra X$ is called an {\em affine-linear bundle} of rank $r$
if there exists an open affine covering
$X= \bigcup_i U_i$ and $U_i$-isomorphisms
$ \psi_i: B|_{U_i} = p^{-1}(U_i) \ra U_i \times \AA^r =\AA^r_{U_i}$
such that for every open affine subset
$U \subset U_i \cap U_j$
the transition mappings
$ \psi_j \circ \psi_i^{-1}: \AA_U^r \ra \AA_U^r$ 
are affine-linear automorphisms, i.e. they are induced by an automorphism
of the polynomial ring $ \Gamma(U, \O_X)[T_1, \ldots, T_r]$
given by $T_i \mapsto \sum_j a_{ij} T_j + v_i$.
\end{definition}

\begin{proposition}
\label{relationaction}
Let $R$ denote a commutative ring and let $f_1, \ldots, f_n \in R$
be elements, set $U=D(I)=D(f_1, \ldots, f_n) \subseteq \Spec \, R$.
Let $f_0 \in R$ denote another element and let
$A=R[T_1, \ldots, T_n]/(f_1 T_1 + \ldots +f_nT_n +f_0)$ denote the forcing
algebra. Set $\tilde{B} = \Spec\, A$
and $B= \Spec\, A|_U$. Then the following hold.

\renewcommand{\labelenumi}{(\roman{enumi})}
\begin{enumerate}

\item 
$B=\Spec \,A|_U$ is an affine-linear bundle of rank
$n-1$ over $U$.

\item
The group scheme $\tilde{V}$ acts on $\tilde{B}$.

\item
This action is simply transitive, i.e.
the morphism
$$ \mu \times p_2:
\tilde{V} \times_X \tilde{B} \lra \tilde{B} \times_X \tilde{B} $$
is an isomorphism.

\end{enumerate}
\end{proposition}

\proof
(i).
Consider the mapping
$\Spec \,A \longrightarrow \Spec \, R$ over the open
subset $U=D(I)=D(f_1,\ldots,f_n)$.
On $D(f_i)$, $i \geq 1$, one can identify
$$ A_{f_i} = (R[T_1,\ldots,T_n]/(f_1T_1+\ldots +f_nT_n +f_0))_{f_i} \cong
R_{f_i}[T_1, \ldots,T_{i-1},T_{i+1},\ldots,T_n] \, .$$
So this mapping looks locally like
$D(f_i) \times \AA ^{n-1} \longrightarrow D(f_i)$.
The transition mapping on $D(f_if_j)$ is given by
$$R_{f_if_j}[T_1,\ldots,T_{i-1},T_{i+1},\ldots,T_n] \longrightarrow
R_{f_if_j}[T_1,\ldots,T_{j-1},T_{j+1},\ldots,T_n] \, , $$
where $ T_k \longmapsto T_k $ for $k \neq i,j $ and
$T_i \longmapsto -1/f_i(\sum_{j \neq i}f_jT_j +f_0) $.
This is an affine-linear mapping,
therefore the forcing bundle $\Spec \, A|_{D(I)}$
is an affine-linear bundle of rank $n-1$.

(ii).
The action of the group scheme
$\tilde{V}= \Spec R[S_1, \ldots, S_n]/(\sum f_iS_i)$
on $\tilde{B} =\Spec R[T_1, \ldots , T_n]/(\sum f_iT_i +f_0)$
is given by the coaction
$$ R[T_1, \ldots , T_n]/(\sum f_iT_i +f_0)
\lra R[T_1, \ldots , T_n, S_1, \ldots, S_n]/
(\sum f_iT_i +f_0, \, \sum f_iS_i)$$
by sending $T_i \mapsto T_i +S_i$.
It is easy to verify that this is indeed a group coaction.

(iii).
The morphism $\mu \times p_2$ corresponds to the ring homomorphism
$$ R[S_1, \ldots ,S_n, T_1, \ldots ,T_n]/ (\sum f_iS_i+f_0,\sum f_iT_i+f_0)
\ \ \ \   \ \  \ \ \ \ $$

$$
\ \ \ \ \ \   \lra
R[S_1, \ldots ,S_n, T_1, \ldots ,T_n]/ (\sum f_iS_i,\sum f_iT_i+f_0)$$
given by $S_i \mapsto S_i +T_i$ and $T_i \mapsto T_i$.
This is clearly an isomorphism.
\qed

\begin{remark}
\label{translations}
The difference of two sections in $\tilde{B}$ (over some open subset)
is a relation for $f_1, \ldots, f_n$.
Therefore the relations act as translations on the affine-linear
bundle given by the forcing algebra.
\end{remark}

\smallskip

\subsection{Cohomology classes, torsors and their geometric realizations}

\label{torsorsection}

\markright{Cohomology classes, torsors and their geometric realizations}

\

\bigskip
The geometric situation arising from forcing data
in a commutative ring has a natural generalization which we consider
in this section and throughout this habilitation thesis.
For notions of the theory of group schemes we refer to
\cite{SGA3}, \cite{GIT}, \cite{waterhouse}.

\begin{definition}
Let $G$ denote a group scheme over a scheme $X$.
A scheme $B \ra X$ together with an operation
$\mu: G \times_X B\ra B $ is called a
{\em geometric {\rm(}Zariski-{\rm)} torsor} for $G$
(or a $G$-{\em principal fiber bundle} or a {\em principal homogeneous space})
if there exists an open
covering $X= \bigcup U_i$ and isomorphisms
$\varphi_i: B|_{U_i} \ra G|_{U_i}$ such that the diagrams
(with $U = U_i$, $\varphi = \varphi_i$)
$$
\begin{CD}
G|_U \times_U B|_U  \, \,   @> \mu >> \, \,  B|_U    \\
\,  @V  id \times \varphi VV   \, \,     @VV  \varphi V \\
G|_U \times _U G|_U  \, \,  @> \mu >>  G|_U
\end{CD}
$$
commute.
\end{definition}

\begin{remark}
This concept is also applicable for other topologies like the flat topology
or the \'{e}tale topology and shows then its real strength
(see  \cite{SGA3}, \cite{giraud} and \cite{milne}).
We will apply this concept only for vector bundles $G=V$ with respect
to the Zariski-topology.
A $V$-torsor $B$ for a vector bundle $V$
is the same as an affine-linear bundle.
\end{remark}

\begin{proposition}
\label{cohotorsor}
Let $X$ denote a Noetherian separated scheme and let
$p:V \ra X$ denote a vector bundle
on $X$ with sheaf of sections $\shS$.
Then there exists a correspondence
between cohomology classes $c \in H^1(X,\shS)$
and geometric $V$-torsors.
\end{proposition}
\proof
We will describe this correspondence, see also
\cite[Proposition III.4.6]{milne}.
Let $B$ denote a $V$-torsor. Then there exists by definition
an open covering $X= \bigcup U_i$ such that $B|_{U_i}$
together with the action of $V|_{U_i}$
is isomorphic to $V|_{U_i}$ acting on itself.
These isomorphisms $\varphi_i : B|_{U_i} \ra V|_{U_i}$
induce isomorphisms
$$ \psi_{ij}= \varphi_j \circ \varphi_i^{-1}:
V|_{U_i \cap U_j} \lra V|_{U_i \cap U_j} \, .$$
These isomorphisms are compatible with the operation
of $V$ on itself,
and this means that they are of the form 
$\psi_{ij} = \id_{V|_{ U_{i} \cap U_j}} + s_{ij} $ with suitable sections
$ s_{ij} \in \Gamma(U_i \cap U_j, \shS)$.
This family defines a $\rm \check{C}$ech-cocycle for the covering
and gives therefore a cohomology class in $H^1(X,\shS)$.

For the reverse procedure, suppose that the cohomology class
$c \in H^1(X,\shS)$ is represented by a
$\rm \check{C}$ech-cocycle $s_{ij} \in \Gamma(U_i \cap U_j, \shS)$
for an open covering $X= \bigcup _i U_i$.
Set $B_i := V|_{U_i}$.
We take the morphisms
$$B_i|_{U_i \cap U_j}= V|_{U_i \cap U_j}
\stackrel{ \psi_{ij} }{\lra}
V |_{U_i \cap U_j}= B_j|_{U_i \cap U_j}$$
given by $\psi_{ij} := \id_{V|_{U_i \cap U_j}} + s_{ij}$
to glue the $B_i$ together to a scheme $B$ over $X$.
This is possible since the cocycle condition
guarantees the glueing condition \cite[0, 4.1.7]{EGAI}.
The action of $V_i $ on $B_i$ by itself glues also together to
give an action on $B$.
\qed

\medskip
This habilitation thesis deals to a large extent
with the following very general problem
which contains the problems of tight closure.

\begin{problem}
\label{cohoclassproblem}
Let $X$ denote a scheme and let $V$ denote a vector bundle
on $X$ with sheaf of sections $\shS$.
Study the cohomological and geometric properties of
a geometric torsor $B \ra X$ corresponding to a cohomology class
$c \in H^1(X, \shS)$ in dependence of $X$, $\shS$ and $c$.
\end{problem}

In order to work in a projective setting and to use the methods
of projective geometry it is very useful to have the following
realization of a torsor.

\begin{proposition}
\label{torsorprojective}
Let $X$ denote a Noetherian separated scheme and let
$V$ be a vector bundle on $X$ with sheaf of sections $\shS$.
Let $c \in H^1(X, \shS)$ be a cohomology class and let
$0 \ra V \ra V' \ra \AA_X \ra 0$ denote
the extension of vector bundles corresponding to 
$c \in H^1(X,\shS) \cong \Ext^1(\O_X,\shS)$.
Then the geometric torsor $B \ra X$ corresponding to $c$
{\rm(}as described in {\rm \ref{cohotorsor})} is isomorphic to
$\PP(V') -\PP(V)$.
\end{proposition}

\proof
The inclusion $e:V \hookrightarrow V'$ induces
a projective subbundle $\PP(V) \subset \PP(V')$.
Let $T : V' \ra \AA_X$ denote the linear form on $V'$ in the sequence,
hence $V$ is the zero set of $T$.
We consider the cone mapping $ V' \dasharrow \PP(V')$.
We claim that the zero set
$Z=V(T-1) \subset V'$ maps isomorphically onto $\PP(V') -\PP(V)$.
This may be checked locally in $X$ and is therefore clear.
Hence we shall show that $Z=V(T-1) \subset V'$ is isomorphic to the torsor
given by $c$.

First we claim that $Z$ is indeed a $V$-torsor by the action
$V \times_X Z \ra Z$ given by the addition in $V'$. This is a well defined
action on $Z$ since $T(v+z)= T(v) + T(z) =1$ (where $v \in V$, $z \in Z$).

Let now $U_i$ denote a covering on $X$ such that there
exist sections $s_i' \in \Gamma(U_i, \shS')$ which map
to $1 \in \Gamma(U_i, \O_X)$.
Consider the mappings $\psi_i := e + s_i' \circ p : V|_{U_i} \ra V'|_{U_i}$.
The image of this mapping lies on $Z$, since $T(s_i')=1$, and hence
we have mappings $\psi_i : V|_{U_i} \ra Z|_{U_i}$, which are isomorphisms.
These identifications commute with the operation,
since (let $p$ denote the projections)
$$ \psi_i (v_1 + v_2) = v_1+v_2 + (s_i' \circ p)(v_1 +v_2)
= v_1+ (v_2 + (s_i' \circ p) (v_2) ) =  v_1+ \psi_i (v_2) \, .$$
This shows that $Z$ is a $V$-torsor.

We can compute the corresponding cohomology class of this torsor
by using these identifications.
The transition mappings are given by
$$ \psi_j^{-1} \circ \psi_i : V|_{U_i \cap U_j} \lra V|_{U_i \cap U_j} $$
which sends $v \mapsto v + (s_i' -s_j')(p(v))$.
Therefore the corresponding cohomology class is given by
$s_{ij}=  s_i' -s_j' \in \Gamma(U_i \cap U_j , \shS)$.
But this is also the image of $1$ under the connecting homomorphism
for the short exact sequence
$0 \ra \shS \ra \shS' \ra \O_X \ra 0$.
\qed

\smallskip

\smallskip

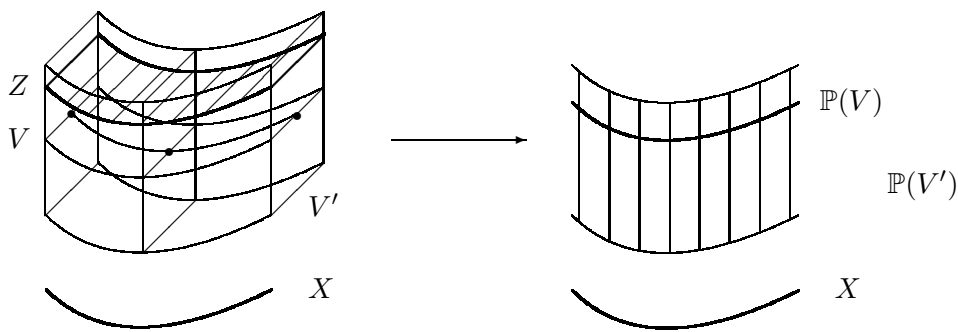
\setlength{\unitlength}{1cm}
\begin{picture}(12,5)

\put(1,0.5){\usebox{\bundleranktwo}}

\put(1, 3){\usebox{\subbundle}}

\put(1, 2.5){\usebox{\bogen}}

\put(1.35, 2.85){\usebox{\bogen}}

\put(1.7,3.2){\usebox{\bogen}}

\put(1.,3.){\line(1,1){0.7}}

\put( 1.35 , 3.35 ){\circle*{0.1}}

\put(2.3, 2.5){\line(1,1){0.7}}

\put( 2.65 , 2.85 ){\circle*{0.1}}

\put(4.,3.){\line(1,1){0.7}}

\put( 4.35 , 3.33 ){\circle*{0.1}}

\put(.5, 2.9){$V$}
\put(.5, 3.6){$Z$}
\put(4.5,2) {$V'$}
\put(4.5, .9){$X$}

\put(5.6, 3){$\vector(1,0){1.8}$}

\put(8,0){\usebox{\ruledsurface}}

\thicklines

\bezier{300}(8.,3.5)(9, 2.5)(11,3.5)

\put(11.5, .9) {$X$}
\put(11.3, 3.4) {$ \PP(V)$}
\put(12.2, 2.3) {$ \PP(V')$}
\end{picture}

\savebox{\subbundle}{}

\savebox{\bundleranktwo}{}

\savebox{\viereckfaser}{}

\begin{remark}
Suppose that we have any projective subbundle
$\PP(V) \subset \PP(V')$ of codimension one.
Then it is easy to show that the complement
$\PP(V') - \PP(V)$ is a $(V \otimes L)$-torsor for some suitable
line bundle $L$.
\end{remark}

We consider now some functorial properties of torsors.

\begin{lemma}
\label{monotonielemma}
Let $X$ denote a Noetherian separated scheme and let $V$ and $W$
be geometric vector bundles on $X$ with sheaf of sections $\shS$ and $\shT$.
Let
$\varphi: \shS \ra \shT$ be a homomorphism of locally free sheaves.
Let $c \in H^1(X,\shS)$ with corresponding extension
$0 \ra V \ra V' \ra \AA_X \ra 0$
and let $\varphi(c) \in H^1(X, \shT)$ be its image
with corresponding extension
$0 \ra W \ra W' \ra \AA_X \ra 0$.
Then there exists an affine morphism
$$\PP(V') - \PP(V) \lra \PP(W') - \PP(W) \, .$$
For the cohomological dimension we
have the inequality
$$ \cd (\PP(V') - \PP(V)) \leq \cd( \PP(W') - \PP(W)) \, .$$
\end{lemma}
\proof
First note that we have a mapping
$\varphi': V' \ra W'$ compatible with the
extensions and with $\varphi$,
see \cite[Ch. 3 Lemma 1.4]{maclanehomology}.
The induced rational mapping
$\PP(V') \dasharrow \PP(W')$ is defined outside the kernel of $\varphi '$.
Locally these mappings on the vector bundles look like
$$
\begin{CD}
\AA^r    @>  >> \AA^r  \times \AA   \\
 \varphi  @VVV     @VVV  \varphi \times id \\
\AA^s    @>  >> \AA^s  \times \AA  \, \, .
\end{CD}
$$
The point on $\PP(V')$ corresponding to the line given by $v=(v,t)$
does not lie on the subbundle $\PP(V)$ for $t \neq 0$. For these points
the rational mapping is defined and the image point $(\varphi(v),t)$
does not lie on the subbundle $\PP(W)$.
Hence we have an affine morphism
$$\PP(V')-\PP(V) \, \lra \,  \PP(W') - \PP(W) \, .$$
The second statement follows from \ref{cdproperties}(ii).
\qed

\begin{lemma}
\label{torsorpullback}
Let $ \varphi: X \ra Y$ denote a morphism between Noetherian
separated schemes,
let $\shS$ denote a locally free sheaf on $Y$,
let $c \in H^1(Y, \shS)$ with corresponding geometric torsor
$B \ra Y$.
Then $\varphi^*(B)= B \times _Y X$ is the geometric torsor
corresponding to $\varphi^{*}(c) \in H^1(X, \varphi^*(\shS))$.
\end{lemma}

\proof
The pull-back of the cohomology class is represented by the pull-back
of the extension of vector bundles
$0 \ra V \ra V' \stackrel{T}{\ra} \AA_Y \ra 0$
corresponding to $c$.
We have $B =V(T-1) \subseteq V'$
due to the proof of \ref{torsorprojective}. Hence
$$\varphi^*(B)= \varphi^{-1} (V(T=1)) = V( (T \circ \varphi) -1) 
\subset \varphi^*(V') \, .$$
\qed

\begin{proposition}
\label{universaltorsor}
Let $ X $ denote a Noetherian separated scheme,
let $\shS$ denote a locally free sheaf on $X$,
let $c \in H^1(X, \shS)$ with corresponding geometric torsor
$B \ra X$.
Then a morphism $\varphi: Z \ra X$ factors through $B$ if and only if
$0=\varphi^*(c) \in H^1(Z,\varphi^*( \shS))$.
\end{proposition}
\proof
The pull-back of the cohomology class vanishes on the torsor
which it defines: a principal fiber bundle has the property that
$\mu \times p_2: V \times _X B \cong B \times_XB$ is an isomorphism.
Hence the pull-back of $B \ra X$ under
the base change $B \ra X$ is isomorphic
to the projection $V \times_X B \ra B$.
This is trivial since it has a section.

For the other direction we consider the diagram
$$
\begin{CD}
B \times_XZ \, \,   @>  >>  \,  B \\
\,  @VVV     @VVV \\
Z  \,  @> \varphi  >> \, X
\end{CD}
$$
Here $B \times _XZ$ is the torsor corresponding to
the cohomology class $\varphi^*(c)$ due to \ref{torsorpullback}.
Since this pull-back is trivial there exists a section
$Z \ra B \times_X Z$ and therefore there exists also a factorization
$Z \ra B$.
\qed

\begin{corollary}
Let $ X $ denote a Noetherian separated scheme,
let $\shS$ denote a locally free sheaf on $X$
which has a decomposition $\shS= \shS_1 \oplus \shS_2$
into locally free sheaves.
Let $c \in H^1(X, \shS)$ denote a cohomology class
$c=(c_1,c_2)$ and let $B$, $B_1$ and $B_2$ denote the corresponding
torsors.
Then $B \cong B_1 \times_X B_2$.
\end{corollary}
\proof
The projections $p_i: \shS_1 \oplus \shS_2 \ra \shS_i$, $i=1,2$
yield due to \ref{monotonielemma} mappings
$B \ra B_i$, $i=1,2$ and therefore a mapping $B \ra B_1 \times_X B_2$.
This mapping is locally an isomorphism, hence an isomorphism.
\qed

\bigskip

\subsection{Projective bundles arising from forcing data}
\label{projectivebundlessection}

\markright{Projective bundles arising from forcing data}

\

\bigskip
Let $f_1, \ldots ,f_n \in R$ be elements in a commutative ring $R$.
Another element $f_0 \in R$ gives rise to the forcing
algebra, which induces an affine-linear bundle on $U=D(f_1, \ldots ,f_n)$.
This is also a torsor on $U$.
On the other hand, the elements $f_1, \ldots ,f_n$ define the
presenting sequence
$0 \ra \shR \ra \O^n_U \ra \O_U \ra 0$
on $U$.
The element $f_0 \in R \ra \Gamma(U, \O_U)$ defines also via
the connecting homomorphism for the presenting sequence
a cohomology class $c = \delta (f_0) \in H^1(U,\shR)$, which defines
itself a torsor on $U$.

\begin{proposition}
\label{forcingbundletorsortorsor}
Let $R$ denote a commutative ring and let $f_1, \ldots, f_n \in R$
be elements. Let $V$ denote the corresponding relation bundle
on $U=D(f_1, \ldots ,f_n)$ and let $\shR$ be the sheaf of relations.
Let $f_0 \in R$ denote another element.
Then the affine-linear bundle
$B= \Spec\, R[T_1, \ldots ,T_n]/(f_1T_1+ \ldots +f_nT_n +f_0)|_U$
given by these forcing data
is a $V$-torsor, which is isomorphic to the $V$-torsor given
by the cohomology class $c= \delta (f_0) \in H^1(U, \shR)$, where
$\delta$ is the connecting homomorphism
for the presenting sequence on $U$.
\end{proposition}
\proof
Set
$B_i = \Spec\, R_{f_i}[T_1, \ldots ,T_n]/(f_1T_1+ \ldots +f_nT_n +f_0)$
and $V_i =R_{f_i}[S_1, \ldots ,S_n]/(f_1S_1+ \ldots +f_nS_n)$
for $i=1, \ldots, n$.
Then the isomorphisms $\varphi_i:B_i \ra V_i$
given by $T_j \mapsto S_j$ for $j \neq i$
and by $T_i \mapsto S_i -f_0/f_i$ are compatible with the action of $V$ on $B$
defined in \ref{relationaction}(ii).
Therefore the forcing bundle is a $V$-torsor.

We shall compute the cohomology class represented by the
forcing bundle $B$, thus we have to look at the mappings
$$\varphi_j \circ \varphi_i^{-1}: V|_{U_i \cap U_j} \ra V|_{U_i \cap U_j}\, .$$
This mapping is given by
$S_k \mapsto S_k$ for $k \neq i,j$,
$S_i \mapsto S_i - f_0 /f_i$ and $S_j \mapsto S_j +f_0/f_j $.
The corresponding relations are
$$ (0, \ldots, 0, f_0/f_i,0, \ldots ,0, -f_0/f_j,0, \ldots ,0)
\in \Gamma(D(f_if_j), \shR) \, .$$
Now look at the connecting homomorphism for the
presenting sequence.
The ele\-ment $f_0 \in \Gamma(U,\O_U)$ (coming from $f_0 \in R$ by restriction)
is locally given by $(0, \ldots, 0, f_0/f_i,0 , \ldots, 0) \in
\Gamma(D(f_i), \O_U^n)$. The difference of two such local
representatives gives $\delta (f_0)$ as a $\rm \check{C}$ech cocycle.
But this is just the same.
\qed

\medskip
The extension of vector bundles $0 \ra V \ra V' \ra \AA_U \ra 0$
on $U=D(f_1, \ldots,f_n)$ corresponding to a cohomology
class $c \in H^1(U,\shR) =\Ext^1( \O_Y, \shR)$
has for $c= \delta (f_0)$ also another natural interpretation.

\begin{proposition}
\label{forcingsequence1}
Let $R$ be a commutative ring and let
$f_1,\ldots, f_n$ and $f_0$ be elements and set
$I=(f_1,\ldots,f_n)$, $U=D(I)$.
The schemes
$$ V= \Spec R[T_1,\ldots, T_n]/(\sum_{i=1}^n f_iT_i)|_U
\, \mbox{ and } \,
V' = \Spec R[T_0,\ldots, T_n]/(\sum_{i=0}^n f_iT_i)|_U $$
are vector bundles on $U$.
They are related by the short exact sequence of vector bundles
$$0 \longrightarrow V \longrightarrow V' 
\stackrel{T_0}{\longrightarrow} \AA^1_U \longrightarrow 0 \, .$$
This extension sequence represents $\delta(f_0) \in H^1(U,\shR)$,
where $\delta$ is the connecting homomorphism induced by
the presenting sequence.
The inclusion $V \subset V'$ yields
a closed embedding
$\PP(V) \hookrightarrow \PP(V')$ of projective bundles over $U$.
Its complement
$\PP(V') -\PP(V)$ is isomorphic to the
forcing affine-linear bundle
$$\Spec\, R[T_1,\ldots,T_n]/(f_1T_1+\ldots +f_nT_n+f_0)|_U \, .$$
\end{proposition}

\proof
The schemes $V$ and $V'$ are vector bundles on $U$
due to \ref{relationbundle}.
The element $T_0 \in \Gamma(V', \O_{V'})$ is a linear function
$T_0 : V' \ra \AA^1_U$. Its zero set is $V$.
Looking at $D(f_i)$, the exactness of the sequence on $U$ is clear.
Therefore $V'$ is indeed an extension of $\AA^1_U$ by $V$,
so we have to compute the cohomology class
defined by the corresponding sequence
$0 \ra \shR \ra \shR' \ra \O_U \ra 0$.

On $D(f_i)$, the unit $1 \in \Gamma(U, \O_U)$
is represented by the relation
$$(0, \ldots ,0, -\frac{f_0}{f_i},0 ,\ldots, 0, 1)
\in \Gamma(D(f_i), \shR')\, .$$
This defines the $\rm \check{C}$ech cocycle
$$(0, \ldots ,0, \frac{f_0}{f_i},0, \ldots, 0 , -\frac{f_0}{f_j},0 ,\ldots, 0)
\in \Gamma(D(f_i f_j), \shR)\, .$$
The other statements follow from
\ref{torsorprojective} and \ref{forcingbundletorsortorsor}
\qed

\begin{definition}
We call the short exact sequence in \ref{forcingsequence1} the
{\em forcing sequence} and we call $\PP(V')$ the {\em projective bundle}
and $\PP(V)$ the {\em forcing projective subbundle}
or the {\em forcing divisor} associated to the elements $f_1, \ldots, f_n;f_0$.
\end{definition}

\begin{remark}
Since the forcing sequence is the extension corresponding
to $\delta (f_0)$ under the presenting sequence it
follows that the presenting sequence itself
is the forcing sequence corresponding to $f_0=1$.
\end{remark}

We may characterize the containment in the solid closure
in the following way.

\begin{corollary}
Let $R$ be a normal excellent domain, let
$I=(f_1, \ldots ,f_n)$ denote an ideal which is primary
to a maximal ideal $\fom$ of height $d \geq 2$.
Let $f_0 \in R$ denote another element
and let $V$ and $V'$ denote the relation bundles on $U=D(\fom)$
according to {\rm \ref{forcingsequence1}}.
Then the following are equivalent.

\renewcommand{\labelenumi}{(\roman{enumi})}
\begin{enumerate}

\item
$f_0 \in (f_1, \ldots ,f_n)^\soclo $

\item
The cohomological dimension of
$\PP(V') -\PP(V)$ is $d-1$.

\item
The cohomological dimension of the geometric torsor
defined by $\delta(f_0) \in H^1(U, \shR)$ is
$d-1$.

\end{enumerate}
\end{corollary}

\proof
This follows from \ref{solidcd},
\ref{forcingbundletorsortorsor} and \ref{forcingsequence1}.
\qed

\medskip
We gather together some characterizations of
$f_0 \in (f_1, \ldots,f_n)$ in terms of
the geometric objects which we consider.

\begin{lemma}
\label{trivial}
Let $R$ be a commutative ring and let
$f_0,f_1, \ldots,f_n \in R$ be elements,
$U=D(f_1,\ldots,f_n)$.
Let $A=R[T_1, \ldots, T_n]/(f_1T_1+ \ldots +f_nT_n+f_0)$
be the forcing algebra. Then the following are equivalent.

\renewcommand{\labelenumi}{(\roman{enumi})}
\begin{enumerate}
\item 
$f_0 \in (f_1,\ldots,f_n)$.

\item
There exists a section $\Spec \, R \ra \Spec \, A$.

\item
The forcing algebra $A$ is isomorphic to
the algebra of relations
$$R[T_1, \ldots,T_n]/ (f_1T_1+ \ldots +f_nT_n) \, .$$
Suppose further that $R =\Gamma(U,\O_X)$ {\rm(}e.g. if $R$ is normal
and ${\rm ht}\, I \geq 2 ${\rm )}.
Then these statements are also equivalent with

\item
The affine-linear bundle $B=\Spec \, A|_U$ has a section over $U$.

\item
There exists a section $U \ra \PP(V')$ which does not meet
$\PP(V)$.

\item
The forcing sequence splits.

\item
The element $f_0$ defines the zero element in $H^1(U, \shR)$.

\end{enumerate}
\end{lemma}

\proof
Suppose that (i) holds, say $-f_0=\sum a_if_i$.
Then $T_i \mapsto T_i+a_i$ is well defined and gives the isomorphism in (iii).
On the other hand, a relation algebra has the zero section, thus the first
three statements are equivalent.

(ii) $\Rightarrow$ (iv) is a restriction, and (iv) $ \Rightarrow $ (ii)
is true under the additional assumption. (iv) and (v) are equivalent
due to \ref{forcingsequence1}.

(i) gives also directly a section for $V' \ra \AA_U \ra 0$, thus we
get (vi), which is equivalent with (vii).
If the sequence splits, then $V' = V \oplus \AA$ on $U$ and the complement
of $\PP(V)$ is the vector bundle $V$, which has the zero section.
\qed

\begin{remark}
What can we say about the topology and the differential structure
of an affine-linear bundle arising from a forcing problem?
Let $K= \CC$ and let $R$ denote a normal $\CC$-domain of finite
type and let $f_1, \ldots, f_n \in R$ be elements generating an
ideal primary to a maximal ideal $\fom$. Suppose that
$D(\fom)$ is smooth and let $M=D(\fom)$, considered as a real differential
manifold.

Then the affine-linear bundle $(\Spec\, A)|_M$, where $A$ denotes the forcing
algebra for $f_1, \ldots , f_n$ and another element $f$
has always a (differential) section and is therefore isomorphic to the
relation bundle.
In particular, the topological and the differential properties
of these affine-linear bundles do not depend on $f$.

This is an easy application of partition of unity,
see \cite[Satz 5.A.3]{storchwiebevier}.
There exists a partition of unity $h_i$ for the open subsets
$M_i =D(f_i)$, i.e. $\sum_i h_i =1$ and
$\overline{\Supp h_i} \subseteq M_i$. Then the function $q_i =h_i/f_i$
may be extended (by zero outside $M_i$) to a function on $M$, since
$q_i |_{M_i \cap (M- \overline{\Supp h_i})} =0$.
Therefore $\sum_i q_if_i =\sum_i h_i =1$ and
$\sum_i (fq_i) f_i =f$. This gives the section.

Another way to prove this is by using that
a short exact sequence of differential vector bundles always splits,
but this fact also rests upon partition of unity.
\end{remark}

\begin{remark}
We adopt the situation of the previous remark, but now
we consider $M=D(\fom)$ and the torsor defined by forcing data
as complex manifolds.
Then the complex structure and the complex coherent cohomology
of the torsor depend on $f$.
The complex cohomological dimension
however may differ from the algebraic cohomological dimension,
see the example below (\ref{steintwo}) of a Stein but non-affine torsor.
\end{remark}

\smallskip

\subsection{The graded case: relation bundles on projective varieties}
\label{gradedsection}

\markright{The graded case: relation bundles on projective varieties}

\

\bigskip
We have shown in the last section how to realize the affine-linear bundle
coming from forcing data inside a projective bundle.
In order to use the full strength of the methods of projective geometry
such as intersection theory
it is necessary that also the base scheme is projective.
Hence we shall consider now the graded situation
and repeat the constructions of the last sections on the
corresponding projective varieties.

Let $K$ be a field and let $R$ be a
standard $\NN$-graded $K$-algebra, i.e. $R_0=K$ and $R$ is generated
by finitely many elements of first degree.
The projective scheme $Y= \Proj\, R$
is a projective variety over $K$.
It carries the very ample invertible sheaf
$\O_Y(1)$.
We denote by $\AA_Y(1)$ the line bundle on $Y$
with sheaf of linear forms $\O_Y(1)$ and sheaf of sections
$\O_Y(-1)$.
Hence $\AA_Y(m) = \Spec \, S(\O_Y(m)) = \Spec \oplus_{k \geq 0} \O_Y(km)$.
Another geometric realization is given by
$$\AA_Y(m) := D_+(R_+) \subset \Proj\, R[T],\, \,  \deg \, (T) =  -m \, ,$$
where the inclusion $R \subset R[T]$ induces the bundle
projection $\AA_Y(1) \ra Y$.
We will denote by $H$ the hyperplane divisor class of
$\O_Y(1)$ on $Y$.

Let $f_1, \ldots ,f_n$ be homogeneous elements of $R$ of degree $d_i$.
We say that the $f_i$ are primary if they are
$R_+$-primary, i.e. $D(R_+) = D(f_1, \ldots,f_n)$.
We may find degrees $e_i$ (possibly negative) for $T_i$ such that
the polynomials
$\, \sum_{i=1}^nf_iT_i \,$, $\, \sum_{i=0}^n f_iT_i \,$ and
$\, \sum_{i=1}^n f_iT_i+f_0$ are homogeneous
(for the last polynomial $e_i=d_0-d_i$ is the only choice).

Let
$A=R[(T_0),T_1, \ldots,T_n]/(P)$, where $P$ is one of these polynomials.
Then $A$ is also graded (but not standard-graded in general)
and we have the following commutative diagram.
$$
\begin{CD}
\Spec \, A \supset D(R_+A) \, \,   @>  >> \, \,  D_+(R_+ A) \subset \Proj\, A \\
\, \,  @VVV     @VVV \, \, \\
\Spec \, R  \supset D(R_+) \, \,  @>  >> \Proj \, R \, \, \, 
\end{CD}
$$
These algebras $A$ may have negative degrees,
but $\Proj \, A$ can be defined as well,
see \cite{brennerschroer}. The open subset
$D_+(R_+A) \subset \Proj \, A$ is the same as
$D_+(R_+A) \subseteq \Proj \, A_{\geq 0}$.

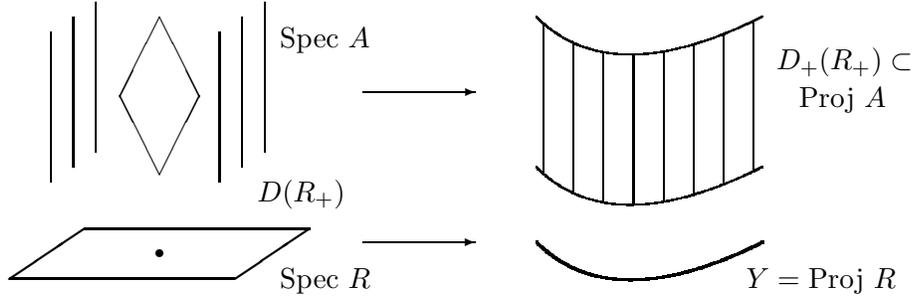
\newsavebox{\raute}

\savebox{\raute}(3,3)[bl]{
\put(0, 0){\line(1,2){0.53}}
\put(0,0){\line(-1,2){0.53}}
\put(0,2.1){\line(1,-2){0.53}}
\put(0,2.1){\line(-1,-2){0.53}}}

\newsavebox{\viereckmitte}

\savebox{\viereckmitte}{
\put(2., .36){\circle*{0.1}}
\thicklines
\put(.0, .0){\usebox{\viereckdick}}}

\newsavebox{\buendelexzept}

\savebox{\buendelexzept}(5,4)[bl]
{\put(.0,.0){\usebox{\viereckmitte}}
\put(2.,1.4){\usebox{\raute}}
\put(0.55,1.3){\usebox{\strichedreistufe}}
\put(2.8, 1.3){\usebox{\strichedreistufe}}}

\setlength{\unitlength}{1cm}

\begin{picture}(12,5)

\put(1,.5){\usebox{\buendelexzept}}

\put(4.3, 1.55){$D(R_+)$}
\put(4.6,.4){$\Spec\, R$}
\put(4.6,3.6){$ \Spec\, A$}

\put(5.7, 3){$\vector(1,0){1.5}$}
\put(5.7, 1){$\vector(1,0){1.5}$}

\put(8,0){\usebox{\ruledsurface}}

\put(10.8, .4) {$ Y= \Proj \,R$}
\put(11.2, 3.3) {$ D_+(R_+) \subset$}
\put(11.5,2.8){$\Proj \, A$}
\end{picture}

\savebox{\buendelexzept}{}

\savebox{\raute}{}

\bigskip
The relation module for the homogeneous elements
$f_1, \ldots ,f_n$ is a graded $R$-module,
where the graded piece is given by
$\Rel (f_1, \ldots ,f_n)_m =$
$$ \{ (r_1, \ldots , r_n) \in R^n
\mbox{ homogeneous }:
\sum_i r_if_i=0  \mbox{ and } \deg (r_i) + \deg (f_i) =m \} \, .$$
We denote the corresponding sheaf of relations on
$Y= \Proj \, R$ by $\shR(m)$.
We call this the {\em sheaf of relations of total degree} $m$.

\begin{proposition}
\label{buendelalsproj}
Let $R$ be a standard-graded $K$-algebra and let
$f_1,\ldots,f_n$ be homogeneous $R_+$-primary elements.
Let $d_i = \deg \, (f_i) $ and
fix a number $m \in \ZZ$ and set $e_i =m-d_i$.
Then the following hold.

\renewcommand{\labelenumi}{(\roman{enumi})}
\begin{enumerate}

\item
Set $\deg \, (T_i) =e_i$. Then
$$ \Proj\, R[T_1,\ldots,T_n]/(\sum_{i=1}^n f_iT_i) \, \, \supset \, \, D_+(R_+)
\longrightarrow \Proj\, R$$
is a vector bundle $V(-m)$ of rank $n-1$ over $Y=\Proj\, R$.
Its sheaf of sections is the sheaf of relations $\shR(m)$ of total degree $m$.

\item
For this vector bundle $V(-m)$ we have the exact sequence
of vector bundles
$$0 \longrightarrow V(-m) \longrightarrow
\AA_Y(-e_1) \times_Y \ldots \times_Y \AA_Y(-e_n)
\stackrel{\sum f_i}{\lra} \AA_Y(-m) \longrightarrow 0 \, .$$

\item
We have $\Det \, V(-m) \cong \AA_Y(-\sum_{i=1}^n e_i +m)
= \AA_Y(\sum_{i=1}^n d_i -(n-1)m) $.

\item
We have $V (-m')= V(-m) \otimes \AA_Y(m-m')$.

\item
The projective bundle
$\PP(V(-m))$ does not depend on the chosen degree $m$.
For the relatively very ample sheaf
$\O_{\PP(V(-m))}(1) $ on $\PP(V(-m))$ we have
$$j^* \O_{\PP(V (-m'))} (1) =  \O_{\PP(V(-m))} (1) \otimes \pi^*\O_Y(m-m')\,,$$
where
$j: \PP(V(-m)) \ra \PP(V(-m) \otimes  \AA_Y(m-m'))$
is the isomorphism and $\pi: \PP(V(-m)) \ra Y$ is the projection.

\end{enumerate}
\end{proposition}

\proof
(i) and (ii).
First note that the natural mapping ($\deg\, T_i=e_i$)
$$\Proj \, R[T_1, \ldots,T_n] \supseteq D_+(R_+) \lra
\AA_Y(-e_1) \times_Y \ldots  \times_Y \AA_Y(-e_n)$$
is an isomorphism.
The ring homomorphism
$R[T] \ra R[T_1,\ldots, T_n], \, T\mapsto \sum f_iT_i$
is homogeneous for $\deg \, (T) =m$.
This gives the epimorphism of vector bundles, since the $D_+(f_i)$
cover $Y$.
Its kernel is given by
$ D_+(R_+) \subset  \Proj \,  R[T_1,\ldots, T_n]/(\sum f_iT_i)$,
thus this is also a vector bundle on $Y$.

(iii) follows from (ii).
If we
tensor the exact sequence for $V(-m)$ with $\AA_Y(m-m')$ we
get the sequence for $V (-m')$, hence (iv) follows.

(v)
$\PP(V)$ does not change when $V$ is tensored
with a line bundle. The relatively very ample sheaves behave like stated due
to \cite[Proposition 4.1.4]{EGAII}.
\qed

\begin{remark}
We call the short exact sequence in \ref{buendelalsproj}
again the presenting sequence.
It yields for the sheaf of relations of total degree $m$
the short exact sequence
$$ 0 \lra \shR (m)  \lra \oplus_i \O_Y(m-d_i)
\stackrel{\sum f_i}{\lra} \O_Y(m) \lra 0 \, .$$
The sheaf of linear forms of total degree $m$
is the dual sheaf $\shF(-m)= \shR(m)^\dual$,
thus $V(-m)= \Spec \, S(\shF(-m))$ and $\PP(V)=\PP(\shF)$.
The corresponding sequence is
$$0 \lra \O_Y(-m) \stackrel{f_1, \ldots,f_n}{\lra}
\oplus_i \O_Y(d_i-m) \lra \shF (-m) \lra 0 \, .$$
The most important choice for $m$ will be $m=d_0$, where $d_0$
is the degree of another homogeneous element $f_0$.
\end{remark}

\begin{remark}
\label{chern}
The sequence in \ref{buendelalsproj} (ii) allows us
to compute inductively the
Chern classes of the vector bundles $V(-m)$ (or of its sheaf of linear forms
$\shF(-m)$).
For the Chern polynomial $c_t(V(-m))= \sum_{i}c_i(V(-m))t^{i}$
we get the relation
(let $H$ denote the hyperplane section of $Y$)
$$c_t(V(-m))(1-mHt)=(1-e_1Ht) \cdots (1-e_nHt) \, .$$
This yields (set $e_i=m-d_i$)
$c_0(V(-m))=1$,
$$ c_1(V(-m))= (-e_1- \ldots -e_n+m)H =(d_1+ \ldots +d_n-(n-1)m)H,\, $$
$c_2(V(-m))=(\sum_{i_1,i_2} e_{i_1}e_{i_2}-(e_1+ \ldots +e_n-m)m) H.H $ etc.
\end{remark}

\smallskip

\subsection{The forcing sequence on projective varieties}

\markright{The forcing sequence on projective varieties}

\label{gradedforcingsection}

\

\bigskip
Let $R$ denote a normal standard-graded $K$-algebra
and let $f_1, \ldots ,f_n$ denote $R_+$-primary homogeneous elements.
We discuss now the forcing sequence defined by another homogeneous element $f_0$
on the projective variety $\Proj \, R$ and how the containment
in the solid closure is expressed in terms of
projective bundles and subbundles.

\begin{proposition}
\label{forcingsequence2}
Let $R$ be a standard-graded $K$-algebra, let
$f_1,\ldots,f_n$ be homogeneous $R_+$-primary elements and let
$f_0 \in R$ be also homogeneous.
Let $d_i = \deg \, (f_i) $ and fix a number $m \in \ZZ$.
Let $\deg \, (T_i) =e_i =m-d_i$.
Let
$$V(-m)=D_+(R_+) \subset \Proj\, R[T_1,\ldots,T_n]/(\sum_{i=1}^n f_iT_i) 
\, \, \mbox{ and }$$
$$V'(-m)= D_+(R_+) \subset \Proj\, R[T_0,\ldots,T_n]/(\sum_{i=0}^n f_iT_i)$$
be the vector bundles on $Y=\Proj\, R$ due to {\rm \ref{buendelalsproj}}.
Then the following hold.

\renewcommand{\labelenumi}{(\roman{enumi})}
\begin{enumerate}

\item 
There is an exact sequence of vector bundles on $Y$,
$$0 \longrightarrow V(-m) \longrightarrow V'(-m)
\stackrel{T_0}{\longrightarrow} \AA_Y(d_0-m) \longrightarrow 0 \, .$$

\item
The embedding $\PP(V) \hookrightarrow \PP(V')$ does not depend
on the degree $m$ {\rm (}and we skip the index $m$ inside $\PP(V)${\rm )}.
The complement of $\PP(V)$ is
$$\PP(V') -\PP(V) \cong D_+(R_+)
\subseteq \Proj \, R[T_1, \ldots,T_n]/(f_1T_1+ \ldots +f_nT_n+f_0) \, .$$

\item
Let $E$ be the Weil divisor {\rm (}the hyperplane section{\rm)} on $\PP(V')$
corresponding to the relatively very ample
invertible sheaf $\O_{\PP(V')} (1)$ {\rm (}depending on the degree{\rm )}.
Then we have the linear equivalence of divisors
$\PP(V) \sim E +(m-d_0) \pi^*H$,
where $H$ is the hyperplane section of $\, Y$.
If $m=d_0$, then the forcing divisor $\PP(V)$ is a hyperplane section.

\item
The normal bundle for
$\PP(V) \hookrightarrow \PP(V')$ on $\PP(V)$ is
$\AA_{\PP(V)}(-1) \otimes \pi^*\AA_Y(d_0-m)$.

\end{enumerate}
\end{proposition}
\proof
(i).
The homogeneous ring homomorphisms
$$R[T_0] \, \lra \, R[T_0, \ldots, T_n]/(\sum_{i=0}^n f_iT_i)
\, \, \, (\deg \,(T_0)= e_0)  \, \mbox{ and }$$
$$ R[T_0, \ldots, T_n]/(\sum_{i=0}^n f_iT_i) \, \lra \,
R[T_1, \ldots, T_n]/(\sum_{i=1}^n f_iT_i), \, \, \, T_0 \longmapsto 0 $$
induce the morphisms on $D_+(R_+)$.
The exactness is clear on $D(f_i)$, $i=1, \ldots,n$, and they cover
$D_+(R_+)$.

The first statement in (ii) is clear due to
\ref{buendelalsproj}(iv), thus we may assume that $m=d_0$.
The homogeneous ring homomorphism
$R[T_0, \ldots, T_n]/(\sum_{i=0}^n f_iT_i) \ra
R[T_0, \ldots, T_n]/(\sum_{i=1}^n f_iT_i +f_0)$ where $T_0 \mapsto 1$
yields the closed embedding
$$ \Proj \,  R[T_0, \ldots, T_n]/(\sum_{i=1}^n f_iT_i +f_0) \supseteq
D_+(R_+) \hookrightarrow V'(-m) \,, $$
where the image is given by $T_0 =1$.
But this closed subset $V_+(T_0-1) \subseteq V'(-m)$
is isomorphic to $\PP(V') -\PP(V)$
under the cone mapping $V'(-m) \dashrightarrow \PP(V')$.

(iii).
The mapping
$T_0: V'(-m) \ra \AA_Y(d_0-m)$ yields via the tautological
morphism $\AA_{\PP(V'(-m))}(1) \ra V'(-m)$ a morphism of line bundles on
$\PP(V'(-m))$,
$$\AA_{\PP(V'(-m))}(1) \ra \pi^*\AA_Y(d_0-m) \, .$$
This corresponds to a section in the line bundle
$\AA_{\PP(V'(-m))}(-1) \otimes \pi^*\AA_Y(d_0-m)$
with zero set $\PP(V(-m))$.
Thus $\PP(V(-m)) \sim E +(m-d_0) \pi^*H$.

(iv).
Let $i:\PP(V) \hookrightarrow \PP(V')$ be the inclusion.
Since $\PP(V)$ is the zero-set of a section in
$ \AA_{\PP(V'(-m))}(-1) \otimes \pi^*\AA_Y(d_0-m) $, we have that
$$i^*(\AA_{\PP(V')}(-1) \otimes \pi^*\AA_Y(d_0-m))
=\AA_{\PP(V)}(-1) \otimes q^*\AA_Y(d_0-m)$$
is the normal bundle on $\PP(V)$.
\qed

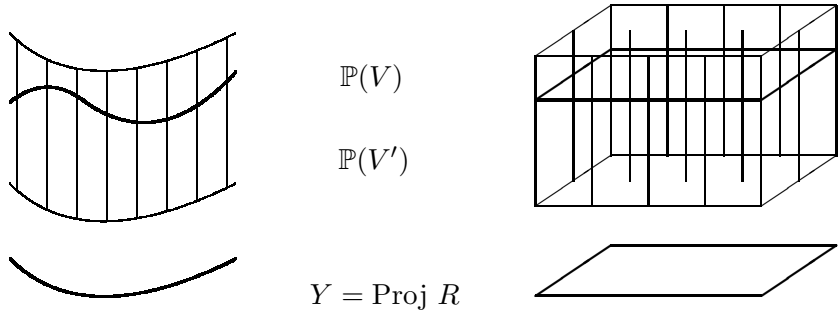
\newsavebox{\viereckbuendel}

\savebox{\viereckbuendel}(5,5)[bl]{
\put(0, 0){\usebox{\viereckdick}}
\put(0, 2.6){\usebox{\viereckdick}}
\put(0, 1.2){\usebox{\viereck}}
\put(0, 3.2){\usebox{\viereck}}
\put(0, 1.2){\usebox{\strichefuenf}}
\put(0.5, 1.533){\usebox{\strichefuenf}}
\put(1, 1.867){\usebox{\strichefuenf}}}

\begin{picture}(12,5)

\put(1,0){\usebox{\ruledsurface}}
\put(1,2.8){\usebox{\krummbogen}}

\put(8, .5){\usebox{\viereckbuendel}}

\put(5,0.4){$Y =\Proj \, R$}
\put(5.38,2.2){$\PP(V')$}
\put(5.4,3.3){$\PP(V)$}

\end{picture}

\savebox{\viereckbuendel}{}

\begin{definition}
\label{forcingdef}
We call the sequence in \ref{forcingsequence2} (i) again
the {\em forcing sequence} and we
denote the situation $\PP(V) \hookrightarrow \PP(V')$
by $\PP(f_1, \ldots, f_n; f_0)$.
This is a projective bundle of rank $n-1$ together
with the forcing divisor $\PP(V)= \PP(f_1, \ldots,f_n)$
over $Y$.
\end{definition}

\begin{remark}
\label{forcingsheaf2}
Corresponding to the forcing sequence of vector bundles in
\ref{forcingsequence2}
we have the exact sequence of the sheaves of relations
$$0 \ra \shR(m) \ra \shR'(m) \ra \O_Y(m-d_0) \ra 0 $$
of total degree $m$.
For $m=d_0$ this extension corresponds to a cohomology
class $c \in H^1(Y, \shR(m))$.
The forcing sequence for the sheaves of linear forms is
$$ 0 \ra \O_Y(d_0-m)  \ra \shF'(-m) \ra \shF (-m)  \ra 0 \, .$$
\end{remark}

\begin{remark}
\label{normalbundle}
The normal sheaf on $\PP(V) \subset \PP(V')$, i.e.
the sheaf of sections in the normal bundle,
is (by \ref{buendelalsproj}(v))
$$\shN=\O_{\PP(V(-m))}(1) \otimes \pi^* \O_Y(m -d_0)
\cong \O_{\PP(V(-d_0))} (1)\, .$$
Its global sections are
$\Gamma(\PP(V), \O_{\PP(V(-d_0))} (1))= \Gamma(Y, \shF(-d_0))$.
Thus the normal sheaf does not depend on $f_0$,
only on its degree $d_0= \deg (f_0)$.
\end{remark}

The next results show that we can express the properties which are of interest
from the tight closure point of view in terms of the projective bundles
on $Y$.

\begin{lemma}
\label{trivialtwo}
In the situation of {\rm \ref{forcingsequence2}} the following
are equivalent.

\renewcommand{\labelenumi}{(\roman{enumi})}
\begin{enumerate}

\item 
$f_0 \in (f_1, \ldots,f_n)$.

\item
There is a section $Y \ra \PP(V')$ disjoined to $\PP(V) \subset \PP(V')$.

\item
The forcing sequence
$0 \ra V(-m) \ra V'(-m) \ra \AA_Y(d_0-m) \ra 0 $ splits.

\item
Let $m=d_0$.
The corresponding cohomological class in $H^1(Y,\shR (m))$ vanishes.
\end{enumerate}
\end{lemma}

\proof
Suppose that (i) holds and write $-f_0 =\sum_{i=1}^n a_if_i$,
where the $a_i$ are homogeneous. Set $m = d_0$.
The $a_i$ define a homogeneous mapping
$$R[T_0,T_1,\ldots,T_n]/(\sum_{i=0}^n f_iT_i) \lra R
\mbox{ by }
T_0 \ra 1,\, \,  T_i \ra a_i \, .$$
The corresponding mapping
$Y \ra V'(-m)$ induces $Y \ra \PP(V')$ and its image
is disjoint to $\PP(V)$.

Suppose that (ii) holds.
A section in $\PP(V'(-m))$ corresponds to a line bundle $L$ on $Y$
and an embedding of vector bundles $L \hookrightarrow V'(-m)$,
see \cite[Proposition 7.12]{haralg}.
Since the section is disjoined to $\PP(V)$, the
morphism $V(-m) \oplus L \ra V'(-m)$ is an isomorphism,
hence the sequence splits.

(iii).
The splitting yields a section
$\AA_Y(d_0-m) \ra V'(-m)$ and this means a homogeneous mapping
$ R[T_0,T_1,\ldots,T_n]/(\sum_{i=0}^n f_iT_i) \ra R[T_0]$.
For $T_0=1$ we get a solution for (i).
(iii) and (iv) are equivalent. 
\qed

\begin{example}
\label{zero}
Let $R$ denote a standard-graded $K$-algebra
and let $f_1, \ldots,f_n $ be homogeneous primary elements of degree $d_i$.
Let $f_0=0$.
Then
$$ R[T_0, \ldots, T_n]/(\sum_{i=0}^n f_iT_i)
= R[T_1, \ldots, T_n]/(\sum_{i=1}^n f_iT_i) [T_0]$$
and we have the splitting forcing sequence
$$0 \lra V(-m) \lra V(-m) \oplus \AA_Y(d_0-m) \lra \AA_Y(d_0-m) \lra 0 \, .$$ 
Then $V(-d_0) \cong \PP(V') - \PP(V)$ and $\PP(V')$ is just the projective
closure of $V(-d_0)$.
\end{example}

\begin{example}
\label{euler}
Let $R =K[X_0, \ldots, X_d]$ denote the polynomial ring
and consider the forcing data $X_0, \ldots , X_d;1$.
Then the forcing sequence (which is here also the presenting
sequence) for the sheaf of relations
on the projective space is just
$$ 0 \lra \shR(m) \lra \oplus_{d+1} \O_{\PP^d}(m-1) \lra \O_{\PP^d}(m) \lra 0 \, .$$
This is for $m=0$ nothing but the Euler-sequence for
the cotangent sheaf. Hence $\shR(0) = \shT^\dual$
and $\shF(0)= \shT$ is the (ample) tangent sheaf on $\PP^d$.
\end{example}

\begin{proposition}
\label{solidgraded}
Let $R$ be a normal standard-graded $K$-algebra of dimension $d \geq 2$,
let $f_1, \ldots,f_n \in R$ be $R_+$-primary homogeneous elements
and let $f_0$ be another homogeneous element.
Let $V$ and $V'$ be as in {\rm \ref{forcingsequence2}}.
Then the following are equivalent.
\renewcommand{\labelenumi}{(\roman{enumi})}
\begin{enumerate}
\item
$f_0 \in (f_1, \ldots,f_n)^\soclo$.

\item
The cohomological dimension of
$\PP(V') - \PP(V)$ is $d -1 = \dim \, Y$.
\end{enumerate}
\end{proposition}
\proof
We have $\PP(V') -\PP(V) \cong D_+(R_+) \subseteq
\Proj \, R[T_1, \ldots,T_n]/(f_1T_1 + \ldots +f_nT_n+f_0)$.
The cohomological dimensions of $D_+(\foa)$ and of $D(\foa)$
are the same due to \ref{cdproperties}(iii).
Hence the cohomological dimensions of $\PP(V') - \PP(V)$ and of
the forcing affine-linear bundle
$D(R_+) \subseteq \Spec R[T_1, \ldots,T_n]/(f_1T_1 + \ldots +f_nT_n+f_0)$
are the same, and the result follows from \ref{solidcd}(ii).
\qed

\medskip
In the presence of a grading it is more natural to compare tight closure
with the graded plus closure, which is defined in the following way.

\begin{definition}
Let $R$ denote an $\NN$-graded Noetherian $K$-domain and let
$I=(f_1, \ldots, f_n)$ denote a homogeneous ideal.
Then
$$I^{\gr} := \{f \in R : f \in IR' \mbox{ for a finite graded extension }
R \subseteq R' \} \, .$$
\end{definition}

$\varphi: R \ra R'$ graded means that there exists $a \in \NN$
such that $\varphi(R_d) \subseteq R'_{ad}$.
If $R \subseteq R'$ is graded and finite, then this induces
a morphism $\Proj R' \ra \Proj R$.
The Frobenius morphism of a graded ring over a field of
positive characteristic is therefore always graded.

If $R \subseteq R'$ is a finite extension such that
$f_0 \in IR'$, then this holds also in $R'' = \sum_k R'_{ak}$,
and by regrading we may find also a finite graded extension
$R \subseteq R''$ such that $f \in R''$ holds, but without any degree shift.
We may characterize the containment
$f_0 \in (f_1, \ldots, f_n)^{\gr}$ in our geometric setting in the following
way.

\begin{proposition}
\label{plus}
Let $R$ be a normal standard-graded $K$-algebra of dimension $d \geq 2$,
let $f_1, \ldots,f_n \in R$ be $R_+$-primary homogeneous elements
and let $f_0$ be another homogeneous element.
Let $V$ and $V'$ be as in {\rm \ref{forcingsequence2}}.
Then the following are equivalent.
\renewcommand{\labelenumi}{(\roman{enumi})}
\begin{enumerate}

\item
$f_0 \in (f_1, \ldots,f_n)^{\gr}$.

\item
There exists a finite surjective morphism $g: Y' \ra Y$ such that
the pull back $g^*\PP(V')= \PP(V') \times_Y Y'$
has a section not meeting $g^*\PP(V)= \PP(V) \times_Y Y'$.

\item
There exists a closed subvariety $Y'' \subset \PP(V')$
not intersecting $\PP(V)$, finite and surjective over $Y$.

\item
There exists a closed subvariety $Y'' \subset \PP(V')$
not intersecting $\PP(V)$ of dimension $d-1 = \dim\, Y$.

\end{enumerate}
\end{proposition}
\proof
First note that the pull-back of a relation bundle $V(-m)$
on $Y$ for homogeneous elements $f_1, \ldots, f_n$
yields the relation bundle $V(-am)$ on $Y'=\Proj \, R'$
for the elements considered in $R'$.

(i) $\Rightarrow $ (ii).
Suppose first that $R \subseteq R'$ is a finite extension
such that $f_0 \in (f_1, \ldots, f_n)R'$ holds.
Then on $Y'= \Proj R' $ there exists due to \ref{trivialtwo}
a section
$s: Y'  \ra \PP(V'_{Y'})$ which
does not meet $\PP(V_{Y'})$, where $V_{Y'}$ denotes the
relation bundle on $Y'$ corresponding to the elements considered
in $R'$.
But $ \PP(V_{Y'})= g^*\PP(V) = \PP(V) \times_Y Y'$,
where $g:Y' \ra Y$ is the finite morphism corresponding
to the ring inclusion $R \subseteq R'$.

(ii) $\Rightarrow $ (i).
Let $g:Y' \ra Y$ be such a morphism.
Then $\O_{Y'}(1)= g^* \O_Y(1)$ is ample
and we consider $R' = \oplus _{k \geq 0} \Gamma(Y' , g^* \O_Y(k))$, which is
finitely generated and finite over $R'$.
We have $Y' \cong \Proj R'$ and we have the commutative diagram
$$
\begin{CD}
D(R_+')   @>  >> D(R_+)  \\
  @VVV     @VVV   \\
Y'    @> g >> Y \, \, .
\end{CD}
$$
Let $B= \PP(V') -\PP(V)$.
The existence of a section 
$s: Y'  \ra \PP(V') \times_Y Y'$ not meeting
$\PP(V) \times_Y Y'$
implies that $g^*(B)$ is trivial.
Then also the pull-back of $B$ to $D(R_+')$ is trivial, but this is the
same as the pull-back of the forcing bundle
$\Spec \, R(T_1, \ldots ,T_n]/(\sum f_iT_i +f_0)|_{D(R_+)}$.
It follows that $f_0 \in (f_1, \ldots, f_n)R''$,
where $R''$ is the normalization of $R'$.

Suppose that (ii) holds. Then the image of the section gives the closed
subvariety $Y''$ finite over $Y$. This gives (iii) and then (iv).
Suppose that (iv) holds.
The mapping $Y'' \hookrightarrow \PP(V') \ra Y$ is projective
and the fibers are zero-dimensional, since $\PP(V_y) \subset \PP(V'_y)$
meets every curve, but $Y'' \cap \PP(V) = \emptyset $.
Hence this mapping is finite and due to the assumption on the dimension
it is surjective.
So suppose that (iii) holds.
The mapping
$g : Y'=Y'' \stackrel{i}{\hookrightarrow} \PP(V') \ra Y$
is finite and surjective,
and the image of the section
$ i \times id_{Y'}  :Y'  \ra \PP(V') \times_Y Y' =g^* \PP(V') $
is disjoined to $g^* \PP(V)= \PP(V) \times_Y Y'$.
\qed

\bigskip

\subsection{Ample and basepoint free forcing divisors}
\label{amplesection}

\markright{Ample and basepoint free forcing divisors}

\

\bigskip
We fix the following situation.

\begin{situation}
\label{gradedsituation}
Let $K$ be an algebraically closed field and let $R$ be a normal
standard-graded $K$-domain of dimension $d$.
Let $f_1, \ldots,f_n$ be homogeneous $R_+$-primary elements 
of degree $d_1, \ldots ,d_n$ and
let $f_0$ be another homogeneous element of degree $d_0$.

Let $V(-m)$ and $V'(-m)$ be the relation bundles on $Y=\Proj \,R$
for $f_1, \ldots, f_n$ and for $f_1, \ldots , f_n,f_0$
of total degree $m$ as described in \ref{forcingsequence2}
with sheaf of sections $\shR(m)$ and $\shR'(m)$ and with sheaf of linear forms
$\shF(-m)$ and $\shF'(-m)$.
Let $Z=\PP(V) \subset \PP(V')$ denote the forcing divisor.
\end{situation}

In this section we examine basic properties of this forcing divisor
$Z=\PP(V) \subset \PP(V')$ corresponding to
homogeneous forcing data $f_1, \ldots,f_n;f_0 \in R$.
When is the forcing divisor
$Z$ ample and when is $Z$ basepoint free?

For $m=d_0$ the forcing divisor is a hyperplane section of
$\O_{\PP(V')}(1)$, and the ampleness of this invertible sheaf is by
definition the ampleness of the locally free sheaf
$\shF'(-d_0)= \pi_*\O_{\PP(V')}(1)$, see \cite[III, \S 1]{haramp} and
chapter \ref{slopechapter} for further ampleness criteria for vector bundles
and applications to tight closure problems.
The following proposition shows that the ampleness property is interesting
mainly in dimension two.

\begin{proposition}
\label{as}
Suppose the situation and notation of {\rm \ref{gradedsituation}}.
Then the following hold.

\renewcommand{\labelenumi}{(\roman{enumi})}
\begin{enumerate}

\item
Suppose that $f_0$ is a unit and that $d_i \geq 1$ for $i=1, \ldots,n$.
Then $Z$ is ample.

\item
If $f_0$ is not a unit, then the cohomological dimension
$cd \, (\PP(V')- Z) \geq d-2 $.

\item
If $f_0$ is not a unit and $d \geq 3$, then $Z$ is not ample.

\end{enumerate}
\end{proposition}
\proof
(i).
We may assume that $f_0=1$.
Then
$$V'(-m) =D_+(R_+) \subset
\Proj\, R[T_0,T_1, \ldots, T_n]/(\sum_{i=1}^nf_iT_i+T_0)
\cong \Proj\, R[T_1,\ldots,T_n] \,,$$
where
$\deg\, T_i=e_i=m-d_i$.
Thus $V'(-m) \cong \AA_Y(d_1-m) \times_Y \ldots \times_Y \AA_Y(d_n-m)$.
For $m=0$ we see that $\shF'(0)$ is a sum of ample invertible sheaves,
hence $\shF'(0)$ is ample due to \cite[III, Corollary 1.8]{haramp}.

(ii).
For $d=0,1$ there is nothing to show, so suppose that $d \geq 2$.
The zero set $V_+ (f_0) \subset Y$ is a closed subset of
dimension $\geq d-2$. There exists a section
$s: V_+(f_0) \ra \PP(V')$ which does not meet $Z$.
Hence $\PP(V')-Z$ contains the
projective subvariety $s(V_+(f_0))$ of
dimension $d-2$, thus the inequality holds for the cohomological dimension
due to \ref{cdproperties}(ii) and (iv).

(iii).
Due to (ii) the complement of $Z$ cannot be affine (it contains
projective curves), hence $Z$ is not ample.
\qed

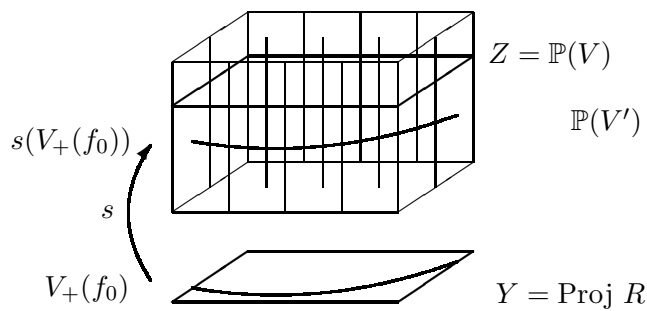
\savebox{\viereckbuendel}(5,5)[bl]{
\put(0, 0){\usebox{\viereckdick}}
\put(0, 2.6){\usebox{\viereckdick}}
\put(0, 1.2){\usebox{\viereck}}
\put(0, 3.2){\usebox{\viereck}}
\put(0, 1.2){\usebox{\strichefuenf}}
\put(0.5, 1.533){\usebox{\strichefuenf}}
\put(1, 1.867){\usebox{\strichefuenf}}}

\newsavebox{\divisorkurvekomplement}

\savebox{\divisorkurvekomplement}(10,5)[bl]{
\put(3 , 0.5){\usebox{\viereckbuendel}}

\thicklines

\bezier{300}(3.25 ,.7 )(5,.4)(6.77 ,1.05)
\bezier{300}(3.25 , 2.65 )(5,2.35)(6.77 ,3)

\put(1.3,.6){$V_+(f_0) $}
\put(7.3,0.5){$Y =\Proj \, R$}
\put(8.3,2.8){$\PP(V')$}
\put(7.2,3.7){$Z=\PP(V)$}
\put(.85,2.55){$s(V_+(f_0))$}

\put(2.7, 2.6){\vector(2,3){.01}}
\put(2.05,1.63){$s$}
\bezier{300}(2.7, 0.8)(2.1, 1.7)(2.7, 2.6)}

\smallskip

\begin{picture}(12,5)

\put(1.5,0){\usebox{\divisorkurvekomplement}}

\end{picture}

\savebox{\divisorkurvekomplement}{}

\savebox{\viereckbuendel}{}

\smallskip
The forcing divisor $Z$ is {\em basepoint free} if and only if
$\O_{\PP(V')}(1)$ is generated by global sections for $m=d_0$.
This is in particular true if
$\pi_* \O_{\PP(V')}(1) = \shF'(-d_0)$ is generated by global sections.

A divisor $Z$ is called {\em semiample}
(\cite[Definition 2.1.14]{lazarsfeldpositive}
if $aZ$ is basepoint free for some $a \geq 1$.
In this case there exists a (projective) morphism
$\varphi: \PP(V') \ra \PP^N$ such that $aZ = \varphi^{-1}(H)$, where $H$
is a hyperplane section in $\PP^N$.
Then $\PP(V') -Z$ is projective over the affine space
$\PP^N -H$. Schemes which are
proper over an affine scheme are called {\em semiaffine} and were
studied in \cite{goodlandman}.

\begin{lemma}
\label{semiaffine}
Suppose the situation and notation of {\rm\ref{gradedsituation}}.
Suppose that the complement
$\, \PP(V') -\PP(V)$ is semiaffine.
Then $f_0 \in (f_1, \ldots,f_n)^\soclo$ if and only if
$f_0 \in (f_1, \ldots, f_n)^{+{\rm gr}}$.
\end{lemma}
\proof
Due to \cite[Corollary 5.8]{goodlandman} the cohomological dimension of
a semiaffine scheme equals the maximal dimension of a closed
proper subvariety.
Thus $f_0 \in (f_1, \ldots, f_n)^\soclo$ implies via \ref{solidcd} that
there exists a projective subvariety $Y' \subset \PP(V')$
of dimension $ \dim \, Y$ which does not meet $\PP(V)$.
Therefore $f_0 \in (f_1, \ldots,f_n)^{+{\rm gr}}$ due to \ref{plus}.
\qed

\medskip
The condition in the following corollary
is usefull mainly for $\dim \, R=2$.
The ampleness of the pull-back $Z|_Z$ is the same as the ampleness of the
bundle $\shF(-d_0)$.

\begin{corollary}
\label{pullbackample}
Suppose the situation of {\rm \ref{gradedsituation}}
and suppose that the pull back $Z|_Z$ is ample.
Then
$f_0 \in (f_1, \ldots,f_n)^\pasoclo$ if and only if
$f_0 \in (f_1, \ldots,f_n)^{+{\rm gr}}$.
\end{corollary}
\proof
The theorem of Zariski-Fujita (see \cite[Remark 2.1.18]{lazarsfeldpositive})
asserts that
$Z$ is semiample. Then the complement of $Z$ is semiaffine
and the result follows from \ref{semiaffine}.
\qed

\begin{corollary}
\label{globalsection}
Suppose the situation of {\rm \ref{gradedsituation}}.
Suppose that the locally free sheaf $\shF'(-d_0)$
is generated by global sections.
Then $f_0 \in (f_1, \ldots,f_n)^\soclo$ if and only if
$f_0 \in (f_1, \ldots, f_n)^{+{\rm gr}}$.
\end{corollary}
\proof
Since $\shF'(-d_0)$ is generated by global sections we know that
the forcing divisor is basepoint free, hence
$\PP(V') - \PP(V)$ is semiaffine.
Hence the result follows from \ref{semiaffine}.
\qed

\begin{remark}
Note that the results \ref{semiaffine} - \ref{globalsection}
yield in characteristic zero 
in fact the stronger result that
$f_0 \in (f_1, \ldots,f_n)^\soclo$ holds if and only if already
$f_0 \in (f_1, \ldots, f_n)$ holds.
\end{remark}

We recover the theorem of Smith (\ref{exclusionbound})
about the exclusion bound within our geometric setting.
Smith proved her result for tight closure
in \cite[Theorem 2.2]{smithgraded}
using differential operators in positive characteristic.
Our version proves the same result for solid closure.

\begin{corollary}
\label{exclusionboundgeo}
Suppose the situation of {\rm\ref{gradedsituation}}.
Suppose that $\deg (f_0) \leq \min_{i} d_i $.
Then $f_0 \in (f_1, \ldots,f_n)^\soclo$ is only possible if
$f_0 \in (f_1, \ldots, f_n)$.
\end{corollary}
\proof
Set $m=d_0$.
Then $e_i=d_0 -d_i \leq 0$ and we have a surjection
$$\O_Y(-e_1) \oplus \ldots \oplus \O_Y(-e_n) \oplus \O_Y
\lra \shF'(-d_0) \ra 0\, .$$
Since the $\O_Y(a)$ for $a \geq 0$ are generated by global sections,
we have also a surjection $\O_Y^k \ra \shF'(-d_0) \ra 0$.
Therefore $\shF'(-d_0)$ is generated by global sections and we have a closed
embedding $V' \hookrightarrow Y \times \AA^k$.

Suppose that $f_0 \in (f_1, \ldots,f_n)^\soclo$. Then by \ref{globalsection}
we know that there exists a subvariety $Y' \subset \PP(V')$
of dimension $\dim \, Y$ not meeting the forcing divisor $Z$.
We may consider $Y' \subset V_+(T_0-1) \subset V'$, since
$V_+(T_0-1)$ is isomorphic to $\PP(V')-\PP(V)$ via the cone mapping
(see the proof of \ref{torsorprojective}).
All together we get a closed embedding
$Y' \hookrightarrow Y \times \AA^k$. Since $Y'$ is a projective variety,
this factors through $Y \times \{P\}$, where $P \in \AA^k$ is a closed point,
and so $Y' \cong Y \times \{P\} \cong Y$, since $K$ is algebraically closed.
Hence we get a section. 
\qed

\medskip
Even if the forcing divisor is not basepoint free, the existence of
linearly equi\-valent effective divisors has consequences on
the existence of closed subvarieties and hence on the existence of
finite solutions (in the sense of \ref{plus} (iii) or (iv))
for the tight closure problem.
See also Propositions \ref{bigaffine} and \ref{nupositive}.

\begin{proposition}
\label{effectivevertreter}
Suppose the situation and notation of {\rm\ref{gradedsituation}},
and suppose furthermore that
$Y=\Proj\, R$ is a smooth variety.
Suppose that there exists a positive {\rm(}effective $\neq 0${\rm )}
divisor $ L \subset Y$ such that for some $a \geq 1$
the divisor $a \PP(V)- \pi^*L $ is linearly equivalent
to an effective divisor.
Then there exists a linearly equivalent effective divisor
$D \sim a\PP(V)$ with the property that
the cohomological dimension of $\PP(V') -{\rm supp}\,  D$ is smaller
than the {\rm(}cohomological{\rm)} dimension of $Y$.
If $Y' \subseteq \PP(V')$ is finite
and surjective over $Y$ and disjoined to $\PP(V)$
{\rm(}as in {\rm \ref{plus}(iii))},
then $Y'$ must lie on the support of $D$.
\end{proposition}
\proof
Let $a\PP(V)- \pi^*L \sim D' $ be effective, hence $a\PP(V) \sim D=D' +\pi^*L$.
The divisor $D'$ cuts out a hyperplane on every fiber,
hence it is also a projective subbundle.
Since a projective bundle minus a dominant effective divisor is relatively
affine over the base we see that
$\PP(V') - {\rm supp}\,D$ is affine over $Y-{\rm supp}\,L$.
But the cohomological dimension of $Y-{\rm supp}\,L$ is smaller than the
dimension of $Y$ due to \ref{cdproperties}(iv),
hence this is also true for $\PP(V') -{\rm supp}\,D$.

Now suppose that $Y'$ is finite and surjective over $Y$
and $Y' \cap \PP(V) = \emptyset$.
Then we have from intersection theory the identities 
$$0=a Y' . \PP(V)= Y'. (D'+\pi^*L) =Y'.D' + Y'.\pi^*L \, .$$
The second summand is a positive cycle, since $Y'$ dominates $Y$.
Hence $Y'.D'$ cannot be effective and the intersection
of $Y'$ and $D'$ must be improper,
so $Y' \subset {\rm supp}\, D' \subset {\rm supp}\, D$.
\qed

\bigskip

\subsection{The two-dimensional situation}
\label{dimensiontwosection}

\markright{The two-dimensional situation}
\

\bigskip
The geometric interpretation of tight closure developed so far
is valid in every dimension, but it is especially
successful if $R$ is a two-dimensional normal standard-graded domain
and $\Proj \, R$ is a smooth projective curve.
There are a lot of reasons which make life more convenient
in the two-dimensional (and projective one-dimensional) situation.
Let us mention the following items.

\smallskip
The notion of affineness is easier to handle
than the notion of cohomological dimension in general, and there
exists various sufficient and necessary conditions for
affineness which we can use to obtain inclusion
and exclusion results for tight closure (see sections
\ref{slopechapter} and \ref{applications}).

Moreover, in dimension $\geq 3$,
the cohomological dimension characterizes tight closure
only in positive characteristic,
whereas in characteristic zero it characterizes solid closure, which has
only in dimension two all the expected properties from tight closure.

The affineness of an open subset is in the projective setting
intimately related to the notion of an ample divisor, which is fundamental in
every part of algebraic geometry.
In fact a lot of results which we shall obtain in the sequel
about affineness and tight closure rest on results about ampleness.
In higher dimensions the notion of ampleness is not well suited
for the study of tight closure problems.
(A conceivable approach in higher dimensions could be achieved
with the help of $k$-ampleness, where $k = \dim\, R-2$,
see \cite[Chapter 2.1]{beltramettisommese} for this notion.)

The top self intersection number of the forcing divisor
can be easily expressed by the degrees of the homogeneous elements
and gives an important indication for ampleness.

The element $f_0 \in R_m$
gives the cohomology class $c= \delta (f_0) \in H^1(Y, \shR(m))$,
which defines the torsor.
It seems that in higher dimension it is necessary
to look also at the higher relation modules
coming from a locally free resolution
$$ \shG_\bullet \lra \oplus_i \O_Y(m-d_i) \lra \O_Y(m) \lra 0 $$
and at the cohomology classes
$c_i \in H^{i}(Y ,\shR_i(m))$, $i=1, \ldots , \dim \, Y$,
defined by the connecting homomorphism
for the short exact sequences
$ 0 \ra \shR_i  \ra \shG_i \ra \shR_{i-1}  \ra 0$
induced by the free resolution.

\medskip
Let us fix the following situation with which we are mainly concerned
in the following chapters.

\begin{situation}
\label{gradedsituationdimensiontwo}
Let $K$ be an algebraically closed field and let $R$ be a normal
standard-graded $K$-domain of dimension two.
Let $Y = \Proj \, R$ denote the corresponding smooth
projective curve over $K$.
Let $f_1, \ldots,f_n$ be homogeneous $R_+$-primary elements 
of degree $d_1, \ldots ,d_n$ and
let $f_0$ be another homogeneous element of degree $d_0$.
Let $V(-m)$ and $V'(-m)$ be the relation bundles on $Y$
for $f_1, \ldots, f_n$ and for $f_1, \ldots , f_n,f_0$
of total degree $m$ (as described in \ref{forcingsequence2})
with sheaf of sections $\shR(m)$ and $\shR'(m)$ and with sheaf of linear forms
$\shF(-m)$ and $\shF'(-m)$.
Let $Z=\PP(V) \subset \PP(V')$ denote the forcing divisor.
\end{situation}

\begin{theorem}
\label{tightaffinecrit}
Suppose the situation and notation of {\rm \ref{gradedsituationdimensiontwo}}.
Then the following are equivalent.

\renewcommand{\labelenumi}{(\roman{enumi})}
\begin{enumerate}
\item
$f_0 \in (f_1, \ldots,f_n)^\soclo$.

\item
The open subscheme $\PP(V')- \PP(V)$ is not affine.

\end{enumerate}
\end{theorem}

\proof
This follows from \ref{solidgraded}.
\qed

\bigskip

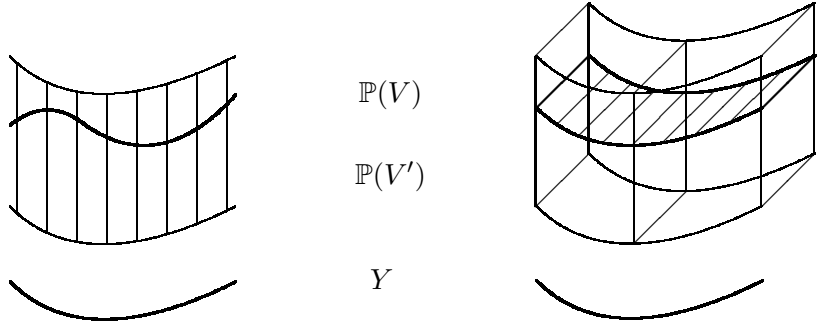
\savebox{\subbundle}(5,5)[bl]{

\put(0, 0.2){\usebox{\bogenfett}}
\put(0.7, 0.9){\usebox{\bogenfett}}

\thicklines

\put(0, 0.7){\line(1,1){0.7}}
\put(3., 0.7){\line(1,1){0.7}}

\thinlines
\put(0.29, .5){\line(1,1){0.7}}
\put(0.57, .34){\line(1,1){0.7}}
\put(0.9, .26){\line(1,1){0.7}}
\put(1.3, .24){\line(1,1){0.7}}
\put(1.71, .27){\line(1,1){0.7}}
\put(2.14, .38){\line(1,1){0.7}}
\put(2.55, .49){\line(1,1){0.7}}
}

\savebox{\viereckfaser}(5,4)[bl]{
\thinlines
\put(0,0){\line(0,1){2}}
\put(0,0){\line(1,1){0.7}}
\put(0.7, .7){\line(0,1){2}}
\put(0,2){\line(1,1){0.7}}
}

\savebox{\bundleranktwo}(5,5)[bl]{

\put(0,0){\usebox{\bogenfett}}

\put(0,1.){\usebox{\bogen}}
\put(0,3.){\usebox{\bogen}}
\put(0.7, 1.7){\usebox{\bogen}}
\put(0.7, 3.7){\usebox{\bogen}}

\put(0, 1.5){\usebox{\viereckfaser}}

\put(1.3, 1.){\usebox{\viereckfaser}}

\put(3, 1.5){\usebox{\viereckfaser}}
}

\begin{picture}(12,5)

\put(1,0){\usebox{\ruledsurface}}

\put(1,2.8){\usebox{\krummbogen}}

\put(5.8, .9){$Y$}
\put(5.65, 3.4) {$ \PP(V)$}
\put(5.6, 2.3){$\PP(V') $}

\put(8, .5){\usebox{\bundleranktwo}}

\put(8, 2.6 ){\usebox{\subbundle}}
\end{picture}

\savebox{\subbundle}{}
\savebox{\bundleranktwo}{}
\savebox{\viereckfaser}{}

\begin{corollary}
\label{ampleaffinesolid}
Suppose the situation and notation of {\rm \ref{gradedsituationdimensiontwo}}.
Suppose that the forcing divisor $Z$ is ample.
Then $f_0 \not\in (f_1, \ldots ,f_n)^\soclo$.
\end{corollary}
\proof
A multiple of the forcing divisor is very ample and defines
an embedding $ i:\PP(V') \ra \PP^N$ such that
$i^{-1} (H)= Z$. Then $\PP(V') -\PP(V)$ is affine and the result follows
from \ref{tightaffinecrit}.
\qed

\begin{remark}
We will see in Theorem \ref{ampleaffineruled}
below that for $n=2$ also the converse
of Corollary \ref{ampleaffinesolid}
is true. For $n \geq 3$ however there exists
non-ample forcing divisor with affine complement, see
example \ref{affinenonample}.
\end{remark}

\begin{theorem}
\label{pluscrit}
Suppose the situation and notation of
{\rm \ref{gradedsituationdimensiontwo}}.
Then the following are equivalent.
\renewcommand{\labelenumi}{(\roman{enumi})}
\begin{enumerate}
\item
$f_0 \in (f_1, \ldots,f_n)^{\gr}$.

\item
There exists a curve in $\PP(V')$ which does not meet
the forcing divisor $Z$.
\end{enumerate}
\end{theorem}

\proof
This follows from \ref{plus}.
\qed

\medskip
An important invariant of a divisor on a projective variety
is its top self intersection number. This may be computed directly
for a forcing divisor.

\begin{theorem}
\label{topselfintersection}
Suppose the situation and notation of {\rm \ref{gradedsituationdimensiontwo}}.
Let $\delta = \deg \,Y$ denote the degree of the curve,
i.e. the degree of the invertible sheaf
$\O_Y(1)$ on $Y$. Let $m=d_0 =\deg (f_0)$.
Then the top self intersection number of the forcing divisor
is
$$ \PP(V)^n = (d_1 + \ldots +d_n - (n-1) m) \delta  \, .$$
\end{theorem}
\proof
Note that $\PP(V')$ is a projective bundle of rank $n-1$, hence its dimension
is $n$ and $\PP(V)^n$ is the top self intersection of the
forcing divisor.
The top self intersection number is also the degree
of the vector bundle
$V(-m)$ (see \cite[Lemma 6.4.10]{lazarsfeldpositive}),
which is the same as the degree of its determinant.
We have computed in \ref{buendelalsproj}(iii)
that
$\Det \, V(-m) \cong  \AA_Y(\sum_{i=1}^n d_i -(n-1)m) $.
Therefore
$\deg V(-m) = (\sum_{i=1}^n d_i -(n-1)m) \delta $.
\qed

\begin{remark}
It is not surprising that the number $\sum_{i=1}^n d_i -(n-1)m$
and the conditions $m \geq (d_1 + \ldots +d_n)/(n-1)$
and $m < (d_1 + \ldots +d_n)/(n-1)$ play a crucial role
throughout this habilitation thesis.
\end{remark}

\begin{remark}
\label{riemannroch}
Suppose the situation and notation of \ref{gradedsituationdimensiontwo}.
A function $ f \in \Gamma( \PP(\shF')- \PP(\shF), \O_{\PP(\shF')})$
has a pole only at the forcing divisor and corresponds therefore
to a section
$$s \in \Gamma( \PP(\shF'), \O( k \PP(\shF))
= \{ q \in Q(\PP(\shF')):\, (q) +k \PP(\shF) \geq 0 \} \, $$
for some $k \in \NN$.
The theorem of Riemann-Roch describes ($m= d_0$) the
Euler-characteristic of $k \PP(\shF)$ in the following way,
\begin{eqnarray*}
\chi (\PP(\shF'), k \PP(\shF)) &= &\chi (Y, S^k(\shF(-m))) \cr
&=& \Gamma(Y, S^k (\shF(-m))) - H^1(Y, S^k (\shF(-m))) \cr
&=& \deg (S^k(\shF(-m)))  + (1-g) \rk (S^k(\shF(-m)))  \cr
&=& (\sum_{i \geq 1} d_i -(n-1)m) \binom{k+n-2}{n-1} \delta
+(1-g) \binom{k+n-2}{n-2}
\end{eqnarray*}
In characteristic zero we may replace $H^1(Y, S^k (\shF(-m)))$ by
$H^0(Y, S^k (\shR(m)))$, but not in positive characteristic.
\end{remark}

We recall the definition of a big divisor,
see \cite{beltramettisommese} or
\cite[Definitions 2.1.3 and 2.2.1]{lazarsfeldpositive}).

\begin{definition}
A Weil divisor $D$ on a normal projective variety $X$ of dimension $d$ is
called {\em big} if
for some multiple the divisor $kD$
defines a rational mapping to some projective space
such that the dimension of the image is
$d$.
\end{definition}

\begin{remark}
\label{bigremark}
A divisor is big if and only if its Iitaka-dimension
is maximal. 
This means that there exists a number $c >0$ such that
$h^0(X, \O_X(kD)) \geq ck^d$ for $k \gg 0$.
An ample divisor is big.
A numerically effective divisor $Z$ is big if and only
if its top self intersection number $Z^n$ is $>0$, see
\cite[Theorem VI.2.15]{kollarrational}
or \cite[Theorem 2.2.14]{lazarsfeldpositive}.
\end{remark}

\begin{remark}
\label{bigaffinremark}
If the complement of an effective Weil divisor $D$ on a normal projective
variety $X$ of dimension $d$ is affine, then
$\Gamma(X-D, \O_X)$ is a finitely generated $K$-algebra of dimension $d$
and there exists a multiple $kD$ such that the corresponding
rational mapping to a projective space induces an isomorphism
on $X-D$.
In particular, an effective divisor which has an affine complement
is big.
\end{remark}

It is sometimes possible to establish that the complement of the forcing
divisor is not affine (and hence that $f_0 \in (f_1, \ldots,f_n)^\soclo$)
by showing that the forcing divisor is not big (see \ref{slopemaxkrit}).
The following proposition deals with equivalent conditions for bigness
in the forcing situation.

\begin{proposition}
\label{bigaffine}
Suppose the situation and notation of {\rm \ref{gradedsituationdimensiontwo}}.
Then the following are equivalent.

\renewcommand{\labelenumi}{(\roman{enumi})}
\begin{enumerate}
\item
There exists an effective divisor $L \subset Y$ of positive degree such that
for some $a \geq 1$ the divisor
$a Z - \pi^*L$ is equivalent
to an effective divisor.

\item
There exists a linearly equivalent effective divisor
$D \sim a Z$ {\rm(}$a \geq 1${\rm )} such that
$\PP(V')- {\rm supp} \, D$ is affine.

\item
The forcing divisor $Z$ is big.

\end{enumerate}

\end{proposition}
\proof
(i) $\Rightarrow$ (ii) follows from \ref{effectivevertreter}.
Suppose that (ii) holds, let $X= \PP(V')$ and let 
$s \in \Gamma(X,\O_X(aZ))$ be a section such that $X_s$ is affine.
The topology of $X_s$ is generated by
subsets $X_t \subseteq X, \, t \in \Gamma(X,\O_X(bZ))$, $ b \geq 1$,
see \cite[Th\'{e}or\`{e}me 4.5.2]{EGAII}.
Therefore the rational mapping defined by $aZ$ is an isomorphism
on $X_s$ and the image has maximal dimension, hence $Z$ is big
(and (ii) $\Rightarrow $ (iii)).
On the other hand, 
if $\emptyset \neq U \subset X_s$ is an affine subset which
does not meet the fiber over a point $P \in Y$,
then there exists also $t \in \Gamma(X,\O_X(bZ))$
such that $ \emptyset \neq X_t \subseteq U$.
Therefore
$bZ+(t) = \sum_i a_i D_i$, $a_i >0$ is an effective divisor
and $\PP(V'_P)$ is one of the $D_i$ (hence (ii) $\Rightarrow $ (i)).
(iii) $\Rightarrow $ (ii).
If $Z$ is big, then for some $a \geq 1$
the multiple $aZ$ defines a mapping which is birational with its image.
Therefore the mapping induces an isomorphism on an open affine subset
$X_s \cong X-V(s)$, $s \in \Gamma(X,\O_X(aZ))$.
\qed

\newsavebox{\rechteck}

\savebox{\rechteck}(5,4)[bl]{
\put(0,0){\line(1,0){4}}
\put(0,0){\line(0,1){3}}
\put(0,3){\line(1,0){4}}
\put(4,0){\line(0,1){3}}}

\newsavebox{\rechteckdreifuenfundzwei}

\savebox{\rechteckdreifuenfundzwei}(5,4)[bl]{
\put(0,0){\line(1,0){3.5}}
\put(0,0){\line(0,1){2}}
\put(0,2){\line(1,0){3.5}}
\put(3.5,0){\line(0,1){2}}}

\newsavebox{\strahl}

\savebox{\strahl}(5,5)[bl]{
\put(0,0){\line(0,1){1.1}}
\put(0,0){\line(0,-1){1.1}}
\put(0,0){\line(3,1){1.}}
\put(0,0){\line(3,5){.55}}
\put(0,0){\line(-3,1){1.}}
\put(0,0){\line(-3,5){.55}}
\put(0,0){\line(3,-1){1.}}
\put(0,0){\line(3,-5){.55}}
\put(0,0){\line(-3,-1){1.}}
\put(0,0){\line(-3,-5){.55}}
}

\newpage

\thispagestyle{empty}

\markboth{3. The two-dimensional parameter case}
{3. The two-dimensional parameter case}

\section{The two-dimensional parameter case}

In this short chapter we apply our theory to
the first non-trivial situation, the case of two
homogeneous parameters $f_1,f_2 \in R$ in a two-dimensional normal
standard-graded domain $R$.

The relation bundle for the parameters is an invertible sheaf
on the projective curve $Y= \Proj R$
and a third homogeneous element $f_0 \in R$ defines an extension of rank two
of two invertible sheaves.
This setting yields a projective bundle of rank one
together with a projective subbundle
of rank $0$ on the projective curve $Y$, i.e a ruled surface together
with a forcing section.
This section is an ample divisor if and only if its complement
is affine, which is equivalent with $f_0 \not\in (f_1,f_2)^\soclo$.

Since in the two-dimensional parameter case
everything is known about the tight closure
we cannot expect new results for tight closure in this chapter.
However, our geometric interpretation yields new geometric proves
for known facts and it yields somewhat surprisingly a new class of
counter-examples to the hypersection problem in complex analysis.

\bigskip

\subsection{Ruled surfaces and forcing sections}

\label{ruledsection}

\markright{Ruled surfaces and forcing sections}

\

\bigskip
We fix the following situation.

\begin{situation}
\label{ruledsituation}
Let $K$ denote an algebraically closed field
and let $R$ denote a normal standard-graded two-dimensional $K$-domain.
Let $Y = \Proj \, R$ denote the corresponding smooth projective curve.
Suppose that $f_1$ and $f_2$ are homogeneous parameters
of degree $d_1$ and $d_2$
and let $f_0$ be another homogeneous element of degree $d_0$.
Fix $m \in \NN$ and let $V(-m)$ and $V'(-m)$ denote the relation bundles
on $Y$ as described in in \ref{forcingsequence2}.
Let $Z=\PP(V) \subset \PP(V')$ denote the forcing divisor.
\end{situation}

\begin{corollary}
\label{forcingrule}
In the situation {\rm\ref{ruledsituation}} the following hold. 

\renewcommand{\labelenumi}{(\roman{enumi})}
\begin{enumerate}

\item
$\PP(V')$ is a ruled surface and $\PP(V) \subset \PP(V')$
is a section, called the forcing section.

\item 
We have
$\Proj \,R[T_1,T_2]/(f_1T_1+f_2T_2) \supset D_+(R_+) =V(-m)
\cong \AA_Y(d_1+d_2-m)$.
In particular, the forcing sequence is
$$ 0 \longrightarrow \AA_Y(d_1+d_2-m) \longrightarrow V'(-m)
\longrightarrow \AA_Y(d_0-m) \longrightarrow 0 \, .$$

\item
We have ${\rm Det} \, V'(-m) \cong \AA_Y( d_1+d_2+d_0 -2m)$.

\item
The normal bundle for the embedding
$Y\cong \PP(V) \subset \PP(V')$
is $\AA_Y( d_0 - d_1 -d_2 )$.

\item
The self intersection number of the forcing section
$Y \cong \PP(V) \hookrightarrow \PP(V')$
is $ (d_1+d_2-d_0 ) \deg\, H$, where
$H$ is the hyperplane section corresponding to $\O_Y(1)$.
\end{enumerate}
\end{corollary}

\proof
(i)
follows from \ref{forcingsequence2}.
(ii).
Let $R[T_1,T_2]/(f_1T_1+f_2T_2)$ be graded
by $ \deg (T_1) = m-d_1$ and $ \deg (T_2) = m-d_2$.
The homomorphism
$R[T_1,T_2]/(f_1T_1+f_2T_2) \ra R[U]$ given by
$T_1 \mapsto f_2U,\, T_2 \mapsto -f_1U$ is then homogeneous for
$\deg\, (U) = m-d_1-d_2$ and induces an isomorphism
on $D_+(R_+)$.
The corresponding line bundle is $\AA_Y(d_1+d_2-m)$. (iii) follows.

The normal bundle for the embedding on $\PP(V) \cong Y$ is
$N = \AA_{\PP(V)}(-1) \otimes \AA_Y(d_0 -m)$ due to \ref{forcingsequence2}(iv).
Furthermore,
$V(-m)=  \AA_Y(d_1+d_2-m)$ on $Y$ and
$V(-m)=\AA_{\PP(V)}(+1)$ is the tautological line bundle
for $\PP(V) \cong Y$. This yields
$N = \AA_Y(m-d_1-d_2) \otimes \AA_Y(d_0-m)=\AA_Y(d_0 - d_1-d_2)$.
Its sheaf of sections is
$\O_Y(d_1+d_2-d_0)$ and its degree is the self intersection number,
hence (iv) follows.
\qed

\begin{remark}
\label{forcingsheaf3}
The corresponding sequence of sheaves are
$$ 0 \longrightarrow \O_Y(d_0-m) \longrightarrow \shF'(-m)
 \longrightarrow \O_Y(d_1+d_2-m) \longrightarrow 0 \, $$
for the linear forms $\shF'(-m)$ and
$$ 0 \longrightarrow \O_Y(m- d_1-d_2) \longrightarrow \shR'(m)
\longrightarrow \O_Y(m-d_0) \longrightarrow 0 \, $$
for the relations $\shR'(m)$ .
These extensions are classified by
$H^1(Y, \O_Y(d_0- d_1-d_2))$ for $m=d_0$, where the elements
$f_1,f_2;f_0$ correspond to the cohomology class
$f_0/f_1f_2$.
\end{remark}

The following theorem gives a criterion for tight closure
in the two-dimen\-sio\-nal parameter case in terms of ampleness of the
forcing divisor.

\begin{theorem}
\label{ampleaffineruled}
Suppose the situation and notation of {\rm\ref{ruledsituation}}.
Then the following are equivalent.

\renewcommand{\labelenumi}{(\roman{enumi})}
\begin{enumerate}

\item 
$f_0 \not\in (f_1,f_2)^\pasoclo$.

\item
$\PP(V')-Z$ is affine.

\item
The forcing divisor $Z$ on $\PP(V')$ is ample.
\end{enumerate}
\end{theorem}

\proof
We know the equivalence (i)$ \Leftrightarrow$ (ii) from
proposition \ref{tightaffinecrit}
so we have to show the equivalence (ii)$ \Leftrightarrow$ (iii).
If $Z$ is ample, then its complement is affine
as mentioned in corollary \ref{ampleaffinesolid}.
If $\PP(V')-Z$ is affine, then it does not
contain projective curves.
Furthermore, there exist global functions on $\PP(V')-Z$ which are
not constant.
Thus $aZ$, $a \geq 1$,
is linearly equivalent with an effective divisor not containing $Z$.
Hence the self intersection number is positive and the
criterion of Nakai yields that $Z$ is ample.
\qed

\unitlength1cm
\begin{picture}(12,5)

\put(5,0){\usebox{\ruledsurface}}

\put(5,2.8){\usebox{\krummbogen}}

\put( 4.7, 2.9){\vector(2,3){.01}}
\put(4,2){$s$}
\bezier{300}(4.7, 1.1)(4, 2.1)(4.7, 2.9)

\put(8.5, .8) {$Y$}
\put(8.2, 3.4) {$ Z=\PP(V) = s(Y) $}
\put(8.8, 2.3){$\PP(V') $}

\end{picture}

\begin{corollary}
\label{degreecrit1}
Let $K$ be an algebraically closed field and let $R$ be a normal
two-dimensional
standard-graded $K$-algebra and let $f_1,f_2$
be homogeneous parameters of degrees
$d_1,d_2$.
Then $R_{ \geq d_1 + d_2} \subseteq (f_1,f_2)^\pasoclo$.

If the characteristic of $K$ is zero,
then
$(f_1,f_2)^\pasoclo= (f_1,f_2)+ R _{\geq d_1 + d_2}$
\end{corollary}
\proof
Let $d_0 = \deg\, f_0 \geq d_1+d_2$.
Then $d_1+d_2-d_0 \leq 0$ and the self intersection of
$Z=\PP(V) \subset \PP(V')$
is not positive. Hence $Z$ is not ampel and $f_0 \in (f_1,f_2)^\pasoclo$
due to theorem \ref{ampleaffineruled}.

Now let $f_0 \in (f_1,f_2)^\pasoclo$, but $f_0 \not\in R_{\geq d_1+d_2}$.
Then the self intersection is positive, but
the forcing divisor $Z$ is not ample.
Thus there must exist a curve $C \subset \PP(V')$ disjoint to $Z$.
By \ref{plus} it follows that $f_0 \in (f_1,f_2)^+$,
so \ref{finiteextension}(ii) gives the result in characteristic $0$.
\qed

\begin{remark}
The second statement in corollary \ref{degreecrit1}
is also true if the characteristic of $K$
is $p \gg 0$. This follows from the theorem of Hartshorne-Mumford
about the ampleness of vector bundles of rank two, which we treat in
the next chapter, see \ref{amplecritranktwo}.
\end{remark}

\begin{remark}
A locally free sheaf $\shE$ of rank two on a projective curve $Y$
is called normalized if $\Gamma(Y, \shE) \neq 0$,
but $\Gamma(Y, \shE \otimes \shL)=0$ for every invertible sheaf
$\shL$ of negative degree.
For every ruled surface $X \ra Y$ there exists a representation
$X = \PP(\shE)$ such that $\shE$ is normalized.
The {\em e-invariant} of a ruled surface is then defined by
$e = - \deg (\shE)$, where $\shE$ is normalized.
In this normalized situation there exists also
a short exact sequence
$0 \ra \O_Y \ra \shE \ra \shM \ra 0$ and the section $C_0$
corresponding to the surjection
is also called normalized. Then $C_0^2= -e$, see \cite[V.2]{haralg}.
\end{remark}

\smallskip 

\subsection{Examples}
\label{examplesparameter}

\markright{Examples}
\

\bigskip

We give some examples of ruled surfaces and their forcing
sections arising from forcing data.
Let $R$ be a standard-graded normal two-dimensional $K$-domain,
$Y= \Proj \, R$,
and let $f,g$ be homogeneous parameters and $h$ homogeneous
of degrees $d_1,d_2,d_0$.

\begin{example}
\label{spaltet}
Suppose that $h \in (f,g)$.
Then the forcing sequence splits and
$$V'(-m) \cong \AA_Y( d_1+d_2-m ) \times \AA_Y(d_0-m)$$
This is normalized for $m = \max (d_0,d_1+d_2)$
and then isomorphic to
$\AA_Y \times \AA_Y(-|d_1+d_2-d_0)$ and
$\PP(V')$ is the projective closure of the line bundle
$\AA_Y(-|d_1+d_2-d_0|)$. The forcing section is either the zero section
of the line bundle or the closure section.
The $e$-invariant is $|d_1+d_2-d_0|\deg \, H$, where
$H$ is the hyperplane section on $Y$.
\end{example}

\begin{example}
\label{zerlegbar}
Suppose that $h=1$. Then
$$\Proj\, R[T_0,T_1,T_2]/(fT_1+gT_2+T_0) \cong \Proj\, R[T_1,T_2] \, .$$
and therefore
$V'(-m)= \AA_Y(d_1-m) \times \AA_Y(d_2-m)$ is decomposable.
The forcing sequence is
$$0 \longrightarrow \AA_Y(d_1+d_2-m) \longrightarrow
\AA_Y(d_1-m)\times \AA_Y(d_2-m) \longrightarrow \AA_Y(-m)
\longrightarrow 0 \, ,$$
where $v \mapsto (gv,-fv)$ and $(a,b) \mapsto -(af+bg)$.
The sequence is normalized for $m= \max(d_1,d_2)$.
Let $d_1 \geq d_2$. Then the normalized sequence is
$$0 \longrightarrow \AA_Y(d_2) \longrightarrow
\AA_Y \times \AA_Y(d_2-d_1) \longrightarrow \AA_Y(-d_1)
\longrightarrow 0 \, .$$
The $e$-Invariant of the ruled surface is $e=|d_2-d_1|\deg\, H$.
\end{example}

\begin{example}
\label{graph}
Let $f$ and $g$ be elements of the same degree $d$ and let $h=1$.
This yields the trivial ruled surface $\PP(V')=Y \times \PP^1$ and the 
forcing section is the graph of the  meromorphic function $q= -f/g$.
For let  $m=d$. Then we have the forcing sequence
$$0 \longrightarrow \AA_Y(d) \stackrel{g,-f} {\longrightarrow}
\AA_Y \times \AA_Y 
\stackrel{fT_1+gT_2}{\longrightarrow} \AA_Y(-d)
\longrightarrow 0 \, ,$$
and $fT_1+gT_2=0$ is equivalent with $T_1/T_2= -g/f$.
If $d \geq 1$, then $h \not\in (f,g)^\pasoclo$, and so
we see via tight closure that the complement of the graph
of a non constant meromorphic function is affine.
\end{example}

\begin{example}
\label{projektivegerade}
Let $R=K[x,y]$, thus $Y=\PP_K^1$.
The ruled surface $\PP(f,g;h)$ must be a Hirzebruch surface
$\PP(\AA_{\PP_K^1} \times \AA_{\PP^1_K} (k))$
and we have to determine the number $k \leq 0$.

Let $(a_1,a_2,a_3)$ and $(b_1,b_2,b_3)$ be a basis of
homogeneous relations for all the relations
for $(f,g,h)$.
Let $e_i$ and $e_i'$ be its degrees and suppose that $e_i \geq e_i'$.
Set $ k= e_i'- e_i$.
The homomorphism
$$K[x,y,T_1,T_2,T_3]/(fT_1+gT_2+hT_3) \stackrel{\psi }{\lra} K[x,y,U,S] $$
given by $T_i \mapsto a_iU+b_iS,\, i=1,2,3$, is well defined
and is homogeneous for
$\deg \, (T_i)=e_i,\, \deg\, (U)=0$, $ \deg\, (S)=e_i -e_i'=-k $.
It induces a mapping
$$ \Proj\, K[x,y,U,S] \supseteq D_+(x,y)
=\AA_{\PP^1_K} \times \AA_{\PP^1_K}(k)  \lra V'(-m) \, .$$
We claim that this is an isomorphism.
Since $(-g,f,0)$ is a relation, there exist $r,s \in R$
such that
$(-g,f,0)= r(a_1,a_2,a_3) + s(b_1,b_2,b_3)$.
Since $(0,-h,g)$ and $(-h,0,f)$ are also relations it follows that
$(g,f) \subseteq (a_3,b_3)$. Hence $(a_3,b_3)$ is $R_+$-primary
and $a_3$ and $b_3$ do not have a common divisor.
Therefore there exists $t \in R$ such that
$r=tb_3,\, s=-ta_3$, hence we may write $f=t(b_3a_2-a_3b_2)$.
The mapping $\psi$ is locally on $D(f)$ given by a linear transformation
$R_f[T_2,T_3] \ra R_f[U,S]$ and its determinant is
$b_3a_2-a_3b_2$, which is a unit in $R_f$.
The same is true on $D(g)$, so the induced mapping is an isomorphism.

The forcing sequence is
$$0 \lra \AA_{\PP_K^1} (d_1+d_2-d_0 -e_3) \lra 
\AA_{\PP^1_K} \times \AA_{\PP^1_K}(e_3'-e_3)
\stackrel{a_3U+b_3S}{\lra} \AA_{\PP_K^1}(-e_3) \lra 0
\, .$$
It follows that $d_1+d_2-d_0=e_3 + e_3'$.

Let $f=x^2,g=y^2$. For $h=xy$ we have $e_3=e_3'=1$, but for
$h=x^2$ we have $e_3=2, \, e_3'=0$.
\end{example}

Let $f,g,h \in R_d$ be homogeneous of the same degree $d$ and
suppose that $f,g$ are parameters and that $h \not\in (f,g)$.
Let $E \subseteq \Gamma(Y, \O_Y(d))$ denote the linear system
spanned by $f,g,h$ and set $m=d$.
Then the presenting sequence for $f,g,h$ from \ref{buendelalsproj}(ii) is
$$0 \lra V'(-d) \lra \AA_Y^3 \lra \AA_Y(-d) \lra 0 \,  .$$
The mapping sends
$(P,t_1,t_2,t_3) \mapsto t_1f(P) +t_2g(P) +t_3h(P)$ and this is
zero if and only if the point $P$ lies on the divisor
defined by the global section $t_1f +t_2g+t_3h$ of 
$\O_Y(d)$.
Therefore
$$ \PP(V')= \{ (P,D): \,\,  P \in D, \, \, D \in E \}  \, $$
and the ruled surface $\PP(V')$ is the incidence variety
for the linear system $E$.

Now let $Y =V_+(F) \subset \PP_K^2 = \Proj\, K[x,y,z]$
be a smooth curve and consider the linear system of lines.
Thus the ruled surface associated to the vector bundle
$$V'(-1)=D_+(x,y,z) \subset \Proj\, K[x,y,z]/(F)[T_1,T_2,T_3]/(xT_1+yT_2+zT_3) \,$$
($\deg \, (T_i)=0$)
consists of the pairs $(P,L)$, where $L$ is a line through $P \in Y$.
Suppose that $x,y$ are parameters for the curve. Then the point $Q=(0,0,1)$
does not belong to $Y$.
The forcing section maps a point $P \in Y$ to the line passing through
$P$ and $Q$, since this is the line given by $T_3=0$.
The self intersection number of the forcing section is
$\deg \,(F)$. If $\deg \,(F) \geq 2$, then $z \not\in (x,y)^\pasoclo$ and
the forcing section is ample.

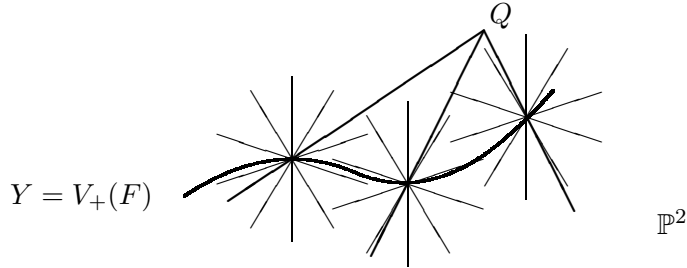
\begin{picture}(12,5)

\put(1.4,1.3){$Y=V_+(F)$}
\thicklines

\put(7.77,3.7){$Q$}
\put(7.7,3.6){\line(-3,-2){3.4}}
\put(7.7,3.6){\line(-1,-2){1.5}}
\put(7.7,3.6){\line(1,-2){1.2}}

\bezier{300}(3.7 , 1.4 )(4.9, 2.2)(6. ,1.7)
\bezier{300}(6, 1.7)(7,1.3)(8 ,2.2)
\bezier{300}(8, 2.2)(8.3,2.45)(8.6 ,2.8)

\put(10,0.9){$\PP^2$}

\thinlines

\put(5.15, 0.79){\usebox{\strahl}}
\put(6.69,0.46){\usebox{\strahl}}
\put(8.26, 1.35){\usebox{\strahl}}
\end{picture}

\savebox{\strahl}{}

\begin{example}
Let $F \in K[x,y,z]$ be a homogeneous polynomial of degree three
such that $Y= \Proj \, K[x,y,z]/(F)$ is an elliptic curve and
suppose that $x,y$ are parameters.
Consider $\PP(x,y;z)$ as above and set $m=1$, and let $\shE= \shF(-1)$
be the sheaf of linear forms for this grading.
The global linear forms $T_1,T_2,T_3$ are a basis for $H^0(Y, \shE) $.
A global linear form $s=a_1T_1+a_2T_2+a_3T_3$ ($a_i \in K$)
belongs to $\Gamma(Y, \shE \otimes \O_Y(-P))$ if and only if
$s|V'_P =0$, and this is the case if and only if
$(a_1,a_2,a_3)$ is a multiple of $(x(P),y(P),z(P))$.
So this can happen at most at one point.
Hence $\dim H^0(Y,\shE \otimes \O_Y(-P)) =1$ and
$H^0(Y, \shE \otimes \O_Y(-P-Q))=0$.
Therefore $\shE \otimes \O_Y(-P)$ is normalized and the $e$-invariant
of $\PP(x,y;z)$ is
$$e= -\deg \, (\shE \otimes \O_Y(-P)) = -3+2=-1 \, .$$
\end{example}

\smallskip

\subsection{Examples over the complex numbers $\CC$}
\label{complex}

\markright{Examples over the complex numbers $\CC$}

\

\bigskip

The ruled surfaces together with their forcing sections arising
from tight closure problems yield analytically interesting
examples over the field of complex numbers $K=\CC$.

\begin{corollary}
\label{steinnonaffine}
Let $F$ be a homogeneous polynomial of degree three
such that $R=\CC[x,y,z]/(F)$ defines an
elliptic curve $Y= \Proj\, R$. Let $f,g,h$ be homogeneous
such that $f$ and $g$ are parameters,
$\deg \, (h) = \deg\, (f) + \deg\, (g)$ and $h \not\in (f,g)$.
Then the complement of the forcing section $\PP(f,g)$ in
the ruled surface $\PP(f,g;h)$ over
$Y$ is not affine, but it is a complex Stein space.
The same is true for the open subset
$D(R_+) \subset \Spec\, R[T_1,T_2]/(fT_1+gT_2+h)$.
\end{corollary}

\proof
The condition on the degrees shows that
$h \in (f,g)^\pasoclo$ and that the open complement is not affine.
Since $h \not\in (f,g)$, the forcing sequence, which is an
extension of $\O_Y$ by $\O_Y$, does not split.
Hence due to \cite{umemura} the complement is Stein.
The corresponding statement for the subset in the cone follows.
\qed

\begin{example}
\label{steintwo}
Let $R=\CC[x,y,z]/(x^3+y^3+z^3)$ and $f=x,g=y,h=z^2$.
Then $z^2 \not\in (x,y)$, but $z^2 \in (x,y)^\pasoclo$.
The open subset
$$D(x,y) \subseteq \Spec\, \CC[x,y,z][T_1,T_2]/(x^3+y^3+z^3,\, xT_1+yT_2+z^2)$$
is not affine, but it is a complex Stein space.
\end{example}

\begin{remark}
The first construction of a non-affine but Stein variety
was given by Serre using non-split extensions
of $\O_Y$ by $\O_Y$ on an elliptic curve, see \cite{umemura},
\cite{bingener} or \cite{neeman}.
Thus we may consider this classical construction
as a construction using forcing algebras.

On the other hand we have to remark that tight closure
takes into account the subtile difference between affine and Stein.
This shows that tight closure is a conception
of commutative algebra and algebraic geometry, not of complex analysis.
\end{remark}

It is also possible to construct new counterexamples
to the hypersection problem. The first counterexample
was given 1997 by Coltoiu and Diederich in \cite{coltoiudiederich}.
For this and related problems in complex analysis
see \cite{diederich}. 
The problem is the following:
Let $X$ denote a complex-analytic space
of dimension $\geq 3$
and let $D \subset X$ be an analytic hypersurface.
Suppose that the intersection
$( X-D) \cap S$ is Stein for every analytic hypersurface
$S \subset X$. Is then $X-D$ itself Stein?

\begin{proposition}
\label{hypersectioncounterexample}
Let $R$ be a standard-graded normal two-dimensional
$\CC$-al\-ge\-bra, let $f,g$ and $h$ be homogeneous elements in $R$
such that $V(f,g)=V(R_+)$, $h \not\in (f,g)$ and
$\deg \, (h) > \deg\, (f) + \deg \, (g)$.
Then
$$W=D(R_+) \subset \Spec \, R[T_1,T_2]/(fT_1+gT_2+h)=X$$
is not Stein {\rm(}considered as a complex space{\rm )},
but it fulfills the assumption in the hypersection problem,
i.e. for every analytic surface $S \subset \Spec\, R[T_1,T_2]/(fT_1+gT_2+h)$
the intersection $W \cap S$ is Stein.
\end{proposition}
\proof
Due to \ref{superhoehe} the superheight of $W$ is one,
and \cite[Theorem 5.1]{brennersuperheight} gives that the assumption of the
hypersection problem holds.
The self intersection number of the forcing section on the
corresponding ruled surface is negative,
thus due to \cite{grauertmod} this section is contractible as a complex space
and therefore its complement is not Stein.
Hence the subset $W$ is also not Stein,
because it is a $\CC^\times$-bundle over this complement.
\qed

\begin{example}
\label{hypersectionexample}
Consider
$R=\CC[x,y,z]/(x^4+y^4+z^4)$ and the forcing algebra for the
elements
$x,y;z^3$, hence
$A=\CC[x,y,z,T_1,T_2]/(x^4+y^4+z^4,xT_1+yT_2+z^3)$.
Then $z^3 \not\in (x,y)$ in $R$, but $z^3 \in (x,y)^\pasoclo$,
for the degree of the self intersection is $-4$.
Therefore
$W=D(R_+) \subseteq \Spec \, A$ is not Stein,
but for every analytic surface $S \subset \Spec A$
the intersection $W \cap S$ is Stein.
\end{example}

\bigskip

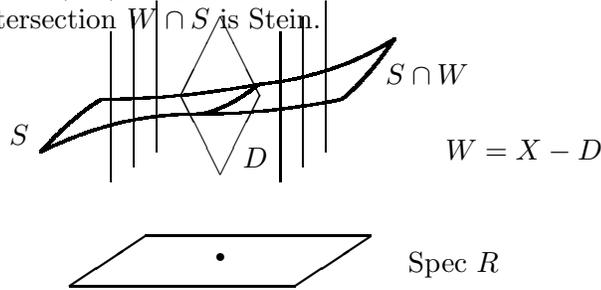
\savebox{\raute}

\savebox{\raute}(3,3)[bl]{
\put(0, 0){\line(1,2){0.53}}
\put(0,0){\line(-1,2){0.53}}
\put(0,2.1){\line(1,-2){0.53}}
\put(0,2.1){\line(-1,-2){0.53}}}

\unitlength1cm
\begin{picture}(12,5)
(0,-1.6)

\put(6, 1.4){\circle*{0.1}}

\put(4,1){\usebox{\viereckdick}}

\put(3.2,2.9){$S$}
\put(8.5,1.2){$\Spec\, R$}
\put(9,2.7){$W =X-D$}
\put(8.2,3.7){$ S\cap W$}
\put(6.3, 2.6){$D$}

\put(6,2.5){\usebox{\raute}}

\thicklines

\bezier{300}(5.6 , 3.3)(6.0 ,3.3 )(6.5 ,3.7)
\bezier{300}(3.6 , 2.8)(4.6 ,3.3 )(5.6 ,3.3)
\bezier{300}(5.6 , 3.3)(6.6 ,3.3 )(7.6 ,3.5)

\bezier{300}(4.4, 3.5)(5.4 ,3.5 )(6.5 ,3.7)
\bezier{300}(6.5 , 3.7)(7.5 ,3.8 )(8.3 ,4.3)

\bezier{300}(3.6 , 2.8)(4.0 ,3.25)(4.4 ,3.5)

\bezier{300}(7.6, 3.5)( 7.95, 3.78 )(8.3 ,4.3)

\thinlines
\put(4.55,2.4){\line(0,1){2}}
\put(4.85,2.6){\line(0,1){2}}
\put(5.15,2.8){\line(0,1){2}}

\put(6.8,2.4){\line(0,1){2}}
\put(7.1,2.6){\line(0,1){2}}
\put(7.4,2.8){\line(0,1){2}}
\end{picture}

\subsection{Plus closure in positive characteristic}
\label{positive}

\markright{Plus closure in positive characteristic}

\

\bigskip
We have mentioned in \ref{tightparameter} the
theorem of Smith that the tight closure and the plus closure of a
parameter ideal are the same.
In this section we will give a proof for this theorem
in the two-dimensional graded case within our geometric setting.
Let $K$ be an algebraically closed field of positive characteristic $p$
and let $f$ and $g$ be homogeneous parameters in
a two-dimensional standard-graded normal $K$-algebra $R$.
Let $h$ be another homogeneous element.
If the complement of the forcing divisor $Z \subset \PP(f,g;h)$
is not affine, then we must find a curve on the ruled surface
disjoined to $Z$.

Suppose first that $ \deg\, (h) > \deg\, (f) + \deg\, (g)$.
Then the cohomology class
$\frac{h}{fg}$ has positive degree
$k \deg\, (h) - \deg\, (f)- \deg\, (g) $ and
applying a Frobenius morphism it becomes
$(\frac{h}{fg})^q \in H^1(Y,\O_Y(qk))$, but
$H^1(Y,\O_Y( qk))=0$ for $q=p^{e}$ sufficently large.
Therefore the forcing sequence splits after a Frobenius morphism
and $h$ belongs to the Frobenius closure of $(f,g)$
and in particular to its tight closure.
Thus we have to consider the case
$\deg\, (h) = \deg\, (f) + \deg\, (g)$.

\begin{proposition}
\label{frobeniusartinschreier}
Let $K$ be an algebraically closed field of characteristic $p >0$.
Let $R$ be a standard-graded normal two-dimensional $K$-domain
and let $f,g,h$ be homogeneous elements such
that $f$ and $g$ are parameters and such that
$\deg \,(h) =\deg\,(f) + \deg \, (g)$.
Then there exists a composition of a Frobenius morphism and an
Artin-Schreier extension of $Y= \Proj \,R$ such that the
image of the cohomology class $c=h/fg \in H^1(Y,\O_Y)$
vanishes.
\end{proposition}

\proof
The Frobenius morphism $\Phi$ acts on
$H^1(Y,\O_Y)$ $p$-linear yielding the so called
Fitting decomposition
$H^1(Y,\O_Y) = V_s \oplus V_n$ such that
$\Phi|V_s$ is bijective and $\Phi|V_n$ is nilpotent
(see \cite[III \S 3]{haramp} or \cite[Lemma III.4.13]{milne}).
Thus we may write $c=c_1+c_2$,
where $c_2$ becomes zero after applying a certain power of the Frobenius.
Thus we may assume that $c=c_1 \in V_s$.
Consider the Artin-Schreier sequence
$$0 \lra \ZZ/ (p) \lra \O_Y \stackrel{\Phi - id}{\lra}  \O_Y \lra 0 \, ,$$
which is exact in the $\acute{\rm e}$tale topology. It yields the
exact sequence
$$ 0 \lra H^1_{et}(Y, \ZZ/(p)) \lra H^1(Y,\O_Y)
\stackrel{\Phi - id}{\lra} H^1(Y,\O_Y) \lra \ldots \, \, .$$
There exists a basis $c_j$ of $V_s$ such that $\Phi(c_j)=c_j$,
see \cite[\S 14]{mumfordabelian}.
Thus we may assume that $\Phi(c)-c=0$ and we consider
$c \in H^1_{et}(Y,\ZZ/(p))$.
Hence $c$ represents an Artin-Schreier extension $Y'$ of $Y$ and
the cohomology class $c$ vanishes on $Y'$, see \cite[III. \S 4]{milne}.
\qed

\begin{remark}
\label{artinschreier}
We describe the Artin-Schreier extension appearing in the
last proof explicitly.
Let $c=h/f_1f_2$ and suppose that $c^p-c=0$ in $H^1(Y,\O_Y)$.
This means that $c^p-c= a_2-a_1$, where
$a_i \in \Gamma(U_i, \O_Y)$, $U_i = D_+(f_i)=\Spec\, R_i$.
Let $U_i'= \Spec\, R_i[T_i]/(T_i^p-T_i+a_i)$, $i=1,2$.
The transition function $T_1 \mapsto T_2+c$ is due to
$$(T_2+c)^p -(T_2+c) +a_1 = T_2^p- T_2 +c^p-c +a_1= T_2^p- T_2+a_2$$
well defined and $U_1'$ and $U_2'$ glue together to a scheme $Y' \ra Y$.
The cohomology class in $Y'$ is $c=h/f_1f_2=T_1-T_2$,
$T_i \in \Gamma(U_i', \O_{Y'})$, so that $c=0$ in $H^1(Y', \O_{Y'})$.
\end{remark}

\begin{remark}
Let $Y = \Proj \, R$ as in {\rm \ref{frobeniusartinschreier}}.
The $p$-rank of $Y$ is the dimension $ \dim V_s \leq g(Y)$.
This is the same as the $p$-rank of the jacobian of $Y$,
see \cite[\S 15]{mumfordabelian}.
The $p$-rank is 0 if and only if the plus closure (=tight closure)
of any homogeneous parameter ideal is
the same as its Frobenius closure.

Suppose now that $Y= \Proj \, K[x,y,z]/(F)$ is an elliptic curve,
thus we have $H^1(Y, \O_Y) \cong K$.
An elliptic curve with $p$-rank $0$ is
called supersingular (or is said to have Hasse invariant $0$).
The criterion \cite[Proposition IV.4.21]{haralg} says that
$Y$ is supersingular if and only if the
coefficient of $(xyz)^{p-1}$ in $F^{p-1}$ is 0,
or equivalently $F^{p-1} \in (x^p,y^p,z^p)$.
On the other hand,
if the Hasse invariant is 1, then $F^{p-1} \not\in (x^p,y^p,z^p)$,
and the criterion of Fedder \cite[Theorem 3.7]{huneketightparameter}
tells us that
$K[x,y,z]/(F)$ is Frobenius pure.
\end{remark}

\newpage

\thispagestyle{empty}

\savebox{\strahl}{}

\markboth{4. Slope criteria for affineness and non-affineness}
{4. Slope criteria for affineness and non-affineness}

\section{Slope criteria for affineness and non-affineness}
\label{slopechapter}

Our interpretation of tight closure in terms of vector bundles and their
extensions has led to the following more general problem:
Suppose that $\shS$ denote a locally free sheaf on a smooth projective
curve $Y$ and let $c \in H^1(Y, \shS)= \Ext^1(\O_Y, \shS)$
denote a given cohomology class (or extension class), and let
$\shS'$ be the locally free sheaf which represents the class.
Let $ 0 \ra \O_Y \ra \shG' \ra \shG \ra 0$ be the dual sequence.
Under which conditions is the open
subscheme $\PP(\shG') - \PP(\shG)$ affine?

The answer depends of course both on $\shS$ and on $c$.
We will study this question in terms of the slope of $\shS$,
which plays a crucial role in this game.
This chapter requires no knowledge about tight closure.
We will derive applications for tight closure in the next chapter.

\bigskip

\subsection{Slope of bundles over a curve}
\label{sectionslope}

\markright{Slope of bundles over a curve}
\

\bigskip
In this section we recall the notion of the slope of a vector
bundle and various related concepts like semistable sheaves
and the Harder-Narasimhan filtration which we will need in the sequel,
our main references are \cite{hardernarasimhan}, \cite{huybrechtslehn},
\cite[Ch. 6.4]{lazarsfeldpositive}, \cite{miyaokachern} and \cite{seshadrifibre}.

Let $\shE$ denote a locally free sheaf of positive rank
on a smooth projective
curve $Y$ over an algebra\-i\-cal\-ly closed field $K$.
The degree of $\shE$ is defined by $\deg \, (\shE) = \deg \, (\bigwedge^r (\shE))$,
where $r$ is the rank of $\shE$.
If $\PP(\shE)= \Proj \, \oplus_{k \geq} S^k(\shE)$
is the corresponding projective bundle of
dimension $r$ and if $\xi$ denotes the divisor class corresponding
to the relatively very ample invertible sheaf
$\O_{\PP(\shE)}(1)$, then also $\deg \, (\shE)= \xi^r$
equals the top self intersection number,
see \cite[Lemma 6.4.10]{lazarsfeldpositive}.
The number $\deg \, (\shE)/ r!$ is also
the coefficient of $k^r$ in the Euler-Hilbert polynomial
$\chi( \O_{\PP(\shE)}(k))$, which equals $\chi( S^k(\shE))$.
The degree is additive on short exact sequences.

The {\em slope} of a locally free sheaf $\shE$ is defined by
$\mu (\shE)= \deg \, (\shE) / \rk (\shE)$.
For two locally free sheaves $\shE_1$ and $\shE_2$
we have that
$ \mu (\shE_1 \otimes \shE_2)= \mu( \shE_1) + \mu ( \shE_2)$
A locally free sheaf $\shE$ is called {\em semistable} if for every
locally free quotient sheaf of positive rank
$\shE \ra \shQ \ra 0$ the inequality
$\mu (\shQ) \geq \mu(\shE)$ holds.
This is equivalent to the property that for every locally free
subsheaf $\shT \subseteq \shE$ the inequality
$\mu (\shT) \leq \mu (\shE)$ holds.

Every locally free sheaf $\shE$
has a unique {\em Harder-Narasimhan filtration}.
This is a filtration of locally free subsheaves
$$ 0=\shE_0 \subset \shE_1 \subset \ldots \subset \shE_s = \shE $$
such that $\shE_i/\shE_{i-1}$ is semistable for every $i=1, \ldots, s$.
$\shE_1$ is called the {\em maximal destabilizing subsheaf}.
The slopes of these semistable quotients
form a decreasing chain
$\mu_1 > \ldots > \mu_s$.
$\mu_{\rm min}(\shE)= \mu_s= \mu(\shE/\shE_{s-1})$ is called the
{\em minimal slope}
and $\mu_{\max}(\shE) =\mu_1 (\shE)$ is called the {\em maximal slope}.
This is the same as
$\mu_{\rm min}(\shE)=\min\{\mu(\shQ):\, \shE \ra \shQ \ra 0 \}$.
For the dual sheaf we have
$\mu_{\rm max} (\shE^\dual )= -\mu_{\rm min}(\shE)$.
If $\mu_{\min} (\shE) > \mu_{\max}(\shF)$, then $\Hom(\shE,\shF)=0$.
In particular, if $\mu_{\rm max} (\shF) < 0$, then
$\Gamma(Y, \shF)= \Hom (\O_Y, \shF)=0$.

If $\varphi: Z \ra Y$ is a finite $K$-morphism between smooth projective
curves over the algebrai\-cal\-ly closed field $K$, then
$\mu(\varphi^*(\shE))= \deg \, (\varphi ) \mu (\shE)$.
If $\varphi $ is separable, then
the pull-back $ \varphi^*(\shE)$ of a semistable sheaf
$\shE$ on $Y$ is again semistable, see \cite[Proposition 3.2]{miyaokachern}.
Hence in the separable case
the Harder-Narasimhan filtration of $\varphi^*(\shE)$ is just the
pull-back of the Harder-Narasimhan filtration of $\shE$,
and also $\mu_{\max}$ and $\mu_{\min}$ transform in the same way
as $\mu$ does.

In the non-separable case however this is not true at all
and the notion of semistability needs to be refined.
A locally free sheaf $\shE$
on $Y$ is called {\em strongly semistable}
if for every finite $K$-morphism $\varphi:Z \ra Y$ the pull-back
$\varphi^*(\shE)$ is again semistable. In characteristic zero,
this is the same as being semistable,
and in positive characteristic it is given by the property that
the pull-back under every $K$-linear Frobenius morphism is semistable,
see \cite[Proposition 5.1]{miyaokachern}.
This difficulty in positive characteristic
is one motivation for the following definition.

\begin{definition}
Let $Y$ denote a smooth projective curve over an algebrai\-cal\-ly closed
field and let
$\shE$ denote a locally free sheaf. Then we define
$$
\bar{\mu}_{\max} (\shE)=
{\rm sup} \{  \frac{\mu_{\rm max} (\varphi^* \shE)}{\deg \, (\varphi) }| \,\,
\varphi: Z \ra Y \mbox{ finite dominant $K$-morphism }      \}
$$
and
$$
\bar{\mu}_{\rm min} (\shE)=
{\rm inf} \{  \frac{\mu_{\rm min} (\varphi^* \shE)}{\deg \, (\varphi)} |\, \,
\varphi: Z \ra Y \mbox{ finite dominant $K$-morphism }      \} \, .
$$
\end{definition}

\begin{remark}
It is enough to consider in the previous definition
only $K$-linear Frobenius morphisms, since every morphism factors
through a separable map and a Frobenius and the maximal and minimal
slope behave well with respect to separable morphisms.
We will see in remark \ref{ampleremark} that $\bar{\mu}_{\min} (\shE)$
is bounded from below,
hence these numbers exist, but it is not clear whether they are obtained.
An equivalent question is whether one may find
a sufficiently high Frobenius power
such that the Harder-Narasimhan filtration of $\varphi^*(\shE)$ consists
of strongly semistable quotients.
\end{remark}

We will also need the following definition,
compare with \cite{ionescutoma}, \cite{langeregel}
and \cite{mukaisakai}
for related invariants of $\shE$ and $\PP(\shE)$.

\begin{definition}
\label{mindegdefinition}
Let $\shE$ denote
a locally free sheaf on a smooth projective curve,
let $1 \leq s \leq r =\rk \, (\shE)$. We set
$$\mindeg_s \, (\shE) = \min \{ \deg \, (\shQ): \, \shE \ra \shQ \ra 0, \, \,\,
\shQ \mbox{ is locally free and } \rk \, (\shQ) =s \} \, $$
and
$$ \maxdeg_s \, (\shE)
= \max \{ \deg(\shT):\, \shT \subseteq \shE
\mbox{ is a locally free subsheaf of rank } s \} \, .$$
\end{definition}

We have $\mindeg_s (\shE)= \deg (\shE) -\maxdeg_{r-s} (\shE)$
and $\maxdeg_s (\shE)= -\mindeg_{s} (\shE^{\dual})$.

\bigskip

\subsection{Ampleness criteria for vector bundles over projective curves}
\label{sectionample}

\markright{Ampleness criteria for vector bundles over projective curves}
\

\bigskip
Recall the following definition
for a locally free sheaf $\shG$ on a scheme $Y$.

\begin{definition}
A locally free sheaf $\shG$ on a scheme $Y$
is called {\em ample}
if the invertible sheaf $\O_{\PP(\shG)}(1)$ on the projective bundle
$\PP(\shG) = \Proj \, \oplus_{k \geq 0} S^k(\shG)$ is ample.
\end{definition}

For this notion we refer to
\cite{hartshorneample}, \cite{haramp} and \cite{lazarsfeldpositive}.
In characteristic zero we have the following
linear criterion of Hartshorne-Miyaoka for ample bundles over a curve.

\begin{theorem}
\label{amplekritnull}
Let $Y$ denote a smooth projective curve over an algebraically closed
field of characteristic $0$.
Let $\shG$ denote a locally free sheaf on $Y$ and
$\PP(\shG)$ the corresponding projective bundle,
and let $\xi$ denote the divisor class corresponding to $\O_{\PP(\shG)}(1)$.
Then the following are equivalent.

\renewcommand{\labelenumi}{(\roman{enumi})}
\begin{enumerate}

\item
The sheaf $\shG$ is ample.

\item
For every projective subbundle $\PP(\shQ) \subseteq \PP(\shG)$
of dimension $s$ we have that $\PP(\shQ) . \xi^s >0 $.

\item
For every locally free quotient sheaf $\shG \ra \shQ \ra 0$, $\shQ \neq 0$,
we have that $\deg \, (\shQ) > 0$.

\item
The minimal slope is $\mu_{\rm min}(\shG) >0$.
\end{enumerate}

\end{theorem}
\proof
The equivalence of (ii), (iii) and (iv) is clear
since we have $\PP(\shQ). \xi^s = (\xi| _{\PP(\shQ)})^s\, = \deg \, (\shQ)$.
(i) $\Rightarrow$ (ii) follows from the Nakai-criterion,
for the other direction see
\cite[Theorem 2.4]{hartshorneamplecurve},
\cite[Corollary 3.5]{miyaokachern} or \cite[Theorem 6.4.15]{lazarsfeldpositive}.
\qed

\bigskip
The minimal degree $\mindeg_1(\shG)$ of a quotient invertible sheaf of
$\shG$ must fulfill a stronger condition to guarantee that
$\shG$ is ample.

\begin{corollary}
\label{amplerankcrit}
Let $Y$ denote a smooth projective curve of genus $g$
over an algebraically closed field of characteristic $0$.
Let $\shG$ denote a locally free sheaf of rank $r$ on $Y$ and
suppose that
$\mindeg_1(\shG) > \frac{r-1}{r} g$.
Then $\shG$ is ample.
\end{corollary}
\proof
This is deduced in \cite[Proposition 2]{ionescutoma}
from the theorem.
\qed

\bigskip
The conditions (ii)-(iv) in Theorem \ref{amplekritnull}
are in positive characteristic not sufficient
for ampleness.
This is due to the fact that the pull back of a semistable sheaf
under a non-separable morphism need not be semistable anymore and the
minimal slope may drop.
It is also related to the failure of the vanishing theorem
in tight closure theory
in small positive characteristic. The following criterion is valid
in every characteristic and is essentially due to Barton
(see \cite[Theorem 2.1]{barton}).

\begin{theorem}
\label{amplekrit}
Let $Y$ denote a smooth projective curve over an algebraically closed
field $K$ of characteristic $p \geq 0$.
Let $\shG$ denote a locally free sheaf of rank $r$ on $Y$.
Then the following are equivalent.

\renewcommand{\labelenumi}{(\roman{enumi})}
\begin{enumerate}

\item
The sheaf $\shG$ is ample.

\item
$\bar{\mu}_{\rm min}(\shG) >0 $.

\item
There exists $\epsilon >0$ such that for every finite $K$-morphism
$\varphi :Y' \ra Y$ and every invertible quotient sheaf
$\varphi^*(\shG) \ra \shM \ra 0$
the inequality $\frac{\deg \, (\shM)}{ \deg \, (\varphi)} \geq \epsilon > 0$
holds.

\end{enumerate}

\end{theorem}

\proof
(i) $\Rightarrow $ (ii).
Suppose that
$\varphi: Y' \ra Y$ is finite, where $Y'$ is another
smooth projective curve,
and let $\varphi ^* (\shG) \ra \shH \ra 0$ be given,
$s= \rk (\shH)$.
We consider first the case $s=1$.
Then $\shM=\shH$ is an invertible sheaf on $Y'$ and
the surjection $\varphi ^* (\shG) \ra \shM \ra 0$
defines a section
$s: Y' \ra \PP(\varphi^*(\shG))$
due to the correspondence described in
\cite[Proposition 4.2.3]{EGAII}
and a curve $Z$ (its image)
in $\PP(\shG)$.
The map $Y' \ra Z$ factors through the normalization of $Z$,
hence $\shM$ is defined already on this normalization. Therefore
we may assume that $Y'$ is the normalization of $Z$.
The numerical class of the curve $Z$ in $\PP(\shG)$
can be written as
$a\xi^{r-1}+ b \xi^{r-2}.f$,
where $\xi$ is again the divisor class corresponding to $\O_{\PP(\shG)}(1)$,
$f$ is the class of a fiber $\PP(\shG) \ra Y$ and $a,b \in \ZZ$.
Furthermore $a$ equals the degree of $\varphi$.
Therefore we have
$$\deg \, (\shM) = Y' . s^*(\O_{\PP(\shG)} (1))
= ( a\xi^{r-1}+ b \xi^{r-2}.f) . \xi=a \xi^r +b
= a \deg \, (\shG) +b \, .$$
Hence
$\deg \, (\shM) / \deg \, (\varphi)  = \deg \, (\shG) + b/a $.
Since $\shG$ and hence $\xi$ is ample, there exists a number $\epsilon >0$
such that $ \xi.Z \geq \epsilon || Z||$ holds for every curve,
where $|| Z||$ denotes any norm on ${\rm Num} \otimes \RR$,
see \cite[Theorem 8.1]{haramp}. Hence
$$ \deg \, (\shG) +\frac{b}{a} = \frac{ \xi. Z}{a}
\geq \frac{\epsilon \sqrt{a^2+b^2}}{a } = \epsilon \sqrt{1+ (\frac{b}{a})^2}
\geq
\epsilon \, .
$$

Now consider the general case.
Since $\shG$ is ample, also its wedge product
$\bigwedge^s \shG$ is ample due to \cite[Corollary 2.6]{hartshorneample}.
The surjection
$\varphi^* (\shG) \ra \shH \ra 0$ yields a surjection
$\bigwedge^s(\varphi^* \shG) \ra \bigwedge^s \shH \ra 0$.
$\shM=\bigwedge^s \shH$ is invertible and since
$\bigwedge^s (\varphi^* \shG)$ is ample, there exists
an $\epsilon_s$ such that
$\deg \, (\bigwedge^s \shH)/ \deg \, (\varphi) \geq \epsilon_s >0$.
Then $\bar{\mu}_{\min} (\shG) \geq \min _s \epsilon_s/s >0$.

\smallskip
(ii) $\Rightarrow $ (iii) is a restriction to invertible
quotient sheaves.

\smallskip
(iii) $\Rightarrow $ (i).
Let $\xi$ denote the hypersection divisor corresponding to
$\O_{\PP(\shG)}(1)$.
Due to the ampleness criterion of Seshadri, see \cite[I\S7]{haramp},
it is enough to show that there exists an
$\epsilon >0$ such that
$\frac{ \xi.Z}{ m(Z)} \geq \epsilon >0$ holds for every
(effective) curve $Z$,
where $m(Z)= {\rm sup} \{ \mult_P(Z)\}$
is the maximal multiplicity of a point on $Z$.
So suppose that $Z$ is an irreducible curve in $\PP(\shG)$.
If $Z$ lies in a fiber $F \cong \PP^{r-1}$, then
$\xi . Z =\deg \, (Z) \geq m(Z)$.
Hence we may assume that $Z$ dominates the base.
Let $Z'$ be the normalization of $Z$, $i: Z' \ra  \PP(\shG)$ the corresponding
mapping and let $\varphi :Z' \ra Y$ be the composition.
Let $\varphi^* \shG \ra \shM \ra 0$ be the corresponding surjection
onto the invertible sheaf $\shM$.

The multiplicity $m(Z)$ is bounded above by $\deg \, (\varphi)$.
Therefore we have
$$ \frac{\xi.Z}{ m(Z)}
=\frac{\deg \, (\shM)}{m(Z)}
\geq \frac{\deg \, (\shM)}{ \deg \, (\varphi)}
\geq \epsilon > 0 \, .
$$
\qed

\begin{remark}
\label{ampleremark}
If $\shG$ is locally free on $Y$ and if $\O_Y(1)$ is an ample invertible sheaf
on $Y$,
then $\shG(n)= \shG \otimes \O_Y(n)$ corresponds to the invertible
sheaf $\O_{\PP(\shG)}(1) \otimes p^*(\O_Y(n))$ on $\PP(\shG)$,
see \cite[Proposition 4.1.4]{EGAII}.
Choosing $n$ high enough we may achieve that
$\shG(n)$ becomes ample. Since the slopes transform like
$\mu_{\rm min}( \shG \otimes \shL^k )
= \mu_{\rm min} (\shG) +k \deg \, (\shL) $ it follows that
$\bar{\mu}_{\min}(\shG)
= \bar{\mu}_{\rm min}(\shG \otimes \O_Y(k)) -k \deg \, (\O_Y(1))$
is bounded from below.
Dually it follows that $\bar{\mu}_{\max} (\shS)$ is bounded from above,
so both numbers exist (but it is not clear whether they are obtained).
\end{remark}

\begin{corollary}
\label{amplekritsemistable}
Let $Y$ denote a smooth projective curve over an algebra\-ical\-ly closed
field $K$.
Let $\shG$ denote a locally free sheaf on $Y$.
Suppose that $\shG$ is strongly semistable.
Then $\shG$ is ample if and only if its degree is positive.
\end{corollary}
\proof
This follows directly from \ref{amplekrit}.
For another proof see \cite[Corollary 3.5 and \S 5]{miyaokachern}.
\qed

\medskip
The ampleness of a locally free sheaf has also the
following consequence on $\mu_{\max}$, which we will use
in section \ref{sectioninclusion}.

\begin{corollary}
\label{mumaxabschaetzung}
Let $\shG$ denote an ample locally free sheaf of rank
$r$ on a smooth projective curve $Y$.
Then we have the estimates (set $\mindeg_0(\shG)=0 $)
$$ \mu_{\rm max}(\shG) \leq
{\rm max}_{s=0, \ldots, r-1}
\frac{ \deg \, (\shG) - \mindeg_s (\shG)}{r-s} 
\leq \deg \, (\shG)
\, .$$
\end{corollary}

\proof
Let $\shT \subset \shG$ denote a locally free subsheaf
of positive rank with a short exact sequence
$ 0 \ra \shT \ra \shG \ra \shQ \ra 0$,
$s= \rk (\shQ)$, $s=0, \ldots, r-1 $.
Then
$$\mu (\shT)= \frac{\deg \, (\shT)}{r-s}
=\frac{ \deg \, (\shG) -\deg \, (\shQ)}{r-s}
\leq \frac{ \deg \, (\shG) -\mindeg_s\, (\shG)}{r-s}
$$
due to the definition of
$\mindeg_s(\shE)$
as the minimum of $\deg \, (\shQ)$, where
$\shQ$ is a quotient sheaf of rank $s$.
Since $\shG$ is ample, it follows that $\mindeg_s (\shG) >0$, hence
the estimate with 
$\deg \, (\shG)$ follows.
(This last estimate holds also for
$\bar{\mu}_{\rm max}$.)
\qed

\medskip
If the bundle has rank two, the following theorem of Hartshorne-Mumford
gives a satisfactory criterion for ampleness also in positive
characteristic.

\begin{theorem}
\label{amplecritranktwo}
Let $Y$ denote a smooth projective curve of genus $g$ defined
over an algebraically closed field of characteristic $p >0$.
Let $\shG$ denote a locally free sheaf of rank two on $Y$.
Suppose that $ \deg \, (\shG) > \frac{2}{p} (g-1)$ and
$\deg \, (\shL) >0$ for every invertible
quotient sheaf $ \shG \ra \shL \ra 0$.
Then $\shG$ is ample.
\end{theorem}

\proof
See \cite[Proposition 7.5 and Corollary 7.7]{hartshorneample}.
\qed

\begin{remark}
Examples of semistable sheaves of positive degree
which are not ample (or which are not strongly semistable)
correspond to examples where the strong vanishing theorem for
tight closure does not hold in
small characteristic.
Look at the examples in \cite[\S3]{hartshorneamplecurve} and
\cite[Lemma 2.4]{giesekerample}) on the one hand and at
\cite[Example 3.11]{hunekesmithkodaira} on the other hand.
\end{remark}

\begin{remark}
There exist more criteria
for ample and very ample vector bundles on curves,
see for example \cite{alzati}, \cite{campanaflenner}, \cite{ionescutoma}.
We omit them, since they
don't seem to have implications on tight closure problems.
\end{remark}

\bigskip 

\subsection{Criteria for affineness}
\label{sectionaffine}

\markright{Criteria for affineness}
\

\bigskip
For ease of reference we fix the following situation,
which is a special case of problem \ref{cohoclassproblem}.

\begin{situation}
\label{extensionsit}
Let $Y$ denote a smooth projective curve over an algebraically closed
field $K$ and let
$\shS$ denote a locally free sheaf on $Y$ and let
$\shG = \shS^\dual$ be its dual sheaf.
A cohomology class $c \in H^1(Y, \shS)= \Ext^1(\O_Y, \shS)$ yields
a short exact sequence
$$ 0 \lra \shS \lra \shS' \lra \O_Y \lra 0 \, .$$
The dual sequence $ 0 \ra \O_Y \ra \shG' \ra \shG \ra 0$
yields a projective subbundle $\PP(\shG) \subset \PP(\shG')$ of codimension
one.
\end{situation}

Such a situation arises in particular from a homogeneous
$R_+$-primary tight closure problem in a homogeneous
coordinate ring $R$ of $Y$.
This tight closure point of view leads to
the question whether the complement
$\PP(\shG') - \PP(\shG)$ is affine or not, in dependence of
$\shS$ and $c \in H^1(Y, \shS)$.
Though we consider in this section mainly the case of characteristic $0$,
subsequent results in positive characteristic are discussed in section
\ref{sectionexclusion} and section \ref{sectionvanishing}.
The ampleness criterion \ref{amplekritnull}
yields at once the following affineness criterion.

\begin{theorem}
\label{ampleaffine}
Let the notation be as in {\rm \ref{extensionsit}} and suppose that
$\Char (K)=0$.
Suppose that $\shG$ is ample
\-- i.e. $\mu_{\rm min}(\shG) > 0$ \--
and that $c \neq 0$.
Then $\PP(\shG')- \PP(\shG)$ is affine.
\end{theorem}

\proof
Since $\shG$ is ample and since $0 \ra \O_Y \ra \shG' \ra \shG \ra 0$
is a non splitting short exact sequence,
the sheaf $\shG'$ is also ample due to
\cite[Theorem 2.2]{giesekerample} (here we use characteristic $0$).
Hence $\PP(\shG) \subset \PP(\shG')$
is an ample divisor and its complement is affine.
\qed

\begin{corollary}
\label{semistablepositiveaffine}
Let the notation be as in {\rm \ref{extensionsit}} and suppose that
$\Char (K)=0$.
Suppose that $\shG$ is semistable of positive degree
and that $c \neq 0$.
Then $\PP(\shG')- \PP(\shG)$ is affine.
\end{corollary}
\proof
This follows from \ref{ampleaffine}.
\qed

\medskip
Even if the divisor $\PP(\shG) \subset \PP(\shG')$
is not ample, the open subset $\PP(\shG') - \PP(\shG)$ may be affine.
We need the following lemma to obtain more general sufficient criteria
for the affineness of the complement.

\begin{lemma}
\label{affinelemma}
Let $Y$ denote a scheme and let $\shS$ and $\shT$
be locally free sheaves on $Y$ and let
$q: \shS \ra \shT$ be a sheaf homomorphism.
Let $c \in H^1(Y,\shS)$ with corresponding extension
$0 \ra \shS \ra \shS' \ra \O_Y \ra 0$
and let $q(c) \in H^1(Y, \shT)$ be its image
with corresponding extension
$0 \ra \shT \ra \shT' \ra \O_Y \ra 0$.

If $\PP((\shT')^\dual) - \PP( \shT^\dual)$ is affine, then also
$\PP((\shS')^\dual)-\PP(\shS^\dual)$ is affine.

\end{lemma}
\proof
This is for $Y$ Noetherian and separated
a special case of lemma \ref{monotonielemma}.
The same proof gives the result without any further condition.
\qed

\begin{theorem}
\label{slopekritaffin}
Let the notation be as in {\rm \ref{extensionsit}} and suppose that
$\Char (K)=0$.
Suppose that there exists a semistable sheaf $\shT$
of negative slope, $\mu (\shT) < 0$,
and a sheaf morphism $q: \shS \ra \shT $
such that $q(c) \neq 0$ in $H^1(Y, \shT)$.
Then the complement $\PP(\shG') - \PP(\shG)$ is affine.
\end{theorem}
\proof
The sheaf $\shH= \shT^\dual$ is
semistable of positive degree, hence ample
due to \ref{amplekritsemistable},
therefore $\PP(\shH')-\PP(\shH)$ is affine due to \ref{ampleaffine}.
Thus $\PP(\shG')-\PP(\shG)$ is also affine due to \ref{affinelemma}.
\qed

\medskip
The first candidates of semistable quotient sheaves to look at 
are $\shS/ \shS_{s-1}$ (the semistable quotient of minimal slope)
and invertible sheaves.

\begin{corollary}
\label{semistablequotientaffine}
Let the notation be as in {\rm \ref{extensionsit}} and suppose that
$\Char (K)=0$.
Suppose that $\mu_{\rm min}(\shS) <0$
{\rm(}equivalently $\mu_{\rm max}(\shG) >0${\rm )}
and that the image of
$c \in H^1(Y,\shS)$ in the semistable quotient $\shQ =\shS/ \shS_{s-1}$
is $\neq 0$.
Then $\PP(\shG') -\PP(\shG)$ is affine.
\end{corollary}

\proof
Since $\mu (\shQ)= \mu_{\min} (\shS) <0 $, the result
follows from \ref{slopekritaffin}.
\qed

\begin{corollary}
\label{quotientnegdegree}
Let the notation be as in {\rm \ref{extensionsit}} and
suppose that $\Char (K)=0$
{\rm(}for $\Char (K) >0$ see the following remark{\rm)}.
Suppose that there exists an invertible
sheaf $\shL$ of negative degree and a morphism
$q: \shS \ra \shL$
such that $0 \neq q(c) \in H^1(Y, \shL)$.
Then $\PP(\shG') -\PP(\shG)$ is affine.
\end{corollary}

\proof
This follows again from \ref{slopekritaffin}.
\qed

\begin{remark}
A sheaf homomorphism $q: \shS \ra \shL$ is equivalent to a section
$ \Gamma(Y, \shS^\dual \otimes \shL)$. Hence \ref{quotientnegdegree}
is only applicable if $\shG= \shS^\dual$
is not normalized, i.e. if there exist
non-trivial global sections of $\shG \otimes \shL$ for an invertible sheaf
of negative degree.
\end{remark}

\begin{remark}
\label{ampleposremark}
In positive characteristic the assumption in \ref{ampleaffine}
(and in the following corollaries \ref{slopekritaffin},
\ref{semistablequotientaffine} and
\ref{quotientnegdegree})
that $c \neq 0$ is to weak to ensure the ampleness of $\shG'$.
We need the stronger condition that
$0 \neq \varphi^*(c) \in H^1(Y', \varphi^*(\shS))$
for every finite morphism $\varphi: Y' \ra Y$.
In the situation coming from a forcing problem in tight closure,
this condition means that $f_0$ does not belong to the plus closure.

However, if $\shG$ is invertible and of positive degree,
then the extension $\shG'$ corresponding to $0 \neq c \in H^1(Y, \shG^\dual)$
is ample also in positive characteristic $p$
under the condition that $p \geq 2 (g-1)$, where $g$ is the genus
of the curve $Y$.
This follows from the ampleness criterion
of Hartshorne-Mumford for bundles of rank two, see
\ref{amplecritranktwo}.

From this it follows also that Corollary \ref{quotientnegdegree}
holds for positive characteristic $p \gg 0$.
For in this case the cohomology class $0 \neq q(c) \in H^1(Y, \shL)$
in the assumption of \ref{quotientnegdegree}
gives rise to an ample sheaf $(\shL')^\dual$ and then
the open subset $\PP( (\shL') ^\dual ) -\PP(\shL^\dual)$ is affine.
The affineness of $\PP(\shG') -\PP(\shG)$ follows then from \ref{affinelemma}.

If $\shG$ is indecomposable and ample of rank $r \geq 2$,
then the degree of
$\shG$ must fulfill stronger conditions to ensure the ampleness
of $\shG'$ defined by $c \in H^1(Y,\shS)$, $\shS= \shG^\dual$.
The result \cite[Theorem 25]{tango} suggests that the
right condition
may be $\deg \, (\shG) > r(r-1)(g-1)+ 2r (g-1)/p$
(or equivalently $\mu(\shG) > (g-1) (r+ p/2 -1)$).
Tango considers only the behavior of a cohomology class
$c \in H^1(Y,\shS)$ under the Frobenius, but along these lines
it should be possible to deduce also an ampleness criterion. 
\end{remark}

\begin{remark}
Another way to obtain results in positive characteristic is via a relative
setting such that $Y \ra \Spec\, D$ is a smooth projective relative curve
over a finitely generated $\ZZ$-domain $D$.
Then the generic fiber has characterstic zero and the special fibers have
positive characteristic.
Since the affineness (and the ampleness) is an open property it follows that
if we have established this property in the generic
situation $Y_{Q(D)}$, then it must also hold on the curves $Y_ \fop$
for almost all $\fop \in \Spec\, D$ or for $p \gg 0$.
Therefore the results in this section hold also in
characteristic $p \gg 0$.
We will come back to this in more detail in \ref{relsituation}.
\end{remark}

\smallskip 

\subsection{Criteria for $\PP(\shG') -\PP(\shG)$ not to be affine}
\label{sectionnonaffine}

\markright{Criteria for $\PP(\shG') -\PP(\shG)$ not to be affine}

\

\bigskip
We look again at the situation of \ref{extensionsit},
but now we look for criteria for $\PP(\shG') -\PP(\shG)$ to be not affine.
We first gather together some trivial but useful criteria.

\begin{proposition}
\label{affinemappingcrit}
Let the notation be as in {\rm \ref{extensionsit}}.

\renewcommand{\labelenumi}{(\roman{enumi})}
\begin{enumerate}

\item
Suppose that $\varphi: X \ra \PP(\shG')$ is an affine mapping such that
$X- \varphi^{-1} (\PP(\shG))$ is not affine.
Then $\PP(\shG') - \PP(\shG)$ is not affine.
This is in particular true for closed subschemes $X \subseteq \PP(\shG')$.

\item
Suppose that $\shG' \ra \shQ \ra 0$ is a surjection of
locally free sheaves. If
$\PP(\shQ) - \PP(\shQ) \cap \PP(\shG)$ is not affine,
then $\PP(\shG') - \PP(\shG)$ is not affine.

\item
Suppose that there exists a curve $X \subseteq \PP(\shG')$
which does not meet $\PP(\shG)$.
Then $\PP(\shG') - \PP(\shG)$ is not affine.

\item
Suppose that there exists a finite morphism
$\varphi: X \ra Y$ such that
$\varphi^*(c)= 0$. Then $\PP(\shG') - \PP(\shG)$ is not affine.

\end{enumerate}
\end{proposition}

\proof
(i) is trivial, since $X- \varphi^{-1}(\PP(\shG)')-
\PP(\shG))$ and since a closed embedding is an affine morphism.
(ii). A locally free quotient sheaf $\shQ$ defines a closed subscheme
$\PP(\shQ) \subseteq \PP(\shG')$, so this follow from (i).
(iii). A curve $X$ in $\PP(\shG')$ is projective. Hence, if the curve $X$
does not meet $\PP(\shG')$ this means that $\PP(\shG') - \PP(\shG)$
contains projective curves, therefore it cannot be affine.
(iv). 
The condition means that the pull back of
$ 0 \ra \O_Y \ra \shG' \ra \shG \ra 0$ splits on $X$.
The splitting surjection $  \varphi^*(\shG') \ra \O_X \ra 0$
yields a section $X \ra \PP(\varphi^*(\shG'))$ disjoined to
$\PP(\varphi^*(\shG))$ and hence there exists a projective curve
inside $\PP(\shG') - \PP(\shG)$.
\qed

\begin{remark}
If $\shG' \ra \shQ \ra 0$ is any quotient sheaf of $\shG'$,
then $\shQ / {\rm torsion}$ is torsion-free and hence locally free
on the smooth curve. Hence this somewhat more general
situation leads also to the situation
in \ref{affinemappingcrit}(ii).
If we have a section $0 \neq s \in \Gamma(\shG' \otimes \shM)$,
where $\shM$ is an invertible sheaf,
then the corresponding non-trivial mapping
$\shM^{-1} \ra \shG'$ induces a surjection
$\shG' \ra \shG'/ \shM^{-1}$, which yields again a projective
subbundle $\PP(\shH)$ of codimension one,
where $\shH = (\shG'/ \shM^{-1})/{\rm torsion}$.
\end{remark}

We will now look for conditions which imply
that the forcing divisor is not big
and hence its complement is not affine (see section
\ref{dimensiontwosection}).

\begin{lemma}
\label{muminpseudo}
Let $Y$ denote a smooth projective curve over an algebraically closed
field and let $\shT$ denote a locally free sheaf on $Y$.
Suppose that $\bar{\mu}_{\min} (\shT) \geq 0$.
Then $\mu_{\min} ( \Gamma^k (\shT)) \geq 0$ for $k \geq 0$, where
$\Gamma^k( \shT) = (S^k(\shT^\dual))^\dual $.
\end{lemma}
\proof
Note that $S^k(\shT) \cong (S^k(\shT^\dual))^\dual$
is only true in characteristic zero.
Assume that there exists a locally free quotient sheaf
$(S^k (\shT^\dual))^\dual \ra \shQ \ra 0$ of negative degree.
We find a finite morphism
$\varphi: Y' \ra Y$ such that we may write
$\varphi^*(\shQ) = \shL^k \otimes \shN $, where
$\shL$ is an invertible sheaf on $Y'$ with $\deg \, (\shL) < 0$ and also
$ \deg \, (\shN) <0$.
Due to our assumption we may assume that this is already true on $Y$.
We tensorize by $\shL^{-k}$ and get
$(S^k (\shT^\dual))^\dual \otimes \shL^{-k} \ra \shN \ra 0$.
But
$$(S^k (\shT^\dual))^\dual \otimes \shL^{-k}
=(S^k(\shT^\dual) \otimes \shL^k)^\dual
=(S^k(\shT^\dual \otimes \shL)   )^\dual
=(S^k( (\shT \otimes \shL^\dual)^\dual)^\dual \,.$$
Now $\deg \, (\shL^\dual) >0$, hence
$\bar{\mu}_{\min} (\shT \otimes \shL^\dual) >0$ and therefore
$\shT \otimes \shL^\dual$ is ample due to \ref{amplekrit}.
Then due to \cite[Theorem 6.6 and Proposition 7.3]{hartshorneample}
it follows that $(S^k( (\shT \otimes \shL^\dual)^\dual)^\dual$ is ample,
buts its quotient sheaf $\shN$ is not, since $\deg \, (\shN) <0$,
which gives a contradiction.
\qed

\begin{theorem}
\label{slopemaxkrit}
Let the notation be as in {\rm \ref{extensionsit}}.
Suppose that
$\bar{\mu}_{\rm max} (\shG ) \leq 0$
{\rm (}equivalently $ \mu_{\rm min}(\shS) \geq 0 ${\rm )}. 
Then $\PP(\shG')- \PP(\shG)$ is not affine.
\end{theorem}

\proof
From the sequence $0 \ra \O_Y \ra \shG' \ra \shG \ra 0$ it follows that
$\bar{\mu}_{\max} (\shG') \leq 0$ holds as well, since the image of
a mapping $\shT \ra \shG'$,
where $\shT$ is a semistable sheaf of positive slope,
must lie inside the kernel of $\shG' \ra \shG$.
Applying Lemma \ref{muminpseudo} to the dual of $\shG'$
it follows that $\mu_{\max} (S^k(\shG')) \leq 0$
for all $k \in \NN$.

If moreover $ \mu_{\rm max} (S^k (\shG')) < 0$,
then $\Gamma(Y, S^k(\shG'))=0$.
So we may suppose that $ \mu_{\rm max} (S^k (\shG')) = 0$.
Then there exists the maximal destabilizing sheaf
$0 \ra \shF_k \ra S^k (\shG') \ra S^k (\shG')/ \shF_k \ra 0$
such that $\shF_k$ is semistable with
$\mu (\shF_k)=0$ and $\mu_{\rm max}(S^k(\shG')/\shF_k) <0$,
so that this sheaf has again no global sections $\neq 0$.

We claim that for a semistable locally free sheaf $\shF$ of degree zero and
rank $s$ the dimension of the global sections is at most $s$:
This is true for invertible sheaves, so we do induction on the rank.
If $\shF$ has a global section $\neq 0$, then $\O_Y \subseteq \shF$
is a subsheaf.
We consider the saturation $\O_Y \subseteq \shM \subseteq \shF$
so that the cokernel $ \shF /\shM$ is torsion free, hence locally free on the
curve (see \cite[1.1]{huybrechtslehn}).
$\shM$ has degree $0$ (since $\O_Y \subseteq \shM$ and
since $\shF$ is semistable) and therefore we may apply the induction hypothesis
to the cokernel.

This gives the estimates (set $r= \rk (\shG)$ and $r+1= \rk (\shG')$)
$$
h^0 (S^k (\shG')) \leq h^0(\shF_k) + h^0(S^k (\shG')/\shF_k)
\leq {\rm rk} (\shF_k) \leq {\rm rk} (S^k (\shG'))= 
\binom{k+r}{r} \, . $$
This is a polynomial in $k$ with leading coefficient
$\frac{1}{r !} k^{r}$, thus it is bounded from above
by $\leq c k^{r}$ for some $c >0$.
Since the dimension of $\PP(\shG')$ is $r+1$ and
since $\pi_* \O_{\PP(\shG')}(k) =S^k(\shG')$ it follows that
$\O_{\PP(\shG')}(1)$ is not big and that $\PP(\shG')-\PP(\shG)$
is not affine by \ref{bigaffine}.
\qed

\begin{corollary}
\label{semistablenotaffine}
Let the notation be as in {\rm \ref{extensionsit}}.
Suppose that $\shG$ is strongly semistable with $\mu(\shG) \leq 0$.
Then $\PP(\shG')- \PP(\shG)$ is not affine.
\end{corollary}

\proof
Since $\shG$ is strongly semistable we have
$\bar{\mu}_{\max} (\shG) = \mu(\shG) \leq 0$, hence
the result follows from \ref{slopemaxkrit}.
\qed

\begin{corollary}
\label{invertibledegreecrit}
Let the notation be as in {\rm \ref{extensionsit}}
and suppose that $\shG$ is invertible.

If $\deg(\shG) \leq 0$, then $\PP(\shG')- \PP(\shG)$ is not affine.

If $\deg (\shG) >0$ and $0 \neq c \in H^1(Y, \shG^\dual)$,
then $\PP(\shG')- \PP(\shG)$ is affine {\rm (}if the cha\-racteristic
is zero or $p \gg 0${\rm )}.
\end{corollary}

\proof
The first statement follows from \ref{semistablenotaffine},
and the second from \ref{semistablepositiveaffine}
and remark \ref{ampleposremark} in positive characteristic.
\qed

\begin{theorem}
\label{splittingsemistable}
Let the notation be as in {\rm \ref{extensionsit}}.
Suppose that we have a decomposition
$$ \shS= \shS_1 \oplus \ldots \oplus \shS_s $$
into strongly semistable locally free sheaves $\shS_j$.
Let $c=(c_1, \ldots ,c_s)$ denote the decomposition of
the cohomology class $c \in H^1(Y, \shS)$.
Then the following are equivalent
{\rm(}for {\rm(i)} $\Rightarrow $ {\rm(ii)}
we have to assume that the
characteristic of $K$ is $0$ or $p \gg 0${\rm)}.

\renewcommand{\labelenumi}{(\roman{enumi})}
\begin{enumerate}

\item
There exists $1 \leq j \leq s$ such that
$\deg (S_j) <0$ and $c_j \neq 0$.

\item
$\PP(\shG')- \PP(\shG)$ is affine.

\end{enumerate}
\end{theorem}
\proof
(i) $\Rightarrow$ (ii).
Suppose that $j$ fulfills the statement in the numerical criterion.
Consider the projection $p_j :  \shS  \ra \shS_j$,
where the corresponding cohomological map sends $c$ to $c_j$.
Due to theorem \ref{semistablepositiveaffine} 
the complement of
$\PP(\shS^\dual_j) \subset \PP((\shS_j')^\dual)$
is affine.
Therefore $\PP(\shG')-\PP(\shG)$
is also affine due to \ref{affinelemma}.

\smallskip
(ii) $ \Rightarrow$ (i).
Assume that the numerical condition of (i) is not fulfilled.
We consider the direct summand
$$\shT = \oplus_{\deg \, (\shS_j) \geq 0}\, \shS_j \subseteq \shS \, .$$
Then all non-zero components of the forcing class $c$
belong to $\shT$, so this class comes from and goes to a cohomology class
in $H^1(Y, \shT)$. Hence we know by Lemma \ref{affinelemma}
that the affineness of $\PP((\shT')^\dual)-\PP(\shT^\dual)$
is equivalent to the affineness of $\PP(\shG') - \PP( \shG)$.

Since the components of $\shT$ are strongly semistable and
of non-negative degree it follows that
$\bar{\mu}_{\rm min} (\shT) \geq 0$
and therefore $\bar{\mu}_{\rm max} (\shT^\dual) \leq 0$, hence
$\PP((\shT')^\dual)-\PP(\shT^\dual)$ is not affine due to
Theorem \ref{slopemaxkrit}.
Then $\PP(\shG')-\PP(\shG)$ is also not affine.
\qed

\begin{corollary}
\label{splittinginvertible}
Let the notation be as in {\rm \ref{extensionsit}}.
Suppose that we have a decomposition
$\shS=  \shL_1 \oplus \ldots \oplus \shL_s$,
where $\shL_j$ are invertible sheaves.
Then $\PP(\shG')- \PP(\shG)$ is affine if and only if
there exists $1 \leq j \leq s$ such that
$\deg (L_j) <0$ and $c_j \neq 0$.
\end{corollary}
\proof
This follows directly from theorem \ref{splittingsemistable},
since invertible sheaves are strongly semistable.
\qed

\begin{corollary}
\label{projectivelineaffinetorsor}
Let the notation be as in {\rm \ref{extensionsit}}
and suppose that $Y= \PP^1$.
Then $\PP(\shG')- \PP(\shG)$ is affine if and only if $c \neq 0$.
\end{corollary}

\proof
Every locally free sheaf on the projective line is the direct sum
of invertible sheaves $\O_{\PP^1}(k)$.
Since $H^1(Y, \O_{\PP^1}(k)) = 0$ for $k \geq 0$, the result follows
from \ref{splittinginvertible}.
\qed

\begin{remark}
The last corollary corresponds to the property tat in the (regular)
two-dimensional polynomial ring every ideal is tightly closed.
\end{remark}

\begin{remark}
\label{remarkelliptic}
Suppose that $Y$ is an elliptic curve.
An indecomposable sheaf on an elliptic curve is strongly semistable,
hence we may always apply theorem \ref{splittingsemistable}
to get an answer whether $\PP(\shG') - \PP(\shG)$
is affine or not, see chapter \ref{tightelliptic}.
\end{remark}

\bigskip

\subsection{Starting an algorithm}
\label{algorithmsection}

\markright{Starting an algorithm}
\

\bigskip
In this section we bring together the results of the previous sections
to describe the first steps of an
``algorithm'' to decide whether an open subset
$\PP(\shG') -\PP(\shG)$ given by a cohomology class
$c \in H^1(Y, \shS)$, $\shS= \shG^\dual$,
is affine or not.
It gives always a complete answer if
$\shS$ is the extension of two semistable sheaves
and in particular if the rank of $\shS$ is two.
This implies that it is possible to decide
whether $f_0 \in (f_1,f_2,f_3)^*$ holds or not
(at least if we are able to compute the
Harder-Narasimhan filtration of the relation sheaf).
We assume that the characteristic of $K$ is zero.
The Harder-Narasimhan filtration of $\shS$,
$$0 = \shS_0 \subset \shS_1 \subset \ldots \subset
\shS_{s-1} \subset \shS_s =\shS \, ,$$
splits into short exact sequences
$$0 \lra \shS_{j-1} \stackrel{i}{\lra} \shS_j
\lra \shS_j/\shS_{j-1} \lra 0 \, ,$$
where the quotients $\shS_j /\shS_{j-1}$ are semistable with
slope $\mu_j (\shS)= \mu(\shS_j /\shS_{j-1})$.
The algorithm uses the fact that for a cohomology class
$c_j \in H^1(Y, \shS_j)$ we have either
$0 \neq \bar{c}_j \in H^1(Y, \shS_j/\shS_{j-1})$
or $c_j = i(c_{j-1})$, where $c_{j-1} \in H^1(Y,\shS_{j-1})$.

If $s=1$, then $\shS$ is semistable and everything is clear
by theorem \ref{semistablepositiveaffine} and
corollary \ref{semistablenotaffine}.
If $s=2$, then we have an exact sequence
$$0 \lra \shS_1 \lra \shS \lra \shQ \lra 0 \, ,$$
where $\shS_1$ and $\shQ$ are semistable of different slope.
In this case the algorithm gives also a complete answer.

We present the algorithm in the following diagram.
Numbers in bracket indicate the proposition on which
the implications in the diagram relies.
Note that for $s=2$ we have
$\mu_{s-1} (\shS)= \mu_1(\shS) = \mu_{\rm max} (\shS)$.
So if this number is $<0$, then we may conclude that
$\PP(\shG') - \PP(\shG)$ is affine by \ref{semistablepositiveaffine}.

\newpage

\small

\thispagestyle{empty}

\setlength{\unitlength}{1cm}

\begin{picture}(12,20.5)

\put(2.5, 20.3){$c \in H^1(Y, \shS) \mbox{ \ \ \ \ \ \ }$
(given by  $f_0$, $\shS= Rel(f_1, \ldots ,f_n)$ )}

\put(3,20.1){\line(0,-1){1}}

\put(0,19){$c=0$}
\put(6,19){$c \neq 0$}
\put(1,19.1){\line(1,0){4.9}}

\put(0.5, 18.8){\line(0,-1){15.8}}

\put(0.6, 18){\rm (\ref{affinemappingcrit})}

\put(6.2 ,18.8){\line(0,-1){1.7}}

\put(1 ,17){ $\mu_s(\shS) = \mu_{\rm min}(\shS) \geq 0 $}

\put(8 ,17){ $\mu_s(\shS) = \mu_{\rm min}(\shS) < 0 $}

\put(4.5 ,17.1){\line(1,0){3.4}}

\put(2.5,16.8){\line(-1,-1){2.0}}

\put(8.3 ,16.8){\line(0,-1){1.7}}

\put(8.4 ,16.3){$ \bar{c}= $}

\put(8.4 , 15.9){$ {\rm im}(c) \in H^1(Y, S/ S_{s-1}) $}

\put(1.5, 15.5){\rm (\ref{slopemaxkrit})}

\put(5,15){$\bar{c} =0$}

\put(11.5 ,15){$ \bar{c} \neq 0$}

\put(5.9, 15.1){\line(1,0){5.4}}

\put(5.3 ,14.8){\line(0,-1){1.7}}

\put(11.8 ,14.8){\line(0,-1){11.8}}

\put(5.4 ,14.2){$ c = {\rm im}(c_{s-1}) \in H^1(Y, S_{s-1}) $}

\put(5.4, 13.8){$ c_{s-1} \neq 0$}

\put(10.6,14.2){\rm (\ref{semistablequotientaffine})}

\put(2,13){$ \mu_{s-1}(\shS) \geq 0$}

\put(6.8,13){$ \mu_{s-1}(\shS) < 0$}

\put(4, 13.1){\line(1,0){2.7}}

\put(2.9,12.8){\line(0,-1){3.3}}

\put(1.5, 11.3){(\rm \ref{slopemaxkrit} on} 

\put(1.9, 10.9){$\shS_{s-1}$)}

\put(7.2 ,12.8){\line(0,-1){1.7}}

\put(8.8 ,13.1){\line(3,-1){3}}

\put(9.5,13){\rm (\ref{ampleaffine} if $s=2$)}

\put(7.2, 11.8){ $ \bar{c}_{s-1} \in H^1(Y,\shS_{s-1}/ \shS_{s-2})$}

\put(4.5,11){$\bar{c}_{s-1} =0$}

\put(8 ,11){$ \bar{c}_{s-1} \neq 0$}

\put(6. ,11.1){\line(1,0){1.9}}

\put( 8.5, 10.8){\line(0,-1){1.4}}

\put(5.3 ,10.8){\line(0,-1){3.3}}

\put(1.3, 9){$ \PP(( \shS_{s-1}')^\dual) - \PP( \shS_{s-1}^\dual)$}

\put(1.8,8.5){is not affine}

\put(7.5, 9){$ \PP(( \shS_{s-1}')^\dual) - \PP( \shS_{s-1}^\dual)$}

\put(8.2,8.5){is affine}

\put(8.5 ,8.3){\line(0,-1){5.3}}

\put(2.9,8.3){\line(-1,-2){2.4}}

\put(1.9,6){\rm (\ref{affinelemma})}

\put(4.2 ,7){$c_{s-1}= \im(c_{s-2})$}

\put(4.3 ,6.6){$ \in H^1(Y,\shS_{s-2})$}

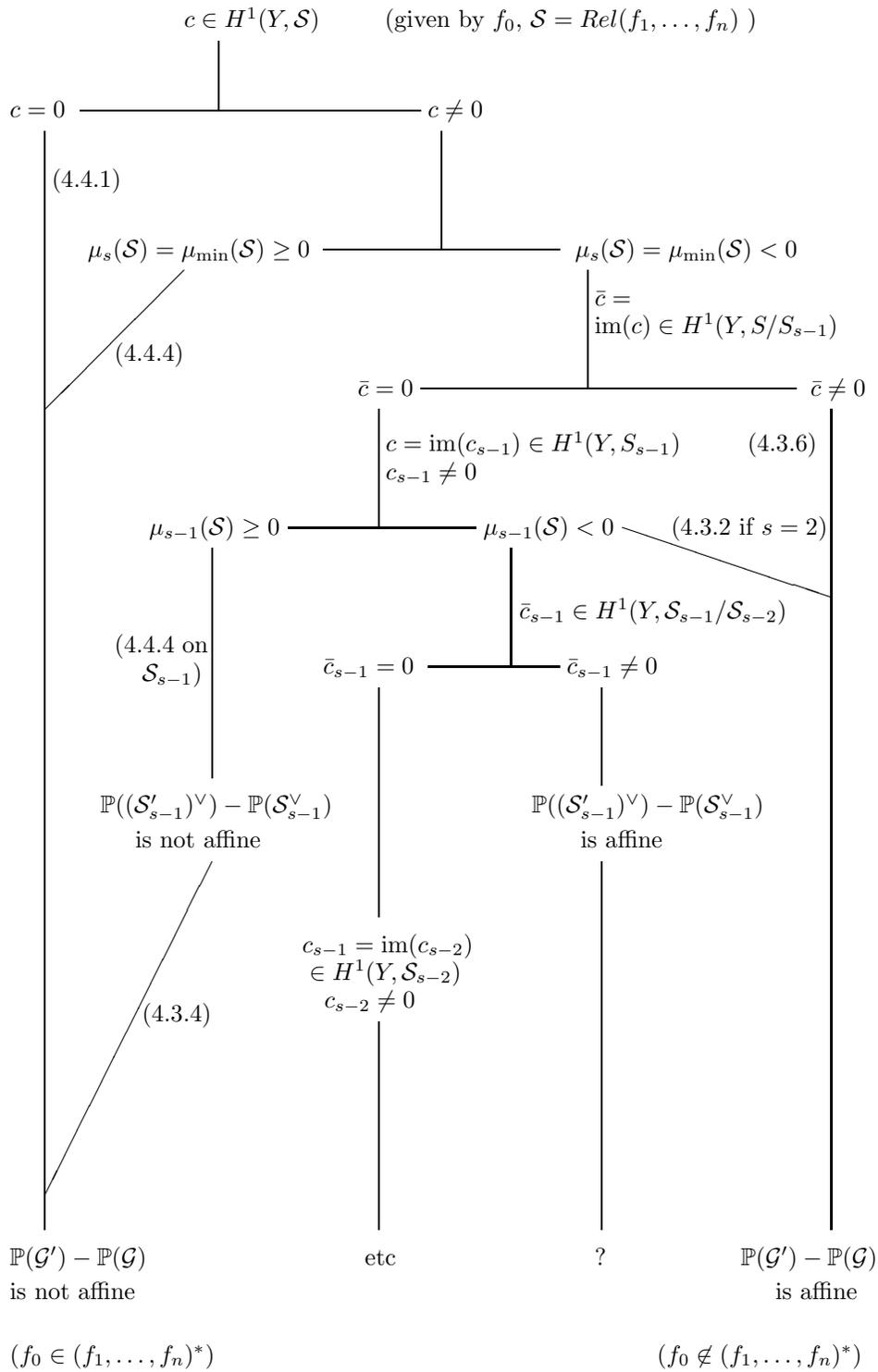
\put(4.4,6.2){  $c_{s-2} \neq 0$ }

\put(5.3, 6.0){\line(0,-1){3}}

\put(0,2.5) {$ \PP(\shG') - \PP(\shG)$}

\put(0,2){is not affine}

\put(5.1 ,2.5){etc}

\put(8.4 ,2.5){?}

\put(10.5,2.5){$ \PP(\shG') - \PP(\shG)$}

\put(11,2){is affine}

\put(0,1.1){($ f_0 \in (f_1, \ldots, f_n)^* $)}
\put(9.3,1.1){($ f_0 \not\in (f_1, \ldots, f_n)^* $)}
\end{picture}

\normalsize

\newpage

\thispagestyle{empty}

\markboth{5. Inclusion and exclusion bounds for tight closure}
{5. Inclusion and exclusion bounds for tight closure}

\newsavebox{\viereckv}

\savebox{\viereckv}(5,4)[bl]{
\thinlines
\put(0,0){\line(1,0){3.5}}
\put(0,0){\line(1,1){0.7}}
\put(.7,.7){\line(1,0){3.5}}
\put(3.5,0){\line(1,1){0.7}}
}

\newsavebox{\viereckfaserdrei}

\savebox{\viereckfaserdrei}(8,4)[bl]{

\put(0,0){\usebox{\viereckfaser}}
\put(1,0){\usebox{\viereckfaser}}
\put(2,0){\usebox{\viereckfaser}}
}

\newsavebox{\viereckfaserrand}

\savebox{\viereckfaserrand}(8,4)[bl]{
\put(0,0){\usebox{\viereckfaser}}
\put(3.5,0){\usebox{\viereckfaser}}
}

\section{Inclusion and exclusion bounds for tight closure}
\label{applications}

\bigskip
We shall now apply the results of the previous chapters to tight
closure problems. In particular we will give inclusion and exclusion
bounds for a homogeneous $R_+$-primary ideal in  a two-dimensional
normal standard-graded $K$-domain.
We derive these bounds from the minimal and the maximal slope of the
relation bundle associated to ideal generators $f_1, \ldots ,f_n$.
If this sheaf of relations is semistable, then the minimal
and the maximal slope coincide and therefore
the inclusion and the exclusion bound coincide,
and we get a strong vanishing theorem. We study in
detail the case of three elements $f_1,f_2,f_3$, which yields a
geometrically rich situation with various new phenomena.

\bigskip

\subsection{Inclusion bounds}

\label{sectioninclusion}

\markright{Inclusion bounds}
\

\bigskip
We fix the following situation.

\begin{situation}
\label{forcingsit}
Let $K$ denote a field with algebraic closure $\bar{K}$ and let
$R$ denote a two-dimensional standard-graded
$K$-algebra such that $R_{\bar{K}}$ is a normal domain.
Let $Y= \Proj\, R_{\bar{K}}$ denote the cor\-responding
smooth projective curve over $\bar{K}$, let $g$ denote its genus
and let $\delta = \deg \, (Y) = \deg \, (\O_Y(1))$ denote the degree
of the very ample invertible sheaf $\O_Y(1)$ on $Y$.
Let $f_1, \ldots ,f_n \in R$ be homogeneous $R_+$-primary elements
of degree $d_i$.
Let $\shR(m)$ be the sheaf of relations
of total degree $m$ on $Y$ and let $\shF(-m)$ be its dual sheaf.

Let $f_0$ denote another homogeneous element of degree $d_0$,
let $\shR'(m)$ be the sheaf of relations for the elements
$f_1, \ldots, f_n,f_0$ of total degree $m$and let
$$ 0 \lra \shR(m) \lra \shR'(m) \lra \O_Y \lra 0 $$
be the corresponding forcing extension and
let $c \in H^1(Y, \shR(m))$ be the corresponding forcing class
defined by $f_0$.
The corresponding surjection
$\shF'(-m) \ra \shF(-m) \ra 0$ yields the embedding
$\PP(\shF(-m)) \subset \PP(\shF'(-m))$.
\end{situation}

\begin{definition}
Let the notation be as in \ref{forcingsit}.
Then we set
$$ \mu_{\rm max}(f_1, \ldots ,f_n) := \mu_{\rm max} (\shF (0)) \, $$
and also for $\mu$, $\mu_{\min}$, $\bar{\mu}_{\max}$ and $\bar{\mu}_{\min}$.
\end{definition}

\begin{remark}
Note that we consider for ideal generators $f_1, \ldots ,f_n \in R$
always the slope of the relation bundle after replacing
$K$ by $\bar{K}$. Since changing the base field does not affect
the affineness of open subsets, it does not affect solid closure.
The slope of $\shF (0)$ is \-- due to the presenting sequence for $\shF(0)$ \--
given by
$$\mu(\shF(0)) = \frac{d_1+ \ldots +d_n}{n-1} \delta  \, ,$$
therefore we get the estimates
$$ \mu_{\rm min}(f_1, \ldots ,f_n) \leq
\frac{d_1+ \ldots +d_n}{n-1} \delta \leq 
\mu_{\rm max}(f_1, \ldots ,f_n) \, .$$
Equality holds if and only if $\shF(-m)$ is semistable.
Furthermore we have
$$\mu(\shF(-m))=\mu (\shF(0) \otimes \O_Y(-m)) =  \mu(\shF(0)) -m \delta 
= ( \frac{d_1+ \ldots +d_n}{n-1} -m )\delta  \, ,$$
and the same rule holds for $\mu_{\rm max}$ etc.
\end{remark}

From the conditions in section \ref{sectionnonaffine}
we derive the following
numerical condition that elements of sufficiently
high degree must belong to the tight closure.

\begin{theorem}
\label{maxin}
Let the situation and notation be as in {\rm \ref{forcingsit}}.
Suppose that
$ \deg \,(f_0) \geq  \frac{1}{\delta} \bar{\mu}_{\rm max} (f_1, \ldots ,f_n)$.
Then $f_0 \in (f_1, \ldots,f_n)^\soclo$.
\end{theorem}
\proof
Let $m = \deg \, (f_0)$.
The condition means that
$$ \bar{\mu}_{\rm max} (\shF (-m))
=\bar{\mu}_{\rm max} (f_1, \ldots, f_n) -m \delta \leq 0 \, .$$
Hence the result follows from \ref{slopemaxkrit}
and \ref{tightaffinecrit}.
\qed

\medskip
To obtain criteria for tight closure membership
we need bounds from above for $\bar{\mu}_{\max} (f_1, \ldots ,f_n)$.
The next proposition gives a general bound for $\bar{\mu}_{\rm max}$.
We will give better bounds under additional
conditions in the next sections.

\begin{proposition}
\label{mumaxabschaetzung2}
Let the notation be as in {\rm \ref{forcingsit}}.
Suppose that the degrees are ordered
$1 \leq d_1 \leq d_2 \leq \ldots \leq d_n$.
Let $\shE= \oplus_i \O_Y(d_i)$ and let
$0 \ra \O_Y \ra \shE \ra \shF(0) \ra 0$ be the presenting sequence
for $\shF(0)$.
Then we have the estimate
$$ \mu_{\rm max}(f_1, \ldots ,f_n)
\leq
\delta \cdot (d_{n-1}+d_n) \, .$$
The same is true for $\bar{\mu}_{\max}$.
\end{proposition}

\proof
Set $\shF=\shF(0)$.
Corollary
\ref{mumaxabschaetzung} together with the inequality
$\mindeg_s(\shF) \geq \mindeg_s(\shE)$
yields
$$
\mu_{\rm max} (\shF)
\leq
\max_{s=0, \ldots ,n-2} \frac{ \deg \, (\shF) -\mindeg_s\, (\shF) }{n-1-s}
\leq
\max_{s=0, \ldots ,n-2}
\frac{ \deg \, (\shE) -\mindeg_s\, (\shE) }{n-1-s} \, .
$$
We claim that
$\mindeg_s \, (\shE) = \delta(d_1 + \ldots +d_s)$.
Since $\O_Y(d_1) \oplus \ldots \oplus \O_Y(d_s)$
is a quotient sheaf of rank $s$, the estimate $\leq$
is clear.
For the other estimate we consider first the case
$s=1$, so suppose that $\shQ$ is an invertible sheaf.
If $\shQ$ is a quotient of $\shE$, then
$\Hom( \O_Y(d_i), \shQ)=H^0(Y, \O_Y(-d_i) \otimes \shQ) \neq 0$ for
at least one $i$. Therefore
$\deg \, (\shQ) \geq \delta {\rm min}_i \, (d_i) = \delta d_1$.

Now suppose that $\shQ$ is a locally free quotient of $\shE$ of rank $s$.
Then we have a surjection
$$\oplus_{i_1 < \ldots < i_s}\,
\O_Y(d_{i_1}) \otimes \ldots \otimes \O_Y(d_{i_s})
\cong
\bigwedge^s \shE
\lra \bigwedge^s \shQ= \det \shQ \, .
$$
Due to the case $s=1$ we know
$\deg \, (\shQ) \geq \delta (d_1 + \ldots + d_s)$, which proves the claim.

Thus we have the estimate
$$
\mu_{\rm max} (\shF)
\leq
\max_{s=0, \ldots ,n-2}
\frac{ \deg \, (\shE) -\mindeg_s\, (\shE) }{n-1-s}
=\max_{s=0, \ldots ,n-2}
\delta \frac{d_{s+1}+ \ldots + d_n}{n-1-s}
 \, .
$$
Here the term for $s=n-2$, which is $\delta ( {d_{n-1} +d_n})$,
is maximal.

If $\varphi: Y' \ra Y$ is a finite morphism,
then the situation is preserved under the pull-back
(even if $\varphi^* \O_Y(1)$ is not very ample anymore).
Then
${\mu}_{\max}(\varphi^*(\shF))
\leq (\delta \cdot \deg \, (\varphi)) (d_{n-1} + d_n) $
and hence the inequality holds also for
$\bar{\mu}_{\max}$.
\qed

\begin{corollary}
\label{inclusionbound}
Let the notation be as in {\rm \ref{forcingsit}}.
Suppose that the degrees are ordered
$1 \leq d_1 \leq d_2 \leq \ldots \leq d_n$.
Then for $m \geq d_{n-1} + d_n$ we have the inclusion
$$ R_{m} \subseteq (f_1, \ldots ,f_n)^\soclo  \, .$$
\end{corollary}
\proof
This follows from \ref{mumaxabschaetzung2} and \ref{maxin}.
\qed

\begin{remark}
This corollary is for the two-dimensional case
a somewhat better estimate than
the estimate $2 \max_i (d_i)$
(and $d_1+ \ldots + d_n$) found by
K. Smith, see \ref{smithinclusionbound}, \ref{exclusionbound}
or \cite[Proposition 3.1 and Proposition 3.3]{smithgraded}.
Smith´s estimate for $R_+$-primary ideals
is however true in every dimension
and holds also for the plus closure
in positive characteristic.
\end{remark}

\begin{remark}
An estimate from below for $\mu_{\max}$ is given by
$\mu_{\rm max} (\shF(0)) \geq \delta \cdot {\rm max}_i (d_i)$.
To see this consider again the sequence
$0 \ra \O_Y \ra \oplus_i\,  \O_Y(d_i) \ra \shF(0) \ra 0$.
If $\O_Y(d_i) \ra \shF(0) $ is the zero map for one $i$, then
$\O_Y \cong \O_Y(d_i)$, $d_i=0$ and $f_i$ is a unit and the statement is
clear from $\shF(0) \cong \oplus_{j \neq i} \O_Y(d_j)$.
Otherwise $\O_Y(d_i) \ra \shF(0) $ is not the zero map, and then
$\mu_{\rm max} (\shF(0)) \geq \mu(\O_Y(d_i)) = \delta d_i$ for all $i$.
\end{remark}

There exists also the following inclusion criterion
to an ideal itself in terms of the slope of the relation sheaf.

\begin{proposition}
Let the notation be as in {\rm \ref{forcingsit}}.
Suppose that
$\deg (f_0) >
\frac{1}{\delta} \mu_{\rm max} (f_1, \ldots, f_n) +\frac{2g-2}{\delta}$.
Then $f_0 \in (f_1, \ldots, f_n)$.
\end{proposition}

\proof
Let $m= \deg \, (f_0)$. It suffices to show that
$H^1(Y, \shR(m))=0$. This is by Serre duality equivalent with
$H^0(Y, \shF(-m) \otimes \omega_Y) =0$. Now
\begin{eqnarray*}
\mu_{\rm max} (\shF(-m) \otimes \omega_Y) &=&
\mu_{\rm max} (\shF(-m)) + \mu(\omega_Y) \cr
&=& \mu_{\rm max} (f_1, \ldots , f_n) - m\delta + 2g-2 \cr
&<& 0
\end{eqnarray*}
by the numerical condition, thus
$\shF(-m) \otimes \omega_Y$ has no global sections $\neq 0$.
\qed

\smallskip 

\subsection{Exclusion bounds for tight closure}
\label{sectionexclusion}

\markright{Exclusion bounds for tight closure}
\

\bigskip
We are now looking for exclusion bounds for tight closure.
A number $a \in \NN$ is an exclusion bound
for $(f_1, \ldots ,f_n)$ if the following holds:
If $\deg (f_0) <a$, 
then $f_0 \in (f_1,...,f_n)^\soclo$ if and only if $f_0 \in (f_1, \ldots, f_n)$.

The theorems in this and the next section hold either in characteristic zero
or in positive characteristic $p$ under the condition that $p  \gg 0$.
To make sense of this statement we have to
suppose that everything is given relatively to a base scheme
such that the generic fiber has characteristic zero and the special fibers
have positive characteristic.
For this we fix the following situation, see
also \cite[Definition 3.3]{hunekesmithkodaira}
and in particular the appendix of Hochster in \cite{hunekeapplications}
for this setting.

\begin{situation}
\label{relsituation}
Let $D$ denote a finitely generated normal $\ZZ$-domain of dimension one.
Let $S$ denote a standard-graded flat $D$-algebra
such that for all $\fop \in \Spec \, D$ the algebras
$S_{\kappa(\fop)} = S \otimes_D \kappa(\fop) $
are two-dimensional geometrically normal standard-graded
$\kappa(\fop)$-algebras (so that $S_{\bar{\kappa}(\fop)}$ are normal domains).
For $\fop =0$ this is an algebra over the quotient field $Q(D)$ of
characteristic zero, and for a maximal ideal $\fop$
the algebra $S_{\kappa(\fop)}$ is an algebra over the finite field
$\kappa(\fop)=D/ \fop$ of positive characteristic.

We suppose that we have $S_+$-primary homogeneous elements
$f_1, \ldots, f_n \in S$ of degree $d_i$ and another
homogeneous element $f_0$.
Let $B$ denote the forcing algebra over $S$ for this data and let
$U= D(S_+) \subseteq \Spec \, B$.
These elements yield homogeneous forcing data for every $S_{\kappa(\fop)}$
and $B_{\kappa(\fop)}=B \otimes_D \kappa(\fop)$
is the corresponding forcing algebra.
For every prime ideal $\fop \in \Spec\, D$ the affineness of
$U_\fop =U_{\kappa(\fop)}= U \cap \Spec \, B_{\kappa(\fop)}$
is equivalent with
$f_0 \not\in (f_1, \ldots ,f_n)^\soclo$ in $S_{\kappa(\fop)}$.

We denote by $Y= \Proj\, S$ the smooth projective (relative) curve
over $\Spec\, D$
and by $\delta$ the common degree of the curves $Y_\fop$, $\fop \in \Spec\, D$.
We denote by $\shR(m)$ the sheaf of relations on $Y$ and by $\shF(-m)$
its dual sheaf and we denote the restrictions to $Y_\fop$
by $\shR_\fop(m)$ and $\shF_\fop(-m)$.
The different notions of slopes and of semistability refer always
to $\shF_{\bar{\kappa}(\fop)}$ on $Y_{\bar{\kappa}(\fop)}$.
The element $f_0$ yields an extension
$0 \ra \shR(m) \ra \shR'(m) \ra \O_Y \ra 0$
and a subbundle $\PP(\shF) \subset \PP(\shF')$, which induces
the projective subbundle $\PP(\shF) \subset \PP(\shF')$
on every curve for every $\fop \in \Spec\, D$.
\end{situation}

\begin{remark}
\label{affineopen}
We will apply several times the following conclusion:
let $D \subseteq B$ denote Noetherian domains and let
$U=D(\foa) \subseteq \Spec\, B$ denote an open subset,
$\foa=(a_1, \ldots ,a_k)$.
Suppose that $U_\eta = U \cap \Spec \, (B \otimes_D \kappa(\eta))$
is affine, where $\eta$ denotes the generic point of $\Spec \, D$.
This means that there exist rational functions
$q_j \in \Gamma(U_\eta , \O_{\eta})$ such that $\sum q_ja_j =1$
(\ref{affinefunctionsremark}).
We find a common denominator $0 \neq g \in D$ such that
these functions $q_j$ are defined on $U \cap D(g) \subseteq \Spec\, B$,
hence also $U \cap D(g)$ is affine.
This means that after shrinking $D$ (i. e. replacing $\Spec D$ by $\Spec D_g$)
we may assume that $U$ itself is affine.
Then for every $P \in \Spec \, D$ the fibers $U_{\kappa(P)}$ are affine.

Suppose in the situation \ref{relsituation}
that $f_0 \not\in (f_1, \ldots ,f_n)^\soclo$ holds over the generic point
$\eta \in \Spec \, D$. This means that the open subset
$U_\eta$ is affine.
Then after shrinking $D$ we may assume that $U$ is affine,
hence that every fiber $ U_{\kappa(\fop)}$ is affine.
This means that $f_0 \not\in (f_1, \ldots ,f_n)^\soclo$ holds in $S_{\kappa(\fop)}$
for all $\fop \in \Spec \, D$ (or for almost all $\fop \in \Spec \, D$
for the old $D$). In this case we say briefly that
$f_0 \not\in (f_1, \ldots ,f_n)^\soclo$ holds for $p \gg 0$.
\end{remark}

\begin{theorem}
\label{minex}
Let the notation be as in {\rm\ref{forcingsit}} and in {\rm \ref{relsituation}}.
Suppose that the characteristic of $K$ is $0$ or $p  \gg  0$.
If $\deg \, (f_0) < \frac{1}{\delta} \mu_{\rm min} (f_1, \ldots ,f_n)$,
then $f_0 \in (f_1, \ldots, f_n)^\soclo$
if and only if $f_0 \in (f_1, \ldots ,f_n)$.
\end{theorem}
\proof
Let $m =\deg \, (f_0)$. Suppose first that the characteristic is zero.
We may assume that $K$ is algebraically closed.
The condition means that the minimal slope
$ \mu_{\rm min}(\shF (-m) ) 
= \mu_{\rm min} (f_1, \ldots, f_n) -m \delta > 0$,
hence $\shF(-m)$ is ample due to \ref{amplekritnull}.
Suppose that $f_0 \not\in (f_1, \ldots ,f_n)$.
This means by \ref{trivialtwo} that the
corresponding forcing class is $c \neq 0$.
Hence $\PP(\shF') -\PP(\shF)$ is affine due to \ref{ampleaffine} and
$f_0 \not\in (f_1, \ldots ,f_n)^\soclo$ due to \ref{tightaffinecrit}.

Now suppose the relative situation \ref{relsituation}.
Note that the slope condition is imposed on the
generic fiber, i.e.
$\deg \, (f_0) < \frac{1}{\delta}(f_1, \ldots, f_n)
=\frac{1}{\delta} \mu_{\rm min}(\shF_\eta (0))$.
We have to show that $f_0 \not\in (f_1, \ldots ,f_n)$
implies $f_0 \not\in (f_1, \ldots ,f_n)^\soclo$ for almost
all $\fop \in \Spec\, D$.
From $f_0 \not\in (f_1, \ldots ,f_n)$ in $S_{\kappa(\fop)}$
it follows $f_0 \not\in (f_1, \ldots ,f_n)$ in $S$
and by shrinking $D$ we may assume that
$f_0 \not\in (f_1, \ldots ,f_n)$ in $S_{Q(D)}$.
From the case of characteristic zero
we know that $U_\eta $ is affine
and the result follows from remark \ref{affineopen}.
\qed

\medskip
From the bounds proved in Theorem \ref{maxin} and in
Theorem \ref{minex} it is easy to derive the following result
of Huneke and Smith (see \cite[Theorem 5.11]{hunekesmithkodaira}).

\begin{corollary}
\label{projdim}
Let the notation be as in {\rm\ref{forcingsit}}
{\rm(}or {\rm \ref{relsituation})}.
Suppose that the characteristic of $K$ is $0$
or $p  \gg  0$.
Suppose that the projective dimension of $R/(f_1, \ldots, f_n)$
is $2$ or equivalently that the sheaf of relations has a splitting
$\shR(0) \cong \O_Y(-a_1) \oplus \ldots \oplus \O_Y(-a_r)$.
Set $a= \max \{ a_i\}$ and $b= \min \{a_i \}$.
Then
$$ R_{\geq a} \subseteq (f_1, \ldots ,f_n)^\soclo 
\mbox{ and }
(f_1, \ldots ,f_n)^\soclo  \subseteq (f_1, \ldots , f_n) + R_{ \geq b}\, .$$
\end{corollary}
\proof
Whenever $\shF(0)$ is a direct sum of invertible sheaves $\shL_j$ we have
that $\bar{\mu}_{\rm max} (\shF(0)) = \max_j \deg \, (\shL_j)$
and $\bar{\mu}_{\min} (\shF(0)) = \min_j \deg \, (\shL_j)$.
So for $\shF(0)= \O_Y(a_1) \oplus \ldots \oplus \O_Y(a_r)$ we find
that $\bar{\mu}_{\max}(\shF(0)) =a \delta$ and
$\bar{\mu}_{\min} (\shF(0))= b \delta $,
so the result follows from \ref{maxin} and \ref{minex}.
\qed

\medskip
If the sheaf of relations on the projective curve splits into
invertible sheaves as in the previous Corollary
\ref{projdim} it is easy to give a numerical criterion for tight closure.

\begin{corollary}
\label{splitting}
Let the notation be as in {\rm\ref{forcingsit}}
{\rm(}or {\rm\ref{relsituation})}.
Suppose that we have a decomposition
$\shR(0)= \shL_1 \oplus \ldots \oplus \shL_r$,
where $\shL_j$ are invertible sheaves.
Let $m= \deg\, (f_0)$ and let $c \in H^1(Y, \shR(m))$
be the forcing class with components
$c_j \in H^1(Y ,\shL_j \otimes \O_Y(m))$.
Suppose that the characteristic is zero or $p \gg 0$.
Then $f_0 \in (f_1,\ldots ,f_n)^\soclo$
if and only if $ \,\deg \, (\shL_j) +m \delta \geq 0 $ or $c_j=0$ holds
for all $1 \leq j \leq r$.
\end{corollary}
\proof
This follows from \ref{splittinginvertible}.
\qed

\begin{remark}
The situation of Theorem \ref{splitting} holds for every primary
homogeneous ideal in $K[x,y]$ (due to the splitting theorem of Grothendieck,
see \cite[Theorem 2.1.1]{okonekvector}),
but for a polynomial ring the computation of tight closure does not make
much problems, so \ref{splitting} gives a help which we do not need in
this case. However, the splitting situation holds also if
$I \subseteq R$ is the extended ideal $I=JR$ of an ideal
$J \subseteq K[x,y] \subset R$, and in this case it is also useful for
computations, see example \ref{indepexample} below.
\end{remark}

\begin{remark}
Let the notation be as in \ref{forcingsit}. From the sequence
$$ 0 \lra \O_Y \lra \O_Y(d_1) \oplus \ldots \oplus \O_Y(d_n) \lra
\shF(0) \lra 0 \, $$ we also get the estimate
$$\bar{\mu}_{\rm min} (\shF(0)) \geq \bar{\mu}_{\rm min} (\oplus_i \O_Y(d_i))
={\rm min}_i\, \mu(\O_Y(d_i))
= \delta {\rm min}_i (d_i) \, .$$
So if $f_0 \neq 0$ and
$ \deg \, (f_0) < {\rm min}_i (d_i)$, then $f_0 \not\in (f_1, \ldots , f_n)$
and due to \ref{minex} also $f_0 \not\in (f_1, \ldots , f_n)^\soclo$.
\end{remark}

The minimal slope in the generic point gives also a bound
for the minimal slope in positive characteristic for $p \gg 0$.
The same is true for the maximal slope.

\begin{proposition}
\label{compare}
Suppose the relative situation {\rm \ref{relsituation}},
let $Y= \Proj  S \ra \Spec  D$ denote the smooth projective
curve of relative dimension one and suppose
that the generic curve $Y_\eta=Y_{Q(D)}$ has at least one
$Q(D)$-rational point.
Let $\mu_{\rm min} (\shF_\eta (-m))$
denote the minimal slope and let $\mu_{\rm max}(\shF_\eta (-m)) $
denote the maximal slope on $Y_{\bar{\kappa}(\eta)}$.
Then for $p \gg 0$ we have the bounds
$$ \bar{\mu}_{\rm min} (\shF(-m))
> \lceil \mu_{\rm min} (\shF_\eta(-m))  \rceil -1 \, $$
and
$$ \bar{\mu}_{\rm max} (\shF(-m))
< \lfloor \mu_{\rm max}(\shF_\eta (-m))  \rfloor +1 \, . $$
\end{proposition}
\proof
First we may assume by shrinking $D$ that there exists a section
for the relative curve $Y$. Hence there exists an invertible sheaf $\shM$
on $Y$ such that the degree of $\shM$
on every fiber $Y_\fop$ is one.
Let $\shL$ denote an invertible sheaf on $Y$ which has on every fiber
$Y_\fop$ the degree
$- \lceil \mu_{\rm min}(\shF_\eta (-m)) \rceil +1$.

Then
$ \mu_{\rm min}(\shF_\eta (-m) \otimes \shL_\eta)
= \mu_{\rm min}(\shF_\eta(-m))
- \lceil \mu_{\rm min}(\shF_\eta (-m))  \rceil +1 >0$
and hence $\shF_\eta (-m) \otimes \shL_\eta $ is ample on $Y_\eta$
due to \ref{amplekritnull}.
Therefore $\shF_\fop (-m) \otimes \shL_\fop $ is ample on
$Y_\fop$ for $p \gg 0$, since ampleness is an open property
(see \cite[Th\'{e}or\`{e}me 4.7.1]{EGAIII}
or \cite[Theorem 1.2.13]{lazarsfeldpositive}; since we need here only
the generic open property, we can also use \ref{affineopen} together with
\cite[Th\'{e}or\`{e}me 4.5.2] {EGAII}).
This means again by \ref{amplekrit} that
$\bar{\mu}_{\rm min}( \shF_\fop(-m) \otimes \shL_\fop) >0$ or that
$$\bar{\mu}_{\rm min}( \shF_\fop(-m)) > - \deg \, (\shL_\fop)
= \lceil \mu_{\rm min}(\shF_\eta (-m))  \rceil - 1 \, .$$

This gives the first result.
The second statement follows by applying the first statement to $\shR(m)$,
\begin{eqnarray*}
\bar{\mu}_{\rm max} (\shF(-m)) &= &-\bar{\mu}_{\rm min} (\shR(m)) \cr
&<& -( \lceil \mu_{\rm min}(\shR_\eta (m))  \rceil -1) \cr
&=& \lfloor -\mu_{\rm min}(\shR_\eta (m)) \rfloor  +1 \cr
&=& \lfloor \mu_{\rm max}(\shF_\eta (-m)) \rfloor +1 \, .
\end{eqnarray*}
\qed

\bigskip  

\subsection{Vanishing theorems for tight closure}

\markright{Vanishing theorems for tight closure}

\label{sectionvanishing}

\

\bigskip
The inclusion bound in Theorem \ref{maxin} and the exclusion bound in
Theorem \ref{minex} coincide if the sheaf of relations is semistable.
This gives a new class of vanishing type theorems (or strong bound theorem)
in dimension two
and generalizes the strong vanishing theorem for para\-meter ideals,
see \cite{hunekesmithkodaira} and Corollary \ref{parametervanishing} below
(the name vanishing is due to the fact
that it is related to the Kodaira vanishing theorem,
see \ref{kodairavanishing}).
We give first the formulation in zero characteristic.

\begin{theorem}
\label{semistablevanishing}
Let the notation be as in {\rm \ref{forcingsit}}
and suppose that the characteristic of $K$ is zero.
Set $k= \lceil \frac{d_1 + \ldots +d_n}{n-1} \rceil $.
Suppose that the sheaf of relations $\shR(m)$ for the elements
$f_1, \ldots ,f_n$ is semistable. Then
$$(f_1, \ldots ,f_n)^\soclo = (f_1, \ldots ,f_n) + R_{\geq k} \, .$$
\end{theorem}

\proof
Let $f_0 \in R$ be homogeneous of degree $m$.
Suppose first that $m \geq k$.
Then
$$m \geq \frac{d_1 + \ldots +d_n}{n-1} =\frac{\mu (\shF(0))}{\delta}
=\frac{\mu_{\max}(f_1, \ldots ,f_n)}{\delta} \, .$$
Hence the numerical condition in \ref{maxin} is fulfilled,
thus $f_0 \in (f_1, \ldots ,f_n)^\soclo$.

Suppose now that $m<k$.
Then $m< \frac{d_1+ \ldots +d_n}{n-1}
=\frac{\mu_{\min}(f_1, \ldots, f_n)}{\delta}$
and \ref{minex} gives the result.
\qed

\medskip
Suppose that in the relative setting (\ref{relsituation})
the sheaf of relations is semistable in the generic point, so that
the vanishing theorem \ref{semistablevanishing} holds
in the generic point.
What can we say about the behavior in positive characteristic?
We know by \cite[\S 5]{miyaokachern} that
$\shF_\fop$ is semistable over an open non-empty subset
of $\Spec\,D$. However, for strongly semistable
we have to take into account the following problem of Miyaoka.

\begin{remark}
Miyaoka states in \cite[Problem 5.4]{miyaokachern} the following problem:
suppose that $C$ is a relative (smooth projective) curve over a
(say) $\ZZ$-algebra $D$ of finite type
and assume that a locally free sheaf $\shF$
is semistable in the generic fiber (characteristic zero).
Let $S$ be the set of points $P \in \Spec \, D$
of positive characteristic such that
$\shF | {C_P}$ is strongly semistable. Is $S$ dense in $\Spec\, D$?
\end{remark}

Therefore we may not expect that semistability in the generic point implies
a vanishing theorem for $p \gg 0$ without any further conditions.
It implies however that the bounds are quite near to the expected number
$(d_1 + \ldots + d_n)/(n-1)$.

\begin{corollary}
\label{comparebound}
Suppose the situation of {\rm \ref{relsituation}} and suppose that
the sheaf of relations is semistable over the generic point.
Let $m =\deg(f_0)$. Then the following hold for $p \gg 0$.

\renewcommand{\labelenumi}{(\roman{enumi})}
\begin{enumerate}

\item
If $m \geq  \frac{d_1 +\ldots + d_n}{n-1} + \frac{1}{\delta} $, then
$f_0 \in (f_1, \ldots ,f_n)^\soclo$.

\item
If $m \leq \frac{d_1 +\ldots + d_n}{n-1} - \frac{1}{\delta}$,
then
$f_0 \in (f_1, \ldots ,f_n)^\soclo$ if and only if $f_0 \in (f_1, \ldots ,f_n)$.
\end{enumerate}
\end{corollary}
\proof
For (i) we have the estimates
\begin{eqnarray*}
\delta m & \geq &  \delta \frac{d_1+\ldots+d_n}{n-1} + 1 \cr
&\geq & \lfloor \delta \frac{d_1+\ldots+d_n}{n-1} \rfloor + 1 \cr
&= & \lfloor \mu_{\rm \max} (\shF_\eta(0)) \rfloor +1 \cr
&> & \bar{\mu}_{\rm max}(\shF(0)) \, ,
\end{eqnarray*}
where the last estimate follows from \ref{compare}.
The statement follows then from \ref{maxin}.

(ii) follows by similar estimates from \ref{minex}.
\qed

\medskip
The previous corollary shows that the inclusion and exclusion bounds
are very near to $(d_1+ \ldots + d_n)/(n-1)$. If this number is not an
integer and if the degree of the curve is big enough, then we get also a
vanishing theorem from \ref{comparebound}.
In general however we get a vanishing theorem in positive characteristic only
for those points $\fop \in \Spec\, D$
for which the sheaf of relations $\shR_\fop$ is strongly semistable.

\begin{theorem}
\label{semistablevanishingp}
Suppose the situation {\rm \ref{relsituation}}.
Set $k= \lceil \frac{d_1 + \ldots +d_n}{n-1} \rceil $.
Suppose that for the generic fiber
the sheaf of relations is semistable. Then for all $p \gg 0$
such that the corresponding sheaf of
relations is strongly semistable we have
$$(f_1, \ldots ,f_n)^\soclo = (f_1, \ldots ,f_n) + R_{\geq k} \, .$$
\end{theorem}

\proof
The inclusion $\supseteq$ follows for every $\fop \in \Spec\, D$ such that
the sheaf $\shR_\fop (m)$ is strongly semistable on $Y_{\bar{\kappa}(\fop)}$
from \ref{maxin}.

For the inclusion $\subseteq $ we do not need the
condition strongly semistable.
From \ref{minex} we know that a single fixed element
$f_0$ with $f_0  \in (f_1, \ldots ,f_n)^\soclo$
belongs also to the right hand side for $p \gg 0$,
but here we state the identity of the two ideals
for $p \gg 0$.
Let $I=(f_1, \ldots , f_n)$ in $S$. We may assume by shrinking
$D$ that $I=S \cap IS_{Q(D)}$.
Let $m <k$ and consider
$S_m /I_m \subseteq H^1(D(S_+),Rel(f_1, \ldots,f_n)_m)$.
We may assume that $S_m/I_m$ is a free $D$-module with a basis
induced by $h_j \in S_m$, $ 1 \leq j \leq t$.

For every field $Q(D) \subseteq L$
the sheaf $\shF_L(-m) = \shR_L(m)^\dual$ is ample on
$\Proj S_{L}$ due to \ref{amplekrit}, hence for every
$h = \sum \lambda_j h_j \neq 0 $
the extension $\shF'_L(-m)(h)$ is also ample
and the open subset $\PP(\shF'_{L,h})- \PP(\shF_L)$
is affine. This is then also true for the open subset
$U_L =D( (S_L)_+) \subseteq
\Spec \, S_L[T_1, \ldots ,T_n]/(f_1T_1 + \ldots +f_nT_n+h)$.

We introduce indeterminates $\Lambda_j, \, 1 \leq j \leq t,$
for the coefficients of an element $h= \sum_j \lambda_jh_j$
and consider the universal forcing algebra
$$C= S[T_i, \Lambda_j] /
(f_1T_1 + \ldots +f_nT_n + \Lambda_1h_1 + \ldots +\Lambda_th_t)$$
over $S[\Lambda_j]$ and over $D[\Lambda_j]$. Set $U=D(S_+) \subseteq \Spec\, C$.

We claim that (after shrinking $D$)
for every point $P \in \Spec \, D[\Lambda_1, \ldots , \Lambda_t]$,
$P \not\in V(\Lambda_1, \ldots , \Lambda_t)$,
the fiber $U_P$ is affine.
We show this by increasing inductively
the open subset where this statement holds.

For the quotient field $L=Q(D)(\Lambda_1, \ldots ,\Lambda_t)
=Q(D[\Lambda_1, \ldots , \Lambda_t])$
we know that $U_L$ is affine.
Therefore also $U \cap D(g)$ is affine (by remark \ref{affineopen}),
where $0 \neq g \in D[\Lambda_1, \ldots , \Lambda_t]$,
and we know then that the fiber $U_P$ is affine for every point
$P \in D(g) \in \Spec \, D[\Lambda_1, \ldots , \Lambda_t]$.
So we know that the claim is true
for a non-empty open subset.

For the induction step suppose that the claim is true for
the open subset $W \subseteq \Spec \, D[\Lambda_1 , \ldots, \Lambda_t]$.
By shrinking $D$ we may assume that all irreducible components
of the complement of $W$ dominate $\Spec \, D$.
Consider such an irreducible component $Z=V(\foq)$ of
maximal dimension
and suppose that $Z \neq V(\Lambda_1, \ldots , \Lambda_t)$.
The generic point $\zeta$ of $Z$ has characteristic zero
and the $\Lambda_j$ are not all zero in $\kappa(\zeta)$,
hence again $U_\zeta$ is affine and we find an open neighborhood
$\zeta \in D(g) \cap Z \subseteq Z$,
$g \in D[\Lambda_1, \ldots , \Lambda_t] $,
such that for every point $P \in D(g) \cap Z$
the fiber $U_P$ is affine.
Hence the claim is now true on a bigger open subset and the number of
components of the complement of maximal dimension has dropped.

So we see that the claim holds
eventually for $D(\Lambda_1, \ldots , \Lambda_t)$.
Now the claim
means in particular that for every prime ideal $\fop \in \Spec\, D$
and every linear combination
$h = \sum_j \lambda_j h_j \neq 0$, $ \lambda_j \in \kappa(\fop)$,
the corresponding
open subset $U_{\fop, h}$ is affine,
hence $h \not\in I^\soclo_{\kappa(\fop)}$ for all $h \not\in I_{\kappa(\fop)}$.

This procedure can be done for every degree $0 \leq m < k$,
hence we find a sufficiently small $\Spec\, D$ such that
the statement holds for all $\fop \in \Spec\, D$.
\qed

\begin{remark}
The Theorems \ref{semistablevanishing}, \ref{comparebound}
and \ref{semistablevanishingp} indicate that the
number $\lceil \frac{d_1 + \ldots +d_n}{n-1} \rceil $
is the ``generic bound'' for the degree of an element
to belong to the tight closure.
It is reasonable to guess that for $R$ of higher dimension the number
$\lceil \frac{ \dim \, R -1}{n-1} (d_1 + \ldots +d_n) \rceil $
should take over this part ($n \geq \dim \, R$).
\end{remark}

The following result of Huneke and Smith
\-- the strong vanishing theorem for para\-meters in dimension two \--
is an easy corollary of \ref{semistablevanishing} and
\ref{semistablevanishingp}, see \cite[Theorem 4.3]{hunekesmithkodaira}.
For characteristic zero we have proved it already in \ref{degreecrit1}.

\begin{corollary}
\label{parametervanishing}
Let the notation be as in {\rm\ref{forcingsit}} and suppose that $n=2$,
so we are concerned with the tight closure of a parameter ideal.
Suppose that the characteristic of $K$ is $0$
or $p >  \frac{2}{\delta} (g-1)$.
Then
$$(f_1,f_2)^\soclo = (f_1,f_2) + R_{\geq d_1+d_2} \, .$$
\end{corollary}
\proof
The statement for characteristic $0$
and characteristic $p \gg 0$
follows at once from \ref{semistablevanishing} and \ref{semistablevanishingp},
since the sheaf of relations
$\shR(m)= \O_Y(m-d_1-d_2)$ is invertible, hence (strongly) semistable.

For the precise statement in positive characteristic we need
the ampleness criterion of Hartshorne-Mumford for bundles of rank two,
see \ref{amplecritranktwo}.
The inclusion $ \supseteq $ is true in any characteristic
by \ref{maxin} (or already \ref{degreecrit1}),
since $\shR(m)= \O_Y(m-d_1-d_2)$.
So suppose that $f_0 \in (f_1,f_2)^\soclo$, but
$f_0 \not\in (f_1,f_2)$ and $m= \deg \, (f_0) < d_1+d_2$.
The corresponding forcing class $c \in H^1(Y, \shR(m))$
gives the forcing sequence
$$ 0   \lra \O_Y \lra \shF'(-m) \lra \shF(-m)  \lra 0 \, .$$
Now $\deg \, (\shF'(-m)) = (d_1+d_2 -m) \delta \geq \delta$,
hence $\deg \, (\shF'(-m)) > \frac{2}{p} (g-1) $.
This gives the first condition in the criterion \ref{amplecritranktwo}.
For the second condition, let
$\shF'(-m) \ra \shM \ra 0$ be an invertible quotient sheaf.
If $\deg \, (\shM) < 0$ or if $ \deg \, (\shM) = 0$ and $\shM \neq \O_Y$,
then the composed mapping $\O_Y \ra \shM$ is zero and $\shM$ is a quotient
of the ample invertible sheaf $\shF(-m)$, which is not possible.
If $\shM= \O_Y$, then the composed mapping
$\O_Y \ra \O_Y$ is not zero,
hence it is an isomorphism and the sequence splits,
which contradicts the assumption $c \neq 0$.
\qed

\begin{remark}
If $Y$ is a smooth plane curve given by a polynomial $F$ of degree $\delta$,
then $g=(\delta -1) (\delta -2) /2$ and the condition in
\ref{parametervanishing} is that
$p > \frac{2}{\delta} (\frac{( \delta -1) (\delta -2)}{2} -1) = \delta -3 $.
The condition in \cite[Theorem 4.3]{hunekesmithkodaira}
for \ref{parametervanishing} to hold in positive characteristic is that
$p$ exceeds the degree of each of a set of
homogeneous generators for the module of $k$-derivations of $R$.
\end{remark}

\begin{remark}
It is not so easy to establish the semistability property
for the relation bundle for given ideal generators
on a given curve.
There exist however restriction theorems saying that this property holds
for a generic curve of sufficiently high degree under the condition
that the bundle is (defined and) semistable on the projective plane,
see \cite{flennerrestrict}, \cite{mehtaramanathansemistable},
\cite{huybrechtslehn} \cite{noma}, \cite{hein}.
\end{remark}

\smallskip 

\subsection{Computing the tight closure of three elements}

\label{n3}

\markright{Computing the tight closure of three elements}
\

\bigskip
Suppose that
we have three homogeneous primary
elements $f_1,f_2,f_3 \in R$ of degree $d_1,d_2,d_3$,
where $R$ is a two-dimensional normal standard-graded
domain over a field.
Then the sheaf of relations $\shR(m)$
on the smooth projective curve $Y= \Proj \, R$ for these three elements
has rank two, and due to section
\ref{algorithmsection} we may decide whether the torsor defined
by the cohomology class $\delta(f_0)=c \in H^1(Y, \shR(m))$
is affine or not.
Therefore we may compute $(f_1,f_2,f_3)^\soclo$,
at least if we can compute the Harder-Narasimhan filtration of $\shR$.

\medskip
Recall that we have defined (\ref{mindegdefinition})
the maximal degree of an invertible subsheaf
by $ \maxdeg_1 (\shS)
= \max \{ \deg(\shL):\, \shL
\mbox{ is an invertible subsheaf of } \shS \} \, .$
For a locally free sheaf $\shS$ of rank two
we have
$\mu_{\rm max} (\shS) = \max \{ \maxdeg_1 (\shS) ,\deg (\shS)/2 \}$
and
$\mu_{\rm min} (\shS)
= \min \{ \deg (\shS)/2 - \maxdeg_1 (\shS) , \deg (\shS)/2 \} $,
and $\shS$ is semistable if and only if
$\maxdeg_1 (\shS) \leq \deg (\shS)/2$.
If we find a subsheaf $\shL \subseteq \shS$
such that $\deg (\shL) \geq \deg (\shS)/2$ and such that
the quotient is itself locally free (i.e. $\shL$ is a subbundle),
then $\deg (\shL) = \maxdeg_1(\shS) = \mu_{\rm max} (\shS)$.

\begin{lemma}
\label{exactsequenceranktwo}
Let $\shS$ denote a locally free sheaf of rank two on
a smooth projective curve $Y$.
Suppose that we have a short exact sequence
$0 \ra \shL \ra \shS \ra \shM \ra 0$, where $\shL$ and $\shM$ are invertible
sheaves.
Then $\maxdeg_1 (\shS) \leq \max (\deg (\shL), \deg (\shM) )$.
If furthermore $\deg (\shL)  \geq \mu (\shS) = \deg (\shS)/2$,
then $\maxdeg_1 (\shS) = \deg (\shL)$
and $\mindeg_1(\shS) = \deg (\shM)$.

If $\deg (\shL) = \mu(\shS)$, then $\shS$ is strongly semistable.
If $\deg (\shL) > \mu(\shS)$, then $\shS$ is not stable
and $\shL$ is the maximal destabilizing subsheaf.
\end{lemma}

\proof
Let $\shN$ denote an invertible sheaf and let
$\varphi: \shN \ra \shS$ be a sheaf morphism.
If $\deg (\shN) > \deg (\shM)$, then the composed morphism
$ \shN \ra \shM$ is zero and $\varphi$ factors through
$\shL$. But then $\deg (\shN) \leq \deg (\shL)$ or $\varphi$
is zero.
So suppose that $\deg (\shL)  \geq \mu (\shS)$.
Then $\deg (\shM) = \deg (\shS) -\deg (\shL)
\leq \deg (\shS) - \deg (\shS)/2 = \deg (\shS)/2$, hence
$ \deg (\shL) \geq \mu ( \shS) \geq \deg (\shM)$.
Thus $\shL$ is a subbundle of maximal degree and $\shM$ is a quotient
invertible sheaf of minimal degree.
The other statements follow.
\qed

\medskip
In particular we have the two alternatives:
the locally free sheaf $\shS$ of rank two on $Y$ is strongly stable:
then $\bar{\mu}_{\rm max}(\shS) = \mu(\shS)= \bar{\mu}_{\rm min}(\shS)$
and the tight closure is easy to compute by the numerical criterion
\ref{semistablevanishing}.

Or $\shS$ is not strongly semistable,
then there exists a finite morphism
$\varphi: Y' \ra Y$ such that there exists a short exact sequence on $Y'$,
$ 0 \ra \shL \ra \shS' \ra \shM \ra 0$,
where $\deg (\shL) \geq \mu (\shS')$. In this case
the pull-back of this sequence for another morphism
$\psi : Y'' \ra Y'$
fulfills also the condition in \ref{exactsequenceranktwo},
and $\bar{\mu}_{\rm min} (\shS) = \deg (\shM) /\deg (\varphi)$.

If we have a short exact sequence for the sheaf of relations
for three elements $f_1,f_2,f_3$, then we can often compute
$(f_1,f_2,f_3)^\soclo$ and $(f_1,f_2,f_3)^+$ according to
the following proposition.

\begin{proposition}
\label{exactsequencecrit}
Let $R$ denote a two-dimensional normal standard-graded $K$-domain
over an algebraically closed field $K$.
Let $f_1, f_2,f_3 \in R$ denote three homogeneous
$R_+$-primary elements.
Let $0 \ra \shL \ra \shR(m) \ra \shM \ra 0$ denote 
a short exact sequence, where $\shL$ and $\shM$ are invertible.
Let $f_0$ denote another homogeneous element of degree $m$
and let $c \in H^1(Y, \shR(m))$ denote its forcing class.
Let $c'$ denote the image of $c$ in $H^1(Y, \shM)$.
Then the following hold
{\rm (}suppose in the first two statements that the characteristic of $K$
is zero or $p \gg 0${\rm )}.

\renewcommand{\labelenumi}{(\roman{enumi})}
\begin{enumerate}

\item
If $\deg (\shL) <0 $ and $\deg (\shM) <0 $ and $c \neq 0$,
then $f_0 \not\in (f_1,f_2,f_3)^\soclo$.

\item
If $\deg (\shM) <0 $ and $c' \neq 0$, then $f_0 \not\in (f_1,f_2,f_3)^\soclo$.

\item
If $c' =0$ and $\deg (\shL) \geq 0$, then
$f_0 \in (f_1,f_2,f_3)^\soclo$.

\item
If $\deg (\shL) \geq \deg (\shM) \geq 0 $,
then $f_0 \in (f_1,f_2,f_3)^\soclo$.

\item
Suppose that the characteristic of $K$ is positive.
Suppose that $c'=0 $ or that $\shM$ is trivial or $\shM$ has positive degree
and that $\shL$ is trivial or has positive degree.
Then $f_0 \in (f_1,f_2,f_3)^+$.

\end{enumerate}
\end{proposition}

\proof
(i). The dual sheaf $\shF(-m)= \shR(m)^\dual$ is ample as an extension
of ample invertible sheaves, hence the result follows from
\ref{ampleaffine}.
(ii) follows from \ref{quotientnegdegree}.

(iii) If $c'=0$, then there exists a cohomology class
$e \in H^1(Y, \shL)$ mapping to $c$.
If $f_0 \not \in (f_1,f_2,f_3)^\soclo$, then
$\PP(\shF') - \PP(\shF))$ would be affine and then
$\PP((\shL') ^\dual) - \PP(\shL^\dual))$ would also be affine due to
\ref{affinelemma}, but this is not true due to
$\deg (\shL) \geq 0$ and \ref{slopemaxkrit}.

(iv). The condition implies in connection with \ref{exactsequenceranktwo} that
$\bar{\mu}_{\rm max}( \shR(m)) \geq 0$, hence the result follows again from
\ref{slopemaxkrit}.

(v). After applying a finite mapping $Y' \ra Y$ we may assume
that $c' =0$. For $\deg (\shM) >0$ this can be done by a Frobenius
power and for $\shM = \O_Y$ this is due to \ref{frobeniusartinschreier}.
Therefore we may assume that $c$ stems from a cohomology class
$e \in H^1(Y', \shL')$.
Due to the assumptions on $\shL$ we can do the same with $e$,
hence there exists altogether a finite mapping $Y'' \ra Y$
such that the pull-back of $c$ is zero.
Therefore $f_0 \in (f_1,f_2,f_3)^+$.
\qed

\medskip
We may apply \ref{exactsequencecrit} to the short
exact sequence given by the Harder-Narasim\-han filtration
to compute the tight closure of $(f_1, f_2, f_3)$
if $\shR$ is not strongly sta\-ble.

\begin{corollary}
\label{exactsequencedecide}
Let $R$ denote a two-dimensional normal standard-graded $K$-domain
over an algebraically closed field.
Let $f_1, f_2,f_3 \in R$ denote three homogeneous
$R_+$-primary elements.
Suppose that the sheaf of relations $\shR$ is not strongly
stable, and let $\varphi: Y' \ra Y$ denote a finite dominant
morphism of smooth projective curves
such that there exists a short exact sequence
$$0 \lra \shL(m) \lra \varphi^*(\shR(m)) \lra \shM(m) \lra 0$$ on $Y'$,
where $\deg (\shL(m)) \geq \mu (\varphi^*(\shR (m)) ) \geq \deg (\shM(m))$.

Let $f_0$ denote another homogeneous element of degree $m$,
let $c$ denote its forcing class
in $H^1(Y',\varphi^*(\shR(m)) )$ and $c'$ its image in $H^1(Y, \shM(m))$.
Then we may decide whether $f_0 \in (f_1,f_2,f_3)^\soclo$
in the following way {\rm(}assume in the first and second statement
that the characteristic is zero or $p \gg 0${\rm )}.

\renewcommand{\labelenumi}{(\roman{enumi})}
\begin{enumerate}

\item
If $\deg (\shL(m)) <0 $ and $c \neq 0$, then $f_0 \not\in (f_1,f_2,f_3)^\soclo$.

\item
If $\deg (\shL(m))  \geq 0 $, $\deg (\shM(m)) <0$
and $c' \neq 0$, then $f_0 \not\in (f_1,f_2,f_3)^\soclo$.

\item
If $\deg (\shL(m))  \geq 0 $, $\deg (\shM(m)) <0$
and $c' = 0$, then $f_0 \in (f_1,f_2,f_3)^\soclo$.

\item
If $\deg (\shM(m)) \geq 0 $, then $f_0 \in (f_1,f_2,f_3)^\soclo$.

\end{enumerate}
\end{corollary}

\proof
This follows from \ref{exactsequencecrit}.
\qed

\begin{remark}
\label{exactsequencepositive}
Suppose the situation of \ref{exactsequencedecide}
and suppose that the characteristic is positive.
Then we need in (iii) and (iv) stronger conditions to
conclude that $f_0 \in (f_1,f_2,f_3)^{\gr}$.
For (iii) we need that $\shL(m)$ is trivial or
that $\deg (\shL(m)) $ is positive.
For (iv) we need that $\shM$ is of positive degree.
This follows from  \ref{exactsequencecrit}(v).
\end{remark}

\medskip
The situation of \ref{exactsequencedecide} holds in particular
if the sheaf of relations is decomposable, i.e. the direct sum of two
invertible sheaves. Corollary \ref{exactsequencedecide}
gives of course the same answer as corollary
\ref{splitting}.

\medskip
The next example gives a negative answer to a question
of Craig Huneke asked at the MSRI (September 2002).
The example gives an ideal which is generated by
$*$-independent elements (meaning that none of them is contained
in the tight closure of the others, see \cite{vraciu} for this notion),
but there does not hold a strong vanishing theorem for it, i.e.
there does not exist a common inclusion and exclusion
bound for tight closure.

\begin{example}
\label{indepexample}
Consider the ideal $I=(x^4,xy,y^2)$ in the Fermat cubic $x^3+y^3+z^3=0$.
These ideal generators come from
the regular polynomial ring $K[x,y] \subset K[x,y,z]/(x^3+y^3+z^3)=R$.
From this it follows that they are $*$-independent in $R$
and that the relation bundle must split, and
in fact
$$\shR(x^4,xy,y^2)(5) = \O_Y(0) \oplus \O_Y(2) \,$$
where
$Y= \Proj \, R$ is the corresponding elliptic curve.
Let $h$ denote a homogeneous element of degree $m$, given rise
to a forcing class
$$c \in H^1(Y, \shR(m))= H^1(Y,\O_Y(m-5)) \oplus H^1(Y,\O_Y(m-3)) \, . $$
From the numerical criterion in \ref{exactsequencedecide}
(or \ref{splitting})
it is easy to deduce the following.

For $m \geq 5$ we have $R_m \subset I^\soclo$.

For $m=4$ an element $h$ belongs to $I^\soclo$ only if it belongs to $I$.

For $m=3$ all possibilities occur.
We find
$ yz^2 \in (x^4,xy,y^2)^\soclo$, but not in the ideal itself,
and $ xz^2 \not\in I^\soclo$.
\end{example}

\smallskip

\subsection{Slope bounds for indecomposable bundles}

\label{sectionindecomposable}

\markright{Slope bounds for indecomposable bundles}

\

\bigskip
In this section we discuss inclusion and exclusion
bounds for tight closure for homogeneous primary ideals generated by
three elements under the condition that the relation bundle
is indecomposable.
We have the following lemma for indecomposable sheaves of rank two.

\begin{lemma}
\label{mumaxindecomposable}
Let $\shG$ denote an indecomposable
locally free sheaf of rank two on a smooth
projective curve $Y$ of genus $g$ over an algebraically closed field.
Then
$\mu_{\rm max}(\shG) \leq \mu (\shG) + g-1 $.
\end{lemma}
\proof
Let $0 \ra \shL  \ra \shG \ra \shM \ra 0$ denote a short exact sequence
with invertible sheaves $\shL$ and $\shM$.
The sequence corresponds to an element
$c \in H^1(Y, \shL \otimes \shM^{-1})$.
Since $\shG$ is indecomposable we have $c \neq 0$.
By Serre duality it follows that
$H^0(Y, \shL^{-1} \otimes \shM \otimes \omega_Y) \neq 0$
and hence that
$\deg (\shL^{-1}) + \deg (\shM) + 2g-2 \geq 0$.
Therefore
\begin{eqnarray*}
\deg (\shL) & \leq & \deg (\shM) +2g-2 \cr
&=& 2 \mu (\shG) - \deg (\shL) + 2g-2
\end{eqnarray*}
and $\deg (\shL) \leq \mu (\shG) +g-1$.
\qed

\begin{remark}
The statement \ref{mumaxindecomposable}
is equivalent to the fact that the $e$-invariant of a ruled
surface is $\leq 2g-2$, see \cite[Theorem V.2.12]{haralg}.
\end{remark}

\begin{theorem}
\label{genusbound}
Let $K$ denote an algebraically closed field of characteristic zero
and let $R$ denote a normal
standard-graded $K$-domain of dimension two. Let $Y= \Proj\, R$
denote the corresponding smooth projective curve of genus $g$
and let $\delta$ denote the degree of $\O_Y(1)$.
Let $f_1,f_2,f_3$ denote $R_+$-primary homogeneous
elements of degree $d_1,d_2,d_3$
and suppose that the sheaf of relations $\shR(m)$
is indecomposable on $Y$. Then
$$\mu_{\rm max} (f_1,f_2,f_3) \leq \delta \frac{d_1+d_2+d_3}{2} + g-1 \, $$
and
$$ \mu_{\rm min} (f_1,f_2,f_3) \geq \delta \frac{d_1+d_2+d_3}{2} - g+1 \, .$$
\end{theorem}

\proof
The first statement follows from
\ref{mumaxindecomposable},
for
$$\mu_{\rm max} (f_1,f_2,f_3) = \mu_{\rm max} (\shF(0))
\leq \mu( \shF(0)) + g-1 =
\delta \frac{d_1+d_2+d_3}{2} + g-1 \, .$$

The bound for $\mu_{\rm min}(f_1,f_2,f_3)$
follows from $\mu_{\rm min}( \shF ) = 2 \mu(\shF) - \mu_{\rm max} (\shF)$.
\qed

\begin{corollary}
\label{genusboundtight}
Let $K$ denote an algebraically closed
field of characteristic zero and let $R$ denote a normal
standard-graded
$K$-domain of dimension two. Let $Y= \Proj\, R$
denote the corresponding smooth projective curve of genus $g$
and let $\delta$ denote the degree of $\O_Y(1)$.
Let $f_1,f_2,f_3$ denote $R_+$-primary homogeneous
elements of degree $d_1,d_2,d_3$
and suppose that the sheaf of relations $\shR(m)$ is indecomposable on $Y$.
Then the following statements hold.

\renewcommand{\labelenumi}{(\roman{enumi})}
\begin{enumerate}

\item
For $m \geq \frac{d_1+d_2+d_3}{2} + \frac{g-1}{\delta} $
we have the inclusion $R_m \subseteq (f_1,f_2,f_3)^\soclo$.

\item
For $m < \frac{d_1+d_2+d_3}{2} - \frac{g-1}{\delta}$
we have $(f_1,f_2,f_3)^\soclo \cap R_m = (f_1,f_2,f_3) \cap R_m$.
\end{enumerate}
\end{corollary}

\proof
The first statement follows from \ref{genusbound}
and from \ref{maxin}.
The second statement follows from \ref{genusbound} and \ref{minex}.
\qed

\begin{corollary}
\label{degreebound}
Let $K$ denote an algebraically closed
field of characteristic zero and let $F \in K[x,y,z]$ denote a
homogeneous polynomial of degree $\delta$
such that $R =K[x,y,z]/(F)$ is a normal domain.
Let $f_1,f_2,f_3 \in R$ denote $R_+$-primary homogeneous
elements of degree $d_1,d_2,d_3$.
Suppose that the sheaf of relations $\shR(m)$ is indecomposable
on the curve $Y= \Proj \, R$. Then the following hold.

\renewcommand{\labelenumi}{(\roman{enumi})}
\begin{enumerate}

\item
$R_m \subseteq (f_1,f_2,f_3)^\soclo$
for $m \geq \frac{d_1+d_2+d_3}{2} + \frac{\delta -3}{2} $.

\item
For $m < \frac{d_1+d_2+d_3}{2} - \frac{\delta -3}{2}$
we have $(f_1,f_2,f_3)^\soclo \cap R_m = (f_1,f_2,f_3) \cap R_m$.
\end{enumerate}

\end{corollary}

\proof
This follows from \ref{genusboundtight} taking into account that
$g=(d-1)(d-2)/2$.
\qed

\medskip
We also have the following result in positive characteristic.

\begin{corollary}
\label{genusboundpositive}
Suppose the relative setting {\rm \ref{relsituation}} and suppose that $n=3$.
Suppose that  $\shR(m)$ is indecomposable
on $Y_{\bar{\eta}}$, where $\eta $ is the generic point.
Then for $p \gg 0$ we have
$$ \bar{\mu}_{\rm max}(f_1,f_2,f_3)
< \lfloor \delta \frac{d_1+d_2+d_3}{2} \rfloor + g $$
and 
$$ \bar{\mu}_{\rm min}(f_1,f_2,f_3)
> \lfloor \delta \frac{d_1+d_2+d_3}{2} \rfloor - g $$
\end{corollary}
\proof
We have from \ref{compare} and \ref{genusbound} the estimates
\begin{eqnarray*}
\bar{\mu}_{\rm max} (f_1,f_2,f_3) &< &
\lfloor \mu_{\rm max}(\shF_\eta (0)) \rfloor  +1 \cr
& \leq & \lfloor  \delta \frac{d_1+d_2+d_3}{2} + g-1 \rfloor +1 \cr
&=& \lfloor  \delta \frac{d_1+d_2+d_3}{2} \rfloor +g 
\end{eqnarray*}
\qed

\begin{remark}
Corollary \ref{genusboundpositive} allows us to formulate
\ref{genusboundtight} and \ref{degreebound} also in positive
characteristic with some modifications.

The pull-back of
an indecomposable locally free sheaf $\shS$ on a projective curve $Y$
under a morphism $Y' \ra Y$ is not indecomposable in general.
Lets call a locally free sheaf $\shS$
{\em strongly indecomposable} if the pull-back $\varphi^*(\shS)$
is again indecomposable for every finite dominant morphism
$\varphi :Y' \ra Y$.
If the rank of $\shS$ is two, then a strongly indecomposable
sheaf is even strongly semistable.
For suppose that
$\shL \subset \shS$ is an invertible subsheaf
such that $\deg (\shL) > \mu (\shS)$ on $Y$.
Let $\varphi: Y' \ra Y$ denote a $K$-linear Frobenius morphism of degree $q$.
Then from \ref{mumaxindecomposable}
we know that
$ \deg ( \varphi^* (\shL)) \leq \mu (\varphi^*( \shS) ) +g(Y')-1
= \deg (\varphi) \mu (\shS) +g(Y')-1$.
Since $g(Y')=g(Y)$, this gives a contradiction if we choose
the degree of $\varphi$ high enough.
\end{remark}

For the plus closure in positive characteristic we get a
somewhat worse inclusion bound.

\begin{proposition}
\label{indecomposableplus}
Let $R$ denote a normal standard-graded $K$-do\-main
over an algebraically closed field $K$ of positive characteristic,
let $Y= \Proj\, R$ be the projective curve of genus $g$
and degree $\delta$ with canonical sheaf $\omega_Y$.
Let $f_1,f_2,f_3$ denote homogeneous primary elements of degree
$d_1,d_2,d_3$ and suppose that their sheaf of relations $\shR(m)$
is indecomposable.

Then $R_m \subseteq (f_1,f_2,f_3)^{\gr}$
for $m \geq \frac{d_1+d_2+d_3}{2} + \frac{2g-2}{\delta}$.
\end{proposition}

\proof
If $H^1(Y, \shR(m))=0$, then there is nothing to prove.
So suppose that $H^1(Y, \shR(m)) \neq 0$.
This means by Serre duality that
$$H^0(\shF(-m) \otimes \omega_Y)
 \neq 0 \, $$
or that there exists a non-trivial morphism
$\shR(m) \ra \omega_Y$.
Since $\mu (\shR(m)) =(m - \frac{d_1+d_2+d_3}{2} )\delta$
and $\deg (\omega_Y)= 2g-2$ it follows from
the numerical condition that $\shR(m)$ is not stable.

On the other hand we have
due to \ref{mumaxindecomposable}
that
$$\mu_{\rm min} (\shR(m))
\geq  \mu( \shR(m)) - (g-1)
= (m- \frac{d_1+d_2+d_3}{2})\delta - (g-1)
\, .$$
This is $>0$ for $g \geq 2$ and the result follows from
\ref{exactsequencepositive}.
For $g=0$ every sheaf is decomposable, so there is nothing to prove.
So suppose that $g=1$ and $2m=d_1+d_2+d_3$.
An indecomposable sheaf on an elliptic curve is
semistable.
The non-trivial mapping $\shR(m) \ra \O_Y= \omega_Y$
is then surjective and we have a short exact sequence
$0 \ra \O_Y \ra \shR(m) \ra \O_Y \ra 0$,
and the result follows from
\ref{exactsequencecrit}(v).
\qed

\begin{example}
\label{150example}
Let $F \in K[x,y,z]$ denote a homogeneous form of degree $\delta =5$
such that $R=K[x,y,z]/(F)$ is a geometrically normal domain.
Suppose that $f_1,f_2,f_3 \in R$ are primary
homogeneous elements of degree $100$ and such that its sheaf of relations
is (generically) indecomposable. Then the bounds of Smith in
\cite{smithgraded} give $R_{\geq 200} \subseteq (f_1,f_2,f_3)^\soclo$ on the one
hand and, on the other hand, that an element of degree $ m \leq 100$ belongs
to $(f_1,f_2,f_3)^\soclo$ if and only if it belongs to the ideal itself.

The bounds in \ref{degreebound} give on the one hand
that $R_{\geq 151} \subseteq (f_1,f_2,f_3)^\soclo$ and on the other
that an element of degree $\leq  148$ belongs
to $(f_1,f_2,f_3)^\soclo$ if and only if it belongs to the ideal itself
(for zero characteristic or $p \gg 0$).
So only the degrees $149$ and $150$ are not covered by \ref{degreebound}.
If furthermore the sheaf of relations is semistable in the generic point,
then \ref{comparebound} shows that also an element of degree $149$
belongs to the tight closure only if it belongs to the ideal itself.
\end{example}

\smallskip

\subsection{The degree of relations}

\label{relations}

\markright{The degree of relations}
\

\bigskip
The notions of semistability and of minimal and maximal degree
of a locally free sheaf $\shS$ on a smooth projective curve $Y$
refer to all locally free subsheaves of $\shS$ (or quotient sheaves).
For a relation sheaf $\shR(m)$ however
defined by homogeneous primary elements $f_1,\ldots ,f_n$
in a two-dimensional normal standard-graded $K$-domain we have
the fixed polarization $\O_Y(1)$ on $Y= \Proj\, R$.
It is then often easier to control the behavior of
$\shR(m) = \shR(0) \otimes \O_Y(m)$
instead of $\shR(0) \otimes \shL$ for all invertible sheaves $\shL$.
The (non-)\-existence of relations $\neq 0$
for $f_1, \ldots ,f_n$ of certain degree has a lot of consequences
on the structure of $\shR(m)$ and hence on the corresponding
tight closure problem.

\begin{lemma}
\label{maxdeglemma}
Suppose that the locally free sheaf $\shS$ on the smooth
projective curve $Y$ of genus $g$ over an algebraically closed field $K$
does not have sections $\neq 0$.
Then $\Gamma(Y,\shS \otimes \shL) =0 $
for every invertible sheaf $\shL$ of degree $\leq -g$.
In particular $\maxdeg_1 (\shS) \leq g-1$.
\end{lemma}
\proof
Suppose the contrary.
Then we have a non-trivial morphism
$\shM \ra \shS$ such that $\deg \shM \geq g$.
But due to the theorem of Riemann-Roch we have
$ h^0(\shM) \geq \deg (\shM) +1-g$, hence the invertible
sheaf $\shM$ must have non-trivial sections,
which gives a contradiction.
\qed

\begin{proposition}
\label{maxinvertbound}
Let $f_1, \ldots, f_n$ denote homogeneous primary elements
in a normal two-dimensional standard-graded $K$-domain $R$ of degree $d_i$,
where $K$ is an algebraically closed field.
Suppose that $Y= \Proj\, R$ has genus $g=g(Y)$ and degree $\delta$.

\renewcommand{\labelenumi}{(\roman{enumi})}

\begin{enumerate}

\item
Suppose that there exists a relation $\neq 0$
for the elements $f_1, \ldots, f_n$ of total
degree $k < (d_1 + \ldots +d_n )/(n-1)$.
Then the sheaf of relations is not semistable
{\rm (}and not stable for ``$ \leq$ ''{\rm )}.

\item
Suppose that there does not exist a relation $\neq 0$
of total degree $k$. Then
$$ \maxdeg_1(\shR(k)) \leq g-1 \, . $$

\item
Let $n=3$. Suppose that there exists a relation $\neq 0$
of total degree $k < \frac{d_1 + d_2 +d_3}{2} - \frac{g-1}{\delta} $.
Then the sheaf of relations is decomposable, i.e. the sum of
two invertible sheaves.

\item
Let $n=3$ and suppose that
there does not exist a relation $\neq 0$
of total degree $k \geq \frac{d_1+d_2+d_3}{2} + \frac{g-1}{ \delta} $.
Then $\shR$ is semistable.
\end{enumerate}
\end{proposition}

\proof
(i). The relation $\neq 0$ induces a non-trivial morphism
$\O_Y \ra \shR(k)$, but the degree
$\deg (\shR(k))= ((n-1)k - d_1- \ldots -d_n) \delta < 0$
is negative, hence $\shR(k)$ is not semistable.

(ii).
The assumption means that $\shS= \shR(k)$ has no global sections $\neq 0$,
hence \ref{maxdeglemma} yields that
$\maxdeg_1 (\shR(k)) \leq g-1$.

(iii).
Since $\shR(k)$ has a section it follows that
$\mu_{\rm max} (\shR(k)) \leq 0$.
On the other hand we have
$\mu (\shR(k)) +g-1 =
(k- \frac{d_1+d_2+d_3}{2}) \delta +g-1 < 0$,
and the result follows from
\ref{mumaxindecomposable}.

(iv).
The numerical condition means that
$g-1 \leq   (k - \frac{1}{2}(d_1+d_2+d_3)) \delta = \deg (\shR(k))/2$,
hence from (ii) we get that
$\maxdeg_1 (\shR(k)) \leq \deg (\shR(k))/2 $
and the sheaf of relations is semistable.
\qed

\begin{remark}
We call a relation of degree $k < (d_1 + \ldots +d_n )/(n-1)$
as in \ref{maxinvertbound}(i) a {\em destabilizing relation}.
\end{remark}

We may derive from Proposition \ref{maxinvertbound}(ii)
the following inclusion bound for tight closure.

\begin{corollary}
\label{relationboundinclusion}
Let $R$ denote a normal two-dimensional standard-graded domain
over an algebraically closed field $K$ of characteristic zero
and let $f_1,f_2,f_3 \in R$ be $R_+$-primary homogeneous elements
of degree $d_1,d_2,d_3$.
Suppose that there does not exist a relation $\neq 0$ of total
degree $k \leq \frac{d_1+d_2+d_3}{2} + \frac{g-1}{\delta}$.
Then $R_m \subseteq (f_1,f_2,f_3)^\soclo$ holds for
$m \geq  d_1+d_2+d_3-k + \frac{g-1}{\delta }$.
\end{corollary}
\proof
If $\shR$ is semistable, then the result follows from
\ref{semistablevanishingp}. If $\shR$ is not semistable, then
with $ \mu_{\rm max} (\shF(-k)) = \deg (\shF(-k))-
\mu_{\rm min} (\shF(-k))= \deg (\shF(-k)) + \mu_{\rm max} (\shR(k))$
we get
\begin{eqnarray*}
\mu_{\rm max} (\shF(-m)) &=& \mu_{\rm max} (\shF(-k)) + (k-m) \delta \cr
&=& \deg (\shF(-k)) - \mu_{\rm min}(\shF(-k)) + (k-m)\delta \cr
&=& (d_1 +d_2+d_3 -2k) \delta + \maxdeg_1 (\shR(k)) +(k-m) \delta \cr
& \leq & (d_1 +d_2+d_3 -2k)\delta + g-1 +(k-m) \delta \cr
&=&  (d_1 +d_2+d_3 -k)\delta + g-1 -m \delta 
\end{eqnarray*}
This is $\leq 0$ if and only if
$m \geq  d_1+d_2+d_3-k + \frac{g-1}{\delta }$.
\qed

\begin{corollary}
\label{semistablecritcor}
Let $R=K[x,y,z]/(F)$ denote a normal two-dimensional
stan\-dard-graded $K$-domain over an algebraically closed field $K$,
where $F$ is a polynomial of degree $\delta$.
Let $f_1,f_2,f_3 \in R$ be $R_+$-primary homogeneous elements
of degree $d_1,d_2,d_3$.
Suppose that there does not
exist a relation $\neq 0$ for $f_1,f_2,f_3$ of total
degree $k$ with $k \leq  \frac{d_1+d_2+d_3}{2} + \frac{\delta -3}{2}$.
Then $\shR$ is semistable.
\end{corollary}
\proof
This follows directly from
\ref{maxinvertbound}(iv).
\qed

\begin{example}
Consider the sheaf of relations for the elements $x^d,y^d,z^d$ on
the Fermat curve given by
$x^\delta+y^\delta+z^\delta=0$.
The equation gives for $\delta \geq d$ at once the relation
$(x^{\delta-d},y^{\delta-d},z^{\delta-d})$ of total degree $\delta$.
\end{example}

\begin{example}
\label{powersexample}
We consider the elements $x^d,y^d,z^d$ on a
smooth projective curve given by an equation $F=0$, where $F$
is a homogeneous polynomial of degree $\delta$.
There exist relations like $(y^d,-x^d,0)$ of total degree $2d$.
Suppose that there does not exist relations of smaller degree.
Then the numerical condition in Corollary \ref{semistablecritcor}
for semistability
is that $2d-1 \geq 3d/2 + (\delta-3)/2$ or equivalently that
$\delta \leq d+1$.
If we want to apply \ref{semistablecritcor}
we have to make sure that the defining polynomial $F$ of degree $\delta$
does not yield relations of degree $< 2d$.

\smallskip
Look at $d=2$ and $\delta =3$.
If the monomial $xyz$ does occur in $F$, then there do not exist relations
of degree $3$ and the relation sheaf is semistable.
This yields however nothing interesting for the tight closure,
since then $xyz \in (x^2,y^2,z^2)$ holds anyway.

\smallskip
Now look at $d=4$ and $\delta =5$.
Under suitable conditions at the coefficients of $F$
there does not exist a relation of degree
$7$ for $x^4,y^4,z^4$.
Write $F=a x^3y^2+ bx^3yz+cx^3z^2+ dx^2y^3 $ etc.
Such a relation is the same as a multiple $FQ$ ($\deg (Q)=2$)
which belongs to $(x^4,y^4,z^4)$.
The six monomials of degree $2$ yield
six linear combinations in the six monomials of degree
$7$ outside $(x^4,y^4,z^4)$, namely
$ x^3y^3z$, $x^3y^2z^2$, $x^3yz^3$, $x^2y^3z^2$, $x^2y^2z^3$ and $xy^3z^3$.
We may choose the coefficients of $F$ in such a way that
these linear combinations are linearly independent.
Then there does not exist a relation of degree $7$.
So in this case the sheaf of relations is semistable.
It follows for $\Char (K)=0$ that $R_6 \subseteq (x^4,y^4,z^4)^\soclo$.
Note that it is not true that
$R_6 \subseteq (x^4,y^4,z^4)$, since there exist
$10$ monomials of degree $6$ outside $(x^4,y^4,z^4)$
in $K[x,y,z]$, hence the dimension of
$R_6/ (x^4,y^4,z^4)$ is at least $10-3$.
\end{example}

\begin{proposition}
\label{amplecritrelation}
Suppose that the characteristic of the algebraically closed field $K$ is zero
and that $R$ is a two-dimensional normal standard-graded
$K$-do\-main. Let $f_1, \ldots ,f_n$ denote homogeneous
$R_+$-primary elements.
Suppose that there does not exist a global relation $\neq 0$ of
total degree $k$.
Then $\shF(-m)$ is ample for
$m < k- \frac{2n-3}{(n-1) \delta} g + \frac{1}{ \delta}$.
In particular, $\shF(-m)$ is ample for $m \leq k - 2g/\delta $.
\end{proposition}

\proof
We know by \ref{maxinvertbound}(ii) that
$\maxdeg_1 (\shR(k)) \leq g-1$ and therefore dually that
$\mindeg_1 (\shF(-k)) \geq -g+1$.
Hence
$$\mindeg_1 (\shF(-m)) = \mindeg_1 (\shF(-k) \otimes \O_Y(k-m))
= \mindeg_1 (\shF(-k)) + (k-m)\delta
\geq -g+1  + (k-m) \delta \, . $$
The numerical condition is equivalent with
$ -g+1  + (k-m) \delta > \frac{n-2}{n-1} g$.
Hence $\mindeg_1(\shF(-m)) > \frac{n-2}{n-1} g$
and the result follows from \ref{amplerankcrit}.
\qed

\begin{corollary}
\label{globalrelationample}
Let $R=K[x,y,z]/(F)$ be a normal standard-graded $K$-domain
over an algebraically closed field $K$ of characteristic zero,
where $F$ is an irreducible polynomial of degree $\delta$.
Let $f_1, \ldots ,f_n$ denote primary homogeneous elements of degree $d_i$.
Suppose that there does not exist a global relation $\neq 0$
for $f_1, \ldots, f_n$ of total degree $k$.
Then $\shF(-m)$ is ample for $m \leq k- \delta + 2$.
An element $f_0 \in R$ of degree
$m \leq k- \delta + 2$ belongs to $(f_1, \ldots ,f_n)^\soclo$
if and only if it belongs to $(f_1, \ldots ,f_n)$.
\end{corollary}
\proof
$m \leq k- \delta +2
\leq k - (\delta -2)(\delta -1)/\delta
= k - 2g/\delta $, hence the result follows from
\ref{amplecritrelation} and \ref{minex}.
\qed

\begin{example}
\label{powersexample2}
We want to apply \ref{amplecritrelation} and \ref{globalrelationample}
to example \ref{powersexample} for $d=4$, $ \delta =5$
under the condition that there does not exist a global
relation of degree $7$.

Then \ref{globalrelationample} shows ($k=7$, $g=6$, $n=3$)
that $\shF(-m)$ is ample only for $m \leq 7-5+2=4$.
The second bound in \ref{amplecritrelation}
gives this for $m \leq 7- \frac{2 \cdot 6}{5} =4.6$.
The first bound however yields ampleness for
$m < 7 - \frac{3}{2 \cdot 5}6 + \frac{1}{5}
=7- \frac{16}{10}= 5.4 $.
\end{example}

\medskip
Suppose further that $\shR$ is the sheaf of relations
on a smooth projective curve
$Y= \Proj R$ for homogeneous primary elements
$f_1, \ldots ,f_n \in R$, where $R$ is a normal standard-graded
$K$-domain over an algebraically closed field $K$.
We say that a relation $r \in \Gamma(Y,\shR(m))$ is a
{\em primary relation} if it has no zero on $Y$, or equivalently,
if $r : \O_Y \ra \shR(m)$ defines a subbundle.

\begin{lemma}
\label{filtrieren}
Let $R$ be a normal standard-graded $K$-domain
over an algebrai\-cal\-ly closed field $K$.
Let $f_1, \ldots ,f_n \in R$ be homogeneous primary elements.
Let $r=(r_i)$ denote a homogeneous primary relation for $f_i$
of total degree $k$.
Then there exists a sequence
$0 \ra \O_Y \ra \shR(k) \ra \shL \ra 0$ such that
$\shL$ is locally free.
\end{lemma}
\proof
The relation is a global element in $\Gamma(Y,\shR (k))$ and yields a
morphism $\O_Y \subseteq \shR (k)$.
This morphism is locally given
by $\O_U \ra \O_U^{n-1}$, $ a \mapsto a(r_1, \ldots, r_n)$,
and since the relation is primary we may assume that
one $r_i$ is a unit.
Then the quotient is also locally free.
\qed

\begin{corollary}
\label{primarysemistable}
Let $K$ denote an algebraically closed field
and let $R$ denote a two-dimensional normal standard-graded $K$-domain.
Let $f_1, f_2,f_3 \in R$ denote three homogeneous
$R_+$-primary elements.
Suppose that there exists a primary relation of
total degree $k$.
Then this relation gives rise to a short exact sequence
$$0 \ra \O_Y \ra \shR(k) \ra \O_Y(2k - d_1-d_2-d_3) \ra 0 \, .$$
If $k \leq (d_1+d_2+d_3)/2$, then $\maxdeg_1(\shR(k))=0$
and $\mindeg_1(\shR(k)) =(2k-d_1-d_2-d_3) \delta \leq 0$.
If moreover $k = (d_1+d_2+d_3)/2$, then the sheaf of relations
is strongly semistable. 
\end{corollary}

\proof
The existence of the short exact sequence follows from
\ref{filtrieren}.
We have $ \deg (\O_Y) = 0 \geq k - (d_1+d_2+d_3)/2 = \mu (\shR(k))$,
hence we are in the situation of \ref{exactsequenceranktwo}.
If $k = (d_1+d_2+d_3)/2$, then $\shR(k)$ is the extension of
the structure sheaf by itself, hence its degree is $0$
and this follows again from \ref{exactsequenceranktwo}.
\qed

\begin{corollary}
\label{primaryexclusion}
Let $R$ denote a two-dimensional normal standard-graded $K$-domain
and let $f_1, f_2,f_3 \in R$ denote three homogeneous
$R_+$-primary elements.
Suppose that there exists a primary relation of
total degree $k \leq (d_1+d_2+d_3)/2$.
Suppose that the characteristic of the algebraically closed field $K$ is zero
or $p \gg 0$.
Let $f_0 \in R$ be a homogeneous element of
degree $ \deg (f_0)= m$. Then the following hold.

If $m <k$, then $f_0 \in (f_1,f_2,f_3)^\soclo $ if and only if $f_0 \in (f_1,f_2,f_3)$.

If $m \geq d_1+d_2+d_3-k$, then $ f_0 \in (f_1,f_2,f_3)^\soclo$.

\end{corollary}
\proof
From \ref{primarysemistable} we
get the short exact sequence
$0 \ra \O_Y(m-k) \ra \shR(m) \ra \O_Y(m+k-d_1-d_2-d_3 ) \ra 0$,
where $ \deg ( \O_Y(m-k) \geq \mu (\shR(m)) \geq (\O_Y(m+k-d_1-d_2-d_3))$.
Thus we are in the situation of \ref{exactsequencedecide} (i) and (iv).
\qed

\medskip
We may also deduce a result about the plus closure.

\begin{corollary}
\label{exactsequence}
Let $R$ denote a two-dimensional normal standard-graded domain
over an algebraically closed field $K$ of positive characteristic.
Let $f_1, f_2,f_3 \in R$ denote three homogeneous
$R_+$-primary elements.
Suppose that there exists a primary relation of
total degree $k \leq (d_1+d_2+d_3)/2$.
Then $(f_1,f_2,f_3)^\soclo = (f_1, f_2,f_3)^{\gr}$.
\end{corollary}

\proof
Let $c \in H^1(Y,\shR(m))$
denote the cohomology class of a homogeneous element
$f_0 \in R$ of degree $m$.
We look at the sequence from \ref{primarysemistable},
$ 0 \ra \O_Y(m-k) \ra \shR(m) \ra \O_Y(m+k-d_1-d_2-d_3) \ra 0$
and run through the cases according to \ref{exactsequencedecide}.
If $m \geq d_1+d_2+d_3 -k$, then $f_0 \in (f_1, f_2,f_3)^+$
follows from \ref{exactsequencecrit} (v).
So suppose that $m < d_1+d_2+d_3 -k $.
If the image of $c$ in $H^1(Y,\O_Y(m+k -d_1-d_2-d_3))$ is $\neq 0$,
then $f_0 \not\in (f_1,f_2,f_3)^\soclo$ by \ref{exactsequencecrit} (i).
So we may assume that
$c$ stems from $e \in H^1(Y, \O_Y(m-k))$.
If $m \geq k$ or $c=0$,
then the pull back of $e$ under a finite mapping is zero,
hence $f_0 \in (f_1, f_2,f_3)^+$.
If $m < k$ and $c \neq 0$, then $f_0 \not\in (f_1, f_2,f_3)^\soclo$
by \ref{exactsequencecrit}(ii).
\qed

\begin{example}
\label{examplerelationindecomposable }
Let $K$ denote an algebraically closed field
and consider
$$R=K[x,y,z]/(x^\delta+ay^\delta+bz^\delta+cxz^{\delta-1}+dyz^{\delta-1}) \,  $$
where $a,b,c,d \neq 0 \, $ are chosen such that $Y= \Proj \, R$ is smooth.
Consider the relation sheaf for the elements
$x^\delta,\, y^\delta,\, z^\delta$.
Then we have a relation of total degree $\delta +1$, given by
$(z,az, bz+cx+dy)$.
Since $\delta +1 < \frac{3 \delta }{2} - \frac{\delta-3}{2}
= \frac{3 \delta}{2} - \frac{g-1}{\delta} $,
it follows from \ref{maxinvertbound}(iii)
that $\shR$ is decomposable.

The relation $(z,az, bz+cx+dy)$ is primary if and only if
$cx+dy$ and $x^\delta + ay^\delta$ have no common
homogeneous zero. This is true if and only if
$(-\frac{d}{c})^\delta \neq -a$. If this is true,
then we have the splitting
$\shR( \delta +1 )= \O_Y \oplus \O_Y(- \delta +2 )$,
where the second summand corresponds
to a relation of total degree $2 \delta -1$.
We can find such a relation in the following way:
There exists a polynomial $P(x,y)$ in $x$ and $y$ of degree $\delta -1$
such that
$(cx+dy) P(x,y) = rx^\delta +sy^\delta$.
Then $(P+rz^{\delta -1},aP+sz^{\delta -1}, bP )$ is a relation
of total degree $2\delta-1$, since
$Px^\delta + aPy^\delta + rx^\delta z^{\delta -1} + sy^\delta z^{\delta -1}
+  bPz^\delta
= Px^\delta +aPy^\delta +Pbz^\delta +Pcxz^{\delta -1} +Pdyz^{\delta -1}=0$.
\end{example}

\begin{corollary}
\label{bastel}
Let $f_1,f_2,f_3 \in K[x,y,z]$ be homogeneous elements of degree $d_1,d_2,d_3$
such that $d_1+d_2+d_3=2k$ is even and $k \geq d_i$ for $i=1,2,3$.
Let $g_1,g_2,g_3 \in K[x,y,z]$
be homogeneous of degree $k-d_i$.
Suppose that $V(f_1,f_2,f_3)=V(g_1,g_2,g_3)=V(x,y,z)$.
Set $F=f_1g_1+f_2g_2+f_3g_3$ and suppose that
$R=K[x,y,z]/(F)$ is a normal domain.
Then the sheaf of relations $\shR(m)$ for $f_1,f_2,f_3$ on
$Y= \Proj R$ is an extension of
the structure sheaf by itself and is strongly semistable.
In particular
$$(f_1,f_2,f_3)^\soclo = R_{\geq k} +(f_1,f_2,f_3) \, .$$
If furthermore the characteristic of $K$ is positive,
then $(f_1,f_2,f_3)^\soclo = (f_1,f_2,f_3)^{+ \rm gr}$.
\end{corollary}

\proof
The relation $(g_1,g_2,g_3)$ is primary
of total degree $k$, thus this follows from
\ref{primarysemistable}.
\qed

\begin{example}
\label{fermatexample}
Consider a Fermat polynomial $x^k+y^k+z^k \in K[x,y,z]$
and let $R=K[x,y,z]/(x^k+y^k+z^k)$.
Let $f=x^{d_1}, \, g=y^{d_2}, h=z^{d_3}$ such that $d_i \leq k$
and $d_1+d_2+d_3 =2k$.
Then we are in the situation of \ref{bastel},
we just may take
$(x^{k-d_1},y^{k-d_2},z^{k-d_3})$ as a primary relation.
Therefore we get
$R_{\geq k} \subseteq (x^{d_1}, y^{d_2}, z^{d_3})^\pasoclo$,
for $d_1+d_2+d_3=2k$ and $d_i \leq k$.
We also get that
$R_{\geq k} \subseteq (x^{d_1}, y^{d_2}, z^{d_3})^{\gr}$
in positive characteristic.

For instance we get
$xyz \in (x^2,y^2,z^2)^\pasoclo$ modulo $x^3+y^3+z^3=0$.
This was stated in \cite{huneketightparameter} as an elementary
example of what is not known
in tight closure theory. The first proof was given in \cite{singhcomputation}.

For $R=K[x,y,z]/(F)$, where $F=x^4+y^4+z^4$,
we get that $R_{\geq 4} \subseteq (x^3,y^3,z^2)^\pasoclo $
($=(x^3,y^3,z^2)^{+ \rm gr}$)
etc.
\end{example}

\begin{example}
Let $F \in K[x,y,z]$ be a homogeneous equation for an elliptic curve.
When is $xyz  \in (x^2,y^2,z^2)^\soclo$ in $R=K[x,y,z]/(F)$?
If the coefficient of $F$ in $xyz$ is not zero, then of course
$xyz \in (x^2,y^2,z^2)$.
Thus we may write $F=Sx^2+Ty^2+Uz^2$ so that
$(S,T,U)$ is a homogeneous relation for $(x^2,y^2,z^2)$
of total degree $3$.
If this relation is primary, i.e. $V(S,T,U)=V(R_+)$,
then $xyz \in (x^2,y^2,z^2)^\soclo$.

If however the relation $(S,T,U)$
is not primary, e. g. for $F=x^3+y^3+(x+y)z^2$,
then we have a decomposition $\shR(3) = \O(P) \oplus \O(-P)$
($P$ a point) and $xyz \not\in (x^2,y^2,z^2)^\soclo$,
since $H^1(Y,\shR(3)) = H^1(Y, \O(-P))$
and \ref{quotientnegdegree}.
\end{example}

\smallskip 

\subsection{Correlations and big forcing divisors}

\label{correlations}

\markright{Correlations and big forcing divisors}

\

\bigskip
Let $K$ denote further an algebraically closed field
and let $R$ denote a two-di\-men\-sio\-nal normal standard-graded
$K$-domain, $Y= \Proj\, R$.
Let $\shR(k)$ denote the sheaf of relations
for primary homogeneous elements $f_1, \ldots ,f_n \in R$
and $\shF(-k)$ its dual sheaf on $Y= \Proj\, R$.
The sections $s \in \Gamma(Y, \shF(-k))$ correspond to sheaf morphisms
$s : \O_Y \ra \shF(-k)$ and to {\em correlations}
$s^\dual: \shR(k) \ra \O_Y$ of total degree $k$.
These sections correspond also to sections
in $\Gamma(\PP({\shF}), \O_{\PP(\shF(-k))} (1) )$.

The existence of correlations of certain degrees has the same
consequences on $\shR(k)$ as the existence of relations of certain degrees.

\begin{proposition}
\label{correlationbound}
Let $f_1, \ldots, f_n$ denote homogeneous primary elements
in a normal two-dimensional $K$-domain $R$ of degree $d_i$.
Suppose that $Y= \Proj\, R$ has genus $g$ and degree $\delta$.

\renewcommand{\labelenumi}{(\roman{enumi})}
\begin{enumerate}

\item
Suppose that there exists a correlation $\neq 0$
for the elements $f_1, \ldots, f_n$ of total
degree $k > (d_1 + \ldots +d_n )/(n-1)$.
Then the sheaf of relations is not semistable.
{\rm (}and not stable for $\geq ${\rm )}.

\item
Let $n=3$. Suppose that there exists a correlation $\neq 0$
of total degree $k > (d_1 + d_2 +d_3 )/2 + (g-1)/\delta  $.
Then the sheaf of relations is decomposable.
\end{enumerate}

\end{proposition}
\proof
Same proof as for \ref{maxinvertbound}(i) and (iii).
\qed

\medskip
Let $f_0$ be another homogeneous element of degree $m$ defining
the forcing sequence
$0 \ra \O_Y \ra \shF'(-m) \ra \shF(-m) \ra 0$.

Suppose that there exists a section
$0 \neq s \in \Gamma(Y, \shF'(-m-1))$.
Then we have seen in \ref{bigaffine} that the forcing divisor
$Z=\PP(\shF) \subset \PP(\shF')$ is big,
since it is linearly equivalent to an effective divisor $D$ with an
affine complement.
Such a section induces also a section
$0 \neq s \in \Gamma(Y,\shF(-m-1))$ which is a correlation
of degree $m+1$.
In this situation it is sometimes possible
to decide whether $f_0 \in (f_1,f_2,f_3)^\soclo$ holds or not
on the divisor $D$.

\begin{proposition}
\label{nupositive}
Let $R$ be a normal two-dimensional standard-graded $K$-do\-main.
Let $f_1, f_2,f_3$ be homogeneous primary elements
and let $f_0$ be another homogeneous element of degree $m$.
Suppose that there exists an effective divisor $L \subset Y$
of positive degree such that
$Z- \pi^* L$ is equivalent to an effective divisor.
Then there exists an effective divisor $D$,
$Z \sim D=H+F$, where $H$ is the horizontal component
and $F$ the fiber components.
Moreover, the following hold.

\renewcommand{\labelenumi}{(\roman{enumi})}
\begin{enumerate}

\item 
If $H-Z \cap H$ is not affine, then $f_0 \in (f_1,f_2,f_3)^\pasoclo$.

\item
If $H-Z \cap H$ is affine
{\rm(}this is fulfilled when the pull back $Z|_H$ is ample or
when $Z \cap H$ contains components which lie in a fiber{\rm)},
then there does not exist a finite graded solution
for the tight closure problem, i.e. $f_0 \not\in (f_1,f_2,f_3)^{+{\rm gr}}$.
\end{enumerate}
\end{proposition}

\proof
Let $D' \sim Z -\pi^*L$ be an effective divisor.
Then $Z \sim D=D' + \pi^*L$ may be written
as $D= H+F$, where $F$ consists of all fiber components.
Then $F \neq 0$ and $H$ is a projective subbundle,
since it intersects every fiber in a line.

So we look at the intersection $Z \cap H \subset H$.
(i). From \ref{affinemappingcrit}(i) we get
that $\PP(V')-Z$ is not affine, hence $f_0 \in (f_1,f_2,f_3)^\pasoclo$.

(ii).
We have to show that the forcing divisor
$Z$ intersects every curve $C \not\subseteq Z$ positively.
$Z.C \sim (H+F). C >0$ for $C \not\subseteq H$
(as in \ref{effectivevertreter}).
If however $C \subset H$, then the assumption (ii)
shows that $C \cap (Z \cap H) \neq \emptyset $.
\qed

\medskip
We close this chapter with some further examples.

\begin{example}
Let $Y= \PP^1_K= \Proj\, K[x,y]$ and consider the projective bundles
corresponding to the forcing data $x,y,1;1$.
Then $Z$ is big, $Z^3=2 >0$ and $Z$ is numerically
effective, but there exists a curve $C$ disjoined to $Z$
(the solution section corresponding to $1 \in (x,y,1)$)
and the complement of $Z$ is not affine and the forcing divisor
is not ample.

The relation algebra is $K[x,y][T_1,T_2,T_3,T_4]/(xT_1+yT_2+T_3+T_4)$.
Eliminating $T_3$ in the forcing
equations yields the splitting forcing sequence
$$0 \lra \AA_Y(1) \times \AA_Y(1) \lra \AA_Y(1) \times \AA_Y(1) \times \AA_Y
\lra \AA_Y \lra 0 \, .$$
The forcing subbundle is $Z \cong \PP^1 \times \PP^1$ given by
the equation $T_4=0$.
We have $\shF'(0) = \O_Y(1) \oplus \O_Y(1) \oplus \O_Y$
and also $\shF'(-1)=\shF'(0) \otimes \O_Y(-1)$ has sections $\neq 0$.
Therefore $Z$ is big.
A section is for example $xT_1$, thus a divisor $D$
linearly equivalent to $Z$ is given by $D=H+F$, where
$H=\{ T_1 =0 \}$ and $F= \{ x=0 \}$.
$H$ is a Hirzebruch surface $\PP( \AA_Y \times \AA_Y(-1))$
(the blowing up of a projective plane).
The intersection
$Z \cap H$ is on $Z \cong \PP^1 \times \PP^1$ a horizontal fiber
and on $H= \PP(\AA_Y \times \AA_Y(-1))$ a line not meeting the exceptional divisor,
which is also the solution section $C$.
The self intersection number of $Z \cap H$ on $H$ is
$$ (Z|_H)^2=  Z^2.H =Z^2(Z-F)=Z^3-Z^2.F= 2-1=1 \, ,$$
hence $Z$ is numerically effective.
(The self intersection of $Z \cap H$ on $Z$ is 0.)
\end{example}

\begin{picture}(12,5)
\thicklines

\put(4.85, 0.8){\line(1,0){3.5}}

\put(3.7,1.2){}
\put(9,.8){$ \PP^1$}
\put(3.9,2.2){$ Z$}
\put(5.305, 4.05){$H$}
\put(8.3, 1.55){$F$}
\put(8.8,2.6){ $Z \cap H$}
\put(3.9 ,3.35){$C$}

\put(4.5,2.3){\usebox{\viereckv}}

\put(4.5,1.5){\usebox{\viereckfaserrand}}

\put(4.85,1.85){\usebox{\rechteckdreifuenfundzwei}}

\put(4.85, 2.65){\line(1,0){3.5}}
\put(4.85, 3.4){\line(1,0){3.5}}

\end{picture}

\begin{example}
Let $K$ denote an algebraically closed field and consider
$$R=K[x,y,z]/(x^4+ay^4+bz^4+cxz^3+dyz^3) \,  $$
where $a,b,c,d \neq 0 \, $ are chosen such that $Y= \Proj \, R$ is smooth.
The elements $y$ and $z$ are para\-meters and
we consider the forcing data $x^4,\, y^4,\, z^4$ and $ xy^2z^3$.
We want to show that both cases described
in \ref{nupositive} do actually occur depending on the coefficients.

First we show that $\Gamma(Y,\shF'(-7)) \neq 0$.
This sheaf is the sheaf of linear forms
for the geometric vector bundle
$$V'(-7) = D_+(x,y) \subseteq \Proj \, A \, ,$$
where 
$$A = R[T_1,T_2,T_3,T_4]/(x^4T_1+ y^4T_2 +z^4T_3 + xy^2z^3T_4)$$
is graded by $\deg (T_1)= \deg (T_2) =\deg (T_3)=3$ and $ \deg (T_4)=1$.
From the curve equation and the relation equation we get
$$ z^3(- (bz+cx+dy)T_1 +zT_3 +xy^2T_4) = y^4(aT_1 -T_2)\, .$$
This gives us the global function $G$ on $V'(-7)$ given by
$$-\frac{bz+cx+dy}{y^4}T_1+\frac{z}{y^4}T_3 + \frac{x}{y^2}T_4 \mbox{ on }
D_+(y)\, \mbox{ and }\,
+\frac{a}{z^3}T_1 - \frac{1}{z^3}T_2 \mbox{ on } D_+(z) \, ,$$
which is a linear form.
We consider the function $yG \in \Gamma(Y, \shF'(-6))$.
It induces a section in
$\Gamma(\PP(\shF'), \O_{\PP(\shF'(-6))}(1))$, and this gives
the linear equivalence
$$ Z=\PP(V) =V_+(T_4) \sim V_+(yG)=V_+(G) + V_+(y)=H+F $$
on $\PP(V')$
and we are in the situation of \ref{nupositive}
($H=V_+(G)$ does not contain fiber components).

When does the intersection $ Z \cap H$ have fiber components?
If $z \neq 0$, then the equation for $G$ on $D_+(z)$ does not vanish,
thus there cannot be fiber components.
So look at $z=0$. Then $x^4+ay^4 =0$ and
the equation for $G$ on $D_+(y)$ becomes just $\frac{cx+dy}{y^4}T_1=0$
(note that $T_4=0$, since we are on $Z$).
Thus there exists a fiber component if and only if
$cx+dy=0=x^4+ay^4$ has a solution, and this means $(d/c)^4=-a$.

Consider the equation $x^4-y^4+z^4+xz^3+yz^3=0$.
This yields a smooth curve and the
intersection has fiber components, hence every curve
intersects the forcing divisor and therefore
$ xy^2z^3 \not\in (x^4,y^4,z^4)^{+{\rm gr}}$.

The equation $x^4+y^4+z^4+xz^3+yz^3=0$ yields also a smooth curve,
and here the intersection does not have a fiber component.
Hence the intersection $Z \cap H$ is irreducible
and its self intersection number
on $H=V_+(G)$ is negative (it is $-4$, as the computation
in the previous example shows). Therefore the complement of
it cannot be affine, hence the complement of the forcing
divisor $Z$ is not affine, thus
$xy^2z^3 \in (x^4,y^4,z^4)^\pasoclo$.
\end{example}

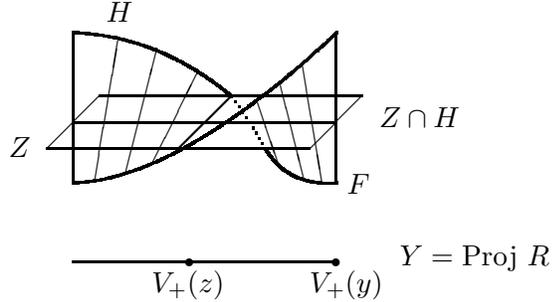
\begin{picture}(12,5)

\thicklines

\put(4.85, .8){\line(1,0){3.5}}

\put(6.4, .8){\circle*{0.1}}
\put(8.35, .8){\circle*{0.1}}

\put(5.9 ,.4){$V_+(z)$}
\put(8 ,.4){$V_+(y)$}

\put(3.7,1.2){}
\put(9.2 ,.75){$Y = \Proj\, R$}
\put(4,2.2){$Z$}
\put(5.3, 4){$H$}
\put(8.5, 1.7){$F$}
\put(8.8,2.6){ $Z \cap H$}

\put(4.5,2.3){\usebox{\viereckv}}

\put(4.5,1.5){\usebox{\viereckfaserrand}}

\put(4.85, 2.65){\line(1,0){3.5}}

\put(4.85, 1.85){\line(0,1){2}}

\put(6.25, 2.3){\line(1,1){.7}}

\put(8.35, 1.85){\line(0,1){2}}

\thinlines

\put(5.14, 1.88){\line(1,6){.31}}
\put(5.55, 1.97){\line(1,4){.42}}
\put(5.9, 2.1){\line(1,2){.62}}

\put(7.3, 2.95){\line(1,-3){.3}}
\put(7.6, 3.2){\line(1,-4){.33}}
\put(7.9, 3.44){\line(1,-6){.26}}

\thicklines

\bezier{300}(4.85,1.85)(6.25, 1.9)(8.35,3.85)

\bezier{300}(4.85, 3.85)(6., 3.8)(6.95,3.)

\bezier{10}(6.95,3)(7.2, 2.65) (7.4, 2.3)

\bezier{300}(7.4,2.3)(7.7, 1.8)(8.35,1.85)

\end{picture}

\smallskip
In the case of a projective bundle of rank two
we cannot characterize the relation $f_4 \not\in (f_1,f_2,f_3)^\pasoclo$
by the ampleness of the forcing divisor, as the following
example shows.
The first example of an affine open subset in a three-dimensional
smooth projective variety with no ample divisor
on the complement was given by Zariski and described in \cite{goodman}.

\begin{example}
\label{affinenonample}
Let $R=K[x,y]$ and consider on $\PP^1=\Proj\, R$
the projective bundle of rank two
defined by the forcing data
$x^4,y^4,x^4;x^3y^3$.
The third self intersection number of the forcing subbundle $Z$
is zero, hence $Z$ is not ample.
But the complement of $Z$ is affine,
since $x^3y^3 \not\in (x^4,y^4) =(x^4,y^4)^\pasoclo 
= (x^4,y^4,x^4)^\pasoclo $ in the regular ring $K[x,y]$.
Therefore $Z$ is also big.

$Z$ is not numerically effective:
for $m=d_0=6$ the forcing divisor $Z$ is a hyperplane section,
i.e. a global section of $\O_{\PP(\shF'(-6))}(1)$.
The pull back of $Z$ on $Z$ yields the hyperplane section
on the ruled surface
$Z= \PP(x^4,y^4,x^4)$ for this grading.
We have $\shR(8) = \O_{\PP^1} \oplus \O_{\PP^1}(4)$
given by the relations $(y^4,-x^4,0)$ and $(1,0,1)$,
hence $Z \cong \PP(\O_{\PP^1} \oplus \O_{\PP^1}(-4))$.
Let $E=Z|_Z$ denote this hyperplane section, let $C =\PP(x^4,y^4) \subset Z$
be the forcing section
and let $L$ be a disjoined section corresponding to $x^4 \in (x^4,y^4)$.
Then we know by \ref{forcingsequence2}(iii) that $C \sim E + 2\pi^*H$ (where
$H$ is the hyperplane section on $\PP^1$)
and therefore
$E.L = C.L - 2\pi^* H.L = - 2\pi^* H.L < 0$.

The divisor $Z$ is also not semiample:
We know that there exists a curve  $L$ such that
$Z.L <0$. Let $P \in L$ and suppose
there exists an effective divisor
$D \sim aZ$ such that $P \not\in D$. Then
$L \not\subset D$ yields a contradiction.

The divisor is also not almost basepoint free in the sense of
\cite[conjecture III.5.2]{haramp}, so this gives also a counter example
to this conjecture.
If we take $x^5,y^5,x^5;x^4y^4$, then the self intersection is negative $-1$,
and the other properties hold again.
\end{example}

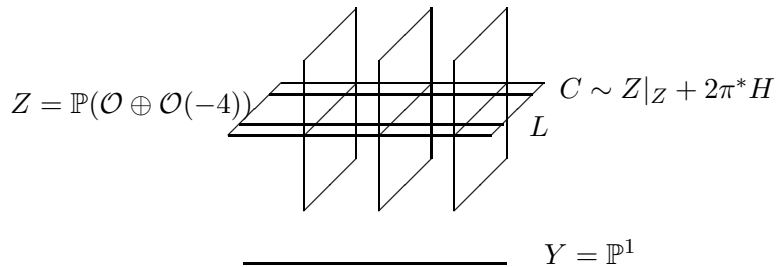
\savebox{\viereckfaser}(5,4)[bl]{
\thinlines
\put(0,0){\line(0,1){2}}
\put(0,0){\line(1,1){0.7}}
\put(0.7, .7){\line(0,1){2}}
\put(0,2){\line(1,1){0.7}}
}

\unitlength1cm
\begin{picture}(12,5)

\thicklines

\put(4.7,.8){\line(1,0){3.5}}
\put(8.7 ,0.8){$Y = \PP^1 $}
\put(10,2.5){$ $}
\put(1.6, 2.8){$ Z=\PP(\O \oplus \O(-4))$}

\put(8.5,2.5){$L$}
\put(8.9,3.){$C \sim Z|_Z +2 \pi^*H$}

\put(4.65, 2.65){\line(1,0){3.5}}

\put(5.05, 3.05){\line(1,0){3.5}}

\put(4.5,2.5){\usebox{\viereckv}}

\thinlines
\put(5.5,1.5){\usebox{\viereckfaser}}
\put(5.5,2.5){\line(1,1){0.7}}
\put(6.5,1.5){\usebox{\viereckfaser}}
\put(6.5,2.5){\line(1,1){0.7}}
\put(7.5,1.5){\usebox{\viereckfaser}}
\put(7.5,2.5){\line(1,1){0.7}}
\end{picture}

\begin{remark}
Suppose that $D \subset X$ is an effective divisor on a smooth
projective curve. Suppose that the complement $X -D$ is affine.
Then the theorem of Goodman tells us that there exists
a modification $X' \ra X$ such that the pull-back of $D$
is an ample divisor on $X'$, see \cite[II Theorem 6.1]{haramp}.
\end{remark}

\newpage

\thispagestyle{empty}

\markboth{6. Tight closure and plus closure
in cones over elliptic curves}
{6. Tight closure and plus closure
in cones over elliptic curves}

\section{Tight closure and plus closure in cones over elliptic curves}

\label{tightelliptic}

\bigskip
In this last chapter we study locally free sheaves
and their geometric torsors over an elliptic curve.
Due to the classification of Atiyah we know much more
about vector bundles on an elliptic curve than
about vector bundles on  a curve of higher genus.
If $c \in H^1(Y,\shS)$ is a cohomology class and if
$\PP(\shG') -\PP(\shG)$ is the corresponding torsor ($\shG = \shS^\dual$),
then we give an easy numerical condition for the affineness of
$\PP(\shG') -\PP(\shG)$ in terms of the
indecomposable summands  of $\shS$.
Moreover, the same numerical condition characterizes
in positive characteristic the non-existence of projective
curves inside $\PP(\shG') -\PP(\shG)$.
From this we derive that for a normal
homogeneous coordinate ring over an elliptic curve over a field
of positive characteristic the tight closure of a primary
homogeneous ideal is the same as its plus closure.

\smallskip

\subsection{Vector bundles over elliptic curves}
\label{elliptic}

\markright{Vector bundles over elliptic curves}

\

\bigskip
Let $Y$ denote an elliptic curve over an algebraically closed field $K$.
We gather together some results on vector bundles over elliptic curves.

\begin{proposition}
\label{ellipticproperties}
Let $Y$ denote an elliptic curve over an algebraically closed field $K$
and let $\shG$ denote a locally free sheaf on $Y$.
Then the following hold.

\renewcommand{\labelenumi}{(\roman{enumi})}
\begin{enumerate}

\item
Suppose that $\shG$ is indecomposable and of positive degree,
$\deg \, (\shG) >0$.
Then $H^1(Y, \shG)=H^0(Y, \shG^\dual)=0$.

\item
If $\shG$ is indecomposable, then $\shG$ is semistable.

\item
If $\shG$ is semistable, then $\shG$ is also strongly semistable.

\end{enumerate}

\end{proposition}

\proof
(i). See \cite[Lemma 1.1]{hartshorneamplecurve}.

(ii). See \cite[Appendix A]{tusemistable}.

(iii). See \cite[\S 5]{miyaokachern}.
\qed

\medskip
The following theorem of Gieseker-Hartshorne
gives an easy numerical criterion for ample bundles over an elliptic curve.

\begin{theorem}
\label{amplekritelliptic}
Let $Y$ denote an elliptic curve over an algebraically closed field $K$
and let $\shG$ denote a locally free sheaf.
Then $\shG$ is ample if and only if the degree of every
indecomposable summand of $\shG$ is positive.
\end{theorem}

\proof
See \cite[Proposition 1.2 and Theorem 1.3]{hartshorneamplecurve} or
\cite[Theorem 2.3.]{giesekerample}.
\qed

\begin{corollary}
\label{ellipticextension}
Let $Y$ denote an elliptic curve over an algebraically closed field $K$
and let $\shS$ denote a locally free sheaf with dual sheaf $\shG=\shS^\dual$.
Suppose that $\shG$ is ample and let
$ 0 \ra \O_Y \ra \shG' \ra \shG \ra 0$ be a non-trivial extension
given by $0 \neq c \in H^1(Y,\shS)$.
Then $\shG'$ is also ample and
$\PP(\shG') - \PP(\shG)$ is affine.
\end{corollary}

\proof
This follows from \ref{amplekritelliptic}
and \cite[Proposition 2.2]{giesekerample}.
\qed

\medskip
The following theorem is a generalization of a theorem of Oda.

\begin{theorem}
\label{odageneral}
Let $Y$ denote an elliptic curve over an algebraically closed field $K$
and let $\shS$ denote an indecomposable locally free sheaf of negative degree.
Let $g: X \ra Y$ be a finite dominant map, where $X$ is another smooth
projective curve.
Then $H^1(Y, \shS) \lra H^1(X, g^* \shS)$ is injective.
\end{theorem}

\proof
Let $0 \neq c \in H^1(Y, \shS)$ be a non zero class and consider
the corresponding extension
$$ 0 \lra \shS \lra \shS' \lra \O_Y \lra 0 \, .$$
Let $\shG$ and $\shG'$ denote the dual sheaves and let
$\PP(\shG) \hookrightarrow \PP(\shG')$ be the corresponding
projective subbundle.
The indecomposable sheaf $\shG$ is of positive degree,
hence ample due to \ref{amplekritelliptic}.
Then also $\shG'$ is ample due to \ref{ellipticextension}.
Hence $\PP(\shG)$ is an ample divisor on $\PP(\shG')$ and
its complement is affine. This property is preserved under the finite
mapping $X \ra Y$, therefore the complement of
$\PP(g^*(\shG)) \subset \PP(g^*(\shG'))$ is affine and there cannot
be projective curves in the complement. Hence the pull-back of the sequence
does not split and $g^*(c) \neq 0$.
\qed

\begin{remark}
Oda proved this statement only for the Frobenius morphism in
positive characteristic, see \cite[Theorem 2.17]{odaelliptic},
and Hartshorne used this theorem to prove the numerical criterion for
ampleness. The proof of this criterion by Gieseker
in \cite{giesekerample}
however is independent
of the theorem of Oda. 
\end{remark}

\medskip
A kind of reverse to \ref{odageneral} is given by the following lemma.

\begin{lemma}
\label{pvanish}
Let $K$ denote an algebraically closed
field of positive characteristic $p$ and let $Y$ be an elliptic
curve.
Let $\shS$ be an indecomposable locally free sheaf on $Y$ of degree
$ \deg (\shS) \geq 0$.
Let $c \in H^1(Y, \shS)$ be a cohomology class.
Then there exists a finite curve
$g: X \ra Y $ such that $g^*(c) \in H^1(X, g^*(\shS))$ is zero.
\end{lemma}

\proof
In fact we will show that the multiplication mappings $[p^{e}]:Y \ra Y,
\, y \mapsto p^{e}y$
have the stated property for $e \gg 0$.
If $\deg \, \shS > 0$, then
$H^1(Y, \shS) \cong H^0 (Y, \shS^\dual)=0$ due to
\ref{ellipticproperties}(i) and there is nothing to prove.
The same is true for an invertible sheaf $\neq \O_Y$ of degree $0$.
For $\shS =\O_Y$, the multiplication map
$[p]: Y \ra Y$ induces the zero map on
$H^1(Y, \O_Y)$, this follows for example
from \cite[Corollary 5.3]{silvermanelliptic} and
\cite[\S 13 Corollary 3]{mumfordabelian}.

Now we do induction on the rank, and
suppose that $r=\rk\, \shS \geq 2$ and $\deg \, \shS =0$.
Due to the classification of Atiyah (see \cite[Theorem 5]{atiyahelliptic})
we may write
$\shS = \shF_r \otimes \shL$, where $\shL$ is an invertible sheaf
of degree $0$ and where
$\shF_r$ is the unique sheaf of rank $r$ and degree zero with
$\Gamma(Y, \shF_r) \neq 0$.
In fact, for these sheaves we know that $H^0(Y, \shF_r)$ and $H^1(Y, \shF_r)$
are one-dimensional and there exists a non-splitting
short exact sequence
$$ 0 \lra \O_Y \lra \shF_{r} \lra \shF_{r-1} \lra  0 \, .$$
This gives the sequence
$ 0 \ra \shL  \ra \shS \ra \shF_{r-1} \otimes \shL \ra 0$.
Let $c \in H^1(Y, \shS)$. Then the image of this class in
$H^1(Y, \shF_{r-1} \otimes \shL)$ is zero after applying
$[p^{e}]$ and comes then from an element in $H^1(Y, \shL)$, which itself is
zero after applying $[p]$ once more.
\qed

\bigskip 

\subsection{A numerical criterion for subbundles to have affine complement}
\label{ellipticaffin}

\markright{A numerical criterion for subbundles
to have affine complement}

\bigskip
In this section we investigate subbundles $\PP(\shG) \subset \PP(\shG')$
of codimension one over an elliptic curve with respect to the properties
which are interesting from the tight closure and plus closure point of view:
Is the complement affine? Does it contain projective curves?

\medskip
Our first main result is the following numerical characterization
for the complement of a projective subbundle to be affine.

\begin{theorem}
\label{numkritaffinelliptic}
Let $K$ denote an algebraically closed field and let
$Y$ denote an elliptic curve.
Let $\shS$ be a locally free sheaf on $Y$ and let $\shS = \shS_1 \oplus \ldots \oplus \shS_s$
be the decomposition into indecomposable locally free sheaves.
Let $c \in H^1(Y, \shS)$ and let
$0 \ra \shS \ra \shS' \ra \O_Y \ra 0$ be the corresponding
extension and let $\PP(\shG) \subset \PP(\shG')$ be the corresponding
projective bundles. Then the following are equivalent.

\renewcommand{\labelenumi}{(\roman{enumi})}
\begin{enumerate}

\item
There exists
$1 \leq j \leq s$ such that $\deg \shS_j < 0$
and $c_j \neq 0$, where $c_j $ denotes the component of $c$
in $H^1(Y, \shS_j)$.

\item
The complement $\PP(\shG') - \PP(\shG)$ is affine.
\end{enumerate}

\end{theorem}
\proof
Proposition \ref{ellipticproperties} shows that
we have a decomposition into strongly
semistable sheaves.
Hence the result follows from
\ref{splittingsemistable}.
For (i) $\Rightarrow$ (ii) in small positive characteristic
we may use \ref{amplekritelliptic} and \ref{ellipticextension}
to show that $\PP( \shG_j ') - \PP(\shG_j)$ is affine
and then $\PP(\shG') - \PP(\shG)$ is affine due to \ref{affinelemma}.
\qed

\medskip
We shall show now that the same numerical criterion
holds in positive characteristic for the (non-)existence of projective curves
inside $\PP(\shG') -\PP(\shG)$.

\begin{theorem}
\label{numkritcurve}
Let $K$ denote an algebraically closed field of positive characteristic and let
$Y$ denote an elliptic curve.
Let $\shS$ be a locally free sheaf on $Y$
and let $\shS = \shS_1 \oplus \ldots \oplus \shS_s$
be the decomposition into indecomposable locally free sheaves.
Let $c \in H^1(Y, \shS)$ and let
$0 \ra \shS \ra \shS' \ra \O_Y \ra 0$ be the corresponding
extension and let $\PP(\shG) \subset \PP(\shG')$ be the corresponding
projective bundles.
Then the following are equivalent.

\renewcommand{\labelenumi}{(\roman{enumi})}
\begin{enumerate}

\item
There exists
$1 \leq j \leq s$ such that $\deg \shS_j < 0$
and $c_j \neq 0$, where $c_j $ denotes the component of $c$
in $H^1(Y, \shS_j)$.

\item
The sequence
$0 \ra \shS \ra \shS' \ra \O_Y \ra 0$ does not split after
a finite dominant morphism $X \ra Y$, where $X$ is another projective
curve.

\item
The subbundle $\PP(\shG) \subset  \PP(\shG')$ intersects every curve
in $\PP(\shG')$ positively.

\end{enumerate}

\end{theorem}
\proof

(i) $\Rightarrow$ (ii).
Suppose that the sequence splits under the finite morphism
$g:X \ra Y$. Then $g^*(c_j)=0$ on $X$ and from
\ref{odageneral}
we see that $\deg \, \shS_j \geq 0$ or $c_j=0$

(ii) $\Rightarrow$ (iii).
If there exists a curve $C$ on
$\PP(\shG')$ not meeting $\PP(\shG)$,
then it dominates the base. Let $X$ be the normalization of $C$ and
let $g: X \ra Y$ the finite dominant mapping. Then
$g^*\PP(\shG') \ra X$ has a section not meeting $g^*\PP(\shG)$
and then the sequence splits on $X$.

(iii) $\Rightarrow $ (i).
Suppose to the contrary that for all the indecomposable components
of $\shS$ with negative degree $c_j =0$ holds.
For every component $\shS_j$ with $\deg \, S_j \geq 0$ there exists
due to \ref{pvanish}
a finite curve $g_j: X_j \ra Y$ such that
$g_j^*(c_j) \in H^1(X_j, g_j^* \shS_j)$
is zero.
Putting these curves together we find a curve $g: X \ra Y$
such that $g^* (c)=0$.
Thus the sequence splits on $X$ and this gives a projective curve
in $\PP(\shG') - \PP(\shG)$.
\qed

\begin{corollary}
\label{geokritaffin}
Let $K$ denote an algebraically closed field of positive characteristic and let
$Y$ denote an elliptic curve.
Let $\shS$ be a locally free sheaf on $Y$,
let $c \in H^1(Y, \shS)$ and let
$0 \ra \shS \ra \shS' \ra \O_Y \ra 0$ be the corresponding
extension.
Let $\PP(\shG) \subset \PP(\shG')$ be the corresponding
projective bundles.
Then $\PP(\shG') - \PP(\shG)$ is affine if and only if
it contains no projective curve.
\end{corollary}
\proof
This follows from \ref{numkritaffinelliptic} and \ref{numkritcurve}, since for
both properties the same numerical criterion holds.
\qed

\bigskip 

\subsection{Numerical criteria for tight closure and plus closure
for cones over elliptic curves}

\label{sectiontightplus}

\markright{Numerical criteria for tight closure and plus closure
for cones over elliptic curves}

\

\bigskip
We are now in the position to draw the consequences to tight closure
and plus closure in a homogeneous coordinate ring of an elliptic curve.

\begin{corollary}
\label{numkrittight}
Let $K$ be an algebraically closed field and
$R=K[x,y,z]/(F)$ a homogeneous coordinate ring over the
elliptic curve $Y= \Proj \, R$, $\deg (F)=3$.
Let $f_1, \ldots ,f_n$ be homogeneous generators of a
$R_+$-primary ideal in $R$ and let $m \in \NN$ be a number.
Let $\shR(m)$ be the corresponding locally free sheaf of relations
of total degree $m$ on $Y$
and let $\shR(m) = \shS_1 \oplus \ldots \oplus \shS_s$
be the decomposition into indecomposable locally free sheaves.
Let $f_0 \in R$ be another homogeneous element of degree $m$
and let $c \in H^1(Y, \shR(m))$ be its forcing class.
Then the following are equivalent.

\renewcommand{\labelenumi}{(\roman{enumi})}
\begin{enumerate}

\item
There exists
$1 \leq j \leq s$ such that $\deg \shS_j < 0$
and $c_j \neq 0$, where $c_j $ denotes the component of $c$
in $H^1(Y, \shS_j)$.

\item
The complement of the forcing divisor is affine.

\item
$f_0 \not\in (f_1, \ldots ,f_n)^\soclo $.

\end{enumerate}
\end{corollary}
\proof
The equivalence (ii) $ \Leftrightarrow $ (iii) was stated in
\ref{pluscrit},
and (i) $\Leftrightarrow$ (ii) is \ref{numkritaffinelliptic}.
\qed

\begin{corollary}
\label{numkritplus}
Let $K$ denote an algebraically closed field of positive characteristic
and let
$R=K[x,y,z]/(F)$ be a homogeneous coordinate ring over the
elliptic curve $Y= \Proj \, R$, $\deg (F)=3$.
Let $f_1, \ldots ,f_n$ be homogeneous generators of
an $R_+$-primary ideal in $R$
and let $m \in \NN$ be a number.
Let $\shR(m)$ be the corresponding locally free sheaf on $Y$
and let $\shR(m) = \shS_1 \oplus \ldots \oplus \shS_s$
be the decomposition into indecomposable locally free sheaves.
Let $f_0 \in R$ be another homogeneous element of degree $m$
and let $c \in H^1(Y, \shR(m))$ be the corresponding class.
Then the following are equivalent.

\renewcommand{\labelenumi}{(\roman{enumi})}
\begin{enumerate}

\item
There exists
$1 \leq j \leq s$ such that $\deg \shS_j < 0$
and $c_j \neq 0$, where $c_j $ denotes the component of $c$
in $H^1(Y, \shS_j)$.

\item
The forcing divisor intersects every curve.

\item
The complement of the forcing divisor is affine.

\item
$f_0 \not\in (f_1, \ldots ,f_n)^\soclo $.

\item
$f_0 \not\in (f_1, \ldots , f_n)^{\gr}$

\end{enumerate}
\end{corollary}
\proof
Follows directly from \ref{numkritcurve}, \ref{numkrittight}
and \ref{pluscrit}.
\qed

\begin{theorem}
\label{tightpluselliptic}
Let $K$ be an algebraically closed field of positive characteristic
and let
$R=K[x,y,z]/(F)$,
where $F$ is a homogeneous polynomial of degree $3$ and $R$ is normal.
Let $I$ be an $R_+$-primary homogeneous ideal in $R$.
Then $I^{\gr} = I^+ = I^\soclo $.
\end{theorem}
\proof
The inclusions $ \subseteq $ are clear. It is known that the tight closure
of a homogeneous ideal is again homogeneous,
see \cite[Theorem 4.2]{hochsterhunekesplitting}.
Hence the statement follows from \ref{numkritplus}.
\qed

\begin{remark}
If the $p$-rank (=Hasse invariant) of the elliptic curve $Y$ is 0,
then the plus closure (=tight closure) of a primary homogeneous ideal is
the same as its Frobenius closure.
\end{remark}

\medskip
The tight closure of a primary homogeneous ideal is
easy to describe by a numerical condition,
if (there exists a system of homogeneous generators such that)
the corresponding
vector bundle is indecomposable.

\begin{corollary}
\label{indecomposable}
Let $K$ be an algebraically closed field and
$R=K[x,y,z]/(F)$ a homogeneous coordinate ring over the
elliptic curve $Y= \Proj \, R$.
Let $f_1, \ldots ,f_n$ be homogeneous generators of
an $R_+$-primary graded ideal in $R$
with $\deg \, f_i = d_i$.
Suppose that the corresponding locally free sheaf $\shR(m)$ on $Y$
is indecomposable.
Let $f_0 \in R$ be a homogeneous element of degree $m$,
defining the cohomology class $c \in H^1(Y,\shR(m))$.
Then the following are equivalent.

\renewcommand{\labelenumi}{(\roman{enumi})}
\begin{enumerate}

\item
$m < \frac{d_1+ \ldots +d_n}{(n-1)}$ and $f_0 \not \in (f_1, \ldots ,f_n)$.

\item
$\deg \shR(m) < 0$ and the cohomology class is $c \neq 0$.

\item
The sheaf $\shF(-m)=\shR(m)^\dual$ is ample and $c \neq 0$.

\item
The sheaf $\shF'(-m)$ is ample.

\item
The complement of the forcing divisor is affine.

\item
$f_0 \not\in (f_1, \ldots ,f_n)^\soclo \,
{\rm(}= (f_1, \ldots ,f_n)^+ $ in positive
characteristic{\rm )}.

\end{enumerate}

\end{corollary}

\proof
The degree of $\shR(m)$ is $\deg \shR(m) =-3(d_1+ \ldots +d_n -(n-1)m)$
due to \ref{buendelalsproj}, hence (i) $\Leftrightarrow$ (ii) is clear.
(ii) $\Leftrightarrow$ (iii) follows from \ref{amplekritelliptic},
(iii) $\Rightarrow $ (iv) follows from \cite[Proposition 2.2]{giesekerample}
and \ref{amplekritelliptic}, (iv) $\Rightarrow $ (v) is clear
and (v) $\Rightarrow $ (ii) is \ref{numkrittight} for indecomposable $\shR(m)$.
\qed

\begin{corollary}
\label{indecomposable2}
Let $K$ be an algebraically closed field and
$R=K[x,y,z]/(F)$ a homogeneous coordinate ring over the
elliptic curve $Y= \Proj \, R$.
Let $f_1, \ldots ,f_n$ be homogeneous generators of an
$R_+$-primary ideal in $R$
with $\deg \, f_i = d_i$ and suppose
that the corresponding locally free sheaf $\shR(m)$ on $Y$
is indecomposable.
Set $k= \lceil \frac{d_1 + \ldots +d_n}{n-1} \rceil $.
Then
$$(f_1, \ldots ,f_n)^\soclo = (f_1, \ldots ,f_n) + R_{\geq k} $$
\end{corollary}

\proof
This follows from \ref{semistablevanishing} and
\ref{semistablevanishingp} in connection with
\ref{ellipticproperties} or from \ref{indecomposable}.
\qed

\begin{example}
The statements in \ref{indecomposable} and \ref{indecomposable2}
do not hold without the condition that $\shR(m)$ is indecomposable.
The easiest way to obtain decomposable relation sheaves
is to look at redundant systems of generators.
Consider the elliptic curve $Y$ given by $x^3+y^3+z^3=0$.

Look at the elements $x^2,y^2,x^2$.
Then the corresponding sheaf is of course decomposable.
For $m=3$ we have that
$\shR(3) \cong \O_Y(1) \oplus \O_Y(-1)$.
Let $f_0 =xyz$. Then the number $d_1+ \ldots +d_n -(n-1)m $ is zero, but
$xyz \not\in (x^2,y^2,x^2)^\pasoclo =(x^2,y^2)^\pasoclo$.
The complement of the forcing divisor in $\PP(x^2,y^2,x^2,xyz)$
is affine, but it is not ample, since its degree is zero (or
since $\O_Y(-1)$ is a quotient invertible sheaf
of $\shF'(-3)$ of negative degree).

Consider now $x,y,z^3$ and $f_0=z^2$.
Then $z^2 \in (x,y,z^3)^\soclo=(x,y)^\soclo$, but $z^2 \not\in (x,y)$ and the
number in \ref{indecomposable}(i) is $1+1+3-2\cdot 2=1 >0$.
The relation sheaf decomposes
$\shR(3)= \O_Y \oplus \O_Y(+1)$.

The ideal $(x^2,y^2,xy)$ provides an example where no generator is
superfluous, but the corresponding sheaf of relations
is anyway decomposable.
\end{example}

\newpage

\thispagestyle{empty}

\end{document}